# Spectral flow calculations for reducible solutions to the massive Vafa-Witten equations

Clifford Henry Taubes

Department of Mathematics
Harvard University
Cambridge, MA 02139

chtaubes@math.harvard.edu

ABSTRACT: The Vafa-Witten equations (with or without a mass term) constitute a non-linear, first order system of differential equations on a given oriented, compact, Riemannian 4-manifold. Because these are the variational equations of a functional, the linearized equations at any given solution can be used to define an elliptic, first order, self-adjoint differential operator. The purpose of this article is to give bounds (upper and lower) for the spectral flow between respective versions of this operator that are defined by the elements in diverging sequences of reducible solutions. (The spectral flow is formally the difference between the respective Morse indices of the solutions when they are viewed as critical points of the functional.) In some cases, the absolute value of the spectral flow is bounded along the sequence, whereas in others it diverges. This is a curious state of affairs. In any event, the analysis introduces localization and excision techniques to calculate spectral flow which may be of independent interest.

## 1. Introduction

The Vafa-Witten equations were introduce many years ago in a paper by Cumrun Vafa and Edward Witten [VW]; they constitute a system of first order differential equations on an oriented Riemannian 4-manifold that are generalizations of the self-dual Yang-Mills equations (solutions to the latter are also solutions to the former). The Vafa-Witten equations ask that a pair, to be denoted by (A, $\omega$), of connection on a principal SU(2) or SO(3) bundle and self-dual 2-form with values in the associated Lie algebra bundle obey

- $F_A{}^+ - \frac{1}{2}[\omega;\omega] = 0$,
- $d_A\omega = 0$,

(1.1)

where the notation is as follows: What is denoted here by $F_A$ is the curvature 2-form of the connection A; and $F_A{}^+$ is the curvature's self-dual part. Meanwhile, $[\omega;\omega]$ is a very specific linear combination of the commutators of the components of $\omega$ (see (1.6) and (1.7) below for the precise definition). Finally, $d_A$ here and in what follows denotes the exterior covariant derivative that is defined by the connection A. The 'massive' Vafa-Witten equations are slight modifications of (1.1) that require the choice of a real number to be denoted by $m$: These ask that (A, $\omega$) obey

- $F_A{}^+ - \frac{1}{2}[\omega;\omega] - m\omega = 0$,
- $d_A\omega = 0$.

(1.2)

The equations in (1.1) and in (1.2) are invariant with respect to the action of the gauge group of the given principal bundle (the group of bundle automorphisms). This is an infinite dimensional group, so if there is but one solution to (1.2), then there are infinitely many. Because of this, the focus is the quotient of the space of smooth solutions to (1.2) by the gauge group action. This quotient space will be denoted by $\mathcal{M}_m$. As it turns out, $\mathcal{M}_m$ has reasonable local structure: A neighborhood of any given point is homeomorphic to the quotient by the action of a subgroup of SU(2) (or SO(3)) on the zero locus of a smooth, equivariant map from a ball in one Euclidean space to another (see for example [Ma], or [DK] for analogous theorems for self-dual Yang-Mills moduli spaces). With regards to the global structure: A paper by the author [T1] examined the limits of non-convergent sequences in $\mathcal{M}_m$, and in particular, limits of sequences where the integral of $|\omega|^2$ diverges along the sequence. For example, such sequences exist for the equations in (1.1) for the case of the product manifold $S^1 \times Y$ when the moduli space of flat, stable $Sl(2;\mathbb{C})$ connections on Y is non-compact. One purpose of this article is to present some examples of $|\omega|$-divergent sequences of solutions to $m \neq 0$ versions of (1.2) on various 4-manifolds. (These will all be 'reducible' in the sense that $\omega$ can be written locally as $\varsigma^+ \otimes \sigma$ with $\varsigma^+$ denoting a self-dual 2-form and with $\sigma$ denoting an A-covariantly constant section of the associated Lie algebra bundle.)

A brief digression before stating the main purpose of this article: The linearization of the equations in (1.2) at a given solution defines a symmetric first order differential operator which can be used in turn to define a (gauge invariant) spectral flow for each solution to (1.2). This operator is depicted below in (1.4). This spectral flow is relevant by virtue of the fact that solutions to (1.2) correspond to critical points of the gauge invariant functional of pairs $(A, \omega)$ that is depicted below in (1.5), and this correspondence identifies the spectral flow with the 'Morse' index of the critical point.

With this last point understood, one might hope that the spectral flow along all $|\omega|$-divergent sequence in $\mathcal{M}_m$ is likewise divergent in which case these divergent sequences will likely not obstruct hypothetical topological applications. The main purpose of this article is to explain why this is not the case for some (but not all) of the reducible, $|\omega|$-divergent solution sequences that are described in this paper. Even so, these reducible, $|\omega|$-divergent solution sequences with bounded spectral flow have to obey delicate cohomological constraints which deserves further study. Instances with bounded and unbounded spectral flow are described in the upcoming Section 1b.

**a) The relevant operator**

The operator of interest is depicted momentarily; notation for the depiction is introduced first. To start, let X denote the given Riemannian 4-manifold. Its bundle of self-dual 2-forms will be denoted by $\Lambda^+$ in what follows (and the bundle of anti-self-dual 2-forms is denoted by $\Lambda^-$). What is denoted below by $\underline{\mathbb{R}}$ signifies the product real line bundle. Supposing that $P \to X$ denotes the given principle SU(2) or SO(3) bundle, the associated $\mathbb{R}^3$ bundle with fiber the Lie-algebra of SU(2) is denoted by ad(P).

The operators of interest in what follows are certain first order, elliptic operators acting as endomorphisms of the space of sections of the vector bundle

$$((\Lambda^+ \oplus \underline{\mathbb{R}}) \oplus T^*X) \otimes \mathrm{ad}(P) . \tag{1.3}$$

Any given instance of the operator requires for its definition a pair $(A, \omega)$ of connection on P and section of $\Lambda^+ \otimes \mathrm{ad}(P)$. The corresponding operator takes a section $((\varsigma, \phi), a)$ of the bundle depicted in (1.3) to the section of the bundle in (1.3) whose respective $\Lambda^+$, $\underline{\mathbb{R}}$ and $T^*X$ components have the schematic form:

- $(d_A a)^+ - [\omega; \varsigma] - [\omega, \phi] - m\varsigma$ ,

- $*d_A*a - *(\omega\wedge\varsigma - \varsigma\wedge\omega) - m\phi$ ,
- $*d_A\varsigma - d_A\phi + *(a\wedge\omega - \omega\wedge a)$ ,

(1.4)

where $*$ denotes here the metric's Hodge star and where $(\cdot)^+$ denotes the $\Lambda^+$ part of the given ad(P) valued 2-form. The endomorphism $\varsigma \to [\omega; \varsigma]$ is specified precisely below in (1.7). The linear map $\phi \to [\omega, \phi]$ is, component-wise, the commutator with $\phi$.

**b) Some iconic examples**

Start with a compact, oriented, Riemannian 4-manifold, X. A simple construction of an $|\omega|$-divergent sequence of reducible solutions (1.2) on this manifold is as follows: Fix a complex, hermitian line bundle on X (this is denoted by L) whose first Chern class in the second cohomology has zero self-cup product (the cup product of this class with it self pairs to zero with the fundamental class of the manifold). Then, take a unitary connection on L (denoted by $\hat{A}$) whose curvature 2-form is harmonic. For each non-negative integer q: The principle SO(3) bundle on the manifold X is the bundle of oriented, orthonormal frames in the underlying real bundle $L^q \oplus \underline{\mathbb{R}}$ with $\underline{\mathbb{R}}$ denoting here the product vector bundle with fiber $\mathbb{R}$. An orthogonal connection (the connection A) on this principle SO(3) bundle is specified by requiring that the $\underline{\mathbb{R}}$-factor be covariantly constant and that the covariant derivative acting on sections of the $L^q$ factor be the covariant derivative defined by $\hat{A}$. The $\underline{\mathbb{R}}$ factor corresponds to a covariantly constant, norm 1 section of the associated Lie algebra bundle which is denoted by $\sigma$. With $\sigma$ in hand, take $\omega$ to be $\frac{1}{m}$ times the self-dual part of the 2-from $-iqF_{\hat{A}}\sigma$. The pair $(A, \omega)$ just described solves the given $m$ version of (1.2). Because $\omega$ is proportional to q, letting q increase without bound leads to an $|\omega|$-divergent sequence of solutions to (1.2).

Some comments: Because the cup product of the first Chern class with itself is zero all the even q versions of this SO(3) bundle are isomorphic to the product principle SO(3) bundle; and all of the odd q versions are isomorphic to the q = 1 version (which will be isomorphic to the product bundle if the first Chern class of L is divisible by 2 in the second cohomology).

The upcoming Theorem 1.1 summarizes what is said in the subsequent sections about the spectral flow for this sequence as $q \to \infty$: To set the stage for the theorem: The theorem assumes that the self-dual part of the real 2-form $-iF_{\hat{A}}$ vanishes transversly which implies that this zero locus (denoted by Z) is a finite union of 1-dimensional submanifolds. (This can be assumed if the metric is generic). The orthogonal complement to the self-dual part of $-iF_{\hat{A}}$ in the bundle $\Lambda^+$ over X–Z is an oriented 2-plane bundle whose Euler class is in $H^2(X;Z)$. As explained in Section 5h, the cup product pairing of the first Chern class of $L^q$ with this Euler class can be defined in a canonical way using the harmonic form $-iF_{\hat{A}}$; and when defined that way, this pairing can be written as $\mathrm{q}n + \mathfrak{e}_{\mathrm{q}}$ with $n$ being a q-independent integer and with the norm of $\mathfrak{e}_{\mathrm{q}}$ having a q-independent upper bound.

**Theorem 1.1**: *Supposing that the self-dual part of* $-iF_A$ *is nowhere zero, then the sequence of spectral flows parametrized by* q *for the* $|\omega|$*-divergent sequence described above is bounded if and only if the integer* $n$ *is zero. Supposing that the self-dual part of* $-iF_A$ *has a non-empty zero locus which is cut out transversally, then this same assertion about the spectral flow (it is bounded as a function of* q *if and only if* $n$ *is zero) still holds if the parameter* $m$ *that appears in (1.2) is not too large, less than* $\frac{1}{\kappa}$ *with* $\kappa > 1$ *being independent of* q.

This theorem is an instance of what is said in Proposition 5.15. What follows momentarily are some illustrative instances of the theorem. More details about some of these examples are in the body of this paper and at the end, in Sections 5i and 5j.

What follows directly is a remark about the upper bound for $m$ in the theorem: If $m$ is large and if the self-dual part if $-iF_{\hat{A}}$ has a non-empty zero locus, then there can be very many zero or nearly zero eigenvalues of the linearized equations (their eigenvectors concentrate near the zero locus); and the number will increase with q. If $m$ is small, then there are at most $\mathcal{O}(1)$ eigenvalues of the linearized equations that are nearly zero.

The first set of examples has X being isometric to the product of a compact oriented 3-manifold and the circle. The line bundle L should be the pull-back via the projection map to the 3-manifold, and the connection Â should be pulled back by this same map. The number $n$ is always zero in this case. (The relevant differential operator is invariant with respect to rotations of the circle factor; and because of this, much of the analysis in Sections 3-5 isn't needed to prove that the spectral flow is bounded.)

The second set of examples has X being a Kähler manifold. If the first Chern class of the line bundle L in the construction is a class of type 1-1 (so L has a holomorphic structure), then $n$ need not be zero. For instance, if the first Chern class of the Kähler manifold has positive self-cup product, then $n$ is not zero when the line bundle L for the construction has type (1,1) first Chern class. However, if the Kähler manifold's first Chern class is not proportional to the class of the Kähler symplectic form, then there are other complex line bundles for the construction whose first Chern class has some type (2,0) part with $n$ being zero.

The third set of examples has X being a symplectic 4-manifold. If the first Chern class of the associated canonical line bundle for the symplectic structure (a complex line bundle) has non-positive self-cup product and if the self-dual and anti-self dual second Betti numbers of X are at least two, then there are complex line bundles for the construction with $n$ being zero and others with $n$ not being zero. On the other hand if the associated first Chern class for the symplectic structure has positive square, then the author doesn't know a universal algorithm for finding line bundles for the construction with the number $n$ being zero. More is said about this case in Section 5j. (The construction that works in the Kähler case doesn't do the job in the symplectic case.)

**c) A look ahead for what is to come**

Section 2 of this paper describes in more detail how to obtain reducible, $|\omega|$-divergent sequences of solutions to (1.2) when $m \neq 0$. These will be the focus for the spectral flow analysis that occupies the rest of the paper paper. (Most of Section 2 was known to Edward Witten.)

One point to note is that the $|\omega|$-divergent sequences of solutions that are constructed here do not lead to the more subtle sort of $\mathbb{Z}/2$ harmonic, self-dual 2-form limits that are characterized in part by a non-trivial real line bundle defined only on the complement of a non-empty, (Hausdorff) dimension 2 subset in X (see [T1], [T2]). (The $\mathbb{Z}/2$ harmonic terminology here refers to a triple $(Z, \mathcal{I}, \varpi)$ with Z being a closed set in X with Hausdorff dimension 2, with $\mathcal{I}$ being a real line bundle defined in the complement of Z with no extension over Z, and with $\varpi$ being an $\mathcal{I}$-valued, self-dual harmonic form on X–Z whose norm extends as a Hölder continuous function over Z that is zero on Z.) I hope at some point to extend the spectral flow analysis from the subsequent sections of this paper to these sorts of limits.

Section 3 of this paper begins the study of the spectral flow for the self-adjoint operator from (1.4). Section 3 is a warm-up of sorts that analyzes the spectral flow for the reducible solutions that are constructed in Section 2 in the case where the form $\omega$ is nowhere zero. Section 3d discusses some of the curious phenomena in the case when X is a Kähler complex manifold.

Section 4 begins the analysis for the case when $\omega$ does have a non-empty zero locus under the assumption that the vanishing of $\omega$ along this locus is transverse in a suitable sense (because $\omega$ is reducible, the locus is then a smooth, 1-dimensional submanifold). The analysis in Section 4 is completed in Section 5. The results are summarized in Proposition 5.15 in Section 5h, and Sections 5i and 5j have remarks regarding this proposition that are relevant to some of the specific cases from Section 2. The analysis in Sections 3-5 introduces novel localization and excision techniques to obtain spectral flow bounds (novel to the author anyway) and these could be of interest in their own right.

To tie up a loose end: As noted above, the equations in (1.2) are the formal equations for the critical point of a functional on the space of pairs $(A, \omega)$ with A denoting a connection on the principal bundle P and with $\omega$ denoting a section of $\Lambda^+ \otimes \mathrm{ad}(P)$. Here is the functional:

$$\int_X \left(\mathrm{trace}(\omega\wedge F_A - \tfrac{1}{6}\omega\wedge[\omega;\omega]) - \tfrac{1}{2} m|\omega|^2\right)$$

(1.5)

(What is denoted by trace(·) is the trace homomorphism from ad(P) ⊗ ad(P) to $\mathbb{R}$ that is induced from an ad-invariant trace on the Lie algebra of SU(2).)

**d) Notation and conventions**

To be clear about the notation in (1.1)-(1.5) and used subsequently: The given 4-manifold, X, is assumed to be oriented and to come with a Riemannian metric. As alluded to above, the metric defines an orthogonal splitting of the bundle of 2-forms, $\wedge^2 T^*X$ as the direct sum of two 3-dimensional vector, $\Lambda^+ \oplus \Lambda^-$ with it understood that $\Lambda^+$ is the +1 eigenspace for the metric's Hodge star operator. Any given 2-form $\varsigma$ will be written with respect to this splitting as $\varsigma^+ + \varsigma^-$. This same ± superscript notation will be used for 2-forms with values in auxilliary vector bundles. Now, supposing that an orthonormal frame for $\Lambda^+$ has been chosen near some given point, then I will use $\{\varsigma^+_a\}_{a=1,2,3}$ for the corresponding components of $\varsigma^+$. By way of a prime example, if (A, ω) is a given pair of connection on a principal SU(2) or SO(3) bundle over X and ω is a section of $\Lambda^+$ with values in the associated lie algebra bundle, then the 3 components of the top equation in (1.2) when written near a point in X using a local, oriented orthonormal frame for $\Lambda^+$ are as follows:

$$F_A{}^+{}_a - \tfrac{1}{2\sqrt{2}}\varepsilon^{abc}[\omega_b, \omega_c] - m\,\omega_a = 0$$

(1.6)

where [ · , ·] denotes the commutator on the Lie algebra of SU(2), and where $\{\varepsilon^{abc}\}_{a.b.c\in\{1.2.3\}}$ are the components of the antisymmetric 3-tensor on $\mathbb{R}^3$ with $\varepsilon_{123} = +1$. (The indices b and c in (1.6) and below in (1.7) are implicitly summed over their allowed set of values.) Lie algebra commutation and the same completely antisymmetric tensor are also used to define the endomorphism $\varsigma \to [\omega;\varsigma]$ in the top bullet of (1.4):

$$[\omega;\varsigma]_a \equiv \tfrac{1}{\sqrt{2}}\varepsilon_{abc}[\omega_b, \varsigma_c]\ .$$

(1.7)

A convention adhered for the most part in what follows will be to write out equations using index notation (as in (1.6) and (1.7)) with respect to some implicit, locally define, oriented orthonormal frame for the tangent bundle of X. Such a frame induces corresponding oriented, orthonormal frames for T*X and $\Lambda^+$ and $\Lambda^-$ which are then used when needed. This is done because anyone checking a calculation will be using a locally defined frame anyway. So this convention cuts to the chase, so-to-speak. When equations are written using local orthormal frames, the Einstein summation convention is then used which is that repeated indices for tensors are implicitly summed.

Some Lie algebra conventions: The Lie algebra of SU(2) and SO(3) is identified with the vector space over $\mathbb{R}$ of $2\times 2$, anti-Hermitian matrices with zero trace. This Lie algebra/vector space is denoted by $\mathfrak{su}(2)$ in what follows. This vector space comes with the SU(2)-invariant inner product that is defined by the rule whereby the inner product of matrices $\sigma$ and $\sigma'$ is $-\frac{1}{2}\,\mathrm{trace}(\sigma\sigma')$. (This is a slightly non-standard view of the Lie algebra of SO(3). By taking this view, the equations in (1.2) can be simultaneously considered for SU(2) and SO(3) bundles.) Supposing that P denotes a principal SU(2) or SO(3) bundle over X, then the associated bundle with fiber $\mathfrak{su}(2)$ is denoted by ad(P). The inner product on $\mathfrak{su}(2)$ is used implicitly to define a fiber inner product on this ad(P) bundle. The latter inner product is denoted by $\langle\cdot,\cdot\rangle$. This same bracket notation is used also for the metric inner product on tensor bundles over X (such as the bundles of 1-forms and 2-forms on X) and for the induced inner product on the the tensor product of the bundle ad(P) with any of these bundles.

An analysis convention: Each appearance of '$c_0$' in what follows signifies a number that is always greater than 1 whose value can be assumed to increase between successive incarnations. But, unless stated to the contrary, the value of $c_0$ is independent of relevant large numbers or small numbers and any given pair $(A, \omega)$--it signifies a 'universal' constant with dependence only on the Riemannian metric and isomorphism class of a given principal bundle.

A different analysis convention: A smooth, non-increasing function on $\mathbb{R}$ to be denoted by $\chi$ is fixed here and forever with value is 1 on $(-\infty, \frac{1}{4}]$ and with value 0 on the ray $[\frac{3}{4}, \infty)$. This function $\chi$ will be used to construct various 'bump' functions with a priori control on the norms of their derivatives.

**e) Green's functions and Hardy's inequality**

Some facts about the Green's functions for certain Laplace-type operators are referred to in later sections. The operators in question act on functions and they have the form $d^\dagger d + r^2$ with d denoting the exterior derivative, $d^\dagger$ denoting its formal adjoint and r denoting a non-negative real number. In particular, the Green's function for this operator with pole at a given point $p \in X$ (denoted $G_p$) obeys bounds of the following sort:

- $c_0^{-1}\frac{1}{\mathrm{dist}(\cdot\,,p)^2}e^{-c_0 r\,\mathrm{dist}(\cdot,p)} \le G_p(\cdot) \le c_0\frac{1}{\mathrm{dist}(\cdot\,,p)^2}e^{-r\,\mathrm{dist}(\cdot,p)/c_0}$.
- $|\nabla G_p(\cdot)| \le c_0\left(\frac{1}{\mathrm{dist}(\cdot\,,p)} + r\right)\frac{1}{\mathrm{dist}(\cdot\,,p)^2}e^{-r\,\mathrm{dist}(\cdot,p)/c_0}$.

(1.8)

This is proved using parametrix built from the Green's function of the Euclidean version of the operator. See for example the Appendix in [DHR] for how to do this when $r = 0$.

A function space inequality known as Hardy's inequality also plays a small role in what follows. This inequality says this: If $f$ is any given smooth function on X, and p is any given point in X, and B denotes a ball in X centered at p, then

$$\int_B \frac{1}{\mathrm{dist}(\cdot,p)^2} f^2 \le c_0 \int_B (|df|^2 + f^2) \ .$$

(1.9)

This inequality implies in particular that the integral over X of the product of $f^2$ with the Green's function $G_p$ for $d^\dagger d + r^2$ is bounded by $c_0$ times the integral over X of $|df|^2 + f^2$. See [HLP]. A short proof in a related context is in [T3].

## 2. Some solutions to the $m \ne 0$ equations

The purpose of this subsection is to describe some $|\omega|$-diverging sequences of solutions to (1.2). The sequences described below will be the focus of the spectral flow computations in Section 3-5. (The very recent paper [Ch] studies the non-compactness of the moduli space of solutions to (1.1) on complex, Kähler 4-manifolds.)

### a) Reducible solutions

For the purpose of this subsection, X denotes a compact, oriented Riemannian 4-manifold. In all of the constructions that follow, the manifold X is assumed to have indefinite cup product pairing on $H^2(X;\mathbb{Z})$.

BUNDLES WITH ZERO FIRST PONTRJAGIN CLASS:

What follows is a more detailed account of the general construction of $|\omega|$–diverging sequences of solutions to (1.2) that was presented at the start of Section 1c. These are all solutions on principle SO(3) or SU(2) bundles $P \to X$ whose ad(P) bundle has zero first Pontrjagin class.

Recall now that the construction started with the choice of a non-zero class in $H^2(X;\mathbb{Z})$ whose cup product with itself is zero. Such a class exists if X has an indefinite intersection pairing on its second homology. Hodge theory will find a harmonic form (to be denoted by $\varsigma$) that represents this class. Since this 2-form has integer periods on 2-cycles in X, it can be viewed as representing the first Chern class of a complex line bundle over X to be denoted by L. In particular, there is a Hermitian metric on this line bundle with a unitary connection whose curvature 2-form is $2\pi i\, \varsigma$. Denote this connection by $\hat{A}$. The line bundle L with its Hermitian metric can be used to construct a principal SO(3) bundle whose associated Lie algebra bundle is isomorphic over $\mathbb{R}$ to the underlying real bundle of $L \oplus \underline{\mathbb{R}}$. (Use the metric to construct the unit circle bundle $S \subset L$, the set of elements with norm 1. This is a principal SO(2) bundle over X. The corresponding SO(3) bundle is $S \times_{SO(2)} SO(3)$.) Meanwhile, the connection $\hat{A}$ induces a connection on this principal bundle SO(3) bundle which will be denoted below by $A_1$. To say more about this connection: The $\underline{\mathbb{R}}$-factor in the Lie algebra bundle has an $A_1$- covariantly constant section with norm 1 (to be denoted by $\sigma$); and $A_1$ induces the connection $\hat{A}$ on the L factor. An important point here is that the curvature of the connection $A_1$ can be written as $\pi\varsigma\sigma$.

By the same token, if q is any given integer, then there is a principal SO(3) bundle over X to be denoted by $P_q$; and also a connection, $A_q$, whose curvature 2-form can be written as $\pi q \varsigma \sigma_q$ with $\sigma_q$ being an $A_q$-covariantly constant section of the associated Lie algebra bundle with norm 1.

To continue with the construction: Fix the positive number $m$ and then choose a second positive number, $r$, so that $rm$ is $\pi$ times an integer. A solution to (1.2) is given by the pair ($A = A_{q=rm/\pi}$, $\omega = r\varsigma^+\sigma_{q=rm}$) where, $\varsigma^+$ denotes the self-dual part of the harmonic 2-form $\varsigma$. An $|\omega|$-diverging sequence of solutions to (1.2) is given by choosing an unbounded sequence of integers for $q = rm/\pi$. Note that if q is even (or if the original cohomology class is even), then all of the principal bundles that appear in this sequence are isomorphic to the product bundle. Otherwise, there is a subsequence with all of the principal bundles pair-wise isomorphic SO(3) principal bundles.

NON-ZERO FIRST PONTRJAGIN CLASS WHEN X IS NOT SPIN:

The bundles constructed above all have zero first Pontrjagin class. What follows next describes $|\omega|$–divergent sequences on bundles with non-zero first Pontrjagin number. The target Pontrjagin number is denoted below by k. Assume first that the cup product pairing on $H^2(X;\mathbb{Z})$ is odd (which means that there is a class whose self-cup product is an odd integer). Assume that $\text{rank}(H^2(X;\mathbb{Z})) \geq 3$. Because the cup product pairing on $H^2(X;\mathbb{Z})$ is diagonalizable over $\mathbb{Z}$ (see [D1]), there are classes $x_1$, $x_2$ and $x_3$ in $H^2(X;\mathbb{Z})$ which are pair-wise cup-product orthogonal and whose self-cup products have absolute value 1 with that of $x_3$ being opposite in sign to those of $x_1$ and $x_2$. Letting $\cdot$ denote the cup product pairing, set $\varepsilon$ to denote the sign of $x_1 \cdot x_1$ and $x_2 \cdot x_2$. Fix integers p and q and note that

$$t \equiv q x_1 + (2p+1) x_2 + (2p-1) x_3 \tag{2.1}$$

has self cup product equal to $(q^2 + 4p)\varepsilon$. As a consequence: If $k \in 4\mathbb{Z}$, then $t \cdot t$ will equal k if q is even and if p is chosen to be $\frac{1}{4}(-q^2 + \varepsilon k)$. If $\varepsilon k = 1 \pmod 4$, then $t \cdot t$ will again equal k if q is odd and if p is again $\frac{1}{4}(-q^2 + \varepsilon k)$. I will consider only these values for k. (Values of k congruent to 2 mod(4) and with $\varepsilon k$ congruent to 3 (mod 4) can be obtained using an analogous construction when the rank of $H^2$ is greater than 3. I will leave these cases for the reader.) The important point in what follows is that with k being fixed, then the integer q and can be chosen to be greater than any given number.

Now having chosen q and then p, there is a harmonic 2-form on X to be denoted by $\varsigma_q$ whose de Rham cohomology class is the class $t$. As before, there is a principal

SO(3) bundle with connection, $A_q$, whose curvature 2-form is $\pi\, \varsigma_q \sigma_q$ with $\sigma_q$ being an $A_q$-covariantly constant section of the associated Lie algebra bundle. A corresponding solution to (1.2) is then ($A = A_q$, $\omega = \pi\frac{1}{m}\, \varsigma_q^{+}\sigma_q$). Taking q ever larger gives an $|\omega|$-divergent sequence of solutions to (1.2) on pairwise isomorphic principal SO(3) bundles with first Pontrjagin class equal to the chosen integer k.

NON-ZERO FIRST PONTRJAGIN CLASS WHEN X IS A SPIN MANIFOLD:

Now suppose that X is a spin manifold with indefinite cup product pairing on $H^2(X;\mathbb{Z})$. In this case, there are classes $x_1$, $x_2$ and $x_3$, $x_4$ each with self-cup product equal to zero, and with the only pairwise non-zero cup products between these classes being $x_1 \cdot x_2 = 1$ and $x_3 \cdot x_4 = 1$ (that this is so follows from Donaldson's celebrated theorem [D2] and from the classification of indefinite, even unimodular forms). Choose integers a, b, c and d and define the class

$$t = a x_1 + b x_2 + c x_3 + d x_4. \tag{2.2}$$

This class has self cup product equal to 2(ab + cd). Any even choice of k can be represented by infinitely many choices of the four integers (a, b, c, d). Having chosen an even integer k, fix a 4-tuple of integers $Q \equiv (a, b, c, d)$ as desired subject to the constraint that 2(ab + cd) = k. Let $\varsigma_Q$ now denote the harmonic 2-form on X that represents the class $t$ depicted by (2.2) in the de Rham cohomology of X. As in the previous cases, there is a principal SO(3) bundle with connection, $A_Q$, whose curvature 2-form is $\pi\, \varsigma_Q \sigma_Q$ with $\sigma_Q$ being an $A_Q$-covariantly constant section of the associated Lie algebra bundle. A corresponding solution to (1.2) is then ($A = A_Q$, $\omega = \pi\frac{1}{m}\, \varsigma_Q^{+}\sigma_Q$). Taking the integers that comprise Q ever larger gives an $|\omega|$ divergent sequence of solutions to (1.2) for principal SO(3) bundles with first Pontrjagin class equal to k. Since there are only finitely many of such SO(3) bundles up to isomorphism, there is a subsequence with all elements being solutions for pair-wise isomorphic SO(3) bundles. (If all the elements in Q are even, then these solutions can be viewed as solutions for pairwise isomorphic principal SU(2) bundles.)

**b) Some less reducible solutions**

The assumption here is that the Riemannian 4-manifold X has non-zero $\mathbb{Z}/2\mathbb{Z}$ first homology. By virtue of this, there is a real line bundle over X that is not isomorphic to the product real line bundle. Denote this bundle by $\mathcal{I}$. Having put a metric on this bundle (which will henceforth be assumed), the set of points in $\mathcal{I}$ with norm 1 define a 2-fold cover of X. Denote this cover by W. The cover W has a fiber preserving free involution which is given by multiplication on $\mathcal{I}$ by -1. This involution is denoted by $\iota$ in what follows. I also assume that there is a non-zero class $x \in H^2(W;\mathbb{Z})$ which is odd under this action (which is to say that $\iota^* x = -x$) and which has square zero (the evaluation of the fundamental class of W on the self cup product $x \cup x$ is zero.)

A simple example: Take $X = (S^1 \times S^3)\#Z$ where Z here is a 4-manifold with indefinite intersection pairing on $H^2(Z;\mathbb{Q})$. (Take $Z = S^2 \times S^2$ for instance.) Then take $\mathcal{I}$ to be the pull-back of the Möbius line bundle over $S^1$ via a map from X to $S^1$ that maps the first homology of the Z summand to zero and identifies the first homology of the $S^1 \times S^3$ summand with $H_1(S^1;\mathbb{Z})$. In this case, W is diffeomorphic to $(S^1 \times S^3)\#_2 Z$; and, moreover, the intersection form (over $\mathbb{Q}$) on the anti-invariant cohomology of W (which has the same dimension as $H^2(Z;\mathbb{Q})$) is isomorphic to the intersection form on $H^2(Z;\mathbb{Q})$. Since Z is assumed to have indefinite intersection pairing there are integer classes in the anti-invariant, 2-dimensional cohomology of W with square zero.

The plan for what follows is to repeat the constructions in Part a) on W but in a suitably $\iota$-equivariant fashion so as to obtain an interesting $|\omega|$-diverging sequence of solutions to (1.2) on X.

As a first step towards this end: Give W the metric that is pulled up from X via the covering projection. Having done that, then use Hodge theory to find a harmonic 2-form on W that represents $2\pi x$. Denote this 2-form by $\varsigma$. Because $\varsigma$ is harmonic and because $\iota$ acts as an isometry, Hodge theory guarantees $\iota^*\varsigma = -\varsigma$.

Now, having fixed $m$ for use in (1.2), chose any positive number $r$ so that $rm$ is $2\pi$ times an integer. Having chosen $r$ subject to this constraint, then the 2-form $\frac{1}{2\pi} rm\varsigma$ represents the integral class $\frac{1}{2\pi} rmx$ in $H^2(W;\mathbb{Z})$. This being the case, $rm\varsigma$ is $\sqrt{-1}$ times the curvature 2-form of a unitary connection on Hermitian complex line bundle over W with first Chern class equal to this class $\frac{1}{2\pi} rm\, x$. Denote this the complex line bundle by L. The upcoming Lemma 2.1 says more about L and a specific connection on L whose curvature 2-form is the 2-form $-i\, rm\varsigma$.

To set the background for Lemma 2.1, remember that complex conjugation defines an $\mathbb{R}$-linear bundle isomorphism from L to a line bundle denoted by $\overline{L}$ and vice-versa. Also: The first Chern class $\overline{L}$ is the class $-\frac{1}{2\pi} rm\, x$.

With regards to connections on L and $\overline{L}$, any unitary connection on L pulls-back via the complex-conjugate isomorphism to a connection on $\overline{L}$. If $A_L$ is a unitary connection on L, then the corresponding connection on $\overline{L}$ will be denoted here by $\overline{A}_L$. With regards to their curvature 2-forms: If $A_L$ has curvature 2-form $-i\, rm\varsigma$, then $\overline{A}_L$ has curvature 2-form $i\, rm\varsigma$.

As noted above, the bundle $\overline{L}$ has first Chern class $-x$, but so does the bundle $\iota^* L$. Therefore, these bundles are isomorphic over $\mathbb{C}$. In addition, because $-\varsigma$ is $\iota^*\varsigma$, if $A_L$ is a unitary connection on L with curvature 2-form $-i\, rm\varsigma$, then a $\mathbb{C}$-linear isomorphism between $\iota^*L$ and $\overline{L}$ identifies the pull-back connection $\iota^*A_L$ with a connection on $\overline{L}$ that has the same curvature 2-form as does $\overline{A}_L$. Lemma 2.1 makes a formal assertion to the effect that $A_L$ and this isomorphism can be chosen so that $\iota^*A_L$ *is* $\overline{A}_L$.

**Lemma 2.1**: *There exists a positive integer (to be denoted by* n*) such that if* $\frac{1}{2\pi} rm$ *is divisible by* n*, then the involution* $\iota$ *on* W *is covered by a* $\mathbb{C}$*-linear map* $\iota^*$: $L \to \overline{L}$ *with complex conjugate* $\overline{\iota}^*: \overline{L} \to L$ *such that* $\overline{\iota}^*\iota^*$ *is the identity bundle automorphism of* L*. Moreover, there is a unitary connection on* L*, to be denoted by* $A_L$*, such that* $\iota^*A_L = \overline{A}_L$.

This lemma is proved momentarily. Accept it for now to continue with the construction of the desired ω-divergent sequence of solutions to (1.2).

Having this lemma in hand, now introduce the $\mathbb{C}^2$ bundle $\overline{L} \oplus L$. This bundle has its unitary connection which is defined by the requirement that the direct sum splitting be covariantly constant and that parallel transport of the respective L and $\overline{L}$ factors is the same as parallel transport by $A_L$ and by $\overline{A}_L$. Denote this connection by $\hat{A}$.

By virtue of what is said in Lemma 2.1, the map $\iota$ is covered by a $\mathbb{C}$-linear map

$$\begin{pmatrix} \overline{\iota}^* & 0 \\ 0 & \iota^* \end{pmatrix} : \overline{L} \oplus L \to L \oplus \overline{L} \tag{2.3}$$

This map is denoted by $\hat{\iota}$. Now let $\tau$ denote the $\mathbb{C}^2$–bundle homomorphism from $L \oplus \overline{L}$ to $\overline{L} \oplus L$ that interchanges the two summands:

$$\tau \equiv \begin{pmatrix} 0 & 1 \\ 1 & 0 \end{pmatrix} . \tag{2.4}$$

This definition is such that the composition $\tau\hat{\iota}$ is a lift of $\iota$ to a $\mathbb{C}$-linear map from $\overline{L} \oplus L$ to itself. Moreover (again by virtue of Lemma 2.1), this lift $\tau\iota^*$ is such that $(\tau\hat{\iota})^*\hat{A} = \hat{A}$. For future reference, the curvature 2-form of $\hat{A}$ can be written as $rm\varsigma\sigma$ with $\sigma$ denoting the endomorphism of $\overline{L} \oplus L$ depicted below:

$$\sigma \equiv \begin{pmatrix} i & 0 \\ 0 & -i \end{pmatrix} . \tag{2.5}$$

This endomorphism $\sigma$ is $\hat{A}$-covariantly constant as it preserves the direct sum splitting of $\overline{L} \oplus L$ and is constant on each summand. But note also that $\tau\hat{\iota}(\sigma(\cdot)) = -\sigma(\tau\hat{\iota}(\cdot))$; this is to say that the action by $\sigma$ anti-commutes with the map $\tau\hat{\iota}$.

With the preceding understood, define a rank 2, complex vector bundle $E \to X$ to be the set of equivalence classes in $\bar{L} \oplus L$ with the equivalence relation being this: Write a point in $\bar{L} \oplus L$ as a pair $(p, e)$ with $p \in W$ and $e \in (\bar{L} \oplus L)|_p$. This point is identified by the equivalence relation with the point that is defined by the pair $(\iota(p), \tau\hat{\iota}(e))$. The definition is such that connection $\hat{A}$ descends to X where it defines a special unitary connection on the bundle E (because $\tau\hat{\iota}$ pulls back $\hat{A}$ to itself). Denote this connection by A; it is the connection of interest for (1.2).

It remains now to define the End(E) valued self-dual 2-form $\omega$ for use in (1.2). To this end, let $\varsigma^+$ denote the self-dual part of $\varsigma$. Then, define a self-dual 2-form $\hat{\omega}$ on W with values in the bundle of anti-Hermitian, traceless endomorphisms of $(\bar{L} \oplus L)$ to be

$$\hat{\omega} \equiv r\varsigma^+\sigma \; . \tag{2.6}$$

Because $\varsigma$ is harmonic, the self-dual 2-form $\varsigma^+$ is closed, and as a consequence $d_{\hat{A}}\hat{\omega} = 0$ (remember that $\sigma$ is $\hat{A}$-covariantly constant). Also: The endomorphism valued 2-form $\hat{\omega}$ on W pulls-back as itself via $\tau\iota^*$ because $\iota^*\varsigma = -\varsigma$ and $\tau\sigma\tau = -\sigma$. As a consequence, $\hat{\omega}$ descends to define a self-dual 2-form on X with values in the bundle of anti-Hermitian, trace zero endomorphisms of E. The latter is $\omega$.

To finish the construction: Because $d_{\hat{A}}\hat{\omega} = 0$, so it is that $d_A\omega = 0$. And, also note that $m\omega$ is the self-dual part of $F_{\hat{A}}$. Thus, the pair $(A, \omega)$ obeys (1.2) with the given value of $m$. Taking $r$ ever larger (subject to the constraint that $rm \in \mathbb{Z}$) produces a $|\omega|$-diverging sequence of solutions to (1.2); and these are on pair-wise isomorphic $\mathbb{C}^2$ bundles because each instance of the bundle E has zero second Chern class (this is because the cup product of the class $x$ with itself is zero).

With regards to the reducibility of these solutions: Any given $r$ version of $(A, \omega)$ is irreducible in the sense that there are no A-covariantly constant sections of the bundle of skew-Hermitian endomorphisms of E. To explain why this is, note first that there are $\hat{A}$-covariantly constant sections on W of the skew-Hermitian endomorphisms of $\bar{L} \oplus L$, these being the constant multiples of the endomorphism $\sigma$ that is depicted in (2.5). Even so, the endomorphism $\sigma$ does descend to X as an endomorphism of E because it *anti-commutes* with $\tau$. Instead, $\sigma$ defines on X an A-covariantly constant homorphism from the real line bundle $\mathcal{I}$ into the bundle of skew-Hermitian endomorphisms of E.

As a final remark: The preceding construction can be modified in a straight forward way to construct a sequence of $|\omega|$-diverging solutions to (1.2) on X for pairwise isomorphic principal bundles with non-zero first Pontrjagin class if there are suitable cohomology classes besides in $H^2(W; \mathbb{Z})$ that change sign when pulled back by $\iota$.

***Proof of Lemma 2.1***: Let $L_1$ denote the bundle for the case when $rm = 2\pi$. It is sufficient to prove that there exists an integer n such that the conclusions of the lemma hold when $rm = 2\pi n$. To this end, let $\iota(L_1)$ denote the pull-back of the line bundle $L_1$ by the diffeomorphism $\iota$. This bundle is isomorphic as a complex line bundle to $\bar{L}_1$, and with this understood, fix an isomorphism $w: \bar{L}_1 \to \iota(L_1)$ and let $w\iota(A)$ denote the image of A via this isomorphism. The connection $w\iota(A)$ can be written as $w\iota(A) = \bar{A} + \nu$ with $\nu$ being a closed, $i\mathbb{R}$ valued 1-form. Moreover, $\nu$ has the following property: Let $\gamma$ denote a loop in W which is such that $\iota$ maps $\gamma$ to itself reversing orientation. Then the path integral of $\nu$ on $\gamma$ is an integer multiple of $2\pi i$.

As for loops $\gamma$ such that $\iota(\gamma) = \gamma$ preserving orientation: Let $u$ denote a closed, $i\mathbb{R}$-valued 1-form and let $A' = A + u$. Then $w\iota(A') = \bar{A}' + u + \iota^*(u) + \nu$. Therefore, A can (and should) be chosen so that the corresponding version of $\nu$ obeys $\nu = -\iota^*(\nu)$. This implies that the path integral of this new version of $\nu$ on $\gamma$ is zero.

By way of a summary: One can assume that $w\iota(A) = \bar{A} + \nu$ with $\nu$ having path integral in $2\pi i\mathbb{Z}$ on any loop in W mapped to itself by $\gamma$. With this point understood, note that the action of $\iota$ on $H_1(W; \mathbb{Q})$ can be diagonalized in the sense that there is a basis of classes that are represented by $\iota$-invariant loops in X. Let $\{\gamma_1, \ldots, \gamma_k\}$ denote such a basis of loops. Any class in $H_1(W; \mathbb{Z})$/torsion can be represented as a linear combination of these loops with rational coefficients. Let N denote the least common multiple of the denominators of these coefficients. The path integral of $N\nu$ along any loop in W (whether $\iota$-invariant or not) will be in $2\pi i\mathbb{Z}$. This implies that $N\nu$ can be written as $u^{-1}du$ with u mapping W to $S^1$.

With that point understood, note that the connection A induces a connection $A_N$ on the line bundle $L = (L_1)^N$ and for that connection, $w\iota(A_N) = \bar{A}_N + u^{-1}du$. Let $(u^{-1}w)\iota(A_N)$ denote the result of mapping $A_N$ first by $\iota$ and then by w and, at the end, by $u^{-1}$ (acting as an automorphism of $\bar{L}$). Then $(u^{-1}w)\iota(A) = \bar{A}$. Thus, there exists and isomorphism from $\bar{L}$ to $\iota(L)$ (denoting this one by y) obeying $y\iota(A_N) = \bar{A}_N$. The complex conjugate isomorphism, $\bar{y}\iota$ then obeys $\bar{y}(y \circ \iota) A_N = A_N$ which implies that $(y \circ \iota)$ is equal to $\eta\,\bar{y}^{-1}$ with $\eta$ denoting a constant, complex number with norm 1. Thus, $y = \eta(\bar{y} \circ \iota)^{-1}$ and hence by taking conjugates, $\bar{y} = \bar{\eta}(y \circ \iota)^{-1}$. Compose the latter with $\iota$ to deduce that $\bar{y} \circ \iota = \bar{\eta}\, y^{-1}$. Taking inverses leads to the identity $y = \eta/\bar{\eta}\; y$ which can hold only in the event that $\eta = \pm 1$. If $\eta = 1$, then the conclusions of the lemma follow.

To finish the proof, it is necessary to show that $\eta = 1$ if L has even first Chern class. To see this, note that if the Chern class of L is even, then L can be written as $\mathcal{L}^2$ and the argument in the preceding paragraph with $\mathcal{L}$ replaced by L produces an isomorphism $y_{\mathcal{L}}$ and a corresponding $\eta_{\mathcal{L}} \in \{\pm 1\}$. The isomorphism $y_{\mathcal{L}}$ induces the isomorphism y for L with corresponding $\eta$ being $\eta_{\mathcal{L}}^2$ which is always equal to 1.

## 3. Bounding the norm of the spectral flow when ω is nowhere zero

This section and the subsequent ones consider the spectral flow for the operator that is depicted in (1.4) when (A, ω) is a pair that is described in Section 2. The first subsection here (Section 3a) gives a precise definition of spectral flow. The remaining subsections of this section analyze the spectral flow when ω is nowhere zero. Sections 4 and 5 study the spectral flow when ω has zeros, albeit only in the case when the zero locus is cut out transverally.

### a) The definition of 'spectral flow'

By way of a reminder about what is meant by 'spectral flow': The spectral flow for a continuous, one parameter family of self adjoint, Fredholm operators (parametrized by $t \in [0, 1]$ via a map $t \to \mathcal{L}_t$) is simplest to describe when the endpoint operators $\mathcal{L}_0$ and $\mathcal{L}_1$ both have trivial kernel. In this case, the spectral flow can be defined by first perturbing the family without changing the end members so that the following generic conditions hold:

- *There are only finitely many parameter values in* $[0,1]$ *where the corresponding* $\mathcal{L}_t$ *has a kernel.*
- *The kernel dimension at any such parameter value is* 1.

(3.1)

Having made this perturbation, then the spectral flow is equal to the difference between the number of parameter values where an eigenvalue crosses zero from below as t increases minus the number of parameter values where an eigenvalue crosses zero from above as t increases. (One can use techniques from [K] to prove that the spectral flow is independent of the precise nature of the perturbations that generate the conditions in (3.1) subject to some mild constraints.)

When the kernel of either or both of $\mathcal{L}_0$ and $\mathcal{L}_1$ is non-trivial, then the preceding definition is ambigous, but the absolute value of the ambiguity is, in any event, no greater than the sum of the dimensions of the respective kernels. Indeed, one can perturb the path ever so slightly so that the new end members have trivial kernel; and then define the spectral flow for this perturbed path. The result depends on the perturbed endpoints, but only up to an ambiguity of the indicated size.

The preceding definition of spectral flow (modulo kernel ambiguity at the endpoints) is used in the following manner to bound the absolute value of the spectral flow in a manner that is independent of the perturbation: First perturb ever so slightly if necessary so that the endpoint operators have trivial kernel; then bound the absolute value of the spectral flow for the perturbed path; and then add the dimensions of the kernels of the unperturbed operators.

In the case at hand, there are paths of self-adjoint, Fredholm operators between any two versions of the operator depicted in (1.4). Indeed, if $(m_0, (A_0, \omega_0))$ and $(m_1, (A_1, \omega_1))$ are the relevant two versions of the data $m$ and $(A, \omega)$ for (1.4), then a path between their corresponding versions of (1.4) is obtained from the formula in (1.4) by using a continuous path $t \to (m_t, (A_t, \omega_t))$ for $t \in [0,1]$ of pairs to use for $(A, \omega)$. One can also consider a more generally family of perturbations by adding to this parametrized version of (1.4) any continously parameterized zero'th order operator (for example, a t-dependent family of endomorphisms of the bundle depicted in (1.3) can be added.) Even more general variations are allowed (for example, the Riemannian metric can be changed along the family.)

With regards to the gauge invariance of the spectral flow: Let $\mathfrak{g}$ denote for the moment an automorphism of the ambient principal bundle P. The spectral flow is zero between the respective versions of (1.4) defined by any given $(A, \omega)$ and it pull-back via $\mathfrak{g}$ (thus, $(A - (d_A\mathfrak{g})\mathfrak{g}^{-1}, \mathfrak{g}\,\omega\,\mathfrak{g}^{-1})$.) This is the case whether or not $\mathfrak{g}$ is homotopically trivial. This assertion of vanishing spectral flow between the operator defined by a pair $(A, \omega)$ and its pull-back via a bundle automorphism follows from Proposition 9.2 of [AS] by using techniques from [APS]. To elaborate just a little, fix a 1-parameter path $(A_t, \omega_t)$ parameterized by [0, 1] that is constant near $t = 0$ where it is $(A, \omega)$, and constant near $t = 1$ where it is the pair $(A - d_A\mathfrak{g}\,\mathfrak{g}^{-1}, \mathfrak{g}\omega\mathfrak{g}^{-1})$. This path can be used to construct a first order, elliptic operator on $S^1 \times X$ whose Fredholm index is the spectral flow along the path (this is the input from [APS].) Proposition 9.2 from [AS] asserts that said index is necessarily zero.

**b) The model spectral flow calculation**

This subsection works through the spectral flow bound for the case when X is the 4-torus $\times_4 S^1$ and with $(A, \omega)$ being an $m \neq 0$ solution to (1.2) of the sort that is described in Section 2a.

To set the stage (and review definitions): Identify $S^1$ with $\mathbb{R}/\mathbb{Z}$ and then fix affine coordinates $(x_1, x_2, x_3, x_4) \in \times_4(\mathbb{R}/\mathbb{Z})$ for X. Now let $\nu = dx_1 \wedge dx_2$. This is a closed 2-form with integer pairing with the cycles that represent the classes in $H_2(X; \mathbb{Z})$. Because of this, if q is any integer, then $q\,(dx_1 \wedge dx_2)$ is the first Chern class of a complex line bundle over X to be denoted by $L_q$. With $L_q$ in hand, let $E_q$ denote the $\mathbb{R}^3$ bundle $L_q \oplus \underline{\mathbb{R}}$. The direct sum depiction of $E_q$ with the complex orientation for $L_q$ gives the bundle $E_q$ an orientation. Meanwhile, a fixed Hermitian metric on $L_q$ and the Euclidean metric on $\mathbb{R}$ gives $E_q$ a fiber metric.

Fix a connection on $L_q$ (to be denoted by $\hat{A}_q$) whose curvature 2-form is $2\pi i\, q\,(dx_1 \wedge dx_2)$. With $\hat{A}_q$ in hand, the bundle $E_q$ has a unique orthogonal connection to be denoted by $A_q$ which is characterized by two requirements: The first is that the direct sum splitting of $E_q$ is $A_q$-covariantly constant; and the second is that parallel transport by $A_q$ on the $L_q$ summand of $E_q$ is the same as parallel transport by $\hat{A}_q$.

The bundle of anti-symmetric endomorphism of $E_q$ has an $A_q$-covariantly constant endomorphism which annihilates the $\underline{\mathbb{R}}$ summand of $E_q$ and acts as multiplication by i on the $L_q$ summand. This endomorphism is denoted $\sigma_q$. Set

$$\omega_q = \frac{\pi q}{m}(dx_1 \wedge dx_2 + dx_3 \wedge dx_4)\sigma_q. \tag{3.2}$$

This is an $A_q$-covariantly constant, self-dual 2-form on X with values in the bundle of skew-symmetric endomorphisms of $E_q$. The pair $(A_q, \omega_q)$ is a solution to (1.2) for the given value of $m$ because $\omega$ is $A_q$-covariantly constant and because the self-dual part of $A_q$'s curvature 2-form is $m\omega_q$.

The principal bundle of oriented, orthonormal frames in $E_q$ has zero for its second Pontrjagin class. If q is an even integer, then that principal bundle is isomorphic to the product SO(3) bundle; and if q is odd, then it is isomorphic to a fixed bundle with non-zero second Stieffel-Whitney class. Denote these bundles by $P_0$ and $P_1$ respectively. As a consequence, if q is even, then $(A_q, \omega_q)$ can be viewed as a pair of connection on $P_0$ and self-dual 2-form with values in $P_0$'s adjoint bundle (the bundle $P_0 \times_{SO(3)} \mathfrak{su}(2)$). By the same token, if q is odd, then $(A_q, \omega_q)$ can be viewed as a pair of connection on $P_1$ and self-dual 2-form with values in $P_1 \times_{SO(3)} \mathfrak{su}(2)$. This view point is taken in the spectral flow proposition that follows directly.

**Proposition 3.1**: *There is a* q*-independent bound on the absolute value of the spectral flow between the* $(A_q, \omega_q)$ *version of (1.4) and the version defined by* $(A_0, \omega_0)$ *or* $(A_1, \omega_1)$ *as the case may be.*

***Proof of Proposition 3.1***: The proof of has seven parts.

*Part 1*: This part of the proof rewrites the operator in (1.4) for use in the later parts of the proof. To start the rewriting: Let $\mathbb{V}$ denote the vector bundle $(\Lambda^+ \oplus \underline{\mathbb{R}}) \oplus T^*X$. In the case relevant to Lemma 4.1, the bundle $\mathbb{V}$ is canonically isometric to the product bundle with fiber $(\mathbb{R}^3 \oplus \mathbb{R}) \oplus \mathbb{R}^4$ which is to say $X \times \mathbb{R}^8$. The isomorphism is given by using the constant orthonormal basis

$$\tfrac{1}{\sqrt{2}}(dx^2 \wedge dx^3 + dx^1 \wedge dx^4),\ \tfrac{1}{\sqrt{2}}(dx^3 \wedge dx^1 + dx^2 \wedge dx^4)\ ,\ \tfrac{1}{\sqrt{2}}(dx^1 \wedge dx^2 + dx^3 \wedge dx^4) \tag{3.3}$$

for $\Lambda^+$ and the basis $\{dx^1, dx^2, dx^3, dx^4\}$ for T*X. Using this basis, the operator depicted in (1.4) for any given pair (A,ω) (which I will denote henceforth by $\mathcal{D}$) can be written schematically using 8 × 8 real matrices $\{\gamma_\alpha, \rho_k\}_{\alpha=1,...4; k=1,2,3}$ and an 8 × 8 real matrix Γ as follows:

$$\mathcal{D} = \gamma_\alpha \nabla_{A\alpha} + \rho_k[\tfrac{1}{\sqrt{2}}\, \omega_k, \cdot\,] - \tfrac{1}{2}\, m\Gamma - \tfrac{1}{2}\, m \tag{3.4}$$

with it understood that there is a summation over the α-indices (from 1 to 4) and a summation over the k-indices (from 1 to 3). The notation in (3.4) is as follows: First, $\{\nabla_{A\alpha}\}_{\alpha=1,2,3,4}$ signifies the directional covariant derivative defined the connection A along the vector field dual to $dx^\alpha$. Second $\omega_k$ are the components of ω when written using the basis in (3.3) for $\Lambda^+$ with k increasing from left to right. With regards to the matrices from the set consisting of the four γ's and the three ρ's and Γ: The γ's and ρ's are anti-symmetric whereas Γ is symmetric. With regards to the original depiction of $\mathbb{V}$ as $(\Lambda^+ \oplus \underline{\mathbb{R}}) \oplus T^*X$, the matrix Γ acts as +1 on the $(\Lambda^+ \oplus \underline{\mathbb{R}})$ summand and as -1 on the T*X summand. Meanwhile, all of the γ's interchange the summands whereas all of the ρ's preserve the summands. Finally, these eight matrices enjoy the anti-commutation relations in the first four bullets below and the commutation relation in the last bullet:

- $\gamma_\alpha\gamma_\beta + \gamma_\beta\gamma_\alpha = -2\mathbb{I}$ .
- $\rho_i\rho_k + \rho_k\rho_i = -2\mathbb{I}$ .
- $\rho_i\gamma_\alpha + \gamma_\alpha\rho_i = 0$ .
- $\Gamma\gamma_\alpha + \gamma_\alpha\Gamma = 0$ .
- $\Gamma\rho_k - \rho_k\Gamma = 0$ .

(3.5)

Here, $\mathbb{I}$ denotes the identity 8 × 8 matrix.

*Part 2*: It proves useful now to introduce a 1-parameter family of operators parameterized by $t \in [0, \infty)$ using the rule

$$\mathcal{D}_t \equiv \gamma_\alpha \nabla_{A\alpha} + \rho_k[\tfrac{1}{\sqrt{2}}\, \omega_k, \cdot\,] - \tfrac{1}{2}\, t\,\Gamma\ . \tag{3.6}$$

Thus $\mathcal{D} = \mathcal{D}_{t=m} - \frac{1}{2}\, m$. By virtue of the anti-commutation relations in (3.5),

$$\mathcal{D}_t^2 = (\gamma_\alpha \nabla_{A\alpha})^2 + \tfrac{1}{\sqrt{2}}\, \gamma_\alpha\rho_k[\nabla_{A\alpha}\omega_k, \cdot\,] + (\rho_k[\tfrac{1}{\sqrt{2}}\, \omega_k, \cdot\,] - \tfrac{1}{2}\, t\,\Gamma)^2\ . \tag{3.7}$$

In particular if ω is $\nabla_A$-covariantly constant (as will henceforth be assumed), then $\mathcal{D}_t^2$ is a sum of squares of self-adjoint operators:

$$\mathcal{D}_t^2 = (\gamma_\alpha\nabla_{A\alpha})^2 + (\rho_k[\tfrac{1}{\sqrt{2}}\,\omega_k,\cdot\,] - \tfrac{1}{2}\,t\,\Gamma)^2 \,. \tag{3.8}$$

Looking ahead, this implies that the kernel of $\mathcal{D}_t$ is annihilated by both the differential operator $\gamma_\alpha\nabla_{A\alpha}$ and the algebraic operator $\rho_k[\frac{1}{\sqrt{2}}\,\omega_k,\cdot\,] - \frac{1}{2}\,t\,\Gamma$.

*Part 3*: The norm of the spectral flow from the $(A_q, \omega_q)$ version of $\mathcal{D}$ to the $(A_0,\omega_0)$ or $(A_1,\omega_1)$ version will be bounded by using a piecewise smooth, continuous, 1-parameter deformation of $\mathcal{D}$ that moves it to an operator $D_\varepsilon$ that has the following form:

$$D_\varepsilon = \gamma_\alpha\nabla^0{}_\alpha - \varepsilon\tfrac{1}{2}\,\Gamma \tag{3.9}$$

where $\varepsilon$ is positive and otherwise can be taken arbitrarily small. Here, $\nabla^0$ is the covariant derivative defined by either $A_0$ or $A_1$ as the case may be. It follows from the $\omega_k \equiv 0$ version of (3.8) that $D_\varepsilon$ is invertible; its spectrum is empty in the interval $(-\frac{1}{2}\,\varepsilon, \frac{1}{2}\,\varepsilon)$. Bounding the norm of the spectral flow from the $(A_q,\omega_q)$ version of $\mathcal{D}$ to $D_\varepsilon$ is sufficient for the purposes at hand because the norm of the spectral flow from $D_\varepsilon$ to the $(A_0,\omega_0)$ or $(A_1,\omega_1)$ version of $\mathcal{D}$ has a q-independent upper bound.

The deformation of $\mathcal{D}$ to $D_\varepsilon$ proceeds in stages.

<u>STAGE 1</u>: Deform $\mathcal{D}$ along the 1-parameter family

$$s \to \gamma_\alpha\nabla_{A\alpha} + \rho_k[\tfrac{1}{\sqrt{2}}\,\omega_k,\cdot\,] - \tfrac{1}{2}\,m\Gamma - \tfrac{1}{2}\,(1-s)\,m \tag{3.10}$$

with s starting at 0 and ending at s = 1. The s = 1 version of the operator in (3.10) is the $t = m$ version of the operator $\mathcal{D}_t$ that is depicted in (3.6) whereas the s = 0 version is the $(A_q, \omega_q)$ version of $\mathcal{D}$.

<u>STAGE 2</u>: Deform the $t = m$ version of $\mathcal{D}_t$ by increasing t so that it is very much greater than $|\omega|$. Let T denote this very large version of t. By virtue of (3.8), the absolute value of the smallest eigenvalue of the operator $\mathcal{D}_T$ will be greater than $\frac{1}{2}\,T(1 - \mathcal{O}(|\omega|/T))$ which is greater than $\frac{1}{4}\,T$ when T is large.

<u>STAGE 3</u>: Deform the $t = T$ version of $\mathcal{D}_t$ by the 1-parameter family

$$s \to \mathcal{D}_{T,s} \equiv \gamma_\alpha\nabla_{A\alpha} + (1-s)\,\rho_k[\tfrac{1}{\sqrt{2}}\,\omega_k,\cdot\,] - \tfrac{1}{2}\,T\,\Gamma \tag{3.11}$$

with s starting from 0 and ending at s = 1.  Since T >> $|\omega|$, no member of this family has an eigenvalue in the range $[-\frac{1}{4}T, \frac{1}{4}T]$.  In particular, there is no spectral flow in this stage. The end member of this stage is the operator $\gamma_\alpha\nabla_{A\alpha} - \frac{1}{2}T\Gamma$.

STAGE 4:  Deform the operator $\gamma_\alpha\nabla_{A\alpha} - \frac{1}{2}T\Gamma$ via a 1-parameter family that moves the connection $A_q$ to $A_0$ or $A_1$.  It is also the case that no member of this family has an eigenvalue in the range $[-\frac{1}{4}T, \frac{1}{4}T]$.  (This follows from the $\omega \equiv 0$ version of (3.8).)  The end member of this family is the operator $D_{\varepsilon=T}$.

STAGE 5:  Decrease T to any positive value to obtain a family of $D_\varepsilon$ operators with $\varepsilon$ starting at T with no member having a zero eigenvalue.

As indicated in the decription above of this deformation process, the only stages where there is spectral flow are STAGES 1 and 2.

*Part 4*:  Consider now the case in Proposition 3.1 where $\omega_1 = \omega_2 = 0$ and $\omega_3 = r\sigma$ with $r$ being a real constant and with $\sigma$ being a $\nabla_A$-covariantly constant, norm 1 map from X to the Lie algebra $\mathfrak{su}(2)$.

Let $\mathfrak{su}_{\mathbb{C}}$ denote $\mathfrak{su}(2)\otimes_{\mathbb{R}}\mathbb{C}$, the complexified Lie algebra of SU(2).  The endomorphism $[\frac{i}{2}\sigma, \cdot\,]$ on $\mathfrak{su}_{\mathbb{C}}$ at any given point in X splits $\mathfrak{su}_{\mathbb{C}}$ into a direct sum of three complex lines which are written below as $L^+\oplus L^0 \oplus L^-$ with $L^+$ being the +1 eigenspace and $L^-$ being the -1 eigenspace; and with $L^0$ being the kernel (the eigenspace with eigenvalue 0).  In this regard:  The Hermitian conjugate of any element in $L^+$ is in $L^-$ and vice-versa.  Each of these lines defines a complex line bundle over X.  The bundle $L^+$ is isomorphic to $L_q^{\,2}$ whereas $L^-$ is isomorphic to $L_q^{-2}$.  Meanwhile, $L^0$ is isomorphic to the product bundle.

Because $\sigma$ is covariantly constant, this splitting is preserved by the covariant derivative $\nabla_A$ which implies that $\nabla_A$ acts on sections of $L^\pm$ and $L^0$.  On the latter bundle, $\nabla_A$ is the product flat connection.

As a consequence of what was just said, the operators $\gamma_\alpha\nabla_{A\alpha}$ and $\rho_3[\frac{1}{\sqrt{2}}\omega_3,\cdot\,]$ map any section of $\mathbb{R}^8\otimes L^+$ to a section of $\mathbb{R}^8\otimes L^+$; they map any section of $\mathbb{R}^8\otimes L^-$ to one of $\mathbb{R}^8\otimes L^-$; and they map any section of $\mathbb{R}^8\otimes L^0$ to a section of $\mathbb{R}^8\otimes L^0$.  Therefore, the same can be said of $\mathcal{D}_t$ and also $\mathcal{D}$.  And, this can also be said for the families of operators that are described in STAGES 1-3 from Part 3.

*Part 5*:  A digression is needed to say more about the endomorphism $\rho_3[\frac{1}{\sqrt{2}}\omega_3,\cdot\,]$ in order to analyze the STAGE 1 and STAGE 2 spectral flow.  To start the digression:  This endomorphism acts on $\mathbb{R}^8$-valued sections of $L^+$ as $-i\sqrt{2}\,x\rho_3$; it acts on $\mathbb{R}^8$-values sections of $L^-$ as $+i\sqrt{2}x\rho_3$; and it acts on $\mathbb{R}^8$-valued sections of $L^0$ as 0. Meanwhile, $i\rho_3$ has eigenvalue ±1 with each eigenvalue occuring with multiplicity two on each summand $\mathbb{R}^4 \otimes_{\mathbb{R}} \mathbb{C}$ of the complexification of $\mathbb{R}^8$ (the first summand corresponds to the complexification $\Lambda^+ \oplus \underline{\mathbb{R}}$ and the second to T*X). Indeed, this follows from the depiction below of $\rho_3$ (entries not depicted are zero):

$$
\begin{pmatrix}
0 & -1 & 0 & 0 & & & & \\
1 & 0 & 0 & 0 & & & & \\
0 & 0 & 0 & -1 & & & & \\
0 & 0 & 1 & 0 & & & & \\
& & & & 0 & 1 & 0 & 0 \\
& & & & -1 & 0 & 0 & 0 \\
& & & & 0 & 0 & 0 & -1 \\
& & & & 0 & 0 & 1 & 0
\end{pmatrix}
\tag{3.12}
$$

By way of a parenthetical remark for now: The endomorphism $\rho_3$ defines a complex structure on X ($= \times_4 \mathbb{R}/\mathbb{Z}$) with complex coordinates $w_1 = x_1 + ix_2$ and $w_2 = x_3 + ix_4$.

Given the observations in the preceding paragraph, then the following is a consequence of (3.8) when $r$ is much greater than the parameter t:

$$\mathcal{D}_t \textit{ has no eigenvalues between } -\sqrt{2}r(1 - \mathcal{O}(\tfrac{t}{r})) \textit{ and } \sqrt{2}r(1 - \mathcal{O}(\tfrac{t}{r})) \textit{ on } \mathbb{R}^8 \otimes L^{\pm}. \tag{3.13}$$

Here are some ramifications in the case when $q >> m^2$ (and so $r >> m$): First, the STAGE 1 spectral flow (from Section 3) can only occur for the family of operators on the restricted space of $\mathbb{R}^8$-valued sections of $L^0$. To say more about the latter, note first that the operator $\gamma_\alpha \nabla_{A\alpha}$ on $\mathbb{R}^8$-valued sections of $L^0$ acts as if the bundle and connection were the product bundle with product connection. As a consequence, its spectrum is that of the operator

$$((u, u_0), a,) \to (((da)^+, d^\dagger a), d^\dagger u + du_0) \tag{3.14}$$

This implies in particular that if $q >> m^2$, then the spectral flow during STAGE 1 is at most equal to the kernel of the operator that is depicted in (3.14) which is 8 dimensional. This is because the STAGE 1 deformation has trivial kernel on the $L^0$ summand of the domain when $s > 0$ since the minimum of the spectrum in this case is $\frac{1}{4} m^2(1 - (1-s)^2)$. (The proof of this invokes (3.8).)

There is a second consequence of (3.14) and (3.8) which is this: There is no spectral flow for the STAGE 2 family of operators on the $L^0$ summand of the domain.

*Part 6*: To see about the STAGE 2 spectral flow on the $L^{\pm}$ summands of the domain, note that (3.8) has the following implication: As t increases from $m$ to some very large value $T \gg r$, the operator $\mathcal{D}_t$ will have an eigenvalue cross 0 from *below* only in the event that one (or both) of the following occur:

- $-\sqrt{2}\, r\, i\rho_3 + \frac{1}{2}t$ *has a kernel on* $TX^* \otimes L^+$.
- $\sqrt{2}\, r\, i\rho_3 + \frac{1}{2}t$ *has a kernel on* $TX^* \otimes L^-$.

(3.15)

In either case, the number of eigenvalues that cross 0 from below is equal to the dimension of the kernel of the operator $a \to (((d_A a)^+, d_A^{\dagger} a))$ when the latter operator is restricted to the kernel in (3.15).

Meanwhile, the operator $\mathcal{D}_t$ can have an eigenvalue cross 0 from *above* when t increases from $m$ to $T \gg x$ only in the event that one (or both) of the following occur:

- $-\sqrt{2}\, r\, i\rho_3 - \frac{1}{2}t$ *has a kernel on* $(\Lambda^+ \oplus \underline{\mathbb{R}}) \otimes_{\mathbb{R}} L^+$.
- $\sqrt{2}\, r\, i\rho_3 - \frac{1}{2}t$ *has a kernel on* $(\Lambda^+ \oplus \underline{\mathbb{R}}) \otimes_{\mathbb{R}} L^-$.

(3.16)

In either case, the number of eigenvalues that cross 0 from above is equal to the dimension of the kernel of the operator $(u, u_0) \to (d_A^{\dagger} u, + d_A u_0)$ restricted to the relevant kernel from (3.16).

Since $i\rho_3$ has eigenvalues $\pm 1$ and since t is increasing from $m$ to some $T \gg r$, it follows non-trivial kernels in the cases depicted by (3.15) and (3.16) can occur only when $t = 2\sqrt{2}\, r$. For this value of t, the requirements in (3.15) and (3.16) are satisfied only for the following eigenspaces:

- *The conditions in (3.15) are met when* $t = 2\sqrt{2}\, r$ *by*
  a) *The* +1 *eigenspace of* $i\rho_3$ *on* $TX^* \otimes_{\mathbb{R}} L^+$.
  b) *The* -1 *eigenspace of* $i\rho_3$ *on* $TX^* \otimes_{\mathbb{R}} L^-$.
- *The conditions in (3.16) are met when* $t = 2\sqrt{2}\, r$ *by*
  a) *The* -1 *eigenspace of* $i\rho_3$ *on* $(\Lambda^+ \oplus \underline{\mathbb{R}}) \otimes_{\mathbb{R}} L^+$.
  b) *The* +1 *eigenspace of* $i\rho_3$ *on* $(\Lambda^+ \oplus \underline{\mathbb{R}}) \otimes_{\mathbb{R}} L^-$.

(3.17)

With the preceding understood, the question at the heart of the matter is whether the operator $\gamma^{\alpha}\nabla_{A\alpha}$ has kernel on any of the subspaces depicted in (3.17) and, supposing there is a kernel, whether the kernel dimension in the top bullet's subspaces is equal to the kernel dimension in the lower bullet's subspaces. This dimension issue is crucial because if these dimensions are equal, then the number of eigenvalues crossing 0 from below is exactly equal to the number crossing 0 from above. As explained in the next part of the proof, the respective kernel dimensions are indeed identical (they do have positive dimension).

*Part 7*: To prove the equality of dimensions, note first that $i\rho_3$ anti-commutes with the operator $\gamma^\alpha\nabla_\alpha$. This implies that the operator which is defined by the rule

$$a \to \mathcal{L}a \equiv ((d_A a)^+, d_A^\dagger a) \tag{3.18}$$

sends any section of the +1 eigenspace of $i\rho_3$ on $T^* \otimes_{\mathbb{R}} L^+$ to a section of the -1 eigenspace of $i\rho_3$ on $(\Lambda^+ \oplus \underline{\mathbb{R}}) \otimes_{\mathbb{R}} L^+$. Likewise, it sends any section of the -1 eigenspace of $i\rho_3$ on $T^* \otimes_{\mathbb{R}} L^-$ to a section of the +1 eigenspace of $i\rho_3$ on $(\Lambda^+ \oplus \underline{\mathbb{R}}) \otimes_{\mathbb{R}} L^-$. (The +1 eigenspace of $i\rho_3$ on $T^* \otimes_{\mathbb{R}} L^+$ can be viewed as the tensor product over $\mathbb{C}$ of $L^+$ with (0,1) summand in the complex cotangent bundle for X with $\rho_3$ defining the complex structure. From this perspective, $\mathcal{L}$ is the operator $(\bar\partial, \bar\partial^*)$ as defined using the product connection on the complex cotangent bundle and the connection induced by $A_q$ on $L^+$.) The Hermitian adjoint operator $\mathcal{L}^\dagger$ maps these eigenspaces in the reverse direction; the operator $\mathcal{L}^\dagger$ sends $(u, u_0)$ to $d_A^\dagger u + d_A u^0$.

With the preceding understood, the key point in all of this is that the operator $\mathcal{L}$ has Fredholm index zero because the eigenspaces of $i\rho_3$ are trivial complex line bundles and the first Chern class of the line bundle $L^+$ has zero square. Meanwhile, the fact that $\mathcal{L}$ has index zero implies that the kernel of $\mathcal{L}$ in the space of sections of the +1 eigenspace of $i\rho_3$ in $T^* \otimes L^+$ has the same dimension as the kernel of $\mathcal{L}^\dagger$ on the space of sections of the -1 eigenspace of $i\rho_3$ in $(\Lambda^+ \oplus \underline{\mathbb{R}}) \otimes L^+$. By the same token, the kernel of $\mathcal{L}$ in the space of sections of the -1 eigenspace of $i\rho_3$ in $T^* \otimes L^-$ has the same dimension as the kernel of $\mathcal{L}^\dagger$ on the space of sections of the +1 eigenspace if $i\rho_3$ when viewed as a subbundle of $(\Lambda^+ \oplus \underline{\mathbb{R}}) \otimes L^-$.

By way of a parenthetical remark: The kernel of $\mathcal{L}$ is non-trivial as can be seen directly by writing out $\mathcal{L}$ explicitly. Let $a_1^+$ denote the $L^+$ part of $a_1$ and let $a_3^+$ denote the $L^+$ part of $a_3$. Then the requirement that $\mathcal{L}a$ should vanish requires the following:

- $(\nabla_{A1} + i\nabla_{A2})a_1^+ = 0$ *and* $\nabla_{A3}a_1^+ = \nabla_{A4}a_1^+ = 0$.
- $(\nabla_{A1} - i\nabla_{A2})a_3^+ = 0$ *and* $\nabla_{A3}a_3^+ = \nabla_{A4}a_3^+ = 0$.

(3.19)

This implies that $a_1^+$ restricts to each constant $(x_3, x_4)$ slice as a holomorphic section of the line bundle $L^+$. There are lots of these if $q > 0$ (the complex dimension of the space of these sections is 2q.) Meanwhile, $a_3^+$ restricts as an anti-holomorphic section; and there are none of these if $q > 0$.

**c) The spectral flow for the cases when ω is nowhere zero**

Section 2 describes certain |ω|-diverging sequences of solutions to (1.2) on a given oriented, Riemannian 4-manifold X with a given principle SO(3) bundle (denoted by P). These solutions are parametrized in part by an unbounded, countable sequence of real numbers. Any given version of the solution (A, ω) has the properties noted below in (3.20). The note refers to a real line bundle $\mathcal{I} \to X$ with a fiber metric.

*The* ad(P)*-valued self-dual 2-form* ω *can be written as* $r\varsigma^+\sigma$ *where r is a positive real number; where* $\varsigma^+$ *is a harmonic section of* $\Lambda^+ \otimes \mathcal{I}$ *with the integral of* $|\varsigma^+|^2$ *over* X *being 1; and where* σ *is an* A*-covariantly constant section of* $\text{ad}(P) \otimes \mathcal{I}$ *with norm* 1.

(3.20)

Some points to note for what follows: Being that $\varsigma^+$ is harmonic and the integral of the square of its norm being 1, the pointwise norm of $\varsigma^+$ has a $c_0$ upper bound as do the of norms of any k'th order covariant derivative of $\varsigma^+$ (the bound depends on k). Note inparticular that these bounds are independent of $r$ and $\varsigma^+$ even in the event that $\varsigma^+$ has some residual $r$ dependence. (Some of the |ω|-divergent sequence of solutions in Section 2 will have versions of $\varsigma^+$ that have the form $\varsigma^+_0 + \mathcal{O}(r^{-\alpha})$ with $\varsigma^+_0$ being independent of $r$ and with α being positive and independent of $r$.) Meanwhile, the curvature tensor of A can be written as $m\, r\varsigma\, \sigma$ with $\varsigma$ denoting an $\mathcal{I}$-valued, harmonic 2-form with a $c_0$ bound on the X-integral of the square of its norm. (This bound comes from a bound on the absolute value of the first Pontrjagin class of ad(P).) This integral bound, and the fact that ς is harmonic leads in turn to a $c_0$ bound on the pointwise norm of ς and a $c_0$ bound on the pointwise norm of its covariant derivatives to any given fixed order.

Some more points to note: In the case when $\mathcal{I}$ is isomorphic to the product $\mathbb{R}$-bundle, the 2-form ς (from $F_A$ which is $m\, r\varsigma\sigma$; and $\varsigma^+$ is ς's self-dual part) defines a 2-dimensional cohomology class on X. This class is denoted [ς]. When $\mathcal{I}$ is not isomorphic to the product line bundle, then the pull-back of ς to the 2-fold covering space defined by $\mathcal{I}$ (which is the set of points in $\mathcal{I}$ with norm 1) is a closed (and coclosed) 2-form. By virtue of it being closed, it defines a 2-dimensional cohomology class on this 2-fold covering space, this denoted by [ς] also. Note that the pull-back of the class [ς] by the generator of the covering space $\mathbb{Z}/2$ action is -[ς].

The proposition that follows momentarily considers the spectral flow for the operator in (1.4) when defined by a solution that is described in (3.20) with the extra assumption that $\varsigma^+$ is nowhere zero. (This is a very strong assumption; but even so, it is an open condition in the space of Riemannian metrics on X.)

Anyway, assuming that $\varsigma^+$ is nowhere zero: In the case when $\mathcal{I}$ is isomorphic to the product real line bundle, then $\Lambda^+$ can be written as $K \oplus \mathbb{R}\varsigma^+$ where K is a real, oriented 2-plane bundle. As such, it has an associate Euler class in $H^2(X; \mathbb{Z})$. This class is denoted by [K]. This is the first Chern class of K when K is viewed as a complex line bundle. If $\mathcal{I}$ is not isomorphic to the product $\mathbb{R}$-bundle, then $\varsigma^+$ can be viewed as an $\mathbb{R}$-valued 2-form on the 2-fold covering space for $\mathcal{I}$. Then, the orthogonal complement to the span of $\mathbb{R}\varsigma^+$ in the bundle of self-dual 2-forms on the 2-fold covering space is a real, oriented 2-plane bundle which has its Euler class in the 2-dimensional cohomology of this 2-fold cover. This class also pulls-back as -1 times itself via the action of the generator of the covering space $\mathbb{Z}/2$ action is also -1 times itself. This Euler class is denoted by [K] also.

**Proposition 3.2**: *Fixed a principle SO(3) bundle* $P \to X$ *and a number* $\Xi > 1$. *Let* $A_0$ *denote a reference connection on* P. *Given this data, there exists* $\kappa > 1$ *with the following significance: Let* $(A, \omega)$ *denote a solution to (1.2) with* A *being a section of* P *and with* $\omega$ *being a section of* $\Lambda^+ \otimes \mathrm{ad}(P)$. *Assume that this pair is described by (3.20) with* $r > \kappa$ *and with the norm of* $\varsigma^+$ *obeying* $|\varsigma^+| > \frac{1}{\Xi}$. *The absolute value of the spectral flow from the* $(A, \omega)$ *version of (1.4) to the* $(A_0, 0)$ *version of* (1.4) *differs by at most* $\kappa$ *from the absolute value of the cup-product pairing between the classes* $[\varsigma]$ *and* [K].

The rest of this subsection gives the proof of Proposition 3.2. Some specific instances of Proposition 3.2 with X being a complex manifold with a Kähler metric are presented in the next subsection.

***Proof of Proposition 3.2***: By way of a look ahead at the proof, the argument follows the same trajectory as the proof of Proposition 3.1. Even so, new difficulties arise due to the fact that the $\mathcal{I}$-valued 2-form $\varsigma^+$ isn't in general covariantly constant; and because of this, $\omega$ is not A-covariantly constant. A convention in the proof is that the values of the incarnations of the number $c_0$ depend on Proposition 3.2's number $\Xi$.

*Part 1*: The first point to make is that any given $(A, \omega)$ version of the operator in (1.4) can be written on a neighborhood of any given point in X using a local oriented, orthonormal frame for T*X and a local, oriented orthonormal frame for $\Lambda^+$ on that neighborhood so as to have the same form as what is depicted in (3.4). In this case, the $\alpha$-indices on $\nabla_A$ designate the directional covariant derivatives along the vector fields dual to the basis for T*X on the relevant neighborhood, and the k-indices on $\omega$ designate its components with respect to the basis for $\Lambda^+$ (each component is a section of ad(P) over the given neighborhood in X). The basic observation in this regards is that the bundle $(\Lambda^+ \oplus \underline{\mathbb{R}}) \oplus T^*X$ is an 8-dimensional Clifford module for a 7-dimensional Clifford algebra that is associated to the oriented orthormal frame bundle of X. To elaborate: Clifford multiplication by elements in $T^*X \oplus \Lambda^+$ acts as an algebra of endomorphisms of the bundle $(\Lambda^+ \oplus \underline{\mathbb{R}}) \oplus T^*X$ that obey the rules depicted in the top three bullets of (3.5). In addition, the Clifford multiplication action by any given section of $T^*X \oplus \Lambda^+$ (call it $p$) on a section $q$ of $(\Lambda^+ \oplus \underline{\mathbb{R}}) \oplus T^*X$ obeys the product rule with regards to the action of the metric's covariant derivative; this is to say that $\nabla(p \bullet q) = (\nabla p) \bullet q + p \bullet \nabla q$ with • used here to denote the Clifford multiplication action on the bundle $(\Lambda^+ \oplus \underline{\mathbb{R}}) \oplus T^*X$. With regards to $\Gamma$ in (3.4) and (3.5): This is the endomorphism that acts as +1 on on the $(\Lambda^+ \oplus \underline{\mathbb{R}})$ summand of $(\Lambda^+ \oplus \underline{\mathbb{R}}) \oplus T^*X$ and as -1 on the T*X summand. In particular, it obeys the algebraic rules that are depicted in the last two bullets of (3.5).

With the preceding understood, then (3.4) will be used henceforth to depict (1.4)'s operator in a suitable neighborhood of any given point using local, oriented orthonormal frames for T*X and $\Lambda^+$ on that neighborhood. As was the case previously, that operator is denoted by $\mathcal{D}$ with it understood that a specific choice for $(A, \omega)$ is needed to define it. By the same token, for any $t \in [0, \infty)$, the formula in (3.6) gives a local depiction of an operator on the space of ad(P)-valued sections of $(\Lambda^+ \oplus \underline{\mathbb{R}}) \oplus T^*X$ which will be denoted by $\mathcal{D}_t$. In the cases under consideration, the local depiction of this operator can be written as

$$\mathcal{D}_t = \gamma_\alpha \nabla_{A\alpha} + \tfrac{1}{\sqrt{2}}\, r \rho_k \varsigma^+_k [\sigma, \cdot\,] - \tfrac{1}{2}\, t\, \Gamma \,. \tag{3.21}$$

Note that $\mathcal{D}_t^2$ is still described by (3.7) with the simplification that the term with the covariant derivative of $\omega$ can be written as

$$\tfrac{1}{\sqrt{2}}\, t\, r\, \gamma_\alpha \rho_k (\nabla_\alpha \varsigma^+)_k [\sigma, \cdot\,], \tag{3.22}$$

with the covariant derivative here coming from the Levi-Civita connection. (Keep in mind that the covariant derivative of $\varsigma^+$ has an $r$-independent upper bound.)

Looking ahead at the upcoming analysis using (3.21) and (3.22): What is denoted above by $[\sigma, \cdot]$ defines an endomorphism of ad(P) only on subsets in X where an a priori product structure for $\mathcal{I}$ has been fixed. But having fixed such a product structure, then it is the case that $\mathcal{D}_t$ commutes with $[\sigma, \cdot]$. In particular, $\mathcal{D}_t$ preserves both the kernel of $[\sigma,\cdot]$ and the orthogonal complement to this kernel.

Viewed globally, commutation with $\sigma$ defines a homomorphism from ad(P) to ad(P)$\otimes \mathcal{I}$ (and vice versa). This fact will be exploited by writing the bundle ad(P) as $\mathfrak{L}\oplus\mathcal{I}$ with the $\mathcal{I}$ summand denoting elements in kernel of $[\sigma, \cdot]$ (which is the image of $\sigma$ when $\sigma$ is viewed in its guise as an isometric homomorphism from $\mathcal{I}$ to ad(P).) Meanwhile, $\mathfrak{L}$ denotes the 2-dimensional subbundle of elements in ad(P) that are orthogonal to the image of $\sigma$.

An important point with regards to the preceding is that $|[\sigma,\nu]| = 2|\nu|$ if $\nu \in \mathfrak{L}$. This is relevant to (3.21) and (3.22), because it implies the following: If $\psi$ is an $\mathfrak{L}$-valued section of $(\Lambda^+\oplus \mathbb{R})\oplus T^*X$, then the $\mathfrak{L}$-valued section of $(\Lambda^+\oplus \mathbb{R})\oplus T^*X)$ given by $\frac{1}{\sqrt{2}}\, r\rho_k\varsigma^+_k[\sigma,\psi]$ obeys

$$|\tfrac{1}{\sqrt{2}}\, r\rho_k\varsigma^+_k[\sigma,\psi]| = \sqrt{2}\, r|\varsigma^+||\psi| \,. \tag{3.23}$$

Of course if $\psi$ is an $\mathcal{I}$-valued section of $(\Lambda^+\oplus \mathbb{R})\oplus T^*X$, then $\frac{1}{\sqrt{2}}\, r\rho_k\varsigma^+_k[\sigma,\psi]$ is zero.

*Part 2*: The first stage of a path from the given $(A,\omega)$ version of $\mathcal{D}$ to an $(A_0,\omega_0)$ version is the path that is depicted in (3.10). As explained directly, there is an $r$-independent upper bound for the spectral flow along this path. To see why this is, suppose for the moment that $s \in [0, 1]$ and that $\psi$ is in the kernel of the parameter s version of (3.10). In this regard, both the $\mathcal{I}$-valued part of $\psi$ and the $\mathfrak{L}$-valued part must separately be in the kernel.

The $\mathcal{I}$-valued part of $\psi$ (written as $\eta\,\sigma$) is killed by the operator at a particular s value if and only if

$$\gamma_\alpha\nabla_\alpha\eta - \tfrac{1}{2}\, m\Gamma - \tfrac{1}{2}(1-s)\, m = 0 \,. \tag{3.24}$$

An important point with regards to (3.24) is that the covariant derivative in (3.24) is defined by the Levi-Civita connection; it is thus indendent of $r$. As a consequence, the spectral flow along this path is independent of the particular choice for (A, ω) (except to the extent that it is described by (3.20). Note in this regard that the s = 1 version of (3.24) has trivial kernel because its square is $(\gamma_\alpha\nabla_\alpha)^2 + \frac{1}{4}\, m^2$.

As for the $\mathfrak{L}$-valued kernel: Suppose for the moment that $\psi$ is a purely $\mathfrak{L}$ valued section of $(\Lambda^+\oplus \mathbb{R})\oplus T^*X$. Use (3.7) on $\psi$ with (3.22) and (3.23) to see that the inequality below holds for any given $t \geq 0$ and $r > c_0(1 + t)$

$$\int_X |\mathcal{D}_t\psi|^2 \geq \int_X |\gamma_\alpha\nabla_{A\alpha}\psi|^2 + \tfrac{1}{c_0}r^2|\psi|^2)$$

(3.25)

In the case at hand (with t = $m$), this inequality implies the following: If $r$ is greater than $c_0(m+1)$, then the kernel of any s ∈ [0, 1] version of the operator $\mathcal{D}_{t=m} - \frac{1}{2}(1-s)m$ is trivial when acting on $\mathfrak{L}$-valued sections of the bundle $(\Lambda^+ \oplus \underline{\mathbb{R}}) \oplus T^*X$. As a consequence, if $r > c_0(m+\tau)$, then there is no $\mathfrak{L}$-valued contribution to the spectral flow along the path $s \to \mathcal{D}_{t=m} - \frac{1}{2}(1-s)m$.

*Part 3*: This part of the proof considers a path of operators acting on ad(P) valued sections of $(\Lambda^+ \oplus \underline{\mathbb{R}}) \oplus T^*X$ that is parametrized by the interval $[r, R]$ for R >> $r$. A given $r \in [r, R]$ operator on this new path of operators is as follows:

$$\mathcal{D}_{t=m,r} \equiv \gamma_\alpha\nabla_{A\alpha} + \tfrac{1}{\sqrt{2}}\, r\,\rho_k\varsigma^+_k[\sigma,\cdot\,] - \tfrac{1}{2}\, m\Gamma\ .$$

(3.26)

Since these operators commute with $[\sigma, \cdot\,]$, the operators on this path can be analyzed by considering their respective restrictions to $\mathfrak{I}$-valued sections of $(\Lambda^+ \oplus \underline{\mathbb{R}}) \oplus T^*X$ and to $\mathfrak{L}$-valued sections. Since the family is constant on the space of $\mathfrak{I}$-valued sections, nothing more needs be said on that score. As for $\mathfrak{L}$-valued sections, a calculation that leads to (3.25) can be repeated to see that

$$\int_X |\mathcal{D}_{t=m,r}\psi|^2 \geq \int_X (|\gamma_\alpha\nabla_{A\alpha}\psi|^2 + \tfrac{1}{c_0}r^2|\psi|^2)$$

(3.27)

for all $r \geq r$ if $r \geq c_0(m+1)$. It follows from this that there is no spectral flow on the path if $r$ is sufficiently large.

*Part 4*: The step after this one will increase t starting with an operator very close to some very large R version of $\mathcal{D}_{t=m,R}$. To this end, two preliminary steps are required. The first preliminary step changes the operator by changing the norm of $\varsigma^+$ so that the resulting section of $\Lambda^+ \otimes \mathfrak{I}$ has norm 1. The modification deforms the r = R version of the operator depicted in (3.26) along a path parametrized by s ∈ [0, 1] that has the form

$$\gamma_\alpha\nabla_{A\alpha} + \tfrac{1}{\sqrt{2}}\, R\,(1 - s + s\tfrac{1}{|\varsigma^+|})\, \rho_k\,\varsigma^+_k[\sigma,\cdot\,] - \tfrac{1}{2}\, m\Gamma\ .$$

(3.28)

The operators along this path will obey the r = R version (3.27) as long as $R \geq c_0(m+\tau)$. Having traversed this path, the notation will be simplified by henceforth using $\varsigma^+$ to denote what would have previously been denoted by $|\varsigma^+|^{-1}\varsigma^+$. To be sure: This new version of $\varsigma^+$ has norm 1 but it is not (in general) a closed 2-form; but that requirement has no bearing in what follows.

Now for the second of the two preliminary steps: What this step does in effect is to replace the Levi-Civita connection with a connection that makes the new, norm 1 version of $\varsigma^+$ covariantly constant. To this end, it proves notationally convenient to choose the local orthonormal frames for $\Lambda^+$ so that $\varsigma^+$ is the third element of the frame. Having done this, then $\varsigma^+{}_k\rho_k$ is now $\rho_3$. Now introduce $\mathrm{T}$ to denote $\nabla\varsigma^+$; this being a section of $\Lambda^+\otimes T^*X$ that has zero metric pairing with $\varsigma^+$ at each point.

With $\mathrm{T}$ in hand, consider the family of operators which is parametrized by $s \in [0,1]$ and given by this rule:

$$\mathcal{D}_{t=m,R,s} \equiv \gamma_\alpha(\nabla_{A\alpha} + \tfrac{1}{2}\, s\, \mathrm{T}^k{}_\alpha\rho_k\rho_3) + \tfrac{1}{\sqrt{2}} R\rho_3[\sigma,\cdot\,] - \tfrac{1}{2}\, m\Gamma\ . \tag{3.29}$$

Because $\mathrm{T}$ is smooth and independent $m$ and R, the r = R version of (3.27) holds for $\mathcal{D}_{t=m,R,s}$ as long as $R \geq c_0(m+1)$. As a consequence, there is no spectral flow along this path of operators.

A key observation with regards to (3.29) is that the s = 1 version of $\mathcal{D}_{t=m,R,s}$ obeys an analog of the identity in (3.8) with any chosen value for $m$ (not just the specific value being considered). This analog of (3.8) is given below with $m$ denoted by t:

$$(\mathcal{D}_{t,R,1})^2 = (\gamma_\alpha(\nabla_{A\alpha} + \tfrac{1}{2}\mathrm{T}^k{}_\alpha\rho_k\rho_3))^2 + (\tfrac{1}{\sqrt{2}} R\rho_3[\sigma,\cdot\,] - \tfrac{1}{2}\, t\Gamma)^2\ . \tag{3.30}$$

This last identity holds because $\nabla_{A\alpha} + \frac{1}{2}\mathrm{T}^k{}_\alpha\rho_k\rho_3$ commutes with $\rho_3$ (since $\mathrm{T}$ has zero pairing with $\varsigma^+$) whereas the $\gamma$'s anti-commute with $\rho_3$. (Keep in mind that $\nabla_A\sigma = 0$ also.)

*Part 5*: With (3.30) understood, the next path is parametrized by $t \in [m, T]$ with T chosen to be much greater than R. The parameter t operator on this path is below:

$$\mathcal{D}_{t,R,1} = \gamma_\alpha(\nabla_{A\alpha} + \tfrac{1}{2}\mathrm{T}^k{}_\alpha\rho_k\rho_3) + \tfrac{1}{\sqrt{2}} R\rho_3[\sigma,\cdot\,] - \tfrac{1}{2}\, t\Gamma\ . \tag{3.31}$$

It follows as a consequence of (3.30) that spectral flow across zero on the path defined by (3.31) can occur only at the values of t where $\frac{1}{\sqrt{2}} R\rho_3[\sigma,\cdot\,] - \frac{1}{2}\, t\Gamma$ has non-zero kernel elements which are annihilated by the operator $\gamma_\alpha(\nabla_{A\alpha} + \frac{1}{2}\mathrm{T}^k{}_\alpha\rho_k\rho_3)$.

With regards to $\frac{1}{\sqrt{2}} R\rho_3[\sigma,\cdot\,] - \frac{1}{2}\, t\Gamma$ having a non-trivial kernel: The endomorphism $\rho_3[\sigma, \cdot]$ has eigenvalues 0 and ±2 at each point, and the corresponding eigenspaces define respective vector bundles on X. (This is so even in the case where $\mathcal{I}$ is not the product line bundle because then both $\rho_3$ and $\sigma$ are $\mathcal{I}$-valued.) Since $\Gamma$ commutes with $\rho_3[\sigma,\cdot\,]$ and since it has eigenvalues $\pm1$, there is just one positive value for t where $\frac{1}{\sqrt{2}}R\rho_3[\sigma,\cdot\,] - \frac{1}{2}t\Gamma$ has a kernel which is $t = 2\sqrt{2}R$. Moreover, for this value of t, any element in the kernel can be written as $\psi_+ + \psi_-$ where

- $\Gamma\psi_+ = \psi_+$ *and* $\rho_3[\sigma, \psi_+] = 2\psi_+$ .
- $\Gamma\psi_- = -\psi_-$ *and* $\rho_3[\sigma, \psi_-] = -2\psi_-$ .

(3.32)

With regards to the kernel of $\gamma_\alpha(\nabla_{A\alpha} + \frac{1}{2}T^k{}_\alpha\rho_k\rho_3)$: Let $n_+$ denote the dimension of the kernel of this operator on the subspace of sections that obey the condition in the top bullet of (3.32). Meanwhile, let $n_-$ denote the dimension of the kernel of this operator on the subspace of sections that obey the condition in the lower bullet of (3.32). The spectral flow along the path defined by (3.31) is equal to $-(n_+ - n_-)$. (This follows using first order perturbation theory given that the t-derivative of $\mathcal{D}_{t,R,1}$ is $-\frac{1}{2}\Gamma$. See for example Chapter 7 of [K].)

With regards to the numbers $n_+ - n_-$: The sections of $((\Lambda^+ \oplus \underline{\mathbb{R}}) \oplus T^*X) \otimes \mathfrak{su}_E$ that obey the top bullet of (3.32) are precisely the sections of subbundle in $(\Lambda^+ \oplus \underline{\mathbb{R}}) \otimes \mathfrak{L}$ where $\rho_3[\sigma,\cdot\,]$ acts as +2. Denote this eigenbundle by $E_+$. (By way of a reminder, $\mathfrak{L}$ is the orthogonal complement in $\mathfrak{su}_E$ to the span of $\sigma$.) By the same token, the sections of $((\Lambda^+ \oplus \underline{\mathbb{R}}) \oplus T^*X) \otimes \mathfrak{su}_E$ that obey the lower bullet of (3.32) are precisely the sections of the subbundle in $T^*X \otimes \mathfrak{L}$ where $\rho_3[\sigma,\cdot\,]$ acts as -2. The latter bundle is denoted by $E_-$. Because the operator $\gamma_\alpha(\nabla_{A\alpha} + \frac{1}{2}T^k{}_\alpha\rho_k\rho_3)$ anti-commutes with both $\Gamma$ and $\rho_3[\sigma,\cdot]$, it sends sections of $E_+$ to sections of $E_-$ and vice-versa. With this understood, let $\eth$ denote the restriction of $\gamma_\alpha(\nabla_{A\alpha} + \frac{1}{2}T^k{}_\alpha\rho_k\rho_3)$ to $C^\infty(X,E_+)$. Then $\eth$ defines a first order, elliptic operator mapping $C^\infty(X;E_+)$ to $C^\infty(X;E_-)$; and as such, it has a Fredholm index which is precisely $n_+ - n_-$. The computation of this index is described in the next part of the proof.

*Part 6*: When $\mathcal{I}$ is the product line bundle, then the Atiyah-Singer index theorem (the Hirzebruch-Riemann-Roch theorem when $\varsigma^+$ is the Kähler symplectic form for a Kähler metric on X) gives the formula below for $n_+ - n_-$ (the notation is explained subsequently):

$$n_+ - n_- = (1 + b_2^+ - b_1) + (t\cdot t - t\cdot[K]) .$$

(3.33)

By way of notation: The integers $b_2^+$ and $b_1$ are the respective self-dual 2'nd Betti number and first Betti number of X. Meanwhile, what is denoted by $t$ is the Euler class of the bundle $\mathfrak{L}$ and [K] denotes the Euler class of the orthogonal complement in $\Lambda^+$ of the span of $\varsigma^+$, an oriented $\mathbb{R}^2$ bundle over X. With regard to $t$: This is $\frac{1}{\pi}\, m\, r[\varsigma]$ with $[\varsigma]$ denoting the class that appears in the statement of Proposition 3.2. With regards to [K]: The orthogonal complement of the span of $\varsigma^+$ is the underlying real bundle of the canonical line bundle (a complex line bundle) for the symplectic structure defined by $\varsigma^+$. Finally, what is denoted in (3.33) by a·b when a and b are 2-dimensional cohomology classes signifies the evaluation of their cup product on the fundamental class of X.

*Part 7*: When $\mathcal{I}$ is a non-trivial line bundle, then the difference $n_+$ - $n_-$ obeys a formula that is analogous to (3.33). In this regard, both $t$ and [K] can be viewed on the 2-fold cover of X (the manifold W from Section 2b) as classes in $H^2(W; \mathbb{Z})$ that change sign when pulled back by the involution ι of W. It follows that their cup product pairings $t{\cdot}t$ and $t{\cdot}$[K] are well defined as classes in $H^4(X;\mathbb{Z})$. With this understanding, (3.33) again follows from the Atiyah-Singer index theorem. This can be seen directly via the heat equation proof of the index theorem which gives a local differential 4-form whose integral is $n_+ - n_-$. To elaborate just slightly: Near any given point in X, $\frac{1}{2}[\sigma, \cdot\,]$ can be viewed as an $\mathbb{R}$-valued endomorphism of $\mathfrak{L}$ that gives $\mathfrak{L}$ the structure of a complex vector bundle. Meanwhile, $\rho_3$ can likewise by viewed near any given point as an endomorphism of $(\Lambda^+\oplus\underline{\mathbb{R}})\oplus T^*X$ giving both $(\Lambda^+\oplus\mathbb{R})$ and $T^*X$ almost complex structures. It follows from this that the global 4-form coming via the heat equation proof of the index theorem must depend on the connection A and σ and $\rho_3$ when $\mathcal{I}$ is non-trivial in the same way it depends on the corresponding data when $\mathcal{I}$ is the product $\mathbb{R}$-bundle.

**d) Instances of Proposition 3.2**

This subsection considers the implications of Proposition 3.2 when X is a complex manifold with a Kähler metric and $\mathcal{I}$ is isomorphic to the product real line bundle. To say more, let w denote the Kähler symplectic form. This is a nowhere vanishing, self-dual harmonic 2 form. Any other self-dual harmonic 2 form can be writen as $\alpha w + c$ with $\alpha \in \mathbb{R}$ and with $c$ being the real part of a holomorphic (2, 0) form. The key point in this case is that w and $c$ are point-wise orthogonal, and that the zero locus of $c$ is a holomorphic complex dimension 1 subvariety. In particular, any version of $\varsigma^+$ has to have this form. If $\varsigma^+$ is to be nowhere vanishing, then the number α is non-zero unless $c$ is nowhere zero which can happen only if the first Chern class of X is zero (which implies that X is the quotient of a 4-torus or the K3 surface by an isometric group action [Y]). Suppose henceforth that this isn't the case for X.

If $\varsigma^+$ is nowhere zero, then the number α can't be zero which implies that the orthogonal complement in $\Lambda^+$ to $\varsigma^+$ is isomorphic to the underlying $\mathbb{R}^2$ bundle of the canonical complex line bundle for the complex structure on X. The Euler class of this $\mathbb{R}^2$ bundle (the class [K] in Proposition 3.2) is then independent of $\varsigma^+$ except for an ambiguity in its sign (±).

With regards to the class [K]: It is a class of type (1,1) which implies that it can be written as $[K] = \mu\,[w] + s$ where $\mu$ is a real number, [w] denotes the cohomology class of the Kähler symplectic form and $s$ denotes the cohomology class of a closed, anti-self dual harmonic 2-form.

As for the bundle $\mathfrak{L}$: As noted previously, it can be viewed as the underlying $\mathbb{R}^2$ bundle of a complex line bundle over X. Since the self-dual part of its curvature 2-form is $2\,i\,m\,r\,\varsigma^+$, the class $t$ must have the form $\frac{1}{\pi}\,m\,r(\alpha\,[w] + [c] + s')$ with $s'$ also denoting the class of an anti-self dual harmonic 2-form. The self-cup product of $t$ is then

$$t{\cdot}t = \tfrac{1}{\pi^2}\,m^2 r^2(\alpha^2 + c{\cdot}c - |s'{\cdot}s'|)\ , \tag{3.34}$$

which is by assumption, $\mathcal{O}(1)$. Meanwhile,

$$t{\cdot}[K] = \tfrac{1}{\pi}\,m\,r\,\mu\,\alpha + s{\cdot}s'\ . \tag{3.35}$$

By way of an example: If X is a minimal complex surface of general type, then it has a Kähler-Einstein metric (Yau's theorem [Y]) which implies that $s$ is zero and $\mu$ is positive. In this case, $t{\cdot}[K]$ is going to be $\mathcal{O}(r)$ and then the spectral flow will be $\mathcal{O}(r)$ unless $\alpha$ is $\mathcal{O}(r^{-1})$ which violates the assumption in Proposition 3.2 that $|\varsigma^+|$ is uniformly bounded away from zero. Thus, in the case when X has its Kähler-Einstein metric, the spectral flow is always $\mathcal{O}(r)$ when the assumptions of Proposition 3.2 are met.

On the other hand, suppose again that X is minimal and of general type; and assume that the self-dual and anti-self dual 2[nd] Betti numbers are no less than 3; and assume that the Kähler metric is not Kähler-Einstein (so that $s \neq 0$). Under these constraints, there are always versions of $\alpha$, $c$ and $s'$ so that the conditions of Proposition 3.2 are met with (3.34) being $\mathcal{O}(1)$ and (3.35) being zero. The spectral flow for these solutions is then $\mathcal{O}(1)$. But there are also versions of $\alpha$, $c$ and $s'$ under these same constraints where the conditions of Proposition 3.2 are satisfied with (3.34) being $\mathcal{O}(1)$ and (3.35) being $\mathcal{O}(r)$. The subsequent paragraphs describe a construction of these versions of $t$, versions with $t{\cdot}t = 0$, with $t{\cdot}[K] = 0$ and with non-zero cup product with the class of the Kähler symplectic form. (The constraints in the constructions that follow that guarantee the vanishing of $t{\cdot}[K]$ can be ignored so as to obtain classes $t$ with $t{\cdot}t = 0$, with $t{\cdot}[K] \neq 0$ and with non-zero cup product against the class of the Kähler form.)

Suppose first that X is not a spin manifold. Then it's intersection form can be diagonalized over $\mathbb{Z}$. Let $\{P_1, \ldots, P_{b^{2+}}\}$ denote an orthogonal basis of elements with positive square, and let $\{Q_1, \ldots Q_{b^{2-}}\}$ denote an orthonormal basis of elements with negative square. Here, $b^{2+}$ and $b^{2-}$ are the respective self-dual and anti-self dual second Betti numbers. Note in this regard that $b^{2+}$ is no less than 3 and $b^{2-}$ is no less than $\frac{1}{2}(b^{2+} - 1 + b^1)$, the latter being the Bogomolov-Miyaoka-Yau inequality. Assume in what follows that $b^{2-} \geq 3$.

The class [K] can be written as $\sum_j n_j P_j + \sum_\alpha m_\alpha Q_\alpha$ with the n's and m's being integers obeying $\sum_j n_j^2 - \sum_\alpha m_\alpha^2 = [K]\cdot[K]$. With this as background, chose three distinct P's (to be denoted by $P_{j(1)}, P_{j(2)}, P_{j(3)}$) and three distinct Q's (to be denoted by $Q_{\alpha(1)}, Q_{\alpha(2)}, Q_{\alpha(3)}$). Having done that, let $c = \{c_1, c_2, c_3\}$ denote a non-zero vector in $\mathbb{R}^3$ with rational entries that is orthogonal to $\{n_{j(1)}, n_{j(2)}, n_{j(3)}\}$ and to $\{m_{\alpha(1)}, m_{\alpha(2)}, m_{\alpha(3)}\}$. There is such a vector because $b^{2+} \geq 3$ and $b^{2-} \geq 3$. The class $t = \sum_{1\leq k\leq 3} c_k(P_{j(k)} + Q_{\alpha(k)})$ will be orthogonal to [K] and have square zero. Meanwhile, the self-dual part of a harmonic 2-form representing $t$ will be nowhere zero if $t$ has non-zero cup product with the class of the Kähler symplectic form. That class can be written as $\sum_j u_j P_j + \sum_\alpha v_\alpha Q_\alpha$. The pairing of $t$ with that class is $\sum_{k=1,2,3} c_k(u_{j(k)} - v_{\alpha(k)})$. One can choose the vector c so that this is non-zero unless $v_{\alpha(k)}/m_{\alpha(k)} = u_{j(k)}/n_{j(k)}$ for all $k \in \{1, 2, 3\}$. With this understood note that the particular choice of the subset $\{P_{j(k)}\}_{k=1,2,3}$ from the set of P's and the subset $\{Q_{\alpha(k)}\}_{k=1,2,3}$ from the set of $Q_\alpha$'s was arbitrary. Choosing any other collection of three P's and three Q's gives a different version of $t$; and if $t\cdot[w]$ is zero for all of these version of $t$, then for all $\alpha$ and j, one must have $v_\alpha/m_\alpha = u_j/n_j$. This implies that [w] is proportional to [K] which was assumed not to be the case. Thus, there exists a choice for the three P's and the three Q's so that the corresponding version of $t$ has non-zero cup product pairing with the class of the Kähler form, zero cup product with [K] and zero self-cup product.

Suppose next that X is a spin manifold. It's intersection form is a direct some of t some number N (the lesser of $b^{2+}$ and $b^{2-}$) of hyperbolic pairs and then an even, definite form which is a direct sum of either $E_8$ lattices or a direct sum of $-E_8$ lattices. Represent the hyperbolic pair classes as $(P_1,Q_1), \ldots, (P_N, Q_N)$ with the P's and Q's having square zero, with each P and Q orthogonal to the P's and Q's with different indices, and with $P_j\cdot Q_j = 1$ for each j. The class [K] can be written as $\sum_j (n_j P_j + m_j Q_j) + T$ where T is orthogonal to all P's and Q's; it is in the definite summand for the intersection pairing. Meanwhile the n's and m's obey $\sum_j n_j m_j + T\cdot T > 0$. By the same token, the class of the symplectic form can be written as $\sum_j (u_j P_j + v_j Q_j) + T'$.

Now choose two elements with distinct indices from the set $\{P_j\} \cup \{Q_j\}$ to be denoted by $U_1$ and $U_2$. Use $x_1$ and $x_2$ to denote the corresponding $n_j$ and $m_j$ as the case may be (either both are n's or both are m's or one of each). Use $y_1$ and $y_2$ denote the corresponding $u_j$ or $v_j$. Let $c = (c_1, c_2)$ denote a vector in $\mathbb{R}^2$ with rational entries such that $c_1x_1 + c_2x_2$ is zero. Then $t = c_1U_1 + c_2U_2$ has square zero and it is orthogonal to [K]. This class will have non-zero cup product pairing with the Kähler class unless it holds that $y_2/x_2 = y_1/x_1$. Since the choice of the U's was arbitrary, this last condition can't hold for all possible choices of $U_1$ and $U_2$ unless $u_j = \alpha n_j$ and $v_j = \alpha m_j$ for all j with $\alpha$ being a real number. Suppose for the sake of argument that this is the case. It then follows that T and T´ must be different vectors in the definite summand of 2'nd cohomology.

There are two cases to consider at this point: The first is when either the vector of n's or the vector of m's (vectors in $\mathbb{R}^N$) are non-zero. The other is when they are both zero. For the first case, assume without loss of generality that the vector of n's is not zero. Let $T_c$ denote for the moment any rational class in the definite summand of the intersection form. Let $(c_1, \dots , c_N)$ denote a non-zero vector in $\mathbb{R}^N$ with rational entries that obeys $\sum_j c_j m_j = 0$, and then let $(x_1, \dots, x_N)$ denote a second vector in $\mathbb{R}^N$. The class

$$t = \sum_j c_j P_j + \sum_j x_j Q_j + T_c \tag{3.36}$$

will have zero cup product with [K] if

$$\sum_j x_j n_j + T_c \cdot T = 0. \tag{3.37}$$

It will have zero self cup product if

$$\sum_j x_j c_j + \tfrac{1}{2} T_c \cdot T_c = 0. \tag{3.38}$$

If the vectors $(n_1, \dots , n_N)$ and $(m_1, \dots, m_N)$ are both non-zero, then these two equations can be solved simultaneouly for any choice of $T_c$ because the vector of n's and the vector of c's are linearly independent (the former having positive inner product with the vector of m's and the latter having zero inner product with that vector). In this case, $T_c$ can be chosen so that it has non-zero inner product with T´ and thus $t$ has non-zero cup product with the Kähler class. If the vector of m's is zero, then the two equations can be solved simultaneously also since the vector of c's can be chosen at will and specifically, to be linearly independent from the vector of n's.

Now assume that both the vector of n's and the vector of m's are zero. In this case, T´ and T can't be colinear (because the metric isn't Kähler-Einstein) in which case one can take $T_c$ to be orthogonal to T which makes the top equation in (3.37) true but not orthogonal to T´. Meanwhile, the vectors of x's is chosen so that (3.38) holds.

Thus in all cases, $t$ can be chosen to have zero self-cup product, have zero cup product with [K] and non-zero cup product with the class of the Kähler symplectic form.

## 4. Model problems for the eigenvalues near zero when $\varsigma^+$ has transveral zeros

An assumption in this section is that (A, ω) are described in (3.20). Two additional assumptions are made here, the first being this: If the section $\varsigma^+$ in (3.20) has non-empty zero locus, then $\varsigma^+$ vanishes transversally along that zero locus. More is said about what this means in the next paragraph. The second assumption in some of the propositions to come is that the number $m$ from (1.2) is small (but still $\mathcal{O}(1)$).

To be explicit about what transversal vanishing implies: It means first that the zero locus is a disjoint union of smoothly embedded circles. This is to say that any given component of Z can be parametrized by a smooth embedding from the circle (parametrized so that the differential has unit norm) into X. Transversal vanishing also implies that a neighborhood of this embedding is diffeomorphic to $S^1 \times B_0$ where $B_0$ is a small radius ball in $\mathbb{R}^3$ centered at the origin. The radius of this ball is denoted by $r_0$. This number is such that these respective neighborhoods for disjoint components of the zero locus of $\varsigma^+$ have disjoint closures. In the applications to come, the diffeomorphism from $S^1 \times B_0$ into X is obtained by first choosing an oriented, orthonormal frame for the normal bundle of each component of Z to write that bundle as a product $\mathbb{R}^3$-bundle and then restricting the metric's exponential map along each component of the zero locus to this normal bundle.

When pulled back to $S^1 \times B_0$ by the diffeomorphism, the form $\varsigma^+$ appears as

$$\varsigma^+ = \mathbb{M}_{ij}(s)\, x^i (ds \wedge dx^j + \tfrac{1}{2}\, \varepsilon_{jkm} dx^k \wedge dx^m) + \mathfrak{e} \quad \textit{with } |\mathfrak{e}| \leq \mathcal{k}|x|^2 . \tag{4.1}$$

where the notation is as follows: First, s is the Euclidean parameter on $S^1$ and $\{x_1, x_2, x_3\}$ are Euclidean coordinates on $\mathbb{R}^3$. Second, $\varepsilon_{jkm}$ are the components of the completely anti-symmetric 3-tensor on $\mathbb{R}^3$ with $\varepsilon_{123} = 1$. Third, $\{\mathbb{M}_{ij}\}_{i,j\in\{1,2,3\}}$ at any given value for the parameter s are the components of a $3\times3$ matrix which is invertible at each parameter value. (The matrix $\mathbb{M}$ is invertible so that $\varsigma^+$ vanishes transversally.) Finally, $\mathcal{k}$ is a non-negative number which has a $c_0$ upper bound. The condition that $\varsigma^+$ be a closed 2-form requires that $\mathbb{M}$ be traceless and symmetric matrix, but these extra constraints are not needed for what is to come; and it proves useful to allow for any invertible $\mathbb{M}$.

With regards the line bundle $\mathcal{I}$ along the zero locus: It must be isomorphic to the product line bundle on each component. This is because the determinant of $\mathbb{M}$ is a section of $\mathcal{I}$ that can't vanish at any parameter value (which is due to the fact that $\mathbb{M}$ is assumed to be invertible at each parameter value.) By way of an aside: The fact that $\mathbb{M}$ is traceless plays no role in the subsequent analysis except for the implication that it is neither positive nor negative definite.

Looking ahead: As was the case before, the bundle ad(P) has the A-covariantly constant direct sum splitting as $\mathfrak{L} \oplus \mathcal{I}$. Since this splitting is respected by the operator $\mathcal{D}$, the spectral flow (as a function of $r$) is the sum of the corresponding spectral flows from $\mathcal{D}$'s restrictions to the respective sections of $\mathcal{I}$ and $\mathfrak{L}$-valued sections of $(\Lambda^+ \oplus \underline{\mathbb{R}}) \oplus T^*X$. The sole focus in what follows is on the restriction of $\mathcal{D}$ to the space of $\mathfrak{L}$-valued sections of $(\Lambda^+ \oplus \underline{\mathbb{R}}) \oplus T^*X$. This is because $\mathcal{D}$'s action on the space of $\mathcal{I}$ valued sections does not depend on the $r$.

**a) A local model case**

The analysis to come is built around a model situation that has a version of the operator $\mathcal{D}$ on the product manifold $S^1 \times \mathbb{R}^3$. The correspondence here is that $S^1 \times \{0\}$ is an idealized model for a component of the zero locus of $\varsigma^+$, and the full $S^1 \times \mathbb{R}^3$ is an idealized model for the rest of X. The connection between the model and the actual situation comes (in part) from the upcoming Proposition 4.1 which says in effect that eigenvectors of $\mathcal{D}$ with small eigenvalues (in absolute value) concentrate very near to the zero locus of $\varsigma^+$ when $r$ is large.

The promised proposition concerns a somewhat more general version of $\mathcal{D}$ which is the operator that is depicted momentarily in (4.2). Just to be sure about assumptions and notation in (4.2) and Proposition 4.1: These take as given a principle SO(3) or SU(2) bundle over X (this is the bundle P) and a real line bundle over X (this is the line bundle $\mathcal{I}$), and an isometric homomorphism $\sigma$ from $\mathcal{I}$ to ad(P). The bundle ad(P) then splits as $\mathcal{I} \oplus \mathfrak{L}$ with $\mathfrak{L}$ being an $\mathbb{R}^2$-bundle; and the commutator with $\sigma$ defines a homomorphism from $\mathfrak{L}$ to $\mathfrak{L} \otimes \mathcal{I}$. What is denoted by $\varsigma^+$ in (4.2) is an $\mathcal{I}$-valued section of $\Lambda^+$ with transverse zero locus (this is the set Z). The generalization of $\mathcal{D}$ depicted (4.2) also refers to a connection on P (this is denoted by A) with $\nabla_A \sigma$ being zero. By virtue of $\nabla_A \sigma$ being zero, the $\nabla_A$-directional covariant derivatives and the operator that is depicted in (4.2) map $\mathfrak{L}$-valued sections of $(\Lambda^+ \oplus \underline{\mathbb{R}}) \oplus T^*X$ to $\mathfrak{L}$-valued sections of this same bundle and that is how they are viewed in Proposition 4.1 and subsequently. Two last notational remarks about (4.2): What (4.2) denotes by $\mathfrak{Q}$ signifies an endomorphism of the bundle $((\Lambda^+ \oplus \underline{\mathbb{R}}) \oplus T^*X) \otimes \mathfrak{L}$; and what (4.2) denotes by R is a non-negative real number.

What follows is promised the generalized version of $\mathcal{D}$:

$$\mathcal{D}_{\ddagger} = \gamma_\alpha \nabla_{A\alpha} + \tfrac{1}{\sqrt{2}} R \varsigma^+_k \rho_k [\sigma, \cdot\,] - \mathfrak{Q} \,. \tag{4.2}$$

And, what follows next is the promised proposition about $\mathcal{D}_{\ddagger}$.

**Proposition 4.1**: *Given $\Xi > 1$ and a positive number $\kappa_{\ddagger}$, there exists $\kappa > 1$ such that the following is true: Fix $R > \kappa$ to define the operator $\mathcal{D}_{\ddagger}$ that is depicted in (4.2) using a pair $(A, \omega)$ as described in (3.20) and an endomorphism $\mathfrak{Q}$ subject to the additional constraints listed below.*

- *The section $\varsigma^{+}$ obeys $|\varsigma^{+}| \geq \frac{1}{\Xi}\mathrm{dis}(\cdot, Z)$*
- *The norm of the curvature of the connection $A$ is bounded by $\kappa_{\ddagger} R$;*
- *The norm of the endomorphism $\mathfrak{Q}$ is bounded by $\frac{1}{\kappa} R$.*

*Now let $\psi$ denote an $\mathfrak{L}$-valued eigenvector for $\mathcal{D}_{\ddagger}$ whose eigenvalue has absolute value at most $\frac{1}{\kappa}\sqrt{R}$ and such that $\int_X |\psi|^2 = 1$. Then $\psi$ obeys*

$$|\psi|(\cdot) \leq \kappa R^2\, e^{-\frac{1}{\kappa} R\, \mathrm{dist}^2(\cdot, Z)}.$$

The proof is in Section 4e. What follows directly is a description of the model problem on $S^1 \times \mathbb{R}^3$.

The model problem requires the choice of a positive number $\ell$ which will be the length of the $S^1$ factor for the product, Euclidean metric on $S^1 \times \mathbb{R}^3$. This is to say that $S^1$ defined as $\mathbb{R}/(\ell\mathbb{Z})$. The $\mathbb{R}/(\ell\mathbb{Z})$ coordinate for the $S^1$ factor is denoted by s; and the metric is such that ds has norm 1. Meanwhile, the Euclidean coordinates for the $\mathbb{R}^3$ factor are $(x_1, x_2, x_3)$ with the metric making their differentials pairwise orthogonal, orthogonal to dt and having norm 1. These coordinates isometrically identify the $\mathbb{R}/(\ell\mathbb{Z}) \times \mathbb{R}^3$ versions of the respective bundles $(\Lambda^{+} \oplus \underline{\mathbb{R}})$ and $T^*X$ with the product $\mathbb{R}^4$ bundles and in so doing, they identify the $\mathbb{R}/(\ell\mathbb{Z}) \times \mathbb{R}^3$ version of the Euclidean metric's Levi-Civita connection with the respective product connections.

The definition of the model operator requires first the choice of an $\mathbb{R}/(\ell\mathbb{Z})$-parametrized family of invertible, $3 \times 3$ matrices. This matrix valued function on $\mathbb{R}/(\ell\mathbb{Z})$ in denoted by $\mathbb{M}$; its components at any given parameter value s are denoted by $\{\mathbb{M}_{jk}(s)\}_{j,k=1,2,3}$. The definition of the model operator also requires the specification of a positive real number to be denoted by R.

With the data $\mathbb{M}$ and R chosen: The model operator is defined initially on the space of smooth compactly supported, $\mathbb{C}$-valued sections over $\mathbb{R}/(\ell\mathbb{Z}) \times \mathbb{R}^3$ of the product $\mathbb{R}^4 \oplus \mathbb{R}^4$ bundle (thus, maps to $\mathbb{C}^4 \oplus \mathbb{C}^4$). It is

$$\mathcal{D}_0 = \gamma_s \frac{\partial}{\partial s} + \mathfrak{D}_0 \tag{4.3}$$

where $\mathfrak{D}_0$ is an s-dependent differential operator on $\mathbb{R}^3$ that has the form

$$\mathfrak{D}_0 = \sum_{k=1}^{3} \gamma_k \frac{\partial}{\partial x_k} + \sqrt{2} i R \sum_{j,k=1}^{3} \mathbb{M}_{jk}(s)\, x_j \rho_k \tag{4.4}$$

The upcoming Lemma 4.2 summarizes what is needed regarding $\mathfrak{D}_0$ at any given value of the parameter s. The subsequent Lemmas 4.3 and 4.4 do the same for $\mathcal{D}_0$.

With regard to notation for Lemma 4.2 and subsequently: Lemma 4.2 and its proof use $\lambda_1, \lambda_2$ and $\lambda_3$ to denote the positive square roots of the eigenvalues at the given value of s of $\mathbb{M}^T\mathbb{M}$. The lemma and the subsequent discussions and lemmas use $\mathbb{L}$ to denote the Hilbert space completion of the vector space of compactly supported, smooth maps from $\mathbb{R}^3$ to $\mathbb{C}^4 \oplus \mathbb{C}^4$ using the norm whose square sends any given map $\psi$ to the $\mathbb{R}^3$-integral of $|\psi|^2$. The inner product on $\mathbb{L}$ is denoted as $\langle\ ,\ \rangle_{\mathbb{L}}$; and the associated norm is called the $\mathbb{L}$-norm; it is denoted by $\|\cdot\|_{\mathbb{L}}$.

**Lemma 4.2**: *The operator $\mathfrak{D}_0$ extends to the Hilbert space $\mathbb{L}$ as an unbounded, self-adjoint operator with discrete spectrum having no accumulation points and finite multiplicities. In this regard,*

- *The kernel of $\mathfrak{D}_0$ is 1-dimensional.*
- *An orthonormal basis of eigenvectors for the nonzero eigenvalues of $\mathfrak{D}_0$ is labeled by the data sets that have the form $(\varepsilon_0, (\delta_1, \delta_2, \delta_3), (n_1, n_2, n_3))$ with $\varepsilon_0 \in \{1, -1\}$, with $\delta_1, \delta_2$, and $\delta_3$ from $\{0, 1\}$, and with $n_1, n_2$ and $n_3$ being non-negative integers; but these are subject to the following constraints:*
  a) *Not all of* $\delta$'s *and* n'*s are zero*
  b) *If only one pair from $\{(\delta_k, n_k)\}_{k=1,2,3}$ is not $(0, 0)$, then the eigenvector is determined by $\varepsilon_0$ and the sum of the non-zero pair (thus $\delta + n$).*

  *The eigenvalue that corresponds to a given data set $(\varepsilon_0, (\delta_1, \delta_2, \delta_3), (n_1, n_2, n_3))$ is the positive $(\varepsilon_0 = 1)$ or negative $(\varepsilon_0 = -1)$ square root of*

  $$2\sqrt{2}R\,((\delta_1 + n_1)\lambda_1 + (\delta_2 + n_2)\lambda_2 + (\delta_3 + n_3)\lambda_3)\ .$$
- *With regards to the parameter $\varepsilon_0$: If $\phi$ is an eigenvector for $\mathfrak{D}_0$ with eigenvalue $\hat{E}$, then $\gamma_s\phi$ is an eigenvalue with eigenvalue $-\hat{E}$.*

This lemma is proved in Section 4b.

By way of a remark: The first and third bullets of the lemma imply that the kernel of $\mathfrak{D}_0$ is a 1-dimensional eigenspace of $\gamma_s$. (Note that $\gamma_s$ has eigenvalues $\pm i$ and the convention in what follows is that the eigenvalue is +i on the kernel of $\mathfrak{D}_0$.)

The simplest case for $\mathcal{D}_0$ is the case when the matrix $\mathbb{M}$ is constant. In this case, the operator $\frac{\partial}{\partial s}$ commutes with $\mathcal{D}_0$, and then the tried and true method of separation of variables can be used to find $\mathcal{D}_0$'s eigenvalues and eigenvectors. The result is the upcoming Lemma 4.3. To set additional notation for Lemma 4.3 and subsequently: What is denoted by $\underline{\mathbb{L}}$ signifies the Hilbert space that is obtained by completeing the space of compactly supported maps from $\mathbb{R}/(\ell\mathbb{Z})\times\mathbb{R}^3$ to $\mathbb{C}^4\times\mathbb{C}^4$ using the norm whose square is the function that sends any given map $\psi$ to the integral on $\mathbb{R}/(\ell\mathbb{Z})\times\mathbb{R}^3$ of the function $|\psi|^2$. (This is the $L^2$-inner product; it is also the $\mathbb{R}/(\ell\mathbb{Z})$-integral of the square of $\psi$'s pointwise $\mathbb{L}$-norm, $\|\psi\|_{\mathbb{L}}$. )

**Lemma 4.3**: *Supposing that $\mathbb{M}$ is constant, then the following is true: The operator $\mathcal{D}_0$ extends to the Hilbert space $\underline{\mathbb{L}}$ as an unbounded, self-adjoint operator with discrete spectrum having no accumulation points and finite multiplicities. In this regard, the eigenvectors and eigenvalues of $\mathcal{D}_0$ have the form* $E = \pm(4\pi^2 n^2/\ell^2 + \hat{E}^2)$ *where* $n \in \mathbb{Z}$ *and* $\hat{E}$ *is a non-negative eigenvalue of* $\mathfrak{D}_0$.

- *If* $\hat{E} = 0$*, then the corresponding eigenvector is* $e^{2\pi i n s/\ell}\,\phi$ *with* $\phi$ *being a non-zero element in the kernel of* $\mathfrak{D}_0$.
- *If* $\hat{E} > 0$*, then the corresponding eigenvector is* $e^{2\pi i n s/\ell}\,((E+\hat{E})\,\phi + 2\pi i n/\ell\,\gamma_s\,\phi)$.

This lemma is also proved in Section 4b.

The upcoming Lemma 4.4 describes (in part) the spectrum of $\mathcal{D}_0$ when $\mathbb{M}$ is not constant. To set the stage for this lemma: The product bundle $\mathbb{R}/(\ell\mathbb{Z})\times\mathbb{L}$ is a Hilbert space bundle over $\mathbb{R}/(\ell\mathbb{Z})$ with a distinguished line subbundle to be denoted by $\mathcal{K}$ whose fiber at any $s \in \mathbb{R}/(\ell\mathbb{Z})$ is the kernel of the parameter s version of $\mathfrak{D}_0$. A section of $\mathcal{K}$ (to be denoted by $\phi_0$) is said to have unit norm means that $\langle\phi_0, \phi_0\rangle_{\mathbb{L}} = 1$ at each $s \in \mathbb{R}/(\ell\mathbb{Z})$.)

What follows is the promised lemma about $\mathcal{D}_0$ when $\mathbb{M}$ is not constant.

**Lemma 4.4**: *When $\mathbb{M}$ is not constant, the operator $\mathcal{D}_0$ still extends to the Hilbert space $\underline{\mathbb{L}}$ as an unbounded, self-adjoint operator with discrete spectrum having no accumulation points and finite multiplicities. Moreover, there exists a number* $\kappa > 1$ *and a number* $\alpha \in [0,1)$ *and a unit normed section* $\phi_0$ *of* $\mathcal{K}$ *with the following significance: Supposing that* $R \geq \kappa$*, then an eigenvalue of* $\mathcal{D}_0$ *with absolute value at most* $\frac{1}{\kappa^2}\sqrt{R}$ *can be written as*

$$E = -(\alpha + n)\frac{2\pi}{\ell} + \mathfrak{r}$$

*with* n *denoting an integer and with* $\mathfrak{r}$ *denoting a real number with norm bounded by* $\frac{1}{\sqrt{R}}\kappa$. *Conversely, if* n *is an integer with norm bounded by* $\frac{1}{\kappa}\sqrt{R}$, *then there is a* *corresponding eigenvalue of $\mathcal{D}_0$ that can be written as above with $|\mathfrak{r}| \leq \frac{1}{\sqrt{R}}\kappa$. The corresponding the eigenvector in all of these cases differs from $e^{2\pi i n s/\ell}\phi_0$ by map from $\mathbb{R}/(\ell\mathbb{Z})\times\mathbb{R}^3$ to $\mathbb{C}^4\oplus\mathbb{C}^4$ with a $\frac{1}{R}\kappa$ bound on the $\mathbb{R}/(\ell\mathbb{Z})\times\mathbb{R}^3$ integral of the square of its norm (which is its norm as an element of $\underline{\mathbb{L}}$). With regards to the number $\kappa$: An upper bound for $\kappa$ is a priori determined given upper bounds for the norms along $\mathbb{R}/(\ell\mathbb{Z})$ of $\mathbb{M}^{-1}$, $\mathbb{M}$ and $\frac{d}{ds}\mathbb{M}$.*

Lemma 4.4 is proved in the next subsection after the proofs of Lemmas 5.2 and 5.3 which follow in turn momentarily.

**b) Analysis for the model problems**

This subsection contains (momentarily) the proofs Lemmas 4.2-4.4. The lemma that follows directly is used in the proofs of all three. To set the notation: The upcoming lemma refers to a smooth Riemannian manifold to be denoted by X which is non-compact but geodesically complete. The lemma concerns a complex or real vector bundle over X with a Hermitian or orthogonal connection as the case may be. The bundle is denoted by V and the connection's covariant derivative is denoted by $\nabla$. Let $\mathcal{L}$ denote the Hilbert space that is obtained by completing the space of compactly supported sections of V using the norm whose square sends a section (denoted by $\mathfrak{s}$) to the X-integral of $|\mathfrak{s}|^2$. (This is the $L^2$-norm.) The lemma also refers to a non-negative function on X (denoted by $\mathfrak{f}$) that is unbounded in the following strong sense: Given $R > 0$, there exists $r > 0$ such that $\mathfrak{f} > R$ where the distance to a specified point is greater than r. (The point is fixed a prior).

**Lemma 4.5**: *Let $\mathbb{D}$ denote a first order, symmetric elliptic operator mapping sections of V to sections of V with the following property: There exists $\kappa > 1$, a number $\mathrm{R} \geq 0$, and a function $\mathfrak{f}$ as described above such that if $\mathfrak{s}$ is any section of V with compact support, then*

$$\int_X |\mathbb{D}\mathfrak{s}|^2 \geq \int_X \left(\tfrac{1}{\kappa}|\nabla\mathfrak{s}|^2 + (\mathfrak{f} - \mathrm{R})|\mathfrak{s}|^2\right) .$$

*Then $\mathbb{D}$ extends to $\mathcal{L}$ as an unbounded, self-adjoint operator with purely discrete spectrum lacking accumulation points and with finite multiplicities.*

Given the growth condition on $\mathfrak{f}$, the proof is only marginally different from the proof of the analogous theorem for compact manifolds. Because of this, the proof is omitted. (See [H], Chapter XIX).

***Proof of Lemma 4.2***: The $3\times3$ matrix $\mathbb{M}$ at any given parameter value s can be written as $U^{-1}V^{-1}DV$ where U and V are special orthogonal matrices and where D is a diagonal matrix with eigenvalues $\varepsilon\lambda_1$, $\varepsilon\lambda_2$, $\varepsilon\lambda_3$ with $\varepsilon$ denoting the sign of the determinant of $\mathbb{M}$. Using V to change coordinates and then using U to rotate the $\rho_k$'s produces the operator below:

$$\sum_{k=1}^{3}\gamma_k \frac{\partial}{\partial x_k} + \sqrt{2}i\,\varepsilon\sum_{k=1}^{3}\lambda_k x_k \rho_k \ . \tag{4.5}$$

To be sure: The $x_k$'s and $\gamma_k$'s in (4.5) and subsequently are obtained from those in (4.4) using V; and the $\rho_k$'s are obtained from those in (4.4) using U and V. A crucial point in this regard is that the new versions of the $\gamma_k$'s and $\rho_k$'s obey the algebra that is depicted in (3.5). This fact is used implicitly in what follows.

To understand $\mathfrak{D}_0$, look first at $\mathfrak{D}_0^{\,2}$:

$$\mathfrak{D}_0^{\,2} = -\sum_{k=1}^{3}\frac{\partial^2}{\partial x_k^2} + \sqrt{2}R\varepsilon\sum_{k=1}^{3}\lambda_k i\gamma_k\rho_k + 2R^2\sum_{k=1}^{3}\lambda_k^2 x_k^2 \ , \tag{4.6}$$

This depiction of $\mathfrak{D}_0^{\,2}$ with an integration by parts can be used to see that the assumptions of Lemma 4.5 are met (with $X = \mathbb{R}^3$ and $\mathfrak{f}(x) = R^2|x|^2$), and so that lemma's conclusions apply to $\mathfrak{D}_0$; which is that $\mathfrak{D}_0$ has purely discrete spectrum with no accumulation points and finite multiplicities.

To say more about the eigenvectors of $\mathfrak{D}_0$, it is useful to introduce a set of three endomorphisms of $((\Lambda^+\oplus\underline{\mathbb{R}})\oplus T^*(S^1\times\mathbb{R}^3))\otimes\mathbb{C}$:

$$\{i\varepsilon\gamma_1\rho_1, i\varepsilon\gamma_2\rho_2, i\varepsilon\gamma_3\rho_3\} \ . \tag{4.7}$$

Each of these has eigenvalues $\pm1$, each eigenspace has dimension 4 (over $\mathbb{C}$), and they pairwise commute. Moreover, these endomorphisms are covariantly constant as they have no $x_k$ dependence. It follows as a consequence of these observations (and the multiplication rules in (3.5)) that the vector space $\mathbb{C}^4\oplus\mathbb{C}^4$ (which is the fiber of $(\Lambda^+\oplus\underline{\mathbb{R}})\oplus T^*(S^1\times\mathbb{R}^3))\otimes\mathbb{C}$ when viewed as a vector bundle over $\mathbb{R}^3$) can be written as a direct sum of 8 product complex lines corresponding to the possible values of the respective eigenvalues $(\pm1, \pm1, \pm1)$ of the three endomorphisms in (4.7). Another point to note which is relevant later is that the endomorphism $\gamma_s$ commutes with each of the three endomorphisms in (4.7) and so it preserves their joint eigenspaces. Likewise for use later: The endomorphism $\Gamma$ anti-commutes with each of the endomorphisms in (4.7) and thus maps +1 eigenspaces to -1 eigenspaces and vice-versa.

To see the relevance of (4.7) to $\mathfrak{D}_0$, look again at the depiction of $\mathfrak{D}_0^2$ in (4.6) and notice in particular that $\mathfrak{D}_0^2$ commutes with each of the endomorphisms in (4.7). As a consequence, $\mathfrak{D}_0^2$ preserves the eigenspaces. Supposing then that $\iota_1, \iota_2, \iota_3 \in \{\pm 1\}$, then $\mathfrak{D}_0^2$ acting on maps to the corresponding eigenspace acts as the operator

$$-\sum_{k=1}^{3} \frac{\partial^2}{\partial x_k^2} + \sqrt{2}R \sum_{k=1}^{3} \lambda_k \iota_k + 2R^2 \sum_{k=1}^{3} \lambda_k^2 x_k^2 \ . \tag{4.8}$$

As is now evident, this operator is a sum of three commuting differential operators that involve pairwise distinct coordinates:

*For each* $k \in \{1, 2, 3\}$*, the operator* $Q_k \equiv -\frac{\partial^2}{\partial x_k^2} + \iota_k \sqrt{2} R \lambda_k + 2R^2 \lambda_k^2 x_k^2$ .

(4.9)

Granted all of this, then the eigenvectors and eigenvalues for $\mathfrak{D}_0^2$ can be deduced using the technique of separation of variables. The eigenvectors are the products of eigenvectors for the three operators $Q_1$, $Q_2$ and $Q_3$ with the corresponding eigenvalue being a sum of the eigenvalues for the three eigenvectors.

To say more about the spectrum and eigenfunctions of the $Q_k$'s: Borrow notation from the quantum physics of the harmonic oscillator and introduce differential operators $\{\mathbf{a}_k\}_{k=1,2,3}$ and $\{\mathbf{a}_k^\dagger\}_{k=1,2,3}$ using the rule

$$\mathbf{a}_k = \frac{\partial}{\partial x_k} + \sqrt{2} R \lambda_k x_k \quad \textit{and} \quad \mathbf{a}_k^\dagger = -\frac{\partial}{\partial x_k} + \sqrt{2} R \lambda_k x_k \ . \tag{4.10}$$

It is important to note that

- $Q_k = \mathbf{a}_k^\dagger \mathbf{a}_k + (1 + \iota_k)\sqrt{2} R \lambda_k$
- $[\mathbf{a}_k, \mathbf{a}_k^\dagger] = 2\sqrt{2} R \lambda_k$.
- *The kernel of* $\mathbf{a}_k$ *is spanned by* $f_k(x_k) = e^{-\frac{1}{\sqrt{2}} R \lambda_k x_k^2}$ .

(4.11)

It follows directly from (4.11) that the eigenfunctions of $Q_k$ are the functions in the set

$$f_k,\ \mathbf{a}_k^\dagger f_k,\ \mathbf{a}_k^\dagger \mathbf{a}_k^\dagger f_k, \ldots\ (\mathbf{a}_k^\dagger)^n f_k, \ldots \tag{4.12}$$

Also from (4.11): The $Q_k$-eigenvalue of $(\mathbf{a}_k^\dagger)^n \phi_k$ is $(1 + \iota_k + 2n)\sqrt{2} R \lambda_k$. By way of an example, the smallest eigenvalue of $Q_k$ is $(1 + \iota_k)\, 2\sqrt{2} R \lambda_k$; and this is zero if and only if $\iota_k$ is equal to -1. It follows as a consequence that the kernel of $\mathfrak{D}_0^2$ is 1-dimensional; it being the span of the product of $f_1 f_2 f_3$ times a non-zero vector in $\mathbb{C}^4 \oplus \mathbb{C}^4$ whose eigenvalues with respect to the endomorphisms in (4.7) are {-1, -1, -1}.

By way of a summary: A basis of eigenvectors of $\mathfrak{D}_0^{\,2}$ can be labeled by data sets $\{(\delta_1, \delta_2, \delta_3), (n_1, n_2, n_3)\}$ with any given $\delta_k$ such that $2\delta_k - 1$ is the eigenvalue of the corresponding $\varepsilon\gamma_k\rho_k$ from (4.7), and with any given $n_k$ being a non-negative integer. The corresponding eigenvalue for the basis vector is $2\sqrt{2}\,R\sum_{k=1,2,3}(\delta_k + n_k)\lambda_k$. This basis for $\mathfrak{D}_0^{\,2}$ will be used directly to analyze $\mathfrak{D}_0$.

With regard to the first bullet of the lemma: Note that $\psi$ is in the kernel of $\mathfrak{D}_0$ if and only if it is in the kernel of $\mathfrak{D}_0^{\,2}$, and thus the kernel($\mathfrak{D}_0$) is 1-dimensional.

With regards to the second and third bullets: Let $\psi$ now denote one of the basis elements described above, and let $\hat{E}^2$ denote its eigenvalue. The element $\psi + |\hat{E}|^{-1}\mathfrak{D}_0\psi$ is the an eigenvector of $\mathfrak{D}_0$ with eigenvalue $|\hat{E}|$; and $\psi - |\hat{E}|^{-1}\mathfrak{D}_0\psi$ is an eigenvector of $\mathfrak{D}_0$ with eigenvalue $-|\hat{E}|$. Conversely, if $\psi$ denotes an eigenvector for $\mathfrak{D}_0$ with non-zero eigenvalue $\hat{E}$, then $\gamma_s\psi$ is an eigenvector of $\mathfrak{D}_0$ with eigenvalue $-\hat{E}$; and both are eigenvector for $\mathfrak{D}_0^{\,2}$ with eigenvalue $\hat{E}^2$.

To complete the proof of Lemma 4.2: It follows from the observations in the preceding paragraphs that the basis of eigenvectors of $\mathfrak{D}_0^{\,2}$ labeled by the sets $\{(\delta_1,\delta_2,\delta_3),(n_1,n_2,n_3)\}$ with not all $\delta$'s and n's equal to zero generate (via the map $\psi \to \psi + \hat{E}^{-1}\mathfrak{D}_0\psi$) a set of eigenvectors for $\mathfrak{D}_0$ that span the eigenspaces for $\mathfrak{D}_0$ with positive eigenvalue. Granted this, (4.11) can be used to see that this map from a basis of eigenvectors of $\mathfrak{D}_0^{\,2}$ with non-zero eigenvalue to a generating set of eigenvectors of $\mathfrak{D}_0$ with positive eigenvalue is 1-1 except in the case where two pair of $(\delta,n)$'s are zero in which case the map is 2-1. For example, the data set with $(\delta_1,n_1) = (0, k)$ for $k > 0$ and $(\delta_2,n_2)$ and $(\delta_3,n_3)$ both zero and the data set with $(\delta_1,n_1) = (1,k-1)$ and $(\delta_2,n_2)$ and $(\delta_3,n_3)$ still both zero label the same eigenvector of $\mathfrak{D}_0$. (The key observation in this regard is that if $\psi$ and $\psi'$ are orthogonal eigenvectors of $\mathfrak{D}_0^{\,2}$ with the same eigenvalue $\hat{E}^2$ and if the corresponding $\mathfrak{D}_0$-eigenvectors $\psi + \hat{E}^{-1}\mathfrak{D}_0\psi$ and $\psi' + \hat{E}^{-1}\mathfrak{D}_0\psi'$ are equal, then $\psi = \hat{E}^{-1}\mathfrak{D}_0\psi'$ and $\psi' = \hat{E}^{-1}\mathfrak{D}_0\psi$. It is only in the case where two from the pairs $\{(\delta_k,n_k)\}_{k=1,2,3}$ are zero that an equality of this sort can hold between basis vectors.)

For subsequent reference (parenthetical for now): A domain for $\mathfrak{D}_0$ can be taken to be the set of elements $\psi$ in the subspace of $\mathbb{L}$ that is obtained by completing the space of smooth, compactly supported maps from $\mathbb{R}^3$ to $\mathbb{C}^4 \oplus \mathbb{C}^4$ using the norm whose square is the function given below

$$\psi \to \int_{\mathbb{R}^3} (|\nabla\psi|^2 + R^2|x|^2|\psi|^2) \ . \tag{4.13}$$

In particular, if $\psi$ is in this domain and if $\psi$ is $\mathbb{L}$-orthogonal to the kernel of $\mathfrak{D}_0$, then

$$\int_{\mathbb{R}^3} |\mathfrak{D}_0\psi|^2 \ \geq c_0^{-1} \int_{\mathbb{R}^3} (|\nabla\psi|^2 + R(1+R|x|^2)|\psi|^2) \ . \tag{4.14}$$

(This inequality follows from the identity in (4.8).)

***Proof of Lemma 4.3***: Lemma 4.5 can be invoked for the case of $\mathcal{D}_0$ because

$$\mathcal{D}_0{}^2 = -\frac{\partial^2}{\partial s^2} + \mathfrak{D}_0{}^2 \tag{4.15}$$

with $\mathfrak{D}_0{}^2$ as depicted in (4.6). (The function $\mathfrak{f}$ in this case is the function $x \to R^2|x|^2$.) That lemma establishes the claim to the effect that $\mathcal{D}_0$ has purely point spectrum with no accumulations and finite multiplicities.

Because $\frac{\partial}{\partial s}$ commutes with $\mathcal{D}_0$, one can use a Fourier decomposition with respect to the s variable to study $\mathcal{D}_0$'s eigenvalues and eigenvectors. To elaborate: Fix a positive integer n and note that $\mathcal{D}_0$ preserves the subspace of its domain that have the form $e^{2\pi i n s/\ell}\,\phi$ with $\phi$ being independent of s. This is an eigenvector of $\mathcal{D}_0$ if and only if

$$(2\pi i n/\ell)\gamma_s\phi + \mathfrak{D}_0\phi = \lambda\phi \quad \textit{with} \;\; \lambda \in \mathbb{R}. \tag{4.16}$$

If $\phi$ is in the kernel of $\mathfrak{D}_0$, then $i\gamma_s\phi$ is equal to $-\phi$ (taking $\gamma_s$ to have eigenvalue i on the kernel of $\mathfrak{D}_0$). Thus, $e^{2\pi i n s/\ell}\,\phi$ is an eigenvector of $\mathcal{D}_0$ with eigenvalue $-2\pi n/\ell$. Now suppose that $\phi$ is $\mathbb{L}$-orthogonal to the kernel of $\mathfrak{D}_0$ and that $e^{2\pi i n s/\ell}\,\phi$ is an eigenvector of $\mathcal{D}_0$. It must also be one for $\mathcal{D}_0{}^2$ which requres that $\phi$ be an eigenvector for $\mathfrak{D}_0{}^2$. If $\hat{E}$ is positive and its square is the corresponding eigenvalue of $\mathfrak{D}_0{}^2$, then $\phi$ must have the form $\alpha\varphi + \beta\gamma_s\varphi$ with $\alpha$ and $\beta$ being complex numbers and with $\varphi$ obeying $\mathfrak{D}_0\varphi = \hat{E}\varphi$. With regards to $\alpha$ and $\beta$: The pair $(\alpha, \beta)$ solves the coupled system of linear equations

$$-(2\pi i n/\ell)\beta + \hat{E}\alpha = E\alpha \quad \textit{and} \quad (2\pi i n/\ell)\alpha - \hat{E}\beta = E\beta \tag{4.17}$$

The eigenvalue $E$ here is $\pm 1$ times the square root of $(4\pi^2 n^2/\ell^2 + \hat{E}^2)$ which is the requirement for (4.15) to have a non-trivial solution. Granted one of these values for $E$, then $\alpha$ and $\beta$ are uniquely determined by (4.15); and they have the form specified by the second bullet of Lemma 4.3.

***Proof of Lemma 4.4***: The first observation concerns $\mathcal{D}_0{}^2$ which is that it has the form

$$-\frac{\partial^2}{\partial s^2} + \mathfrak{D}_0{}^2 + \sqrt{2}\, i R \gamma_s \dot{\mathbb{M}}(x,\rho) \tag{4.18}$$

where $\dot{\mathbb{M}}(x,\rho)$ denotes here (and subsequently) the sum $\sum_{j,k=1}^{3} (\frac{d}{ds}\mathbb{M}_{jk}(s))\, x_j\rho_k$

It follows from (4.18) and (4.6) that the prerequisites for Lemma 4.5 are met using $\mathcal{D}_0$ with X being $\mathbb{R}/(\ell\mathbb{Z})\times\mathbb{R}^3$ and with $\mathfrak{f}$ denoting the function $(s, x) \to R^2|x|^2$ The conclusions of Lemma 4.5 are the assertions in Lemma 4.4 to the effect that $\mathcal{D}_0$ is an unbounded, self adjoint operator having purely point spectrum with no accumulations and finite multiplicities.

To say more about the spectrum of $\mathcal{D}_0$, digress for the moment for some remarks about the line bundle $\mathcal{K}$: The fiber-wise $\mathbb{L}$-orthogonal projection from $\mathbb{R}/(\ell\mathbb{Z})\times\mathbb{L}$ to $\mathcal{K}$ is the linear map sending any given section $\psi$ of the product bundle to $\langle\phi_0,\psi\rangle_{\mathbb{L}}\,\phi_0$. This endomorphism of the product $\mathbb{L}$-bundle over $\mathbb{R}/(\ell\mathbb{Z})$ is denoted by $\Pi$. The fiberwise projection in $\mathbb{R}/(\ell\mathbb{Z})\times\mathbb{L}$ orthogonal to $\mathcal{K}$ is $(1-\Pi)$. By way of an application: Differentiate the identity $\mathfrak{D}_0\phi_0 = 0$ with respect to the parameter s to see that the s-derivative of $\phi_0$ obeys the equation

$$\mathfrak{D}_0(\tfrac{\partial}{\partial s}\phi_0) + \sqrt{2}\,i\,R\,\dot{\mathbb{M}}(x,\rho)\,\phi_0 = 0 \tag{4.19}$$

As explained below, $\Pi(x_j\rho_k\phi_0) = 0$ which implies that $\frac{\partial}{\partial s}\phi_0$ can be written as

$$\tfrac{\partial}{\partial s}\phi_0 = i\,\alpha_*\,\phi_0 - \sqrt{2}iR\big(\mathfrak{D}_0^{-1}(\dot{\mathbb{M}}(x,\rho)\,\phi_0)\big) \tag{4.20}$$

with $\alpha_*$ denoting a real valued function on $\mathbb{R}/(\ell\mathbb{Z})$. In the preceding formula and henceforth in this proof, $\mathfrak{D}_0^{-1}$ should be viewed as a bounded, linear map from $(1-\Pi)\mathbb{L}$ to $(1-\Pi)\mathbb{L}$ (remember that $(1-\Pi)L$ is the orthogonal complement in $\mathbb{L}$ of the kernel of $\mathfrak{D}_0$). With regards to this function $\alpha_*$: No generality is lost at this point by assuming that $\alpha_*$ is constant with its value in the interval $[0, 2\pi/\ell)$. This is because the section $\phi_0$ can be changed by multipling it by a map from $\mathbb{R}/(\ell\mathbb{Z})$ to U(1) to obtain this outcome.

Now suppose that $\psi$ is an eigenvector for $\mathcal{D}_0$ with eigenvalue $E$ obeying the upper bound $E \le \frac{1}{c}\sqrt{R}$. Here, $c$ is a number greater than 1 with a lower bound to be determined momentarily (the lower bound is ultimately determined by the infimum over s of the minumum of the absolute values of eigenvalues of $\mathbb{M}|_s$, and by the supremum over s of the norm of the s-derivative of $\mathbb{M}|_s$.) Write $\psi = \mu\,\phi_0 + \psi^\perp$ with $\mu$ being a $\mathbb{C}$-valued function on $\mathbb{R}/(\ell\mathbb{Z})$ and with $\psi^\perp = (1-\Pi)\psi$. The section $\psi^\perp$ obeys

$$(1-\Pi)\,\mathcal{D}_0\psi^\perp + \mu\,\gamma_s(1-\Pi)\tfrac{\partial}{\partial s}\phi_0 = E\,\psi^\perp \tag{4.21}$$

which is the $(1-\Pi)$ projection of the eigenvalue equation $\mathcal{D}_0\psi = E\psi$. What with (4.20) and with $\gamma_s$ acting as i on $\phi_0$, the preceding equation asserts that

$$(1-\Pi)(\mathcal{D}_0-\text{E})\psi^{\perp} = -\sqrt{2}\text{R}\mu\,\mathfrak{D}_0^{-1}(\dot{\mathbb{M}}(x,\rho)\,\phi_0)\ . \tag{4.22}$$

As explained below, if $c > c_0$, then this equation uniquetly determines $\psi^{\perp}$ from $\mu$. Briefly for now: This is because the absolute value of the smallest eigenvalue of $(1-\Pi)\,\mathcal{D}_0(1-\Pi)$ acting on sections of $(1-\Pi)\mathbb{L}$ over $\mathbb{R}/(\ell\mathbb{Z})$ is greater than $c_0^{-1}\sqrt{R}$. As a consequence, if $c > c_0$, then the operator $(1-\Pi)(\mathcal{D}_0-\text{E})(1-\Pi)$ is invertible on the space of sections over $\mathbb{R}/(\ell\mathbb{Z})$ of the vector bundle $(1-\Pi)\mathbb{L}$. Supposing that this is so, then (4.22) says that

$$\psi^{\perp} = -\sqrt{2}\,\text{R}\big((1-\Pi)(\mathcal{D}_0-\text{E})(1-\Pi)\big)^{-1}\big(\mu\,\mathfrak{D}_0^{-1}(\dot{\mathbb{M}}(x,\rho)\,\phi_0)\big) \tag{4.23}$$

Meanwhile, $\mu$ must obey the equation below (which is the projection to $\mathcal{K}$ of the eigenvalue equation $\mathcal{D}_0\psi = \text{E}\psi$):

$$\tfrac{d}{ds}\mu + i\alpha_*\mu + \langle\phi_0, \tfrac{\partial}{\partial s}\psi^{\perp}\rangle_{\mathbb{L}} = -i\text{E}\mu \tag{4.24}$$

(Remember that the eigenvalue of $\gamma_5$ on $\mathcal{K}$ is taken here to be +i.) The chain of identities below is used momentarily rewrite (4.24). What follows is the promised chain of identities.

- $\langle\phi_0, \tfrac{\partial}{\partial s}\psi^{\perp}\rangle_{\mathbb{L}} = -\langle\tfrac{\partial}{\partial s}\phi_0, \psi^{\perp}\rangle_{\mathbb{L}}$ ,
- $\langle\phi_0, \tfrac{\partial}{\partial s}\psi^{\perp}\rangle_{\mathbb{L}} = -\sqrt{2}\,i\text{R}\langle\mathfrak{D}_0^{-1}(\dot{\mathbb{M}}(x,\rho)\,\phi_0), \psi^{\perp}\rangle_{\mathbb{L}}$
- $\langle\phi_0, \tfrac{\partial}{\partial s}\psi^{\perp}\rangle_{\mathbb{L}} = 2i\text{R}^2\big\langle(\mathfrak{D}_0^{-1}(\dot{\mathbb{M}}(x,\rho)\,\phi_0), ((1-\Pi)(\mathcal{D}_0-\text{E})(1-\Pi))^{-1}\big(\mu\,\mathfrak{D}_0^{-1}(\dot{\mathbb{M}}(x,\rho)\,\phi_0)\big)\big\rangle_{\mathbb{L}}$

(4.25)

By way of an explanation: The first identity follows because $\langle\phi_0,\psi^{\perp}\rangle_{\mathbb{L}} = 0$; the second identity follows from the first with the help of (4.20); and the third follows from the second with the help of (4.23).

Now a key point with regards to the last identity in (4.25) (as explained below) is that the $\mathbb{R}/(\ell\mathbb{Z})$-integral of the norm of what is depicted on the right hand side of the third bullet of (4.25) is bounded by $c_0\,\frac{1}{\sqrt{R}}\,(\int_{\mathbb{R}/\ell\mathbb{Z}} |\mu|^2\,)^{1/2}$ (supposing that $\text{E} < c^{-1}\sqrt{R}$ and supposing that $c > c$). Granted this last bound, use it in (4.24) to see that

$$\int_{\mathbb{R}/\ell\mathbb{Z}} |\tfrac{d}{ds}(e^{-i(\alpha_*+\text{E})s}\mu)| \le c_0\,\frac{1}{\sqrt{R}}\,(\int_{\mathbb{R}/\ell\mathbb{Z}} |\mu|^2\,)^{1/2}\ . \tag{4.26}$$

A first implication of (4.26) is that the function on $\mathbb{R}/\ell\mathbb{Z}$ given by $|\mu|$ differs from a constant by at most $c_0 \frac{1}{\sqrt{R}}$. (That constant is at most 1 because the $\mathbb{R}/(\ell\mathbb{Z})$-integral of $|\mu|^2$ is no greater than the integral of $|\psi|^2$ over $\mathbb{R}/(\ell\mathbb{Z}) \times \mathbb{R}^3$which is 1.) With this first implication in hand, then a second implication is that $e^{-i(\alpha_*+E)s}\mu$ differs from a constant (to be denoted by $\mu_0$) by at most $c_0 \frac{1}{\sqrt{R}}$. And that bound implies in turn that $\mu$ differs from $\mu_0 e^{-i(\alpha_*+E)s}$ by at most $c_0 \frac{1}{\sqrt{R}}$. Meanwhile, the bound

$$|\mu - \mu_0 e^{-i(\alpha_*+E)s}| \le c_0 \frac{1}{\sqrt{R}} \tag{4.27}$$

can hold only if $E + \alpha_*$ differs from an integer multiple of $2\pi/\ell$ by at most $c_0 \frac{1}{\sqrt{R}}$. This is one of the assertions of Lemma 4.4 if the number $\alpha$ is taken to be $\frac{\ell}{2\pi}\alpha_*$.

Another key point (also explained below) is that the square of the $\mathbb{L}$-norm of $\psi^\perp$, which is the $\mathbb{R}/(\ell\mathbb{Z})$-integral of the square of the fiberwise $\mathbb{L}$ norm of $\psi^\perp$, is bounded by $c_0 \frac{1}{R} \int_{\mathbb{R}/\ell\mathbb{Z}} |\mu|^2$ when $c > c_0$. This last point with (4.27) complete the proof of Lemma 4.4 except for the promised justifications for the following three claims: The first is the claim that the operator $(1-\Pi)(\mathcal{D}_0 - E)(1-\Pi)$ is invertible; the second claim is the asserted bound for the $\mathbb{R}/(\ell\mathbb{Z})$-integral of the right hand side of the third bullet in (4.25); and the third claim is the asserted bound for the $\mathbb{R}/(\ell\mathbb{Z})$-integral of the square of the $\mathbb{L}$-norm of $\psi^\perp$. These justifications are in the subsequent paragraphs.

To start the justification: It follows from Lemma 4.2 that any given incarnation of the operator $\mathfrak{D}_0$ is invertible on the complement of its kernel. More to the point, it follows from this lemma (or more directly, (4.14)) that for any $q \in (1-\Pi)\mathbb{L}$,

$$\|\mathfrak{D}_0^{-1} q\|_{\mathbb{L}} \le c_0^{-1} \frac{1}{\sqrt{R}} \|q\|_{\mathbb{L}} \tag{4.28}$$

with $\|\cdot\|_{\mathbb{L}}$ denoting here the norm on $\mathbb{L}$. In the instance when $q = \dot{\mathbb{M}}(x,\rho)\,\phi_0$, this gives

$$\|\mathfrak{D}_0^{-1}(\dot{\mathbb{M}}(x,\rho)\,\phi_0)\|_{\mathbb{L}} \le c_0 \frac{1}{R}\,. \tag{4.29}$$

To exploit this bound, note next that if $u$ is a smooth section over $\mathbb{R}/(\ell\mathbb{Z})$ of the subbundle in the product $\mathbb{L}$ bundle with fiber $(1-\Pi)\mathbb{L}$, then

$$(1-\Pi)\mathcal{D}_0 u = \mathcal{D}_0 u - \phi_0\langle\phi_0, \gamma_s \tfrac{\partial}{\partial s} u\rangle\,; \tag{4.30}$$

which is to say that

$$(1-\Pi)\mathcal{D}_0 u = \mathcal{D}_0 u + i\,\phi_0\langle \tfrac{\partial}{\partial s}\phi_0, u\rangle \,.$$

(4.31)

What with (4.29), the preceding identity leads directly to the bound

$$\int_{\mathbb{R}/(\ell\mathbb{Z})} \|((1-\Pi)\mathcal{D}_0 u\|_{\mathbb{L}}^2 \geq c_0^{-1} \int_{\mathbb{R}/(\ell\mathbb{Z})} \|\mathcal{D}_0 u\|_{\mathbb{L}}^2 - c_0 \int_{\mathbb{R}/(\ell\mathbb{Z})} \|u\|_{\mathbb{L}}^2 \,.$$

(4.32)

Meanwhile, the Bochner-Weitzenboch formula in (4.16) and (4.14) imply that

$$\int_{\mathbb{R}/(\ell\mathbb{Z})} \|\mathcal{D}_0 u\|_{\mathbb{L}}^2 \geq c_0^{-1} \int_{\mathbb{R}/(\ell\mathbb{Z})} (\|\tfrac{\partial}{\partial s} u\|_{\mathbb{L}}^2 + \|\nabla u\|_{\mathbb{L}}^2 + R(1+R|\cdot|^2)\|u\|_{\mathbb{L}}^2)$$

(4.33)

when $u$ is a smooth section over $\mathbb{R}/\ell\mathbb{Z}$ of the bundle with fiber $(1-\Pi)\mathbb{L}$. Together with (4.32), this last bound implies the following: If $R \geq c_0$, then the operator $(1-\Pi)\mathcal{D}_0$ is invertible on the space of sections of $(1-\Pi)\mathbb{L}$ with its inverse obeying

$$\int_{\mathbb{R}/(\ell\mathbb{Z})} \|((1-\Pi)\mathcal{D}_0(1-\Pi))^{-1} q\|_{\mathbb{L}}^2 \leq c_0^{-1}\frac{1}{R} \int_{\mathbb{R}/(\ell\mathbb{Z})} \|q\|_{\mathbb{L}}^2 \,.$$

(4.34)

The preceding bound leads directly to the following obervation: If $E$ is a number with norm at most $c_0^{-1}\sqrt{R}$, then $(1-\Pi)(\mathcal{D}_0 - E)$ is also invertible on the space of sections of $(1-\Pi)\mathbb{L}$ with its inverse obeying the bound in (4.32) albeit with a larger value for $c_0$. This last observation with (4.29) leads to the $c_0\frac{1}{\sqrt{R}}\big(\int_{\mathbb{R}/\ell\mathbb{Z}} |\mu|^2\big)^{1/2}$ bound that was asserted for the $\mathbb{R}/(\ell\mathbb{Z})$-integral of the expression on the right hand side of the inequality in the third bullet of (4.25). It also leads directly to a $c_0\frac{1}{\sqrt{R}}\big(\int_{\mathbb{R}/\ell\mathbb{Z}} |\mu|^2\big)^{1/2}$ bound for the $\mathbb{R}/(\ell\mathbb{Z})$-integral of $\|\psi^{\perp}\|_{\mathbb{L}}^2$.

The arguments above explain why the eigenvalues of $\mathcal{D}_0$ with norm at most $\frac{1}{c_0}\sqrt{R}$ can be written in the manner asserted by the lemma, as $(-\alpha+n)/\frac{2\pi}{\ell} + \mathfrak{r}$ with $\alpha$ from $[0,2\pi)$, with n being an integer and with $\mathfrak{r}$ obeying $|\mathfrak{r}| \leq c_0\frac{1}{\sqrt{R}}$. The paragraphs that follow give the argument for the converse which is that if $R \geq c_0$, then any given integer (call it n) with norm at most $\frac{1}{c_0}\sqrt{R}$ corresponds to an eigenvalue of $\mathcal{D}_0$ that differs by at most $c_0\frac{1}{\sqrt{R}}$ from $-(\alpha+n)\frac{2\pi}{\ell}$. (The argument is an application of standard perturbative constructions (see e.g [K] applied to the context here. Looking ahead, the key input for the proof are (4.23)–(4.25) and what is said subsequently about these equations.)

To start the argument, fix a number $c$ to be greater than $c_0$ and then fix an integer (to be denoted by n) with norm less than $c^{-1}\sqrt{R}$. The plan is to construct simultaneously an eigenvalue for $\mathcal{D}_0$ having the form $-(\alpha+n)\frac{2\pi}{\ell}+\varepsilon$ with $\varepsilon$ small and a corresponding eigenvector, $\psi = \mu\phi_0+\psi^{\perp}$ with $\mu = e^{2\pi i n s/\ell}+w$ where w also has small norm.

To give the details, fix for the moment a number $\varepsilon$ with norm at most $c^{-1}\sqrt{R}$ and a funtion w with $|w| < c^{-1}$ and with the $\mathbb{R}/\ell\mathbb{Z}$ integral of its derivative also less than $c^{-1}$. Assume in addition that the $\mathbb{R}/(\ell\mathbb{Z})$ integral of $e^{-2\pi i n s/\ell}w$ is zero. Having done all of this, then define a section (denoted by $\psi^{\perp}$) of $Y\times\mathbb{L}$ using the formula in (4.23) with $\mu = e^{2\pi i n s/\ell}+w$. With $\psi^{\perp}$ so defined, then the section $\psi = (e^{2\pi i n s/\ell}+w)\phi_0+\psi^{\perp}$ will be an eigenvector of $\mathcal{D}_0$ and E will be its eigenvalue if (4.24) is obeyed; thus $(w, \varepsilon)$ are such that

$$\frac{d}{ds}w - i\frac{2\pi}{\ell}nw + \langle\phi_0, \frac{\partial}{\partial s}\psi^{\perp}\rangle_{\mathbb{L}} = -i\varepsilon w - i\varepsilon e^{2\pi i n s/\ell}\,. \tag{4.35}$$

To see about a solution to (4.35) with small w and $\varepsilon$, let $L^2(\mathbb{R}/(\ell\mathbb{Z}))$ denote the Hilbert space completion of the space of $\mathbb{C}$-valued functions on $\mathbb{R}/(\ell\mathbb{Z})$ with respect to the norm that assigns to a function (call it $f$) the $\mathbb{R}/(\ell\mathbb{Z})$ integral of $|f|^2$. Also: Let $P_n$ denote the orthogonal projection in $L^2(\mathbb{R}/(\ell\mathbb{Z}))$ to the span of the function $e^{2\pi i n s/\ell}$; and let $D_n$ denote the operator $\frac{d}{ds} - i\frac{2\pi}{\ell}n$. The key point to keep in mind is that $D_n$ is invertible on the image of $(1-P_n)$ in $L^2(\mathbb{R}/(\ell\mathbb{Z}))$ and if $P_n y = 0$, then

$$\int_0^{\ell}\Big(\frac{1}{(1+n^2)}\big|\frac{d}{ds}(D_n^{-1}y)\big|^2 + |D_n^{-1}y|^2\Big) \le \frac{2\pi}{\ell}\int_0^{\ell}|y|^2\,. \tag{4.36}$$

The operator $D_n$ is introduced so as to view (4.35) as a fixed point equation for the pair $(w, \varepsilon)$:

- $w = -D_n^{-1}\big((1-P_n)\langle\phi_0, \frac{\partial}{\partial s}\psi^{\perp}\rangle_{\mathbb{L}} + i\varepsilon w\big)$.
- $\varepsilon = i\frac{1}{\ell}\int_0^{\ell} e^{-2\pi n s/\ell}\langle\phi_0, \frac{\partial}{\partial s}\psi^{\perp}\rangle$.

(4.37)

A Banach space contraction mapping argument will be used to prove that (4.37) has a unique, small norm fixed point when $c > c_0$. The relevant Banach space is the completion of the product of $C^{\infty}(\mathbb{R}/(\ell\mathbb{Z}))\times\mathbb{R}$ using the norm whose square on any given $(w, \varepsilon)$ is the number

$$\int_0^{\ell} (|\tfrac{d}{ds} w|^2 + |w|^2) + \varepsilon^2 \,.$$

(4.38)

This norm is denoted below by $\|\cdot\|_{\ddagger}$. The Banach space is denoted by $\mathbb{B}_{\ddagger}$. As explained below, the right hand side of (4.37) (when viewed as a function of inputs w and ε) defines a contraction mapping from a small radius ball in $\mathbb{B}$ to itself when $c > c_0$.

With the fixed point plan in mind, return to (4.37) and its right hand side. With regards to the dependence of $\langle \phi_0, \frac{\partial}{\partial s} \psi^{\perp}\rangle_{\mathbb{L}}$ on w: This is via $\psi^{\perp}$ in (4.23); and the first thing to note is that $\psi^{\perp}$ has linear dependence on $\mu$ which implies in this case that it is a sum of the $\mu = e^{2\pi i n s/\ell}$ version of (4.23) and the $\mu$ = w version; and thus, so is $\langle \phi_0, \frac{\partial}{\partial s} \psi^{\perp}\rangle_{\mathbb{L}}$ in (4.37). In this regard: The $\mathbb{R}/((\ell\mathbb{Z})$-integral of the $\mu = e^{2\pi i n s/\ell}$ version of $|\langle \phi_0, \frac{\partial}{\partial s} \psi^{\perp}\rangle_{\mathbb{L}}|^2$ is bounded by $c_0 \frac{1}{R}$. Meanwhile, the $\mu$ = w version of $\langle \phi_0, \frac{\partial}{\partial s} \psi^{\perp}\rangle$ has linear dependence on w; and, for any given w, the $\mathbb{R}/((\ell\mathbb{Z})$-integral of its square is at most

$$c_0 \frac{1}{R} \int_{\mathbb{R}/\ell\mathbb{Z}} |w|^2 \,.$$

(4.39)

(All of these bounds follow from what is said in and directly after (4.25).)

With regards to the ε dependence of $\langle \phi_0, \frac{\partial}{\partial s} \psi^{\perp}\rangle_{\mathbb{L}}$: Let $E(\varepsilon)$ denote $-(\alpha+n)\frac{2\pi}{\ell} + \varepsilon$. The key observation is that

$$(1-\Pi)\big((\mathcal{D}_0 - E(\varepsilon))^{-1} - (\mathcal{D}_0 - E(\varepsilon'))^{-1}\big)(1-\Pi) = (\varepsilon - \varepsilon')(1-\Pi)(\mathcal{D}_0 - E(\varepsilon))^{-1}(\mathcal{D}_0 - E(\varepsilon'))^{-1}(1-\Pi) \,,$$

(4.40)

which has two implications which are these: Let $\psi_{\varepsilon}^{\perp}$ denote the $E = E(\varepsilon)$ version of $\psi^{\perp}$. Then the inequalities below in (4.41) hold for any pair ε and ε´ such that both are bounded by $\frac{1}{c}\sqrt{R}$ with it understood that a version of w is used whose square has $\mathbb{R}/(\ell\mathbb{Z})$ bounded by 1 (any a priori bound for the norm of w is acceptable at the expense of changing the version of $c_0$ that appears below).

- $\int_{\mathbb{R}/(\ell\mathbb{Z})} \langle \phi_0, \frac{\partial}{\partial s} \psi_0^{\perp}\rangle_{\mathbb{L}}^2 \le c_0 \frac{1}{R}$.
- $\int_{\mathbb{R}/(\ell\mathbb{Z})} (\langle \phi_0, \frac{\partial}{\partial s} \psi_{\varepsilon}^{\perp}\rangle_{\mathbb{L}} - \langle \phi_0, \frac{\partial}{\partial s} \psi_{\varepsilon'}^{\perp}\rangle_{\mathbb{L}})^2 \le c_0 \frac{1}{R^2} |\varepsilon - \varepsilon'|$.

(4.41)

These bounds with (4.39), with what is said about the w-dependence of $\langle \phi_0, \frac{\partial}{\partial s} \psi^{\perp}\rangle_{\mathbb{L}}$ in the paragraph preceding (4.39), and with (4.36) have the following two implications when $c > c_0$: The first is that the right hand side of (4.37) is a contraction mapping from the radius 1 ball about the origin in $\mathbb{B}_{\ddagger}$ to itself when $c > c_0$. Thus, if $c > c_0$, then there is a unique fixed point in this ball. The second implication concerns the norm of this fixed point, which is that

$$|\varepsilon| \le c_0 \frac{1}{\sqrt{R}} \quad and \quad \int_0^{\ell} \left(\frac{1}{(1+n^2)} |\frac{d}{ds} w|^2 + |w|^2\right) \le c_0 \frac{1}{R} . \tag{4.42}$$

As noted previously, this fixed point gives the desired eigenvalue and eigenvector to complete the proof of the lemma.

**c) Some generalizations of the model case**

This subsection considers two generalizations of the model $\mathcal{D}_0$ operator that is considered by Lemma 4.4. The descripion the first generalization requires the specification of a map $s \to M(s)$ from $\mathbb{R}/(\ell\mathbb{Z})$ to $\mathbb{R}^3$ and a map $s \to W(s)$ from $\mathbb{R}/(\ell\mathbb{Z})$ to the vector space of anti-symmetric tensors on $\mathbb{R}^3$. The components of $M(\cdot)$ and $W(\cdot)$ are functions on $\mathbb{R}/(\ell\mathbb{Z})$ that are denoted respectively by $\{M_j\}_{j=1,2,3}$ and $\{W_{jk}\}_{j,k=1,2,3}$. Also needed is the choice of a small positive number to be denoted by $r_0$, chosen so that $r_0(|M|+|W|) < \frac{1}{10{,}000}$. With this number chosen, introduce by way of notation $\chi_0$ to denote the function on $\mathbb{R}^3$ given by $\chi(\frac{|x|}{r_0} - 1)$. To be sure: This function is equal to 1 where $|x| \le r_0$ and it is equal to zero where $|x| \ge 2r_0$. Given this data, the new model operator is

$$\mathcal{D}_1 \equiv \gamma_s\left((1 - \chi_0 x_j M_j)\frac{\partial}{\partial s} - \tfrac{1}{2}\chi_0 x_j \dot{M}_j - \chi_0 W_{ij} x_j \frac{\partial}{\partial x_i}\right) + \mathfrak{D}_0 \tag{4.43}$$

with $\mathfrak{D}_0$ as before. Here, $\dot{M}_j$ signifies the s-derivative of $M_j$. Also, as before, repeated indices are summed over the set $\{1, 2, 3\}$. The use of $\chi_0$ with this bound on $r_0$ is to guarantee that $\mathcal{D}_1$ is an elliptic operator. (This operator is formally self-adjoint; this is due to the appearance of the term with $\dot{M}$ and by the fact that the matrix $W$ is anti-symmetric.)

The lemma that follows describes the spectrum of $\mathcal{D}_1$ near zero; it says in effect that the conclusions of Lemma 4.4 hold for $\mathcal{D}_1$ also.

**Lemma 4.6**: *The operator $\mathcal{D}_1$ extends to the Hilbert space $\underline{\mathbb{L}}$ as an unbounded, self-adjoint operator with discrete spectrum having no accumulation points and finite multiplicities. Moreover, there exists $\kappa > 1$, a number $\alpha \in [0, 1)$ and a unit normed* *section $\phi_0$ of $\mathcal{K}$ such that if* $R \geq \kappa$*, then an eigenvalue of* $\mathcal{D}_1$ *with absolute value at most* $\frac{1}{\kappa}\sqrt{R}$ *can be written as*

$$E = -(\alpha + n)\tfrac{2\pi}{\ell} + \mathfrak{r}$$

*with* n *denoting an integer and with* $\mathfrak{r}$ *denoting a real number with norm bounded by* $\frac{1}{\sqrt{R}}\kappa$*. Conversely, if* n *is an integer with norm bounded by* $\frac{1}{\kappa}\sqrt{R}$*, then there is a corresponding eigenvalue of* $\mathcal{D}_1$ *that can be written as above with* $|\mathfrak{r}| \leq \frac{1}{\sqrt{R}}\kappa$*. In this case, the corresponding eigenvector differs from* $e^{2\pi i n s/\ell}\phi_0$ *by map from* $\mathbb{R}/(\ell\mathbb{Z})\times\mathbb{R}^3$ *to* $\mathbb{C}^4\oplus\mathbb{C}^4$ *with a* $\frac{1}{R}\kappa$ *bound on the* $\mathbb{R}/(\ell\mathbb{Z})\times\mathbb{R}^3$ *integral of the square of its norm. As for the number* $\kappa$*: An upper bound for* $\kappa$ *is a priori determined given upper bounds for the norms along* $\mathbb{R}/(\ell\mathbb{Z})$ *of* $\mathbb{M}^{-1}$ *and also* $\mathbb{M}$, M, W *and their* s*-derivatives.*

This lemma is proved momentarily.

The second generalization requires the specification of two more maps from $\mathbb{R}/(\ell\mathbb{Z})$ to $\mathbb{R}^3$, these denoted by B and C with respective components by $\{B_j\}_{j=1,2,3}$ and $\{C_j\}_{j=1,2,3}$. The second new operator is

$$\mathcal{D}_2 = \gamma_s\left((1 - \chi_0 x_j M_j)\tfrac{\partial}{\partial s} - \tfrac{1}{2}\chi_0 x_j \dot{M}_j - \chi_0 W_{ij} x_j \tfrac{\partial}{\partial x_i} + 2 i x_i B_i\right) + \mathfrak{D}_C \tag{4.44}$$

where $\mathfrak{D}_C$ at any given $s \in \mathbb{R}/(\ell\mathbb{Z})$ designates the following operator:

$$\mathfrak{D}_C = \gamma_i\left(\tfrac{\partial}{\partial x_i} + i x_j \varepsilon_{jik} C_k\right) + \sqrt{2}\, i R \mathbb{M}_{ik} x_i \rho_k\,. \tag{4.45}$$

Said differently: $\mathfrak{D}_C$ is $\mathfrak{D}_0 + i\gamma_i x_j \varepsilon_{jik} C_k$. (Here again, repeated indices are summed over the set {1, 2, 3}. Also, $\{\varepsilon_{ijk}\}_{i,j,k\in\{1,2,3\}}$ are as before, the components of the completely anti-symmetric 3-tensor with $\varepsilon_{123} = 1$.)

**Lemma 4.7**: *The operator* $\mathcal{D}_2$ *extends to the Hilbert space* $\underline{\mathbb{L}}$ *as an unbounded, self-adjoint operator with discrete spectrum having no accumulation points and finite multiplicities. Moreover, there exists* $\kappa > 1$*, a number* $\alpha \in [0, 1)$ *and a unit normed section* $\phi_0$ *of* $\mathcal{K}$ *such that if* $R \geq \kappa$*, and if* $\sup_{\mathbb{R}/\ell\mathbb{Z}}(|B| + |C|) < \frac{1}{\kappa}R$*, then any eigenvalue of* $\mathcal{D}_2$ *with absolute value at most* $\frac{1}{\kappa^2}\sqrt{R}$ *can be written as*

$$E = -(\alpha + n)\tfrac{2\pi}{\ell} + \mathfrak{r}$$

*with* n *denoting an integer and with* $\mathfrak{r}$ *denoting a real number with norm bounded by* $\frac{1}{\sqrt{R}}\kappa$*. Conversely, if* n *is an integer with norm bounded by* $\frac{1}{\kappa}\sqrt{R}$*, then there is a* *corresponding eigenvalue of* $\mathcal{D}_2$ *that can be written as above with* $|\mathfrak{r}| \le \frac{1}{\sqrt{R}}\,\kappa$. *The corresponding eigenvector differs from* $e^{2\pi i n s/\ell}\,\phi_0$ *by map from* $\mathbb{R}/(\ell\mathbb{Z})\times\mathbb{R}^3$ *to* $\mathbb{C}^4\oplus\mathbb{C}^4$ *with a* $\frac{1}{R}\,\kappa$ *bound on the* $\mathbb{R}/(\ell\mathbb{Z})\times\mathbb{R}^3$ *integral of the square of its norm. As for the number* $\kappa$: *An upper bound for* $\kappa$ *is a priori determined given upper bounds for the norms along* $\mathbb{R}/(\ell\mathbb{Z})$ *of* $\mathbb{M}^{-1}$, $\mathbb{M}$, $M$, $W$ *and their respective derivative with respect to the parameter* $s$.

The rest of this subsection contains the proofs of these two lemmas.

***Proof of Lemma 4.6***: The proof that follows is much like the proof of Lemma 4.4. The first part of the proof explains why the eigenvalues with norm at most $c_0^{-1}\sqrt{R}$ and their eigenvectors have the required form. The second part of the proof constructs an eigenvalue of the form $-(\alpha+n)\frac{2\pi}{\ell} + \mathfrak{r}$ with $|\mathfrak{r}| \le c_0\frac{1}{\sqrt{R}}$ from any given integer $n$ with a norm less than $c_0^{-1}\sqrt{R}$.

To start the first part of the proof: The operator $\mathcal{D}_1^{\,2}$ has the schematic depiction

$$\mathcal{D}_1^{\,2} = -\nabla_s^{\,2} + \mathfrak{D}_0^{\,2} + \mathfrak{p}_s\frac{\partial}{\partial s} + \mathfrak{p}_j\frac{\partial}{\partial x_j} + R\,x_j\,\mathfrak{q}_j \tag{4.46}$$

where the notation uses $\nabla_s$ as shorthand for $(1 - \chi_0 x_j M_j)\frac{\partial}{\partial s} - \frac{1}{2}\chi_0 x_j \dot{M}_j - \chi_0 W_{ij}\, x_j\frac{\partial}{\partial x_i}$, and where it uses $\mathfrak{p}_s$ and each $\mathfrak{p}_j$ and $\mathfrak{q}_j$ to denote endomorphisms of the product $\mathbb{C}^4\oplus\mathbb{C}^4$ bundle over $\mathbb{R}/(\ell\mathbb{Z})$ with norms bounded by $c_0$. It follows as a consequence of (4.46) and (4.6) that the operator $\mathcal{D}_1$ meets the prerequisites for an appeal to Lemma 4.5 (using again $\mathbb{R}/(\ell\mathbb{Z})\times\mathbb{R}^3$ for X and the function $(s, x) \to R^2|x|^2$ for $\mathfrak{f}$). The consequences of that appeal are the conclusions that $\mathcal{D}_1$ has purely point spectrum with no accumulations and finite multiplicities.

To say more about the eigenvalues of $\mathcal{D}_1$: Supposing that $\psi$ is an eigenvector for $\mathcal{D}_1$ and $E$ is its eigenvalue, write $\psi$ as $\mu\phi_0 + \psi^\perp$ with $\psi^\perp = (1-\Pi)\psi$. Here as previously, $\phi_0$ is a section of the line bundle $\mathcal{K}$ (whose fiber is the kernel of $\mathfrak{D}_0$) whose $\mathbb{L}$ norm on each fiber is 1. Meanwhile, $\Pi$ is the fiberwise $\mathbb{L}$-orthogonal projection to $\mathcal{K}$. The $(1-\Pi)$ projection of the equation $\mathcal{D}_1\psi = E\psi$ has the following form

$$(1-\Pi)\,\mathcal{D}_1\psi^\perp + i\,\mu(1-\Pi)\Big((1 - \chi_0 x_j M_j)\frac{\partial}{\partial s}\phi_0 - \frac{1}{2}\chi_0 x_j \dot{M}_j\phi_0 - \chi_0 W_{ij}\, x_j\frac{\partial}{\partial x_i}\phi_0\Big)$$
$$- i\big(\tfrac{d}{ds}\mu\big)(1-\Pi)\,\chi_0\, x_j M_j\,\phi_0 = E\psi^\perp. \tag{4.47}$$

The argument used in the proof of Lemma 4.4 to show that $(1-\Pi)\mathcal{D}_0(1-\Pi)$ is invertible can be repeated with only marginal changes (using (4.47) in lieu of (4.21)) to show that the operator $(1-\Pi)\,\mathcal{D}_1\,(1-\Pi)$ is invertible (if $R > c_0$) on the Hilbert subspace in $\underline{\mathbb{L}}$ consisting of the sections of $\mathbb{R}/(\ell\mathbb{Z}) \times \mathbb{L}$ that fiberwise orthogonal to $\mathcal{K}$. There are two essential points in this regard, the first being this: If $u$ is a smooth section of the product $\mathbb{L}$ bundle that obeys $(1-\Pi)u=0$ at each point in $\mathbb{R}/(\ell\mathbb{Z})$, then

$$\int_{\mathbb{R}/(\ell\mathbb{Z})} \|\mathcal{D}_1 u\|_{\mathbb{L}}^2 \geq c_0^{-1} \int_{\mathbb{R}/(\ell\mathbb{Z})} (\|\tfrac{\partial}{\partial s}u\|_{\mathbb{L}}^2 + \|\nabla u\|_{\mathbb{L}}^2 + R(1+R|\cdot|^2)\|u\|_{\mathbb{L}}^2)\ . \tag{4.48}$$

The second essential point is that for this same $u$,

$$(1-\Pi)\mathcal{D}_1 u = \mathcal{D}_1 u - \phi_0\big(\langle\phi_0,(1-\chi_0\, x_j M_j)\tfrac{\partial}{\partial s}\, u\rangle_{\mathbb{L}} - \tfrac{1}{2}\langle\phi_0 w_{ij} x_i \tfrac{\partial}{\partial x_i}\, u\rangle_{\mathbb{L}}\big) \tag{4.49}$$

which implies (with the help of (4.19) and (4.29)) that

$$\int_{\mathbb{R}/(\ell\mathbb{Z})} \|(1-\Pi)\mathcal{D}_1 u\|_{\mathbb{L}}^2 \geq c_0^{-1} \int_{\mathbb{R}/(\ell\mathbb{Z})} \|\mathcal{D}_1 u\|_{\mathbb{L}}^2 - c_0\Big( \int_{\mathbb{R}/(\ell\mathbb{Z})} \|u\|_{\mathbb{L}}^2 + \tfrac{1}{R} \int_{\mathbb{R}/(\ell\mathbb{Z})} \|\nabla u\|_{\mathbb{L}}^2 \Big)\ . \tag{4.50}$$

The two inequalities (4.48) and (4.50) lead directly to the following conclusion: If $u$ is a smooth section of the product $\mathbb{L}$ bundle with $(1-\Pi)u \equiv 0$, then

$$\int_{\mathbb{R}/(\ell\mathbb{Z})} \|(1-\Pi)\mathcal{D}_1 u\|_{\mathbb{L}}^2 \geq c_0^{-1} \int_{\mathbb{R}/(\ell\mathbb{Z})} (\|\tfrac{\partial}{\partial s}u\|_{\mathbb{L}}^2 + \|\nabla u\|_{\mathbb{L}}^2 + R(1+R|\cdot|^2)\|u\|_{\mathbb{L}}^2)\ . \tag{4.51}$$

The preceding inequality implies a $\mathcal{D}_1$ analog of (4.34) when $R > c_0$ that says this: Supposing that $q$ is a smooth section of the product $\mathbb{L}$-bundle with $(1-\Pi)q = 0$ at each point in $\mathbb{R}/(\ell\mathbb{Z})$, then there exists a unique section $u$ of this same bundle with $(1-\Pi)u = 0$ at each point such that $(1-\Pi)\mathcal{D}_1 u = q$. Moreover,

$$\int_{\mathbb{R}/(\ell\mathbb{Z})} (\|\tfrac{\partial}{\partial s}u\|_{\mathbb{L}}^2 + \|\nabla u\|_{\mathbb{L}}^2 + R\|u\|_{\mathbb{L}}^2) \leq c_0 \int_{\mathbb{R}/(\ell\mathbb{Z})} \|q\|_{\mathbb{L}}^2\ . \tag{4.52}$$

Granted (4.52), standard perturbation theoretic arguments can be used to prove that $(1-\Pi)(\mathcal{D}_1 - E)(1-\Pi)$ is invertible when $E \leq c_0^{-1}\sqrt{R}$ and $R > c_0$; and that (4.51) is obeyed with $\mathcal{D}_1 - E$ replacing $\mathcal{D}_1$ on the left hand side (with a larger version of $c_0$ on the right hand side). Likewise, if $(1-\Pi)(\mathcal{D}_1 - E)u = q$ with both $q$ and $u$ in the kernel of $(1-\Pi)$ at each point, then (4.52) is obeyed (assuming a priori that $E \le c_0^{-1}\sqrt{R}$ and $R > c_0$).

The observations in the preceding paragraph are applicable to (4.47) with $u = \psi^\perp$; in which case the resulting version of (4.52) says in effect that

$$\int_{\mathbb{R}/(\ell\mathbb{Z})} (\|\tfrac{\partial}{\partial s}\psi^\perp\|_{\mathbb{L}}^2 + \|\nabla\psi^\perp\|_{\mathbb{L}}^2 + R\|\psi^\perp\|_{\mathbb{L}}^2) \le c_0 \int_{\mathbb{R}/\ell\mathbb{Z}} |\mu|^2 + c_0 \tfrac{1}{R} \int_{\mathbb{R}/\ell\mathbb{Z}} |\tfrac{d}{ds}\mu|^2 \tag{4.53}$$

given the a priori assumption that $E \le c_0^{-1}\sqrt{R}$ and that $R > c_0$.

The analysis now turns its focus on the $\Pi$ projection of the equation $\mathcal{D}_1\psi = E\psi$. This projection leads to an equation for $\mu$ which is the $\mathcal{D}_1$ analog of (4.24):

$$\tfrac{d}{ds}\mu + i\alpha_*\mu - \mu\langle\phi_0, x_jM_j\tfrac{\partial}{\partial s}\phi_0\rangle + \langle\phi_0,(1-\chi_0 x_jM_j)\tfrac{\partial}{\partial s}\psi^\perp\rangle_{\mathbb{L}} - \tfrac{1}{2}\langle\phi_0 w_{ij}x_i\tfrac{\partial}{\partial x_i}\psi^\perp\rangle_{\mathbb{L}} = -iE\mu\ . \tag{4.54}$$

Here, as before, $\alpha_*$ is a constant from the interval $[0, 2\pi/\ell)$. Note with regards to (4.54): Terms that are proportional to $\langle\phi_0,\chi_0 x_jM_j\phi_0\rangle_{\mathbb{L}}\tfrac{d}{ds}\mu$ and to $\langle\phi_0 w_{ij}x_j\tfrac{\partial}{\partial x_i}\phi_0\rangle_{\mathbb{L}}\mu$ are absent because both of the $\mathbb{L}$-inner products are zero. (To see that these are zero: Use the coordinates where $\mathbb{M}$ is diagonal and then note from the third bullet of (4.11) that neither $\phi_0$ nor $\chi_0$ change when any of the $x_j$'s is replaced by $-x_j$. It follows from this fact that $\langle\phi_0,\chi_0 x_jM_j\phi_0\rangle_{\mathbb{L}}$ and $\langle\phi_0 w_{ij}x_j\tfrac{\partial}{\partial x_i}\phi_0\rangle$ are sums of $\mathbb{R}^3$-integrals whose integrands change sign when some particular $x_j$ is changed to $-x_j$. This implies that these integrals are zero). The next paragraphs derive some useful bounds for the terms that do appear in (4.54).

With regards to the term $\mu\langle\phi_0, x_jM_j\tfrac{\partial}{\partial s}\phi_0\rangle$ that appears in (4.54): It follows from (4.20) and (4.29) that this term is bounded at each $s \in \mathbb{R}/\ell\mathbb{Z}$ by $c_0|\mu|\tfrac{1}{\sqrt{R}}$. With regards to the two terms in (4.54) that involve derivatives of $\psi^\perp$: There is first the term $\langle\phi_0, \tfrac{\partial}{\partial s}\psi^\perp\rangle_{\mathbb{L}}$ and then there are the terms $\langle\phi_0,\chi_0 x_jM_j)\tfrac{\partial}{\partial s}\psi^\perp\rangle_{\mathbb{L}}$ and $\tfrac{1}{2}\langle\phi_0,\chi_0 w_{ij}x_i\tfrac{\partial}{\partial x_i}\psi^\perp\rangle_{\mathbb{L}}$. The latter two are pointwise bounded in $\mathbb{R}/(\ell\mathbb{Z})$ by $c_0\tfrac{1}{\sqrt{R}}\|\tfrac{\partial}{\partial s}\psi^\perp\|_{\mathbb{L}}$ and $c_0\tfrac{1}{\sqrt{R}}\|\nabla\psi^\perp\|_{\mathbb{L}}$. It then follows as a consequence of (4.54) that the $\mathbb{R}/(\ell\mathbb{Z})$ integral of the sum of the squares of these two terms is bounded by

$$c_0\Big(\tfrac{1}{R}\int_{\mathbb{R}/\ell\mathbb{Z}} |\mu|^2 + \tfrac{1}{R^2}\int_{\mathbb{R}/\ell\mathbb{Z}} |\tfrac{d}{ds}\mu|^2\Big)\ . \tag{4.55}$$

When used with (4.54), the bounds in (4.55) and in the preceding paragraph lead to the following bound (keeping in mind that $E \le c_0^{-1}\sqrt{R}$ and $R \ge c_0$):

$$\int_{\mathbb{R}/\ell\mathbb{Z}} |\tfrac{d}{ds}(e^{-i(\alpha_*+E)s}\mu)|^2 \le c_0\Big(\tfrac{1}{R}\int_{\mathbb{R}/\ell\mathbb{Z}} |\mu|^2 + \tfrac{1}{R^2}\int_{\mathbb{R}/\ell\mathbb{Z}} |\tfrac{d}{ds}(e^{-i(\alpha_*+E)s}\mu)|^2\Big)\,. \tag{4.56}$$

Since the $\mathbb{R}/(\ell\mathbb{Z})$ integral of $\mu^2$ is bounded by 1, this last bound can hold only in the event that E differs from $-\alpha_*$ plus some integer multiple of $2\pi/\ell$ by at most $c_0\frac{1}{\sqrt{R}}$. Moreover, if n denotes that integer, than $\mu$ differs from $e^{2\pi i n s/\ell}$ by map from $\mathbb{R}/(\ell\mathbb{Z})$ to $\mathbb{C}$ with a $\frac{1}{R}\kappa$ bound on the $\mathbb{R}/(\ell\mathbb{Z})$ integral of the square of its norm.

The observations in the paragraph above with the bounds in (4.45) complete the proof of the assertion in Lemma 4.6 to the effect that eigenvalues of $\mathcal{D}_1$ with norm at most $c_0^{-1}\sqrt{R}$ and their eigenvectors have the required form.

The proof of Lemma 4.6 ends with an explanation for the claim that any given integer (call it n) with norm at most $c_0^{-1}\sqrt{R}$ corresponds to an eigenvalue of $\mathcal{D}_1$ that can be written as $-(\alpha+n)\frac{2\pi}{\ell} + \tau$ with $|\tau| \le c_0\frac{1}{\sqrt{R}}$. The explanation for why this is so follows the same contraction mapping route as the explanation in the proof of Lemma 4.4 for the analogous claim in that lemma. The only substantive differences are the appearance of extra terms in the fixed point equation of (4.37) to account for the terms that changed (4.24) to (4.54). In this regard, the map sends a pair (w, ε) to the pair whose first component (a function on $\mathbb{R}/(\ell\mathbb{Z})$) and second component (a number) are

$$D_n^{-1}(1-P_n)\mathfrak{Y} \quad \textit{and} \quad i\tfrac{1}{\ell}\int_0^{\ell} e^{-2\pi n s/\ell}\mathfrak{Y} \tag{4.57}$$

where $\mathfrak{Y}$ is given below

$$\mathfrak{Y} = (e^{2\pi i n s/\ell}+w)\langle\phi_0, x_j M_j \tfrac{\partial}{\partial s}\phi_0\rangle - \langle\phi_0,(1-\chi_0 x_j M_j)\tfrac{\partial}{\partial s}\psi^{\perp}\rangle_{\mathbb{L}} + \tfrac{1}{2}\langle\phi_0 w_{ij} x_i \tfrac{\partial}{\partial x_i}\psi^{\perp}\rangle_{\mathbb{L}} - i\varepsilon w\,. \tag{4.58}$$

As was the case before, $\mathfrak{Y}$ is an affine function of w. In this regard, the $\mathbb{R}/(\ell\mathbb{Z})$ integral of the square of the norm of the w = 0 version of $\mathfrak{Y}$ is bounded by the $\mu = e^{2\pi i n s/\ell}$ version of (4.55) (thus by $c_0(\frac{1}{R} + \frac{n^2}{R^2})$); and the square of the norm of the part that is linear in w is bounded by the $\mu$ = w version of (4.55). Meanwhile, the dependence of $\mathfrak{Y}$ on ε is straightforward to analyze using (4.40), the result being the analog of the bounds in (4.41) with the respective integrands replaced by the ε = 0 version of $\mathfrak{Y}^2$ and the square of the difference between given ε and ε´ version of $\mathfrak{Y}$. The version of $c_0$ will also be different (and the implicit assumption is that the $\mu$ = w version of (4.55) is at most 1; or any given fixed number at the expense of changing the version of $c_0$ that must be used).

It is a direct consequence of the preceding observations on $\mathfrak{Y}$ that (4.57) defines a contraction mapping of the radius 1 ball in $\mathbb{B}_{\ddagger}$ to itself if $|n| < c_0^{-1}\sqrt{R}$ and $R > c_0$.

Therefore, it has a unique fixed point in this same ball; and this fixed point is a pair $(w, \varepsilon)$ with the $\mathbb{R}/(\ell\mathbb{Z})$ integral of $w^2$ bounded by $c_0 \frac{1}{R}$ and with $|\varepsilon| \le c_0 \frac{1}{\sqrt{R}}$. The existence of such a pair (and $\varepsilon$ in particular) confirms the claim that $\mathcal{D}_1$ has an eigenvalue that differs from $-(\alpha + n)\frac{2\pi}{\ell}$ by at most $c_0 \frac{1}{\sqrt{R}}$.

***Proof of Lemma 4.7***: Except for cosmetics, the proof is much the same as the proof of Lemma 4.6. For this reason, the details are omitted with the hope that the remarks that follow are sufficient for the reader. The first remark is this: If the norms of B and/or C are not bounded by a small multiple of R (a multiple bounded by $c_0^{-1}$), then the terms in (4.44) and (4.45) with B and C will not be small with respect to the $\sqrt{2}\, i R \mathbb{M}_{ik} x_i \rho_k$ term in $\mathfrak{D}_0$ at any given $s \in \mathbb{R}/(\ell\mathbb{Z})$. (The terms with B and C are linear with respect to the coordinate x as is $\sqrt{2}\, i R \mathbb{M}_{ik} x_i \rho_k$.) In particular, the terms with B at any given $s \in \mathbb{R}/(\ell\mathbb{Z})$ won't be a small perturbation of $\mathfrak{D}_0$ and thus not a small perturbation to $\mathcal{D}_0$. Indeed, if $|B|$ and/or $|C|$ are very much greater than R, then one should consider the $\sqrt{2}\, i R \mathbb{M}_{ik} x_i \rho_k$ term at any given $s \in \mathbb{R}/(\ell\mathbb{Z})$ as a perturbation of the operator $\gamma_k(\frac{\partial}{\partial x_i} + \frac{1}{2} i x_j \varepsilon_{jik} C_k) + \gamma_s i x_j B_j$. In this case, the analysis changes significantly.

The second remark is technical for working through the changes to prove Lemma 4.7 using just minor modifications of the arguments for Lemma 4.6. To state these remarks, let $x\cdot\mathfrak{L}$ now denote the sum of the terms with B and C in (4.44) and (4.45). If $\phi$ is in kernel($\mathfrak{D}_0$) at any given $s \in \mathbb{R}/(\ell\mathbb{Z})$, then $x\cdot\mathfrak{L}\phi$ is orthogonal to the kernel of $\mathfrak{D}_0$ at that value of s (this is because $\mathfrak{L}$ is linear in x). In general if $\psi$ at any give $s \in \mathbb{R}/(\ell\mathbb{Z})$ is an eigenvector of $\mathfrak{D}_0$, then

$$\| x\cdot\mathfrak{L}\psi \|_{\mathbb{L}} \le c_0 \frac{1+\hat{E}}{\sqrt{R}} (|B| + |C|) \|\psi\|_{\mathbb{L}} \tag{4.59}$$

where $\hat{E}$ denotes the corresponding $\mathfrak{D}_0$ eigenvalue.

## d) Eigenvalue bounds from the model case

The model operators in the previous subsection are used in this subsection to analyze the spectrum of the operator $\mathcal{D}_{\ddagger}$ that is depicted in (4.2) with $(A, \omega)$ described by (3.20). In particular, this subsection first states and then proves a proposition to the effect that the operator $\mathcal{D}_{\ddagger}$ has only a few eigenvalues with small absolute value when three conditions hold: The number R is large, the norm of the curvature tensor of A is small relative to R, and the endomorphism Q has compact support in the complement of Z. By way of notation, this upcoming proposition uses p to denote the number of components of Z and it implicitly labels the components of Z as $\{Z_1, \ldots, Z_p\}$. It then uses $\ell_k$ for any givne $k \in \{1,\ldots, p\}$ to denote the length of the component of $Z_k$.

**Proposition 4.8**: *Given* $\Xi \geq 1$ *and* $r_\lozenge \in (0,1]$ *and* $c_\ddagger > 0$*, there exists* $\kappa > 1$ *such that what follows is true. Fix* $R \geq \kappa$ *and use* R *when defining the operator* $\mathcal{D}_\ddagger$ *that is depicted in (4.2). With regards to other data used in (4.2), assume the additional constraints below.*

- *The section* $\varsigma^+$ *obeys* $|\varsigma^+| \geq \frac{1}{\Xi}\text{dis}(\cdot,Z)$
- *The norm of the curvature of the connection* A *is bounded by* $\frac{1}{\kappa^2}R$ *and that of its covariant derivative is bounded by* $c_\ddagger R$*;*
- *The endomorphism* $\mathfrak{Q}$ *can be written as* $\mathfrak{Q}_0 - q\Gamma$ *where* $\mathfrak{Q}_0$ *denotes an endomorphism with norm bounded by* $\frac{1}{\kappa^2}\sqrt{R}$ *and with support where the distance to* Z *is greater than* $r_\lozenge$*; and where q is constant with norm bounded by* 1.

*If* $\mathcal{D}_\ddagger$ *is defined using data obeying these constraints, then there exist* p *numbers* $\{\alpha_k \in [0,1)\}_{k\in\{1,\ldots,p\}}$ *with the following significance: A given real number is an eigenvalue of* $\mathcal{D}_\ddagger$ *with absolute value at most* $\frac{1}{\kappa}\sqrt{R}$ *only if it can be written as*

$$\sum_{k=1}^{p} -(\alpha_k + n_k)\frac{2\pi}{\ell_k} + \mathfrak{r}$$

*with each* $k \in \{1, \ldots, p\}$ *version of* $n_k$ *being an integer; and with* $\mathfrak{r}$ *obeying* $|\mathfrak{r}| \leq \kappa\frac{1}{\sqrt{R}}$.
*Moreover, each eigenvalue in the indicated range has multiplicity at most* p*. Conversely, if* $\{n_k\}_{1\leq k\leq p}$ *is a set of integers with each having norm bounded by* $\frac{1}{\kappa}\sqrt{R}$*, then there is a corresponding eigenvalue of* $\mathcal{D}_\ddagger$ *that has the form depicted above.*

***Proof of Proposition 4.8***: The various incarnations of $c_0$ that appear below in the proof depend implicitly on the numbers $\Xi$, $r_\lozenge$ and $c_\ddagger$ in the statement of the theorem.

A remark is in order with regards to $\Xi$: The uniform lower bound for the norm $|\varsigma^+|$ (it is greater than $\frac{1}{\Xi}\text{dis}(\cdot,Z)$) has various implications, the first being that the norm of the inverse of the matrix $\mathbb{M}$ from (4.2) is bounded by $c_0\frac{1}{\Xi}$. That bound and the fact that $\varsigma^+$ is harmonic with $|\varsigma^+|$ bounded by $c_0$ leads to both upper and lower bounds for the total length of Z. The bounds for $|\mathbb{M}^{-1}|$ and $\varsigma^+$ also lead to upper bounds for the norm of Z's extrinsic curvature (as a submanifold in X), to uniform lower bounds for the radii of various tubular neighborhoods of the components of Z that are used in the subsequent analysis, and to upper bounds for various norms of difference between the Riemanian metric and a fiducial metric on these neighborhood. Various incarnations of $c_0$ in the proof incorporate the implicit $\Xi$ dependence of these bounds.

The proof has three parts.

*Part 1*: Proposition 4.1 says in effect that the support of eigenvectors of $\mathcal{D}_{\ddagger}$ with eigenvalue between $-c_0^{-1}\sqrt{R}$ and $c_0^{-1}\sqrt{R}$ is mostly contained in a small radius tubular neighborhood of Z. This part of the proof sets up a coordinate system on the radius $r_0$ tubular neighborhood of Z to exploit this concentration phenomenon. To start the set-up, fix attention on a given component of Z (which henceforth will be denoted by this same letter Z) and then choose a unit speed parametrization of Z by $\mathbb{R}/(\ell\mathbb{Z})$ with $\ell$ denoting the length of Z.

To continue the set-up, fix an oriented, orthonormal frame for the normal bundle to Z in X to depict this bundle isometrically as the product bundle over $\mathbb{R}/(\ell\mathbb{Z})\times\mathbb{R}^3$. Use the metric's exponential map from the normal bundle to Z in X to define a diffeomorphism from the product of $\mathbb{R}/(\ell\mathbb{Z})$ and a small radius ball in $\mathbb{R}^3$ centered at the origin onto a neighborhood of Z in X. Write that radius as $100\,r_0$ and let $B_0$ denote the ball in $\mathbb{R}^3$ centered at the origin with radius $r_0$. (Choose $r_0$ so that the closures of the radius $100\,r_0$ tubular neighborhoods of the components of Z are pairwise disjoint.) This diffeomorphism will be used below implicitly to identify the neighborhood of Z with $\mathbb{R}/(\ell\mathbb{Z})\times B_0$.

To depict the operator $\mathcal{D}_{\ddagger}$ using the preceding identification: Let s denote the Euclidean affine parameter for $\mathbb{R}/(\ell\mathbb{Z})$ and let $(x_1, x_2, x_3)$ denote Euclidean coordinates for $\mathbb{R}^3$. The pull-back of the metric on X to $\mathbb{R}/(\ell\mathbb{Z})\times B_0$ has the schematic form

$$(1+x_iM_i)^2\, ds\otimes ds + (dx_i + w_{ij}x_j ds)\otimes(dx_i + w_{ik}x_k ds) + \mathfrak{g} \tag{4.60}$$

where the notation is as follows: First, repeated indices are summed over the set $\{1, 2, 3\}$. Meanwhile, $\{M_i\}_{i=1,2,3}$ denote the components of a smooth map from $\mathbb{R}/(\ell\mathbb{Z})$ to $\mathbb{R}^3$, and $\{w_{ij}\}_{i,j\in\{1,2,3\}}$ denote the components of a $3\times 3$, anti-symmetric matrix valued map from $\mathbb{R}/(\ell\mathbb{Z})$. Their norms have a priori $c_0$ bounds. (Just to be sure: The number $c_0$ depends on the Riemannian metric and on the number $\Xi$ that appears in the statement of the proposition.) Finally, what is denoted by $\mathfrak{g}$ signifies a symmetric, bilinear form with norm bounded by $c_0|x|^2$, with derivative norm bounded by $c_0|x|$, and with second derivative norm bounded by $c_0$.

With regards to the line bundle $\mathcal{I}$: The restriction of $\mathcal{I}$ to Z must be isomorphic to the product real line bundle if $\varsigma^+$ vanishes transversally along Z. To see why, note first that transversal vanishing of $\varsigma^+$ along Z is the requirement that $\nabla\varsigma^+$ define an isomorphism along Z from Z's normal bundle in X (which is isomorphic to the product $\mathbb{R}^3$ bundle) to $\Lambda^+|_Z\otimes\mathcal{I}$. Meanwhile, Z's normal bundle and $\Lambda^+$ are both oriented vector bundles along Z so both are isomorphic to the product $\mathbb{R}^3$ bundle along Z. As a consequence, the determinant of the isomorphism $\nabla\varsigma^+$ defines an isomorphism on Z between the product $\mathbb{R}$ bundle and its tensor product with $\mathcal{I}$, which is an isomorphism between the product bundle and $\mathcal{I}$.

An isomorphism along $\mathbb{R}/(\ell\mathbb{Z})\times\{0\}$ between $\mathcal{I}$ and the product $\mathbb{R}$ bundle identifies $\sigma$ along along $\mathbb{R}/(\ell\mathbb{Z})\times\{0\}$ as section of $\mathrm{ad}(P)$ along $\mathbb{R}/(\ell\mathbb{Z})\times\{0\}$. To say more about this section, note that the restriction of P to $\mathbb{R}/(\ell\mathbb{Z})\times\{0\}$ is necessarily isomorphic to the product bundle; and having chosen an isomorphism, it can be viewed as such. This isomorphism with the product bundle can then be extended to the whole of $\mathbb{R}/(\ell\mathbb{Z})\times B_0$ using parallel transport by the connection A along the rays from the origin in $B_0$ at each fixed parameter value in $\mathbb{R}/(\ell\mathbb{Z})$. When viewed in this light, the homomorphism $\sigma$ on $\mathbb{R}/(\ell\mathbb{Z})\times B_0$ defines a map from $\mathbb{R}/(\ell\mathbb{Z})\times B_0$ to the unit sphere in $\mathfrak{su}(2)$ and is thus homotopic to a constant map (because $\pi_1(S^2)$ is trivial). Any such homotopy can be lifted to a map from $\mathbb{R}/(\ell\mathbb{Z})\times B_0$ to SU(2), and such a lift can then be used to construct an automorphism of the product principal bundle over $\mathbb{R}/(\ell\mathbb{Z})\times B_0$ with two salient features: The automorphism pulls back $\sigma$ to a constant element in $\mathfrak{su}(2)$ on part of $\mathbb{R}/(\ell\mathbb{Z})\times B_0$ where $|x| \le \frac{1}{4} r_0$; and the automorphism extends to the whole of X as the identity outside of the $|x| \le \frac{1}{2} r_0$ part of $\mathbb{R}/(\ell\mathbb{Z}) \times \mathbb{R}^3$. Let $B_1$ now denote the $|x| \le \frac{1}{4} r_0$ part of $B_0$. The depiction of $\sigma$ as being constant on $\mathbb{R}/(\ell\mathbb{Z})\times B_1$ is used henceforth.

With regards to the $\mathbb{R}^2$ bundle $\mathfrak{L}$: Since $\sigma$ is constant on $\mathbb{R}/(\ell\mathbb{Z})\times B_1$, the $\mathbb{R}^2$ bundle $\mathfrak{L}$ appears over $\mathbb{R}/(\ell\mathbb{Z})\times B_1$ as the product $\mathbb{R}^2$ bundle with fiber being the orthogonal complement to $\sigma$ in $\mathfrak{su}(2)$. This bundle is henceforth viewed as the product $\mathbb{C}$ bundle by using the endomorphism $\frac{1}{2}[\sigma, \cdot\,]$ to define what is meant by multiplication by i.

With regards to the connection A: The isomorphism described above between P and the product bundle over $\mathbb{R}/(\ell\mathbb{Z})\times B_1$ can be used to depict the connection A on the part of P over $\mathbb{R}/(\ell\mathbb{Z})\times B_1$ as the sum of the product connection (denoted by $\theta_0$) and an $\mathfrak{su}(2)$ valued 1-form that is proportional to $\sigma$: Thus

$$A = \theta_0 + \big((b_0 + x_i C_i)ds + \tfrac{1}{2} x_i \varepsilon_{ijk} B_k dx_j + \mathfrak{z}\big)\sigma \tag{4.61}$$

with the notation as follows: A priori, $b_0$ is a function on $\mathbb{R}/(\ell\mathbb{Z})$, but more is said about it momentarily. Meanwhile, each element from the sets $\{B_i\}_{i=1,2,3}$ and $\{C_i\}_{i=1,2,3}$ denotes a function on $\mathbb{R}/(\ell\mathbb{Z})$ whose norm is bounded by $c_{\ddagger 0}R$ with $c_{\ddagger 0}$ denoting a positive number which is at most $c_0 \frac{1}{R} \sup_X|F_A|$. As for $\mathfrak{z}$, it denotes a 1-form on $\mathbb{R}/(\ell\mathbb{Z})\times B_1$ whose norm is bounded by $c_{\ddagger 1} R|x|^2$ with $c_{\ddagger 1}$ being another positive number, this one with norm at most $c_0 \frac{1}{R} \sup_X|\nabla_A F_A|$. Thus, the enforcement of an upper bound on the norms of $F_A$ and $\nabla_A F_A$ over X enforces an upper bound on $c_{\ddagger 0}$ and $c_{\ddagger 1}$. (Upper bounds on the norm of $F_A$ is enforced by the upcoming Lemma 4.9.)

With regards to $b_0$: The product structure for ad(P) over $\mathbb{R}/(\ell\mathbb{Z})$ can be changed if necessary keeping $\sigma$ constant so that $b_0$ is a constant from the interval $[0, 2\pi/\ell)$. This constant version of $b_0$ is assumed in what follows. An important point to keep in mind for what comes later is that the chosen product structure that makes $b_0$ constant and in the interval $[0, 2\pi/\ell)$ is tailored specifically to the connection A. If the connection A is varied along a path of connections, then the corresponding tailored product structures need not vary continously along the path; and if they don't, then the corresponding path of $b_0$'s won't change continuously either (the versions of $b_0$ along the path can jump from just above zero to just below $2\pi/\ell$ and vice-versa). This last fact is a manifestation of the fact (also important to note) that this isomorphism of the product principal bundle over $\mathbb{R}/(\ell\mathbb{Z})$ that changes $b_0$ to lie in the prescribed interval might not be homotopic to the identity through isomorphisms that fix $\sigma$.

*Part 2:* No generality is lost by taking the number $r_\Diamond$ in the statement of Proposition 4.8 to be less than $\frac{1}{50} r_0$. With this understood, let $B_\Diamond$ denote the concentric ball inside $B_1$ with radius $\frac{1}{100} r_\Diamond$. The observations in Part 1 are used now to write $\mathcal{D}_\ddagger$ on $\mathbb{R}/(\ell\mathbb{Z})\times B_\Diamond$ as

$$\mathcal{D}_\ddagger = \gamma_s\big((1 - x_i M_i)\tfrac{\partial}{\partial s} - \tfrac{1}{2}\, x_j \dot{M}_j - W_{ij} x_j \tfrac{\partial}{\partial x_i} + 2i(b_0 + x_i C_i)\big) \\ + \gamma_i\big(\tfrac{\partial}{\partial x_i} + i\, x_j \varepsilon_{jik} B_k\big) + \sqrt{2}\, i R \mathbb{M}_{ik} x_i \rho_k - q\Gamma + \mathfrak{P} \tag{4.62}$$

with $\mathfrak{P}$ denoting a first order differential operator acting on maps from $\mathbb{R}/(\ell\mathbb{Z})\times B_\Diamond$ to $\mathbb{C}^4 \oplus \mathbb{C}^4$ that can be written as $\mathfrak{p}_s \frac{\partial}{\partial s} + \mathfrak{p}_j \frac{\partial}{\partial x_i} + \mathfrak{p}_0$ with the endomorphisms $\mathfrak{p}_s, \{\mathfrak{p}_j\}_{j=1,2,3}$ and $\mathfrak{p}_0$ obeying the bounds below ($c_1$ denotes a positive constant).

- $|\mathfrak{p}_0| \le c_1(1 + |x| + R|x|^2)$ ,
- $|\mathfrak{p}_s| + |\mathfrak{p}_1| + |\mathfrak{p}_2| + |\mathfrak{p}_3| \le c_1 |x|^2$ .

(4.63)

With regards to $c_1$: This depends in part on part on the number $c_{\ddagger 1}$ that is used to bound the norm of the $\mathfrak{z}$ term in (4.61). It also depends in part on the sup-norm of $|\nabla\varsigma^+|$, on the extrinsic curvature of Z and on the Riemannian metric near Z.

To exploit this depiction of $\mathcal{D}_\ddagger$, fix $\kappa_\ddagger$ so that $|F_A| \le \kappa_\ddagger R$. Having done that, let $\kappa_\Diamond$ denote the version of the number $\kappa$ that appears in Proposition 5.1. Suppose that $\psi$ is an eigenvector for the operator $\mathcal{D}_\ddagger$ with the absolute value of its eigenvalue, E, being at most $\frac{1}{100\kappa_\Diamond}\sqrt{R}$, and with the integral of $|\psi|^2$ on X equal to 1. For any $r \in (0, \frac{1}{100} r_\Diamond)$, Proposition 4.1 can be invoked to see that the contribution to the integral of $|\psi|^2$ from the part of X where $\mathrm{dist}(\cdot, R) \ge r$ is very small when R is large. Keeping this in mind, introduce $\chi_r$ to denote the function $\chi(\frac{\mathrm{dist}(\cdot,Z)}{r} - 1)$ and then set $\psi_r$ to denote $\chi_r\psi$. (Thus, $\psi_r = \psi$ where $\mathrm{dist}(\cdot,Z) \le r$ and $\psi_t = 0$ where $\mathrm{dist}(\cdot,Z) \ge 2r$.) According to Proposition 5.1,

$$|\mathcal{D}_{\ddagger}\psi_r - E\psi_r| \le c_0 R^2 e^{-\frac{1}{\kappa_\Diamond} R r^2} , \tag{4.64}$$

the right hand side having support where $\mathrm{dist}(\cdot,Z) \in [r, 2r]$. In any event, the right hand side is very small if r is $\mathcal{O}(1)$ and R is large.

The section $\psi_r$ can now be pulled back to $S^1\times B_\Diamond$ when $r < \frac{1}{4} r_\Diamond$, and having compact support there (because of $\chi_r$), it can be viewed there as an element in the domain of a version of the operator $\mathcal{D}_2$ that is depicted in (4.44) and (4.45). And, it follows from what is said in Proposition 4.1, and from (4.64) and (4.62) that

$$|\mathcal{D}_2\psi_r - E\psi_r - q\Gamma\psi_r + 2i\gamma_5 b_0\psi_r + \mathfrak{P}\psi_r| \le c_0 R^2 e^{-\frac{1}{\kappa_\Diamond} R r^2} . \tag{4.65}$$

The next part of the proof exploits this last inequality.

*Part 3*: Granted (4.65), then the analysis for the proofs of Lemmas 4.6 and 4.7 can be repeated with only marginal changes (which include an appeal to (4.63) and some facts about $\Gamma$) to deduce the following assertion:

**Lemma 4.9**: *Given the number $c_1$ from (4.63), there exists $\kappa > 1$ such that if $|E| \le \frac{1}{\kappa}\sqrt{R}$ and $R \ge \kappa$ and $|q| \le \frac{1}{\kappa} R^{1/4}$, and if r is between $\kappa\frac{1}{\sqrt{R}}\ln R$ and $\frac{1}{\kappa}$, then the eigenvalue E can be written as*

$$E = (-\alpha + n)\frac{2\pi}{\ell} + \mathfrak{r}_1$$

*with n being an integer, with $\alpha \in [0, 2\pi)$, and with $\mathfrak{r}_1$ obeying $|\mathfrak{r}_1| \le \kappa(1+t^2)\frac{1}{\sqrt{R}}$. Meanwhile, $\psi_r$ has the form $\psi_r = e^{2\pi i n s/\ell}\phi_0 + \mathfrak{k}$ with $\phi_0$ denoting a particular n-independent, unit norm section of $\mathcal{K}$, and with $\mathfrak{k}$ obeying $\int_{\mathbb{R}/\ell\mathbb{Z}} \|\mathfrak{k}\|_{\mathbb{L}}^2 \le \kappa(1+|q|^4)\frac{1}{R}$.*

*Conversely, given an integer n with $|n| \le \frac{1}{\kappa}\sqrt{R}$, there is an eigenvalue of $\mathcal{D}_{\ddagger}$ that can be written in the manner just described.*

This lemma with (4.64) and (4.65) lead directly to what is asserted by Proposition 5.8.

***Proof of Lemma 4.9***: As noted, the argument for Lemma 4.9 differs little from the arguments used to prove Lemmas 4.6 and 4.7. For this reason, the account that follows will focus only the new terms $-q\Gamma$ and $\mathfrak{P}_0$, $\mathfrak{P}_1$ and $2ib_0\gamma_5$. In any event, the argument starts just as before by writing $\psi_r$ as $\mu\,\phi_0 + \psi^\perp$ with $\psi^\perp = (1-\Pi)\psi_r$.

The key point with regards to the term $q\Gamma$ is that $\langle\phi_0, \Gamma\phi_0\rangle = 0$. As a consequence, this term adds a term proportional to $q\mu\Gamma\phi_0$ to the left hand side of (4.47) which adds a multiplicative factor of $c_0(1+|q|)^2$ before the integral of $|\mu|^2$ that appears on the right hand side of (4.53). The $q\Gamma$ term also adds a term to the left hand side of (4.54) that has the form $-q\langle\phi_0, \Gamma\psi^\perp\rangle_{\mathbb{L}}$. This term modifies (4.55) and (4.56) by adding a multiplicative factor of $c_0(1+|q|)^4$ before the integral of $|\mu|^2$ that appears in each of these equations. There are also the corresponding modifications to $\mathfrak{Y}$ in (4.58) and to the bounds on $\mathfrak{Y}$ that are given in the paragraphs that follow (4.58).

With regards to the $(2i\gamma_5 b_0 + \mathfrak{P}_0 + \mathfrak{P}_1)\psi_r$ part of (4.65): The overall effect of this is to add a term to the left hand side of (4.54) and to (4.58). These then contribute to the versions of $c_0$ that appear in (4.55) and (4.56) and to the bounds that are derived in the paragraphs after (4.58). The new term on the left hand side of (4.54) changes (4.54) to the equation depicted schematically below.

$$(1+Q_s)\tfrac{d}{ds}\mu + i(\alpha_* + 2b_0 - Q)\mu - \langle\phi_0, x_j M_j \tfrac{\partial}{\partial s}\phi_0\rangle\mu \\ + \langle\phi_0, (1 - \chi_0 x_j M_j)\tfrac{\partial}{\partial s}\psi^\perp\rangle_{\mathbb{L}} - \tfrac{1}{2}\langle\phi_0 w_{ij} x_i \tfrac{\partial}{\partial x_i}\psi^\perp\rangle_{\mathbb{L}} - i\langle\phi_0, (\mathfrak{P}_0 + \mathfrak{P}_1)\psi^\perp\rangle_{\mathbb{L}} = -iE\mu\,, \tag{4.66}$$

with the function $Q_s$ being function $-i\langle\phi_0, \mathfrak{p}_s\phi_0\rangle_{\mathbb{L}}$ and with $Q_0$ being $\langle\phi_0, \mathfrak{P}\phi_0\rangle_{\mathbb{L}}$. By virtue of (4.63), these functions obey the bounds

$$|Q_s| \le c_0 c_1 \tfrac{1}{R} \quad \textit{and} \quad |Q| \le c_0 c_1\,, \tag{4.67}$$

with $c_1$ coming from (4.63).

It proves useful now to write $Q$ as $\underline{Q} + \frac{d}{ds} f$ with $\underline{Q}$ denoting the average of $Q$ over $\mathbb{R}/(\ell\mathbb{Z})$ and with $f$ denoting a function on $\mathbb{R}/(\ell\mathbb{Z})$. If $\mu$ is written as $e^{if}\hat{\mu}$ then (4.67) holds for $\hat{\mu}$ with $Q$ replaced by $\underline{Q} + \mathfrak{e}$ with $\mathfrak{e}$ being a function on $\mathbb{R}/(\ell\mathbb{Z})$ with norm at most $c_0 c_1 \frac{1}{\sqrt{R}}$. If $r \le c_0^{-1}$, then this $\hat{\mu}$ equation can be plugged directly into the arguments used for Lemmas 4.6 and 4.7 and those arguments then lead directly to the conclusions of Lemma 4.9.

With regards to writing $\mu$ as $e^{if}\hat{\mu}$ so as to view (4.67) as an equation of $\hat{\mu}$ with $Q$ replaced by $\underline{Q} + \mathfrak{e}$. This is equivalent to using the original version of $\mu$ but defined using a new section of $\mathcal{K}$; if the original section is denoted by $\phi_0$, then the new one is $e^{if}\phi_0$.

**e) The proof of Proposition 4.1**

This subsection is dedicated to proving Proposition 4.1. The various instances of $c_0$ in this subsection depend implicitly on the numbers $\Xi$ and $\kappa_\ddagger$. The proof has three parts.

*Part 1*: With no $\mathfrak{Q}$ term, the operator $\mathcal{D}_{\ddagger}^2$ acting on $\mathfrak{L}$-valued sections of the bundle $(\Lambda^+ \oplus \underline{\mathbb{R}}) \oplus T^*X$ has a Bochner-Wietzenboch formula that can be written as

$$\mathcal{D}_{\ddagger}^2 = \nabla^\dagger \nabla + \mathfrak{R} + \mathfrak{F}_A + \tfrac{1}{\sqrt{2}} R(\nabla_\alpha \varsigma^+{}_k)\gamma_\alpha \rho_k[\sigma, \cdot\,] + (\tfrac{1}{\sqrt{2}} R \varsigma^+{}_k \rho_k[\sigma, \cdot])^2 \tag{4.68}$$

with $\mathfrak{R}$ and $\mathfrak{F}_A$ as follows: What is denoted here by $\mathfrak{R}$ signifies an endomorphism of the bundle $(\Lambda^+ \oplus \underline{\mathbb{R}}) \oplus T^*X$ that is determined by the Riemannian metric; and what is denoted here by $\mathfrak{F}_A$ signifies an endmorphism of $((\Lambda^+ \oplus \underline{\mathbb{R}}) \oplus T^*X) \otimes \mathfrak{L}$ that is linear in its dependence on $F_A$. Of note here is that $|\mathfrak{F}_A| \le c_0 \kappa_{\ddagger} R$. As a consequence of this Bochner-Weitzenboch formula: If the $\mathfrak{Q}$ term in (4.2) obeys $|\mathfrak{Q}| \le c_0 \sqrt{R}$ and if $\psi$ is an eigenvector of the corresponding version of $\mathcal{D}_{\ddagger}$ with eigenvalue between $-c_0^{-1}\sqrt{R}$ and $c_0^{-1}\sqrt{R}$, then the square of the norm of $\psi$ obeys the inequality

$$\tfrac{1}{2} d^\dagger d|\psi|^2 + |\nabla_A \psi|^2 + c_0^{-1} R^2 \operatorname{dist}(\cdot, Z)^2 |\psi|^2 - c_0(1 + \kappa_{\ddagger}) R |\psi|^2 \le 0 \ . \tag{4.69}$$

Some immediate consequences of this inequality are described in the next lemma. The lemma uses $\kappa_\Diamond$ to denote the version of $c_0$ that appears in (4.69)

**Lemma 4.10**: *Given positive numbers $\kappa_{\ddagger}$ and $\kappa_\Diamond$, there exists $\kappa > 0$ with the following significance: Fix $R \ge r$. Suppose that the version of $c_0$ that appears in (4.58) is bounded by $\kappa_\Diamond$; and suppose that $\psi$ obeys the inequality in (4.58) with*

$$\int_X |\psi|^2 = 1.$$

*Then the inequalities below hold.*

- $\int_X (|\nabla_A \psi|^2 + R^2 \operatorname{dist}(\cdot, Z)^2 |\psi|^2) \le \kappa(1 + \kappa_{\ddagger}) R$ .
- $\int_X \operatorname{dist}(\cdot, Z)^2 |\nabla_A \psi|^2 \le \kappa(1 + \kappa_{\ddagger})$ .
- $|\psi| \le \kappa(1 + \kappa_{\ddagger}) R$ .

***Proof of Lemma 4.10***: The top bullet's inequality follows by integrating both sides of (4.69) over X. To obtain the second bullet's inequality: On the domain in X where

$$\operatorname{dist}(\cdot, Z) \ge c_0(1 + \kappa_{\ddagger})^{1/2} \tfrac{1}{\sqrt{R}} \ , \tag{4.70}$$

the inequality in (4.69) implies directly that

$$\tfrac{1}{2}\, d^\dagger d|\psi|^2 + |\nabla_A\psi|^2 + c_0^{-1}R^2\,\mathrm{dist}(\cdot,Z)^2\,|\psi|^2 \le 0\ .$$

(4.71)

With the latter point understood: Multiply both sides of (4.71) by $|\varsigma^+|^2$ (which is greater than $c_0^{-1}\mathrm{dist}(\cdot,Z)^2$), and then integrate over X. Then, integrate by parts to remove the operator $d^\dagger d$ from $|\psi|^2$. What with (4.70), the resulting inequality leads directly to the bound in the second bullet.

To see about the third bullet: Fix $p \in X$ and let $G_p$ now denote the Green's function for the operator $d^\dagger d + 1$ on X. Multiply both sides of (4.69) by $G_p$ and then integrate by parts to remove $d^\dagger d$ from $|\psi|^2$. The result is a bound that has the form

$$\tfrac{1}{2}\,|\psi|^2(p) \le c_0(1+\kappa_\ddagger)(1+R)\int_X \frac{1}{\mathrm{dist}(\cdot,p)^2}|\psi|^2\ .$$

(4.72)

which follows from the bound in (1.8) on $G_p$. This last inequality with the Hardy inequality bound in (1.9) and the bound in the top bullet of the lemma lead directly to the bound that is claimed by the third bullet of the lemma.

*Part 2*: To further exploit (4.59)-(4.61), let $r_0$ again denote a small positive number such that the set of vectors in $TX|_Z$ with norm less than $100r_0$ is well inside the fiber subbundle that is mapped diffeomorphically onto a tubular neighborhood of Z by the metric's exponential map. (This number has no R dependence.) Reintroduce the standard bump function $\chi$ from Section 1d and then let $\chi_0$ denote the function $\chi(\frac{2}{r_0}\mathrm{dist}(\cdot,Z) - 1)$. This function is equal to 1 where the distance to Z is less than $\frac{1}{2}r_0$, and it is equal to 0 where the distance to Z is greater than $r_0$. Meanwhile, given a point $p \in X$, let $G_p$ denote the Green's function on X with pole at p for the operator $d^\dagger d + c_0^{-1}R^2 r_0^2$. Take p to have distance greater than $2r_0$ from Z and then multiply both sides of (4.59) by $(1-\chi_0)G_p$. Having done that, integrate both sides over X and then integrate by parts to see that

$$|\psi|(p) \le c_0\frac{1}{r_0^2}$$

(4.73)

where the distance to Z is greater than $2r_0$. (Use (1.8) for the bound on $G_p$ and on its derivatives.)

The preceding pointwise bound will now be boot-strapped using the maximum principle: Let $u$ denote the solution to the equation $d^\dagger d\,u + c_0^{-1}R^2 r_0^2 u = 0$ on the domain in X where the distance to Z is greater than $2r_0$ subject to the constraint that $u = c_0 r_0^{-2}$ on the boundary of this domain. It then follows from (4.70)) and (4.71) and (4.73) using the maximum principle that $|\psi| \le u$ in this same domain. Meanwhile, the function $u$ obeys bounds of the form

$$|u| \le c_0 \frac{1}{r_0^2} \exp(-c_0^{-1} R\, r_0^2) \tag{4.74}$$

where the distance to Z is greater than $4r_0$ (use the Green's function bounds in (1.8) to prove this). Thus, $|\psi|$ also obeys the same bound, which is to say that

$$|\psi| \le c_0 \exp(-c_0^{-1} R) \tag{4.75}$$

where the distance to Z is greater than $4\,r_0$ (keep in mind that $r_0 \ge c_0^{-1}$).

*Part 3*: To consider $\psi$ where the distance to Z is at most $4\,r_0$, fix for the moment a positive number greater than 1 to be denoted by $c$ and use it to define a function (to be denoted by $f$) by the rule

$$f = e^{-\frac{1}{c} R |\varsigma^+|^2} \tag{4.76}$$

Since $|\varsigma^+| \ge c_0^{-1}\mathrm{dist}(\cdot, Z)$ and $|d|\varsigma^+|| \ge c_0^{-1}$ on and near Z, this function obeys a bound that has the form

$$d^\dagger d f + c_0 \frac{1}{c^2} R^2 \mathrm{dist}(\cdot, Z)^2 f \ge 0 \tag{4.77}$$

where the distance to Z is between $c_0\, c^{1/2} \frac{1}{\sqrt{R}}$ and $100\,r_0$.

To continue, fix a positive number to be denoted by $z$ so that the inequality in (4.70) holds where $\mathrm{dist}(\cdot,Z) \ge z \frac{1}{\sqrt{R}}$. Let X denote the supremum of $\psi$ on locus where $\mathrm{dist}(\cdot, Z) = z \frac{1}{\sqrt{R}}$, which is the boundary of this region. (The third bullet of Lemma 4.10 says that $X \le c_0\, \kappa_\diamond (1 + \kappa_\ddagger) R$.) Now if $c > c_0\, \kappa_\diamond$ and if $z > c_0\, c^{1/2}$, then (4.67) and (4.60) imply that the function $u \equiv |\psi| - c_0 X f$ is negative where $\mathrm{dist}(\cdot,Z) = z \frac{1}{\sqrt{R}}$ and that it obeys the inequality $d^\dagger d u < 0$ where $\mathrm{dist}(\cdot,Z)$ is between $z \frac{1}{\sqrt{R}}$ and $100 r_0$. Therefore, the maximum principle can be invoked to see that $u$ can not have a positive maximum in this region unless it is taken on where $\mathrm{dist}(\cdot,Z) = 100 r_0$. What with (4.75), this fact via the maximum principle implies that

$$|\psi| \le c_0 X\, e^{-\frac{1}{c_0 c} R |\varsigma^+|^2} \tag{4.78}$$

over the whole of the $\mathrm{dist}(\cdot,Z) \ge z \frac{1}{\sqrt{R}}$ part of X.

## 5. Estimating spectral flow when $\varsigma^+$ has transverse zeros

The operator $\mathcal{D}$ in the case at hand has the form that is depicted in (3.4) with (A, ω) as described in (3.20). In particular $\omega = r\varsigma^+\sigma$ and with the connection A such that $\nabla_A\sigma = 0$. What follows is the version of $\mathcal{D}$ in the case at hand:

$$\mathcal{D} = \gamma_\alpha\nabla_{A\alpha} + \tfrac{1}{\sqrt{2}}\, r\varsigma^+{}_k\rho_k[\sigma,\cdot\,] - \tfrac{1}{2}\, m\Gamma - \tfrac{1}{2}\, m\ . \tag{5.1}$$

This section derives necessary and sufficient conditions (when $m$ is small) for the existence of an $r$-independent bound on the norm of the spectral flow for $\mathcal{D}$. Implicit in what follows are that the curvature tensor of A obeys a bound $|F_A| \le \frac{1}{c_*} r$ with $c_*$ large enough to invoke the $R \ge r$ versions of Propostions 4.1 and 4.8 with the number Ξ determined by $\varsigma^+$ by requiring that $|\varsigma^+| \ge \frac{1}{\Xi}\mathrm{dist}(\cdot, Z)$. The values for $c_*$ and Ξ are fixed throughout. Various incarnations of $c_0$ depend implicitly on these values.

By way of a look ahead: The analysis that follows analyzes the spectral flow for $\mathcal{D}$ on stages of a path that form a sequential series of 1-parameter families of deformations of (5.1). These stages are described below. The results of the analysis are summarized in Proposition 5.15 which is in Subsection 5h.

<u>STAGE</u> <u>1</u>: The first stage in this series is to change (5.1) along the one-parameter family of operators with parameter t ∈ [0, 1] that is defined by the rule

$$t \to \gamma_\alpha\nabla_{A\alpha} + \tfrac{1}{\sqrt{2}}\, r\varsigma^+{}_k\rho_k[\sigma,\cdot\,] - (1-t)\tfrac{1}{2}\, m(\Gamma - 1)\ . \tag{5.2}$$

<u>STAGE</u> <u>2</u>: This stage moves the t = 1 version of (5.2) along a 1-parameter family of operators (parametrized by [0, 1] also) that is obtained by changing the Riemannian metric, the connection A (with σ being A-covariantly constant) and the $\mathcal{I}$-valued self-dual 2-form $\varsigma^+$ in a small radius tubular neighborhood of Z. The end member of this family is a t = 1 version of (5.2) which has a simple form near Z. The end member's data, consisting of a Riemannian metric, an orthogonal connection on $\mathfrak{L}$ and a version of $\varsigma^+$ are used implicitly in the subsequent steps without notational distinctions.

<u>STAGE</u> <u>3</u>: Having fixed $R_* >> r$, the third stage changes the end member operator from the second stage along the one-parameter family of operators (with parameter R from the interval $[r, R_*]$) that is defined by the rule

$$R \to \gamma_\alpha\nabla_{A\alpha} + \tfrac{1}{\sqrt{2}}\, R\cdot\varsigma^+{}_k\rho_k[\sigma,\cdot\,]\ . \tag{5.3}$$

Any sufficiently large version of $R_*$ will suffice for the applications to come. More is said in what follows about a lower bound for $R_*$

<u>STAGE 4</u>: This stage adds an endomorphism to the $R = R_*$ version of (5.3) that plays the role here of the term with $\mathrm{T}$ that appears in (3.29)–(3.31). The 1-parameter family that implements this addition is parametrized by the interval [0, 1] according to the rule below.

$$t \to \gamma_\alpha(\nabla_{A\alpha} + \tfrac{1}{2}\, t\, \mathrm{T}_\alpha) + \tfrac{1}{\sqrt{2}} R\cdot\varsigma^+{}_k\rho_k[\sigma,\cdot\,] \ . \tag{5.4}$$

The endomorphism $\mathrm{T}$ is described momentarily. It is defined using only the Riemannian metric and $\varsigma^+$; as such, its norm and those of its derivatives have an $r$ and R independent bound.

<u>STAGE 5</u>: Having fixed $T >> R_*$, the fourth stage in the series moves the $R = R_*$ version of (5.3) along the 1-parameter family of operators with parameter $t \in [0,T]$ that is defined using the rule

$$t \to \gamma_\alpha(\nabla_{A\alpha} + \tfrac{1}{2}\, \mathrm{T}_\alpha) + \tfrac{1}{\sqrt{2}} R_*\varsigma^+{}_k\rho_k[\sigma,\cdot\,] - \tfrac{1}{2}\, t\,\Gamma \ . \tag{5.5}$$

<u>STAGE 6</u>: This stage modifies the $t = T_*$ version of (5.5) along a 1-parameter family that decrease $R_*$ to zero, removes the $\mathrm{T}$ term and moves the connection A to a chosen fiducial, connection on P.

<u>STAGE 7</u>: This last stage moves the $t = T_*$ version of (5.5) by decreasing the parameter t to a very small but still positive value.

As explained below, if $r > c_0$ and $m < c_0^{-1}$, and if $R_*$ and subsequently $T_*$ are sufficiently large, then only the fifth stage (where t is increased in (5.5)) can have significant spectral flow, the absolute value of the spectral flow in the other stage being a priori bounded. See Proposition 5.15 in Section 5h for a summary of the subsequent spectral flow analysis.

### a) STAGES 1-5 for $\mathcal{I}$-valued sections

By way of a reminder: The splitting of ad(P) as $\mathcal{L} \oplus \mathcal{I}$ is A-covariantly constant and the operator $\mathcal{D}$ preserves this splitting. Moreover, the restriction of $\mathcal{D}$ to the space of $\mathcal{I}$-valued sections of $(\Lambda^+ \oplus \mathbb{R}) \oplus T^*X$ is independent of $r$ and the connection A. In addition, STAGES 1-5 use families of operators that each map $\mathcal{I}$-valued sections of the bundle $(\Lambda^+\oplus \underline{\mathbb{R}})\oplus T^*X$ to $\mathcal{I}$-values sections of this same bundle; and each likewise maps $\mathfrak{L}$-valued sections to $\mathfrak{L}$-valued sections. And, the restriction of each operator in these stages to the space of $\mathcal{I}$-valued sections of $(\Lambda^+\oplus \underline{\mathbb{R}})\oplus T^*X$ is independent of $r$ and $\sigma$ and the connection A (and the parameter R in STAGE 3) since these restrictions have the form $\gamma_\alpha\nabla_\alpha - \frac{1}{2}(t\Gamma+t')$ with $\nabla$ defined solely by the metric's Levi-Civita connection and with t and t′ being real numbers.

Because of this, there is, given $m$, a $c_0$ bound on the absolute value of the spectral flow along the [0, 1]-parameterized path of operators on the space of $\mathcal{I}$-valued sections of $(\Lambda^+\oplus \underline{\mathbb{R}})\oplus T^*X$ that moves $m$ to zero:

$$t \to \gamma_\alpha\nabla_\alpha - \tfrac{1}{2}(1-t)\,m(\Gamma+1)) \,. \tag{5.6}$$

By the same token, having specified a [0, 1]-parameterized path of Riemannian metrics on X, there is a corresponding $c_0$ bound on the absolute value of the spectral flow of the resulting path of operators (the parameter t version being the parameter t metric's version of $\gamma_\alpha\nabla_\alpha$) acting on this same space of $\mathcal{I}$-valued sections of $(\Lambda^+\oplus \underline{\mathbb{R}})\oplus T^*X$.

Thus, the absolute value of the spectral flow in STAGES 1 and 2 have at most a $c_0$ contribution to the total spectral flow. Meanwhile, there is no spectral flow for the operators in STAGE 3 acting on $\mathcal{I}$-values sections of $(\Lambda^+\oplus \underline{\mathbb{R}})\oplus T^*X$. As for the STAGE 4 spectral flow, this is bounded by $c_0$ if $|\mathrm{T}|$ and $|\nabla\mathrm{T}|$ have a $c_0$ upper bound (which will be the case).

As for STAGE 5: The absolute value of the spectral flow in STAGE 5 is bounded for the STAGE 5 operators acting on $\mathcal{I}$-valued sections of $(\Lambda^+\oplus \underline{\mathbb{R}})\oplus T^*X$ because each positive t operator in this family has trivial kernel on the space of $\mathcal{I}$-valued sections of $(\Lambda^+\oplus\underline{\mathbb{R}})\oplus T^*X$. (The square of the parameter t operator in this step on the space of $\mathcal{I}$-valued sections of $(\Lambda^+\oplus \underline{\mathbb{R}})\oplus T^*X$ is $(\gamma_\alpha\nabla_\alpha)^2 + t^2$ which is strictly positive when $t > 0$.)

**b) Using Proposition 4.8 for spectral flow bounds**

Proposition 4.8 will be used to bound the absolute value of the spectral flow along various paths of operators with each operator along the path acting on the space of $\mathfrak{L}$-valued sections of $(\Lambda^+\oplus \underline{\mathbb{R}})\oplus T^*X$ and with each such operator having form

$$\mathcal{D}_\Diamond = \gamma_\alpha\nabla_{A\alpha} + \tfrac{1}{\sqrt{2}}\mathrm{R}\varsigma^+{}_k\rho_k[\sigma,\cdot\,] + \mathfrak{Q}_0 - q_0\Gamma - q_1 \tag{5.7}$$

with $\mathfrak{Q}_0$ denoting an endomorphism with compact support in X–Z and with $q_0$ and $q_1$ being constants with norm bounded by 1. In this regard, the metric on X, the connection on P, the number R, the section $\varsigma^+$ of $\Lambda^+\otimes\mathcal{I}$, the unit norm section $\sigma$ of $\mathrm{ad}(P)\otimes\mathcal{I}$, the endomorphism $\mathfrak{Q}_0$ and both $q_0$ and $q_1$ can depend on the path parameter.

Two sets of constraints are imposed on the variation of this data with respect to the path parameter. Here is the first set: A fixed (parameter independent) version of Z (which is a finite, disjoint union of embedded curves) is the zero locus of $\varsigma^+$; and $\varsigma^+$ vanishes transversally along each component of its zero locus. Moreover, each version of $\varsigma^+$ obeys $|\varsigma^+| \geq \frac{1}{\Xi}\,\mathrm{dist}(\cdot,Z)$ with a fixed version of $\Xi$. This implies in particular that $\varsigma^+$ can be depicted near each component of Z as in (4.1) using a non-degenerate matrix valued function on that component. This matrix valued function can depend on the parameter, but there is a parameter independent bound on the pointwise norm of its inverse.

Here is the second set of constraints: Each path parameter version of $\mathcal{D}_\Diamond$ in (5.7) satisfies the assumptions of Proposition 4.8 and thus each of these versions of $\mathcal{D}_\Diamond$ satisfies the conclusions of Proposition 4.8.

Granted the preceding assumptions, it follows from Proposition 4.8's conclusions that the spectral flow along the 1-parameter family is determined (up to a $c_0$ error) solely by the path dependence of the numbers $\{\alpha_k\}_{k=1,\ldots,p}$ that appear in Proposition 4.8. To elaborate: The contribution of any given $\alpha_k$ to the spectral flow along the path of operators in (5.7) is within ±2 of the difference between the number of times along the path that $\alpha_k$ jumps from just below 1 to 0 minus the number of times $\alpha_k$ jumps from 0 to just below 1. The preceding implies in particular that the absolute value of the spectral flow along the path of operators is at most $c_0$ plus the $\frac{1}{2\pi}$ times the sum of the total variation of each $\{\alpha_k\}_{k=1,\ldots,p}$ along the path.

With the preceding understood, the challenge is to obtain a useful bound for the total variation of any given $\alpha_k$ along the path of operators. To this end, the key observation (which follows from (4.66)) is that Proposition 4.8's version of $\alpha$ can be written as follows:

$$\alpha = \frac{\ell}{2\pi}(\alpha_* - 2b_0) - \frac{1}{2\pi}\int_0^{\ell}(\langle\varphi_0, \mathfrak{P}\varphi_0\rangle_{\mathbb{L}} + q)ds + \mathfrak{t} \quad mod\ \mathbb{Z}. \tag{5.8}$$

where the notation is as follows: What is denoted by $\alpha_*$ is from in (4.20), and what is denoted $b_0$ comes from (4.61). Meanwhile, $\mathfrak{P}$ come from (4.62). Finally, what is denoted by $\mathfrak{t}$ has norm bounded by $c_0\frac{1}{\sqrt{R}}$, this term accounting for the sum of the contributions to $\alpha$ from the three right most terms on the left hand side of (4.66) and the appearance of $Q_s$ in (4.66).

A crucial point now with regards to (5.8) is that $q$ and $\mathfrak{t}$ and $\langle\varphi_0, \mathfrak{P}\varphi_0\rangle_{\mathbb{L}}$ are a priori $\mathbb{R}$-valued, not $\mathbb{R}/2\pi\mathbb{Z}$ valued; and this implies that the sum of the jumps of $\alpha$ along the path from 0 to just below 1 minus the sum of the jumps from just below 1 to 0 differ by at most $c_0$ from the analogous sum for the $\mathbb{R}/2\pi\mathbb{Z}$ valued number defined below

$$\alpha_{\ddagger} = \tfrac{\ell}{2\pi}(\alpha_{*} - 2b_0)\ mod\ \mathbb{Z}. \tag{5.9}$$

With regards to $b_0$: Remember that $b_0$ is constrained to lie in the interval $[0,2\pi/\ell)$. As a consequence of this, it need not vary continuously when the connection is changed along a continuous path of connections. In particular, the number $\frac{\ell}{2\pi} b_0$ can jump from just below 1 to 0 and vice versa if the corresponding path of tailored product structures for $\mathfrak{L}$ along $\mathbb{R}/(\ell\mathbb{Z})\times\{0\}$ changes discontinuously at any given point along the path of connections. (A change in the product structure can add an integer multiple of $2\pi/\ell$ to the original version of $b_0$) However, if the tailored product structures along $\mathbb{R}/(\ell\mathbb{Z})\times\{0\}$ do change continously as the connection is varied, then the $b_0$ term can be discarded from (5.9) at the expense of adding $c_0$ to a spectral flow bound that algebraically counts only jumps in the number $\alpha_*$.

With regards to variations in the number $\alpha_*$: The number $\alpha_*$ comes via (4.20). In particular, it follows from the definition of $\phi_0$ in Section 5b that $\alpha_*$ is not affected by variations in the $\lambda_k$'s that appear in (4.5), these being the square roots of the eigenvalues of $\mathbb{M}^T\mathbb{M}$. By the same token, $\alpha_*$ is not affected by any change of the number R as long as R stays positive. However, $\alpha_*$ can change due to variations in the orthogonal matrices U and V that write $\mathbb{M}$ as $U^{-1}V^{-1}DV$ because these will change the endomorphisms of $\mathbb{C}^4\oplus\mathbb{C}^4$ that appear in (4.7). It follows as a consequence of the preceding that the total variation in $\alpha_*$ has an a priori upper bound that is independent of R which is determined solely by the path in the space of smooth maps from $\mathbb{R}/(\ell\mathbb{Z})$ to $Gl(3;\mathbb{R})$) between the initial version of $\mathbb{M}$ and the final version of $\mathbb{M}$.

**c) Stage 1-4 spectral flow for $\mathfrak{L}$-valued sections**

With regards to STAGE 1: What is said in the previous section about the $q$ and t terms that appear in (5.8) implies that the absolute value of the spectral flow in this stage is bounded by $c_0 m$.

With regards to STAGE 2: Fix a positive number to be denoted by $r_1$ that is less than $\frac{1}{1000} r_\Diamond$ with $r_\Diamond$ specified in the statement of Proposition 4.8. This number $r_1$ must also be less than $\frac{1}{100} r_0$ with $r_0$ being the radius of the tubular neighborhood of Z (see Part 1 of the proof of Proposition 4.8). A smaller upper bound will be required momentarily, less than an upcoming version of $c_0^{-1}$, but in any event the value of $r_1$ is chosen once and for all (independent of parameters). The transition from the new metric, new connection and new $\varsigma^+$ to the given one will occur in the radius $2r_1$ tubular neighborhood of Z. Because of this, the STAGE 2 path of operators change only inside this tubular neighborhood. This localization is implemented using a cut-off function on X to be denoted by $\chi_1$ which is defined by the rule whereby $\chi_1(\cdot) = \chi(2 - \frac{\mathrm{dist}(\cdot,Z)}{r_1})$. When defined in this way, the function $\chi_1$ is equal to 1 where the distance to Z is greater than $2r_1$ and equal to 0 where the distance is less than $r_1$.

The specification of the new versions of metric, endomorphism $\sigma$, connection and $\varsigma^+$ is done momentarily as is the associated path of operators. To set the stage, focus attention on a given component of Z (to be denoted in what follows as Z also); and having done that, fix an oriented, orthonormal frame for the normal bundle of Z. Use this frame to identify the normal bundle with $\mathbb{R}/(\ell\mathbb{Z}) \times \mathbb{R}^3$. Let $B_0$ again denote the radius $r_0$ ball in $\mathbb{R}^3$; and, as before, use the metric's exponential map to obtain a diffeomorphism between $\mathbb{R}/\ell\mathbb{Z}\times B_0$ and the radius $r_0$ tubular neighborhood of Z in X.

Granted that the original metric on $\mathbb{R}/(\ell\mathbb{Z})\times B_0$ is depicted in (4.60), then the new metric on this same domain is given below:

$$(1+\chi_1 x_i M_i)^2 ds\otimes ds + (dx_i+\chi_1 w_{ij}x_j ds)\otimes(dx_i+\chi_1 w_{ik}x_k ds) + \chi_1 \mathfrak{g} \ . \tag{5.10}$$

The use of $\chi_1$ here makes this metric Euclidean (thus $ds\otimes ds + dx_i\otimes dx_i$) where $|x|$ is less than $r_1$, and equal to the original metrix from X where $|x|$ is greater than $2r_1$. The path from the original metric to the new metric is parameterized by $t \in [0, 1]$ and any given parameter t metric on this path is defined by replacing $\chi_1$ in (5.10) by $((1-t)+t\chi_1)$. (This will be a path of metrics on X if $r_1$ is chosen less than $c_0^{-1}$ which will be assumed.) For use momentarily, let $\wp^+(t)$ for $t \in [0,1]$ denote the orthogonal projection (as defined by the parameter t metric) from $\wedge^2 T^*X$ to the parameter t metric's version of $\Lambda^+$.

Granted that the original connection on $\mathbb{R}/(\ell\mathbb{Z}) \times B_0$ is depicted in (4.61), then the new connection on this same domain is given below:

$$A = \theta_0 + \chi_1\big((b_0 + x_i B_i)ds + \tfrac{1}{2} x_i \varepsilon_{ijk} B_k dx_j + \mathfrak{z}\big)\sigma \ . \tag{5.11}$$

Thus, the new connection is the product connection where $|x|$ is less than $r_1$, and it is the original connection from X where $|x|$ is greater than $r_1$. The path from the original connection to the new connection is parameterized by $t \in [0, 1]$ and any given parameter t connection on this path is defined by replacing $\chi_1$ in (5.11) by $((1-t)+t\chi_1)$.

The definition of the new version of $\varsigma^+$ and the path from the original version to the new one requires a brief digression which starts here with the reminder that the pull-back to $\mathbb{R}/(\ell\mathbb{Z}) \times B_0$ of the original version of $\varsigma^+$ is depicted schematically in (4.1). The new version of $\varsigma^+$ is defined with the help of a homotopy of the the $\mathbb{R}/(\ell\mathbb{Z})$-parametrized family of matrices $\mathbb{M}(\cdot)$ that appear in (4.1). This homotopy is parameterized by [0, 1] with each parameter $\tau$ member of this family being an invertible matrix. The $\tau = 1$, end-member member of the family has one of two forms: Either the end-member is constant with respect to the $\mathbb{R}/(\ell\mathbb{Z})$ parameter with all eigenvalues either +1 or all either -1, or the resulting $\mathbb{R}/(\ell\mathbb{Z})$-parametrized family is

$$\mathbb{M}(s) = \pm \begin{pmatrix} \cos(2\pi s/\ell) & -\sin(2\pi s/\ell) & 0 \\ \sin(2\pi s/\ell) & \cos(2\pi s/\ell) & 0 \\ 0 & 0 & 1 \end{pmatrix} \tag{5.12}$$

(The ± here is the sign of the determinant of the original version of $\mathbb{M}$.) The end-member in (5.12) is used when the original $\mathbb{R}/(\ell\mathbb{Z})$-parametrized family of invertible matrices from (4.1) defines the non-trivial element in its component of $\pi_1(Gl(3;\mathbb{R})$. (There are two path components $Gl(3;\mathbb{R})$ with the determinant being positive in one and negative in the other; and each component has fundamental group $\mathbb{Z}/(2\mathbb{Z})$.) Let $\underline{\mathbb{M}}$ denote smooth map from $[0,1] \times \mathbb{R}/(\ell\mathbb{Z})$ to the space of $3 \times 3$ invertible matrices that implements the desired homotopy. Thus, $\underline{\mathbb{M}}(0,\cdot)$ is the original version of $\mathbb{M}$ from (4.1); and $\underline{\mathbb{M}}(1,\cdot)$ is the final version which is either constant, or the map depicted in (5.12).

To complete the definition of the new version of $\varsigma^+$ and the path from the original version, let $\mathfrak{w}$ now denote the terms that contribute to the $\mathcal{O}(|x|^2)$ part of (4.1). The path from the original version of $\varsigma^+$ to the new version is parametrized by $t \in [0, 1]$ with the $t=0$ end-member being the original and with the general parameter t member given by the rule below

$$\varsigma^+|_t = \underline{\mathbb{M}}_{ij}((1-\chi_1)t, \cdot)\, x^i\, \wp^+(t)\cdot\big((ds \wedge dx^j + \tfrac{1}{2}\,\varepsilon_{jkm}\, dx^k \wedge dx^m) + ((1-t)+t\chi_1)\mathfrak{w}\big). \tag{5.13}$$

This will be a path of sections of $\Lambda^+ \otimes \mathcal{I}$ that vanish only on Z and vanish transversally there if $r_1$ is chosen less than $c_0^{-1}$. This upper bound for $r_1$ will be assumed.

The path of metrics, connections and section $\varsigma^+$ defined as just now provides, for each $t \in [0, 1]$, a corresponding version of the operator $\mathcal{D}_{\ddagger}$; and if $r > c_0$ and $|F_A| \le c_0^{-1} r$, then Proposition 5.8 can be invoked for each such $\mathcal{D}_{\ddagger}$. Doing so, it follows from what is said in the previous section (Section 5b) that these modifications can result in some spectral flow; but it also follows from what is said in Section 5b that the absolute value of this spectral flow will be at most $c_0$ in any event.

With regards to STAGE 3: The family of STAGE 3 operators are parametrized by the value of R in (5.3). Each of these operators is described by Proposition 4.8; any two differ only with respect to their respective R values. As a consequence (see the previous subsection), changing the value of R results in a spectral flow that is bounded by $c_0$ no matter how large the final R value (this is because changing R does not change the value of $\alpha_*$ in (5.8).)

With regards to STAGE 4: To set up the spectral flow for this stage, introduce by way of notation $\chi_\bullet$ to denote the function $\chi_1(\cdot) = \chi(2 - \frac{100\,\mathrm{dist}(\cdot,Z)}{r_1})$. When defined in this way, the function $\chi_\bullet$ is equal to 1 where the distance to Z is greater than $\frac{1}{50} r_1$ and equal to 0 where the distance is less than $\frac{1}{100} r_1$. To continue the set-up, let $\hat{\varsigma}$ denote $\frac{\varsigma^+}{|\varsigma^+|}$, this being an $\mathcal{I}$-valued section of $\Lambda^+$ over X–Z with length equal to 1. What is denoted by $T_\alpha$ in (5.4) and (5.5) is the following $\mathcal{I}$-valued endomorphism of $(\Lambda^+ \oplus \underline{\mathbb{R}}) \oplus T^*X$:

$$T_\alpha = \chi_\bullet(\nabla_\alpha \hat{\varsigma})^k \hat{\varsigma}^j \rho_k \rho_j \ . \tag{5.14}$$

As with the previous stages, it follows from what is said in Section 5b that the spectral flow along this STAGE 4 family is bounded by $c_0$ (this is independent of $R_*$).

There are two aspects of the end-member operator for the STAGE 4 family that should be kept in mind: The first is that $\hat{\varsigma}^j \rho_j$ anti-commutes with $\gamma_\alpha(\nabla_{A\alpha} + \frac{1}{2} T_\alpha)$ where the distance to Z is greater than $\frac{1}{50} r_1$. (This is because $\hat{\varsigma}^j \rho_j$ commutes with each of $\{\nabla_{A\alpha} + \frac{1}{2} T_\alpha\}_{\alpha=1,2,3,4}$.) The second is that the connection A that appears in the end-member version of (5.4) (actually all t versions) is the product connection where the distance to Z is less than $r_1$ and thus where the distance to Z is between $\frac{1}{50} r_1$ and $r_1$. Also, the Riemannian metric is the product metric from $\mathbb{R}/(\ell\mathbb{Z}) \times B_0$ on this same region; and $\varsigma^+$ here has the canonical form with $\mathbb{M}$ either the identity matrix or as depicted in (5.12).

**d) Localization of the STAGE 5 spectral flow on $\mathcal{L}$-valued sections**

The spectral flow in STAGE 5 concerns the 1-parameter family of operator that is depicted in (5.5) acting on the space of $\mathcal{L}$-valued sections of $(\Lambda^+ \oplus \underline{\mathbb{R}}) \oplus T^*X$. Keep in mind with regards to this family that the connection, the metric, and $\varsigma^+$ are fixed to be the end-members of their respective 1-parameter families that are described at the end of the previous section. Meanwhile, the number $R_*$ that appears in (5.5) is very large. Minimum allowed values for $R_*$ will be increased as needed. In any event, $R_*$ is much greater than the number $r$. With the preceding understood: The proposition that follows directly in this subsection is a localization assertion for eigenvectors akin to Proposition 4.1 which allows for arbitrarilly large values of the parameter t in (5.5)'s 1-parameter family.

**Proposition 5.1**: *There exists* $\kappa > 1$ *such that if* $R_*$ *is sufficiently large (given the Riemannian metric, the connection A, and the section* $\varsigma^+$*), then the following is true: Fix* $t \in [0, \infty)$ *to define the operator in (5.5), and suppose that* $\psi$ *is an eigenvector of this operator with eigenvalue* E *obeying* $|E| \le \frac{1}{\kappa^2} R_*$. *For any given positive number* r,

$$|\psi| \le \kappa R_* e^{-\sqrt{R_*}\, r/\kappa} \int_X |\psi|^2 \quad \textit{on the part of } X \textit{ where } \quad \Big| |\varsigma^+| - \tfrac{1}{2\sqrt{2}} \tfrac{t}{R_*} \Big| \ge r.$$

This proposition is proved momentarily.

Looking ahead: The crucial point is that Proposition 5.1 localizes the spectral flow for values of t on the order of $R_*$ to level sets of the norm of $\varsigma^+$. This localization phenomena is exploited in the subsequent subsections to construct a version of 'excision' for the spectral flow in (5.5)'s family. This excision is then used to compare the spectral flow for three versions of (5.5)'s family, the first being the given version and the second being a version that replaces the bundle $\mathfrak{L}$ with the bundle $\mathbb{R} \oplus \mathcal{I}$. The third version is on a new Riemannian manifold with a nowhere zero version of $\varsigma^+$. That last version is analyzed using the technology from Section 3.

***Proof of Proposition 5.1***: The proof is along the same lines as the proof of Proposition 4.1 in Section 5e. The operator that is depicted in (5.5) is denoted by $\mathcal{D}_\bullet$; and the starting point for the proof is the observation that $\mathcal{D}_\bullet{}^2$ has a Bochner-Weitzenboch formula that is akin to the one in (4.58):

$$\mathcal{D}_\bullet^2 = \nabla_A{}^\dagger \nabla_A + \mathfrak{R}' + \mathfrak{F}_A + \tfrac{1}{\sqrt{2}} R_* (\nabla_\alpha \varsigma^+{}_k) \gamma_\alpha \rho_k[\sigma, \cdot\,] + (\tfrac{1}{\sqrt{2}} R_* \varsigma^+{}_k \rho_k[\sigma, \cdot] - \tfrac{1}{2} t\, \Gamma)^2 \ , \tag{5.15}$$

where the notation here uses $\mathfrak{R}'$ to denote an endomorphism of $((\Lambda^+ \oplus \underline{\mathbb{R}}) \oplus T^*X) \otimes \mathfrak{L}$ whose norm is bounded by $c_0$. (This endomorphism is determined by the metric, $\varsigma^+$ and the choice for $r_1$.)

With (5.15) in mind, fix for the moment a number $c > 1$ so as to define the subset $X_{t,c} \subset X$ to be the domain where the inequality below holds

$$\Big| |\varsigma^+| - \tfrac{1}{2\sqrt{2}} \tfrac{t}{R_*} \Big| \ge \frac{c}{\sqrt{R_*}} \ . \tag{5.16}$$

If $\psi$ is an eigenvector of $\mathcal{D}_\bullet$ with eigenvalue E, then (5.15) on $X_{t,c}$ implies that

$$\tfrac{1}{2} d^\dagger d |\psi|^2 + |\nabla_A \psi|^2 + c_0{}^{-1} c^2 R_* |\psi|^2 - c_0 (1 + r + R_* + E^2) |\psi|^2 \le 0 \ . \tag{5.17}$$

And, if (5.17) holds on $X_{t,c}$ with $c > c_0$ and $E \le c_0 c^{-1} \sqrt{R_*}$, then the following inequality must also hold on $X_{t,c}$:

$$d^\dagger d|\psi| + R_*|\psi| \le 0.$$
(5.18)

To exploit this bound, let $\chi_*$ denote the non-negative function on X that is given by the rule $\chi_*(\cdot) = \chi(3 - \sqrt{R_*}c^{-1}||\varsigma^+| - \frac{1}{2\sqrt{2}}\frac{t}{R_*}|)$. This function is equal to zero on $X–X_{t,c}$ and it is equal to one where the left hand side of (5.16) is greater than $\frac{4c}{\sqrt{R_*}}$. Having fixed point p where $\chi_*(\cdot) = 1$, let $G_p$ denote for the moment the Greens function for the operator $d^\dagger d + R_*$ with pole at p. Multiply both sides of (5.18) by $\chi_* G_p$, integrate the result over X and then integrate by parts. Doing this gives the inequality that is claimed by the proposition for a suitable κ for values of r greater than $c_0 \frac{c}{\sqrt{R_*}}$.

For values of r less than this, the assertion of the proposition follows from the claim that $|\psi| \le c_0 R_*$ on X. To prove that claim, note first that if $R_*$ is sufficiently large, then (5.15) leads to the inequality $d^\dagger d|\psi|^2 \le c_0 R_*|\psi|^2$ which holds on all of X. With this understood: Let p denote a point in X where $|\psi|^2$ is largest; and then let $G_p$ denote the Green's function of the operator $d^\dagger d + 1$ with pole at p. Multiply both sides of the inequality $d^\dagger d|\psi|^2 \le c_0 R_*|\psi|^2$ by $G_p$ and then integrate by parts. Doing that bounds $\sup_x|\psi|^2$ by $c_0 R_*$ times the integral over X of $\frac{1}{\mathrm{dist}(\cdot,p)^2}|\psi|^2$ (this is because $G_p \le c_0 \frac{1}{\mathrm{dist}(\cdot,p)^2}$ (see (1.8)). To see about that integral, fix $\varepsilon \in (0, 1)$ for the moment and break the integral of $\frac{1}{\mathrm{dist}(\cdot,p)^2}|\psi|^2$ into the part where the distance to p is greater than $\varepsilon \frac{1}{\sqrt{R_*}}$ and the distance where it is less than that. The contribution of the former part to the bound for $\sup_X|\psi|^2$ is at most $c_0 \frac{1}{\varepsilon^2} R_*^2$; and the contribution to the latter is at most $c_0 \varepsilon^2 \sup_X|\psi|^2$. Take $\varepsilon = c_0^{-1}$ with $c_0$ chosen so that this bound is $\frac{1}{2}\sup_X|\psi|^2$. The resulting inequality asserts that $\sup_X|\psi|^2 \le \frac{1}{2}\sup_X|\psi|^2 + c_0 R_*^2$ which can't hold unless $\sup_X|\psi|^2$ is less than $c_0 R_*^2$.

**e) Comparing spectral flows**

This subsection exploits the localization results in Proposition 5.1 with the forthcoming Proposition 5.2 that compares the respective spectral flows for versions of the operator in (5.5) on $\mathfrak{L}$-valued sections of $(\Lambda^+ \oplus \underline{\mathbb{R}}) \oplus T^*X$ for different versions of X, $\mathfrak{L}$, orthogonal connections on $\mathfrak{L}$ and section $\varsigma^+$.

To set the stage for the proposition (and for subsequent developments), it proves useful to reintroduce some notation by using $\mathcal{D}_\bullet$ to again denote the operator in (5.5) when acting on $\mathfrak{L}$-valued sections of $(\Lambda^+ \oplus \underline{\mathbb{R}}) \oplus T^*X$, and to write it on X–Z as

$$\mathcal{D}_\bullet = \gamma^\alpha(\nabla_\alpha + \tfrac{1}{2}\chi_\bullet(\nabla_\alpha\hat{\rho})\hat{\rho}) + \sqrt{2}R_*|\varsigma^+|\hat{\rho} - \tfrac{1}{2}t\Gamma .$$
(5.19)

Here, $\hat{\rho}$ denotes the endomorphism $\frac{1}{2}|\varsigma^+|^{-1}\varsigma^+{}_k\rho_k[\sigma,\cdot]$ of the bundle $((\Lambda^+\oplus\mathbb{R})\oplus T^*X)\otimes\mathfrak{L}$ and $\chi_\bullet$ denotes $\chi(2-\frac{10^6\,\mathrm{dist}(\cdot,Z)}{r_1})$. With regards to $\hat{\rho}$: This endomorphism $\hat{\rho}$ is symmetric, it obeys $\hat{\rho}\,\hat{\rho}=1$, it commutes with $\Gamma$, and it anti-commutes with each $\gamma_k$.

Any given version of $\mathcal{D}_\bullet$ is defined from a data set $\mathrm{D}=(X, P, \mathcal{I}, A, \sigma, \varsigma^+)$ whose constituents are as follows: First, there is the compact, oriented Riemannian 4-manifold X, the principal SO(3) bundle $P\to X$ and the real line bundle $\mathcal{I}\to X$. Then, A denotes a connection on P, and $\sigma$ denotes an A-covariantly constant isometry from $\mathcal{I}$ to P's associated Lie algebra bundle. Finally, $\varsigma^+$ denotes an $\mathcal{I}$-valued section of $\Lambda^+$ with transversal zero locus

Proposition 5.2 considers two of these data sets, D and D´, subject to the constraint that there exist a non-negative number $c_0$ and a number $c_1>c_0$ such that the data sets are isomorphic on the respective domains in X and X´ where the norms of $\varsigma^+$ and $\varsigma'^+$ are in the interval $[c_0, c_1]$. Isomorphic in this context means that there exists an isometry from a neighborhood of the X version of the domain to a neighborhood of the X´ version, and then there exist respective bundle isomorphisms over the domain in X from P and $\mathcal{I}$ to the pull-backs of P´ and $\mathcal{I}'$ that identify A with the pull-back of A´ and $\sigma$ with the pull-back of $\sigma'$ and $\varsigma^+$ with the pull-back of $\varsigma'^+$.

**Proposition 5.2**: *Suppose that* D *and* D´ *are data sets as described above and that numbers* $c_0>0$ *and* $c_1>c_0$ *are given such that* D *and* D´ *are isomorphic on neighborhoods of the respective domains where* $|\varsigma^+|$ *and* $|\varsigma'^+|$ *are in* $[c_0, c_1]$. *There exists* $\kappa>0$ *such that if* $R_*>\kappa$, *then the following hold:*

- *There exists a positive number (to be denoted by* $E_0$*) which is at most* $e^{-\sqrt{R_*}/\kappa}$ *and such that neither* $E_0$ *nor* $-E_0$ *is an eigenvalue of the* D *or* D´ *version of* $\mathcal{D}_\bullet$ *when* t *is either* $2\sqrt{2}R_*c_0$ *or* $2\sqrt{2}R_*c_1$. *Moreover, the respective dimensions of the span of the eigenvectors of the* D *and* D´ *versions of* $\mathcal{D}_\bullet$ *with eigenvalue between* $-E_0$ *and* $E_0$ *are equal when* t *is* $2\sqrt{2}R_*c_0$*; and these dimensions are also equal when* t *is* $2\sqrt{2}R_*c_1$.
- *The respective spectral flows for the* D *and* D´ *versions of* $\mathcal{D}_\bullet$ *across* $E_0$ *as* t *varies from* $2\sqrt{2}R_*c_0$ *to* $2\sqrt{2}R_*c_1$ *are the same, as are the respective flows across* $-E_0$ *for these same two paths of operators.*

The rest of this section is dedicated to the proof. By way of notation, the proof refers to the $L^2$ norm and its associated inner product for the space of sections of the vector bundle $((\Lambda^+\oplus\mathbb{R})\oplus T^*X)\otimes\mathfrak{L}$. The square of this norm sends any given section $\psi$ to the integral over X of the function $|\psi|^2$. Of course, there is the analogous $L^2$ norm for the primed version of this vector bundle.

***Proof of Proposition 5.2***: The proof has five parts which start momentarily. What follows directly are some preliminary observations.

By assumption, there exists $\varepsilon > 0$ such that the data sets D and D´ are identified on the respective domains in X and X´ where $|\varsigma^+|$ and $|\varsigma'^+|$ are between $c_0 - 2\varepsilon$ and $c_1 + 2\varepsilon$. With $\varepsilon$ so chosen, this identication is assumed henceforth. The common domain is denoted by U. With regards to the specific choice of $\varepsilon$: The number $\varepsilon$ can and should be chosen so that that both $c_0 - \varepsilon$ and $c_1 + \varepsilon$ are regular values for the function $|\varsigma^+|$. (This is Sard's theorem.) The $|\varsigma^+| = c_0 - \varepsilon$ level set is denoted in what follows by $Y_0$, and the level set where $|\varsigma^+|$ is $c_1 + \varepsilon$ is denoted by $Y_1$.

With $\varepsilon$ chosen once and for all (specifically, it is independent of $R_*$ and any eigenvectors of either version of $\mathcal{D}_*$), the convention in what follow is that the various incarnations of the number $c_0$ can depend on $\varepsilon$.

*Part 1*: The plan of the proof is to consider respective D and D´ paths of operators from the $t = 2\sqrt{2}R_* c_0$ version of (5.19) to the $t = 2\sqrt{2}R_* c_1$ version that do more than just move the parameter t in the D and D´ versions of $\mathcal{D}_\bullet$. This new path is the concatentation of five segments. The first two segments (they are described momentarily in Parts 2 and 3 of the proof) keep t fixed at $2\sqrt{2}R_* c_0$ but modify the respective D and D´ versions of $\mathcal{D}_\bullet$ on the complement of the common domain where $|\varsigma^+|$ is between $c_0 - \frac{1}{2}\varepsilon$ and $c_1 + \frac{1}{2}\varepsilon$ so that the resulting two versions of $\mathcal{D}_\bullet$ for each t between $2\sqrt{2}R_* c_0$ and $2\sqrt{2}R_* c_1$ have the *same* eigenvectors and eigenvalues for eigenvalues between -1 and 1. The third segment of the deformation then moves t for these two modified version of $\mathcal{D}_\bullet$. Of course, the respective spectral flow for this third segment of the D and D´operator paths are identical. The fourth and fifth segments undo the modifications in the first two segments while keeping t fixed at $2\sqrt{2}R_* c_1$ to obtain the respective D and D´ versions of this operator that again have the precise form as the original versions as depicted in (5.19) with the larger value of t. Thus, this new path has all of its spectral flow in the first, second and then fourth and fifth segments of the path.

*Part 2*: This part of the proof describes the first segment of the new, concatenated path for the D and D´ versions of $\mathcal{D}_\bullet$. The first segment of this path keeps t fixed at $2\sqrt{2}R_* c_0$ while changing both D and D´ data sets in a small radius tubular neighborhood of the submanifolds $Y_0$ and $Y_1$ (which are the respective $c_0 - \varepsilon$ and $c_1 + \varepsilon$ level sets of $|\varsigma^+|$ from the starting versions of D and D´). The identical change is done for both data sets so that the data sets remain identified on their common domain U. The two data sets are changed along a 1-parameter family by changing the metric and $\varsigma^+$ so that the end metric in the family is a product metric in a small radius tubular neigborhood of these two level sets. (The data set element $\varsigma^+$ has to change with the metric so that each parameter's version is self-dual for the corresponding metric.) An elaboration with regards to the modifications to D (the modifications to D´ are identical) is given below for the changes near the submanifold $Y_0$ only; this is because the elaboration for what happens near $Y_1$ is identical but for some notation.

To start the elaboration: The exponential map from the original metric near $Y_0$ when restricted to the normal bundle to $Y_0$ (oriented by $d|\varsigma^+|$) identifies a neighborhood of $Y_0$ with a neighborhood of $Y_0 \times \{0\}$ in $Y_0 \times \mathbb{R}$; this identification being an isometry along $Y_0$ with the product metric on $Y_0 \times \mathbb{R}$. Fix a small positive number to be denoted by $r_*$ so that this tubular neighborhood of $Y_0 \times \{0\}$ contains the domain $Y_0 \times (-r_*, r_*)$ and so that $|\varsigma_+|$ on the latter domain is between $c_0 - \frac{101}{100}\varepsilon$ and $c_0 - \frac{99}{100}\varepsilon$. The path of data sets that are described momentarily are such that each metric along the path differs from the original metric only on $Y_0 \times (-\frac{3}{4}r_*, \frac{3}{4}r_*)$. And, whereas the original metric is isometric to the product metric only along $Y_0$, the end member metric of the family *is* the product metric on the whole of $Y_0 \times (-\frac{1}{2}r_*, \frac{1}{2}r_*)$.

The desired family of metrics is parametrized by $[0,1]$. To describe this family, let x denote the Euclidean coordinate on the product manifold $Y_0 \times \mathbb{R}$. Letting $\mathfrak{g}_0$ denote the metric on $Y_0$, then the sought after product metric is $\mathfrak{g}_0 + dx \otimes dx$. By comparison, the pull back of the metric from X to $Y_0 \times (-r_*, r_*)$ by the exponential map has the form $\mathfrak{g}_0 + dx \otimes dx + x\,\mathfrak{h}$ where $\mathfrak{h}$ denotes here a symmetric, x-dependent tensor on Y. With this understood, then the parameter $s \in [0,1]$ metric along the 1-parameter family is the metric

$$\mathfrak{g}_0 + dx \otimes dx + (1 - s\chi(\tfrac{3x}{2r_*} - 1))\, x\, \mathfrak{h} \ . \tag{5.20}$$

Some points to keep in mind with regards to this family: Each member of this family of metrics is a symmetric tensor that differs by at most $c_0 r_*$ from the original (as measured using the original metric) and the covariant derivative of each member (with covariant derivative defined using the original metric) is bounded by $c_0$.

Since $\varsigma^+$ won't necessarily be self-dual for the metrics in the family, this $\mathcal{I}$-valued 1-form must be changed along the family also. This change replaces $\varsigma^+$ at each parameter value along the family by the $\mathcal{I}$-valued, self-dual 2-form at the given parameter value that is obtained from the original version by first taking its self-dual projection as defined by the metric at the given parameter value, then dividing by the norm of that projection (as defined by the metric at the given parameter value), and then multipying by the norm of the original version. (This multiplication and then division makes all of the $s \in [0, 1]$ versions of $\varsigma^+$ have the same norm if the metric at the given value of s is used to define the norm.)

A completely analagous deformation to a product metric and change in $\varsigma^+$ on a small radius tubular neighborhood $Y_1$ (which is the $c_1 + \varepsilon$ level set for the original version of $|\varsigma^+|$) can and should be made simultaneously. Note that by taking $r_*$ smaller if necessary, the same value of $r_*$ can be used for these constructions at both $Y_0$ and $Y_1$. Do this to keep the notation under control.

*Part 3*: As noted at the outset, the second segment in the path of operators does not change D or D´. Instead, this second segment of the path of operators modifies boundary conditions for $\mathcal{D}_\bullet$ along $Y_0$ and $Y_1$ (which are the respective $|\varsigma^+| = c_0 - \varepsilon$ and $|\varsigma^+| = c_1 + \varepsilon$ level sets of Part 1's various versions of $|\varsigma^+|$).

To set the background for the upcoming analysis: The metric used below in this Part 2 is the end member of Part 1's deformation (this will be the metric throughout the paths in Segments 2-4 of the five segment path of deformations). This metric was constructed so that $Y_0$ has an isometric diffeomorphism with the product of $Y_0$ and the interval $[-\frac{1}{2} r_*, \frac{1}{2} r_*]$, and likewise a neighborhood of $Y_1$ has an isometric diffeomorphism with the product of $Y_1$ and the interval $[-\frac{1}{2} r_*, \frac{1}{2} r_*]$. With regards to the focus below: The focus for the most part in what follows will be on the neighborhood of $Y_0$ because the analysis for the neighborhood of $Y_1$ is identical but for notation.

To set notation: With regards to the interval factor: The notation again uses x to denote the Euclidean coordinate for the $[-\frac{1}{2} r_*, \frac{1}{2} r_*]$. With regards to the $Y_0$ factor: Supposing that $\{e^k\}_{k=1.2.3}$ denotes an orthonormal frame for $T^*Y_0$, then this frame with dx constitute an orthonormal frame for $T^*X$ on $Y_0 \times [-\frac{1}{2} r_*, \frac{1}{2} r_*]$. Given this orthonormal frame for $T^*Y_0$, the notation below uses $\{\nabla_k\}_{k=1,2,3}$ to denote the corresponding components of the covariant derivative $\nabla$ along the $Y_0$ factor of $Y_0 \times [-\frac{1}{2} r_*, \frac{1}{2} r_*]$; it acts on sections $((\Lambda^+ \oplus \underline{\mathbb{R}}) \oplus T^*X) \otimes \mathfrak{L}$ over $Y_0 \times [-\frac{1}{2} r_*, \frac{1}{2} r_*]$ (The dual basis to $\{e^k\}_{k=1,2,3}$ define a basis for the tangent space to the $Y_0$ factor.)

The operator $\mathcal{D}_\bullet$ can be depicted on $Y_0 \times [-\frac{1}{2} r_*, \frac{1}{2} r_*]$ using the preceding notation:

$$\mathcal{D}_\bullet = \gamma_k(\nabla_k + \tfrac{1}{2}(\nabla_k \hat{\rho})\hat{\rho}) + \gamma_x(\nabla_x + \tfrac{1}{2}(\nabla_x \hat{\rho})\hat{\rho}) + \sqrt{2}R_*|\varsigma^+|\hat{\rho} - \tfrac{1}{2} t\Gamma . \tag{5.21}$$

A crucial player in what comes next is $\gamma_x$ times a part, of (5.21), the operator

$$D \equiv \gamma_x\big(\gamma_k(\nabla_k + \tfrac{1}{2}(\nabla_k \hat{\rho})\hat{\rho}) + \sqrt{2}R_*|\varsigma^+|\hat{\rho} - \tfrac{1}{2} t\Gamma\big) . \tag{5.22}$$

A key point about D is that its derivatives are along the $Y_0$ factor; and in this regard, its restriction to any constant x hypersurface is self-adjoint with respect to the $L^2$ inner product on the space of sections over that hypersurface of $((\Lambda^+ \oplus \underline{\mathbb{R}}) \oplus T^*X) \otimes \mathfrak{L}$. (The $L^2$ inner product here is the integral of the pointwise inner product over $Y_0 \times \{x\}$.)

Moreover, being that D is first order and elliptic, it has purely real, discrete spectrum with finite multiplicities and no accumulations. The spectrum is also unbounded from above and from below.

With regards to the spectrum: The absolute value of an eigenvalue of D on any given $x \in [-\frac{1}{2} r_*, \frac{1}{2} r_*]$ version of $Y_0 \times \{x\}$ is greater than $c_0^{-1} R_*$ when $R_* > c_0$. This is so because there is a Bochner-Weitzenboch formula for D which has the same form as (5.15); and that formula on any $x \in [-\frac{1}{2} r_*, \frac{1}{2} r_*]$ version of $Y_0 \times \{x\}$ leads to the inequality

$$\int_{Y_0 \times \{x\}} |D\psi|^2 \geq c_0^{-1} \int_{Y_0 \times \{x\}} (|\nabla^Y \psi|^2 + R_*^2 |\psi|^2) \tag{5.23}$$

when $R_* > c_0$. (This last inequality follows because the left hand side of (5.16) is greater than $c_0^{-1}$ on $Y_0 \times \{x\}$.)

There is one more key point to be made about D which is that the endomorphism $\gamma_x$ anti-commutes with D. As a consequence, if $\phi$ is an eigenvector of D on $Y_0 \times \{x\}$ with eigenvalue $\hat{E}$, then $\gamma_x \phi$ is an eigenvector of D on $Y_0 \times \{x\}$ with eigenvalue $-\hat{E}$. (If D had a kernel (which it doesn't if $R_* > c_0$), then $\gamma_x$ would map the kernel of D on $Y_0 \times \{x\}$ to itself and, since $\gamma_x^2 = -1$, its restriction to the kernel of D would give the kernel of D a symplectic form which is the pairing that sends kernel elements $\phi$ and $\phi'$ to the $L^2$ inner product of $\phi$ with $\gamma_x \phi'$.)

Let $\Pi_0$ denote the $L^2$ orthogonal projection on $Y_0 \times \{0\}$ to the span of the eigenvectors of D with positve eigenvalue (as noted twice now, there are no zero eigenvalues when $R_* > c_0$). An important point here is that $\gamma_x \Pi_0 = (1 - \Pi_0)\gamma_x$ which is to say that the span of $\Pi_0$ is a Lagrangian subspace for the symplectic form defined by $\gamma_x$ on the space of sections over $Y_0 \times \{0\}$ of the bundle $((\Lambda^+ \oplus \mathbb{R}) \oplus T^*X) \otimes \mathfrak{L}$.

An analogous $L^2$ orthogonal projection is defined on $Y_1$ (which is the $|\varsigma^+| = c_1 + \varepsilon$ level set). This projection is defined using the $Y_1$ version of D, it projects to the span of the latter's negative eigenvalues on $Y_1$. (The inequality in (5.23) also holds with $Y_0$ replaced by $Y_1$, and for the same reason.) The $Y_1$ projection is denoted by $\Pi_1$ in what follows.

The projections $\Pi_0$ and $\Pi_1$ will be used directly to define a 1-parameter family of self-adjoint boundary conditions for the direct sum of $\mathcal{D}_\bullet$'s restrictions to the following codimension 0 submanifolds with boundary in X:

$$X_- = \{p \in X: |\varsigma^+| \leq c_0 - \varepsilon\} \ , \quad X_0 = \{p \in X: c_0 - \varepsilon \leq |\varsigma^+| \leq c_1 + \varepsilon\} \ , \quad X_+ = \{p \in X: |\varsigma^+| \geq c_1 + \varepsilon\} \ . \tag{5.24}$$

(A completely analogous 1-parameter family will also be defined on the corresponding three domains in X´.) To simultaneously elaborate and avoid notational clutter, subsequent identities in this part of the proof use $\mathcal{V}$ to denote the vector bundle $((\Lambda^+\oplus\mathbb{R})\oplus T^*X)\otimes\mathfrak{L}$. With this shorthand understood, use $\mathcal{C}$ to denote the direct sum of the respective sections of $\mathcal{V}$ over $X_-$, $X_0$ and $X_+$, thus

$$\mathcal{C} \equiv C^\infty(X_-;\mathcal{V}|_{X_-}) \oplus C^\infty(X_-;\mathcal{V}|_{X_0}) \oplus C^\infty(X_-;\mathcal{V}|_{X_+}) \ . \tag{5.25}$$

Let $\mathcal{L}$ denote the Hilbert space that is obtained by completing of $\mathcal{C}$ using the norm whose square sends any given element $\Psi = (\psi_-, \psi_0, \psi_+)$ to the sum of the respective $X_-$, $X_0$ and $X_+$ integrals of $|\psi_-|^2$, $|\psi_0|^2$ and $|\psi_+|^2$. (This norm on $\Psi$ is denoted below by $\|\Psi\|_{\mathcal{L}}$.)

The operator $\mathcal{D}_{\mathcal{C}}$ defined below

$$\mathcal{D}_{\mathcal{C}} \equiv \mathcal{D}_{\bullet}|_{X_-} \oplus \mathcal{D}_{\bullet}|_{X_0} \oplus \mathcal{D}_{\bullet}|_{X_+} \tag{5.26}$$

is an unbounded operator on $\mathcal{L}$ with dense domain $\mathcal{C}$.

Each parameter from the one parameter family (the parameter space is $[0,1]$) determines a subspace in $\mathcal{C}$ by specifying the behavior of sections along the $Y_0$ and $Y_1$ boundaries of the spaces in (5.24). To be sure, the boundary space here is the direct sum of the boundary spaces for each summand in (5.25), thus the space below:

$$C^\infty(Y_0;\,\mathcal{V}|_{Y_0}) \oplus \big(C^\infty(Y_0;\,\mathcal{V}|_{Y_0}) \oplus C^\infty(Y_0;\,\mathcal{V}|_{Y_1})\big) \oplus C^\infty(Y_0;\,\mathcal{V}|_{Y_1}) \ . \tag{5.27}$$

The upcoming Lemma 5.3 asserts that the restriction of the operator $\mathcal{D}_{\mathcal{C}}$ to each parameter's subspace in $\mathcal{C}$ completes to a self-adjoint (unbounded) operator on $\mathcal{L}$ with pure point spectrum having finite multiplicities and no accumulation points. The parameter 0 version of this self-adjoint operator is the original version of $\mathcal{D}_{\bullet}$ acting on the space of sections of $\mathcal{V}$ over the whole of X, whereas the parameter 1 version has no coupling between the three summands in (5.25); it is a direct sum of three self-adjoint operators, these acting on their respective summands in (5.25) with no boundary constraints from the other summands.

To present the boundary conditions for a given parameter from $[0, 1]$, let $\psi_-$, $\psi_0$ and $\psi_+$ denote for the moment sections of $((\Lambda^+\oplus\mathbb{R})\oplus T^*X)\otimes\mathfrak{L}$ over the respective subsets of X depicted in (5.24). The boundary conditions for $s \in [0, 1]$ are as follows:

- *On* $Y_0$ *(which is the* $|\varsigma^+| = c_0 - \varepsilon$ *level set)*:
  a) $\Pi_0\psi_- = (1-s)\cdot\Pi_0\psi_0$
  b) $(1-\Pi_0)\psi_0 = (1-s)(1-\Pi_0)\psi_-$ .
- *On* $Y_1$ *(which is the* $|\varsigma^+| = c_1 + \varepsilon$ *level set)*:

a) $\Pi_1\psi_0 = (1 - s)\cdot\Pi_1\psi_+$

b) $(1 - \Pi_1)\psi_+ = (1 - s)(1 - \Pi_1)\psi_0$ .

$$(5.28)$$

As noted above, Lemma 5.3 makes the formal statement to the effect that any $s \in [0,1]$ version of the boundary conditions in (5.28) result in a self-adjoint version of $\mathcal{D}_C$. The lemma refers to a subspace of $\mathcal{L}$, denoted by $\mathcal{H}$, which is the completion of $C$ using the norm whose square on any given element $\Psi = (\psi_-, \psi_0, \psi_+)$ from $C$ is the sum

$$\int_{X_-} (|\nabla\psi_-|^2 + |\psi_-|^2) + \int_{X_0} (|\nabla\psi_0|^2 + |\psi_0|^2) + \int_{X_+} (|\nabla\psi_+|^2 + |\psi_+|^2) \ . \tag{5.29}$$

This norm for $\mathcal{H}$ is denoted by $\|\Psi\|_{\mathcal{H}}$.

**Lemma 5.3**: *There exists* $\kappa > 1$ *with the following significance: Fix* $R_* > \kappa$ *and a value of* $t$ *from the interval* $[2\sqrt{2}R_*c_0, 2\sqrt{2}R_*c_1]$ *to define the operator* $\mathcal{D}_C$. *Given* $s \in [0,1]$, *let* $\mathcal{H}_s$ *denote the subspace in* $\mathcal{H}$ *whose elements obey the corresponding boundary conditions in* (5.28). *This space* $\mathcal{H}_s$ *is a closed domain for* $\mathcal{D}_C$ *acting on* $\mathcal{L}$, *and* $\mathcal{D}_C$ *on this domain is self-adjoint, with purely point spectrum having finite multiplicities and no accumulation points. Moreover, supposing that* $\Psi \in \mathcal{H}_s$, *then*

$$\|\mathcal{D}_C\Psi\|_{\mathcal{L}}^2 \geq \|\Psi\|_{\mathcal{H}}^2 - \kappa R_* \|\Psi\|_{\mathcal{L}}^2 \ ;$$

*and if, in addition,* $\Psi$ *has support only where* $|\varsigma^+| \notin [c_0 - \frac{1}{1000}\varepsilon, c_1 + \frac{1}{1000}\varepsilon]$, *then*

$$\|\mathcal{D}_C\Psi\|_{\mathcal{L}}^2 \geq \kappa^{-1}R_*^2\|\Psi\|_{\mathcal{L}}^2.$$

This lemma is proved at the end of the subsection.

A point to remember about the (5.28) family of boundary conditions: The $s = 0$ version says in effect that $\psi_- = \psi_0$ on the common $Y_0$ boundary of their respective domains of definition, and that $\psi_+ = \psi_0$ on the common $Y_1$ boundary of their respective domains of definition. This is to say that $(\psi_-, \psi_0, \psi_+)$ define an honest section over the whole of X of the bundle $((\Lambda^+ \oplus \mathbb{R}) \oplus T^*X) \otimes \mathfrak{L}$. This being the case, then the self adjoint version of $\mathcal{D}_C$ for the $s = 0$ boundary conditions is just the operator $\mathcal{D}_\bullet$ acting as it should on the space of sections of this bundle over the whole of X. By comparison, the $s = 1$ boundary conditions decouple the three summands in (5.25) which implies that the corresponding self-adjoint version of $\mathcal{D}_C$ is the direct sum of three self-adjoint versions of $\mathcal{D}_\bullet$ each acting on it's own summand in (5.25) without knowledge of the other summands.

And a final point to stress: The analogous 1-parameter family of self-adjoint operators can be constructed using the D´ data on X´; and the key observation in this regard is that the s = 1 version of this family and the s = 1 version of the D family have the same self-adjoint operator acting on the middle summand in (5.25); and its domain and action have no knowledge of the respective D and D´ operators acting on the other summands in the respective D and D´ versions of (5.25).

*Part 4*: This part of the proof supplies upper bounds for the norm of the spectral flow for the respective D and D´ paths of operators that are described in Parts 2 and 3 (which are the first and second segments of the five segment concatenation). This is done by comparing the operators along these paths with a fixed operator, the latter being the s = 0 version of the path of operators from Part 1 which is the starting operator.

The easiest segment to deal with is the first segment, it is easiest because Proposition 5.1 can be used with out modification to analyze things. To set the stage for this, introduce by way of notation the function $\chi_\varepsilon$ which is defined by the rule whereby $\chi_\varepsilon(\cdot) = \chi(\frac{5}{4\varepsilon}(c_0 - |\varsigma^+|)\chi(\frac{5}{4\varepsilon}(|\varsigma^+| - c_1)$, this being a function with two crucial properties: It has compact support where $|\varsigma^+|$ is between $c_0 - \frac{4}{5}\varepsilon$ and $c_1 + \frac{4}{5}\varepsilon$ which is in the interior of the part of what is denoted by $X_0$ in (5.24) where the metric does not change along the Part 1 metric deformation. This function is also equal to one where $|\varsigma^+|$ is greater than or equal to $c_0$ but less than or equal to $c_1$. By way of notation, the region where $|\varsigma^+|$ is between $c_0 - \frac{4}{5}\varepsilon$ and $c_1 + \frac{4}{5}\varepsilon$ is denoted by $U_*$ in what follows.

A key point now is that the metric deformation in Part 1 and the corresponding modification to $\varsigma^+$ does not depend on the value of $R_*$. As a consequence of this, the conclusions of Proposition 5.1 hold uniformly along the Part 1 path segment if $R_*$ is sufficiently large and t is set equal to $2\sqrt{2}R_* c$ with $c$ any value in the interval $[c_0, c_1]$. For any such t and for any s ∈ [0, 1] for Part 1's family, Proposition 5.1 implies this: Let $\psi$ denote an eigenvector for the corresponding version of $\mathcal{D}_\bullet$ with eigenvalue E between -1 and 1 and with $L^2$ norm equal to 1. Then

$$|(\mathcal{D}_\bullet - E)(\chi_\varepsilon\psi)| \le c_0 e^{-\sqrt{R_*}/c_0} \quad \textit{and} \quad \Big|1 - \int_{U_*} |\chi_\varepsilon\psi|^2\Big| \le c_0 e^{-\sqrt{R_*}/c_0} . \tag{5.30}$$

A crucial consequence of these two inequalities is this: Suppose again that $\psi$ is an eigenvector of some parameter s ∈ [0,1] version of $\mathcal{D}_\bullet$ with eigenvalue E from the interval between -1 and 1. Fix $\delta > 0$ and introduce $\wp_{E,\delta}$ to denote the parameter s = 0 version of the $L^2$ orthogonal projection to the span of the eigenvectors of the s = 0 version $\mathcal{D}_\bullet$ with eigenvalues that differ by *at most* $\delta$ from E. Also, let $c_1$ denote the version of the number $c_0$ that appears in (5.30). Then, by virtue of (5.30),

$$\int_{X'} |(1-\wp_{E,\varepsilon})(\chi_\varepsilon\psi)|^2 \le c_0 \tfrac{1}{\delta}\, e^{-\sqrt{R_*}/c_1}$$

(5.31)

which is an inequality whose right hand sides is at most $c_0 e^{-\sqrt{R_*}/2c_1}$ if $\delta > e^{-\sqrt{R_*}/2c_1}$. Of course, the analogous observation hold with the roles of the parameter s and the parameter 0 switched so that $\psi$ is an eigevectors of the parameter 0 version of $\mathcal{D}_\bullet$ with eigenvalue E between -1 and 1, and $\wp_{E,\delta}$ is the $L^2$ orthogonal projection to the span of the eigenvectors of the parameter s version of $\mathcal{D}_\bullet$ with eigenvalue within $\delta$ of E.

What follows is a second consequence of (5.30): If $\psi_1$ and $\psi_2$ are both eigenvectors of a parameter s version of $\mathcal{D}_\bullet$ with respective eigenvalues between -1 and 1 which are $L^2$-orthogonal on X, and which have $L^2$ norms equal to 1, then

$$\left|\int_{U_*} (\chi_\varepsilon\psi_1)^\dagger(\chi_\varepsilon\psi_1)\right| \le c_0\, e^{-\sqrt{R_*}/c_0} .$$

(5.32)

The lemma below is used momentarily to obtain the promised spectral flow bounds for the Part 1 path of operators.

**Lemma 5.4**: *Given a compact family of data sets (versions of* D*), there exists* $\kappa > 1$ *with the following significance: Fix* $R_* > \kappa$ *and* $t \ge 0$ *and a data set from the family to define the operator* $\mathcal{D}_\bullet$ *as depicted in (5.20). The dimension of the span of the eigenvectors of this* $\mathcal{D}_\bullet$ *with eigenvalue between* $-\frac{1}{\kappa}\sqrt{R_*}$ *and* $\frac{1}{\kappa}\sqrt{R_*}$ *is at most* $\kappa R_*^2$.

Assume this lemma for now. It is proved at the end of this subsection.

Keeping this lemma in mind, what follows directly are a construction and some observations to ponder: Fix a cover of the interval $[-\frac{1}{2},\frac{1}{2}]$ by intervals of length $e^{-\sqrt{R_*}/c_1}$ such that at most two intervals overlap at any given point. Each eigenvalue of the s = 0 version of $\mathcal{D}_\bullet$ between $-\frac{1}{2}$ and $\frac{1}{2}$ lands in one such interval. Now, supposing that $R_* > c_0$, then at most $c_0R_*^2$ intervals are occupied (Lemma 5.4); and because of this fact, the original, length $e^{-\sqrt{R_*}/c_1}$ intervals can be concatenated to obtain a cover of $[-\frac{1}{2},\frac{1}{2}]$ by intervals of length at most $c_0\, e^{-\sqrt{R_*}/5c_1}$ such that the endpoint of any of these longer intervals has distance greater than $e^{-\sqrt{R_*}/4c_1}$ from every eigenvalue of the s = 0 version of $\mathcal{D}_\bullet$. But this implies (via (5.31)) that the endpoints of these longer intervals have distance greater than $e^{-\sqrt{R_*}/3c_1}$ from any eigenvalue of any $s \in [0, 1]$ version of $\mathcal{D}_\bullet$. In particular, one of these longer intervals will contain the number 0 and thus there is an interval of length at most $e^{-\sqrt{R_*}/5c_1}$ containing 0 whose endpoints are not eigenvalues of any operator from the Part 1 family of operators. This implies in particular that there is no spectral flow across the endpoints of this interval (actually, this is true for any of these longer intervals). This leads to a bound on the absolute value of the spectral flow for the Part 1 family given an upper bound on the number of eigenvalues of the s = 0 version of $\mathcal{D}_{\bullet}$ in this one interval of length at most $e^{-\sqrt{R_*}/5c_1}$ that contains 0. (Such a bound is contingent on that particular operator; but a particularly relevant case is analyzed in Section 5f)

So much for the spectral flow for the Part 1 path of operators. The argument given above, copied almost verbatim, provides the identical spectral flow bounds for the Part 2 path of operators provided that (5.30) holds for the operators on the Part 2 path. Meanwhile, (5.30) will hold given the upcoming Lemma 5.6 which is an analog of Proposition 5.1 for operators on the Part 2 path. The statement of the lemma introduces the bump function $\chi_{\varepsilon'}$ to denote $\chi(\frac{100}{\varepsilon}(c_0 - |\varsigma^+|)\chi(\frac{100}{\varepsilon}(|\varsigma^+| - c_1)$ which is zero where $|\varsigma^+|$ is less than $c_0 - \frac{1}{100}\varepsilon$ and where $|\varsigma^+|$ is greater than $c_1 + \frac{1}{100}\varepsilon$. (In particular, it is zero on the support of $d\chi_\varepsilon$)

**Lemma 5.5**: *There exists* $\kappa > 1$ *such that if* $R_* \geq \kappa$*, then the following is true: Let* $\Psi = (\psi_-, \psi_0, \psi_+)$ *denote an eigenvector for the self-adjoint operator described in Lemma 5.4 with any given* $s \in [0,1]$ *version of the boundary conditions in (5.28). If the eigenvalue of* $\Psi$ *has norm less than* $\frac{1}{\kappa^2}\sqrt{R_*}$*, then*

$$\int_{X_-} |\psi_-|^2 + \int_{X_+} |\psi_+|^2 + \int_{X_0} |(1-\chi_{\varepsilon'})\psi_0|^2 \leq \kappa e^{-\sqrt{R_*}/\kappa} .$$

***Proof of Lemma 5.5***: The arguments for the proof of Proposition 5.1 that uses the Green's function for $d^\dagger d + R_*$ can be repeated to see that $|\psi_0|$ obeys the bound

$$|\psi_0| \leq c_0 e^{-\sqrt{R_*}/c_0} \tag{5.33}$$

on the domain where the distances to $Y_0$ and to $Y_1$ are greater than $\frac{1}{100} r_*$. Essentially the same argument with the same Green's function leads directly to a $c_0 e^{-\sqrt{R_*}/c_0}$ pointwise bound for $|\psi_-|$ and $|\psi_+|$ on the respective parts of $X_-$ and $X_+$ where the respective distances to $Y_0$ and $Y_1$ are greater than $\frac{1}{100} r_*$. Now let $\Psi$ denote $(\psi_-, (1-\chi_{\varepsilon'})\psi_0, \psi_+)$. This is an element in $\mathcal{H}_s$ whose support avoids the region where $|\varsigma^+| \in [c_0 - \frac{1}{1000}\varepsilon, c_1 + \frac{1}{1000}\varepsilon]$. Therefore, the last assertion in Lemma 5.3 can be invoked which will happen momentarily. To do that, let E denote the eigenvalue of $\Psi$. Then $\mathcal{D}_C\Psi = E\Psi - (\gamma_\alpha \nabla_\alpha \chi_{\varepsilon'})\psi_0$, and therefore, by virtue of (5.33),

$$\| \mathcal{D}_C \Psi \|_{\mathcal{L}} \leq c_0 E \|\Psi\|_{\mathcal{L}} + c_0 e^{-\sqrt{R_*}/c_0} . \tag{5.34}$$

Granted this bound, and supposing that $E \leq c_0^{-1} R_*$, then the claim in Lemma 5.5 follows directly from the last assertion in Lemma 5.3.

*Part 5*: This part of the proof of Proposition 5.2 supplies the proofs to Lemmas 5.3 and 5.4.

***Proof of Lemma 5.3***: To check that these boundary conditions lead to a self-adjoint version of $\mathcal{D}_C$, suppose that $(\eta_-, \eta_0, \eta_+)$ defines an element from $C$ that obey the conditions in (5.28). One must first show that the sum of the respective $X_-, X_0$ and $X_+$ integrals of $\langle\eta_-, \mathcal{D}_\bullet\psi_-\rangle$, $\langle\eta_0, \mathcal{D}_\bullet\psi_0\rangle$ and $\langle\eta_+, \mathcal{D}_\bullet\psi_+\rangle$ is the same as the sum of the corresponding integrals of $\langle\mathcal{D}_\bullet\eta_-, \psi_-\rangle$, $\langle\mathcal{D}_\bullet\eta_0, \psi_0\rangle$ and $\langle\mathcal{D}_\bullet\eta_+, \psi_+\rangle$. In general, the difference between these two sums is a sum of integrals over $Y_0$ and $Y_1$ (use integration by parts); it is the sum depicted below with $\psi_{0-}$ and $\eta_{0-}$ denoting the respective restrictions of $\psi_0$ and $\eta_0$ to $Y_0$, and with $\psi_{0+}$ and $\eta_{0+}$ denoting their restrictions to $Y_1$.

$$\int_{Y_0}\langle\eta_-, \gamma_x\psi_-\rangle - \int_{Y_0}\langle\eta_{0-}, \gamma_x\psi_{0-}\rangle + \int_{Y_1}\langle\eta_{0+}, \gamma_x\psi_{0+}\rangle - \int_{Y_0}\langle\eta_+, \gamma_x\psi_+\rangle \tag{5.35}$$

The key point is that the conditions in (5.28) guarantee that the sum of the two integrals in (5.35) over $Y_0$ is zero as is the sum of the two integrals over $Y_1$.

One must also show that the boundary conditions in (5.28) are such that the domain of the adjoint of $\mathcal{D}_C$ with a given boundary condition from (5.28) is the same as that of $\mathcal{D}_C$. Formally, this is equivalent to showing that the boundary conditions in (5.28) define a Lagrangian subspace for the symplectic form on the space in (5.27) which is defined so that when $(\psi_-, (\psi_{0-}, \psi_{0+}), \psi_+)$ and $(\eta_-, (\eta_{0-}, \eta_{0+}), \eta_+)$ denote two elements in the boundary space of (5.27), then their symplectic pairing is the sum of integrals in (5.35). The task of showing this amounts to some straightforward algebra and it is left to the reader.

With regards to the Hörmander inequality $\|\mathcal{D}_C\Psi\|_{\mathcal{L}}^2 \geq \|\Psi\|_{\mathcal{H}}^2 - c_0 R_* \|\Psi\|_{\mathcal{L}}^2$: It is sufficient to prove this for elements in $C$ that obey (5.28) for the given value of s. To see that this holds, one need only confirm this ineequality for elements with compact support where the distance to $Y_0$ and $Y_1$ is at most $\frac{1}{100} r_*$. This because the general case can be proved from this special case by using appropriately chosen bump functions to write any given $\Psi$ as a sum of two elements, one with support as above and the other with compact support in the complement of $Y_0 \cup Y_1$. The Hörmander inequality for the second term follows because it holds for $\mathcal{D}_\bullet$ as can be seen using the Bochner-Weitzenboch formula. To see about a version of $\Psi$ with support very near to $Y_0$, note first that $\mathcal{D}_\bullet$ near to $Y_0$ can be written as $\mathcal{D}_\bullet = \gamma_x(\hat{\nabla}_x - D)$ with $\hat{\nabla}_x$ denoting $\nabla_x + \frac{1}{2}(\nabla_x\hat{\rho})\hat{\rho}$. Granted this identity for $\mathcal{D}_\bullet$ and granted that $\Psi$ has the form $(\psi_-, \psi_0, 0)$ with $\psi_-$ and $\psi_0$ having support very near to

$Y_0$, where this for $\mathcal{D}_*$ is valid, then integration by parts writes $\|\mathcal{D}_C\Psi\|_{\mathcal{L}}{}^2$ as a sum of integrals over the domains in $X_0$ and $X_-$ plus boundary integrals over $Y_0$ as done below:

$$\int_0^{\frac{1}{2}r_*}\int_{Y_0\times\{x\}} (|\hat{\nabla}_x\psi_0|^2 + |D\psi_0|^2 + \langle\psi_0,[\hat{\nabla}_x,D]\psi_0\rangle) + \int_{-\frac{1}{2}r_*}^{0}\int_{Y_0\times\{x\}} (|(\hat{\nabla}_x\psi_-|^2 + |D\psi_-|^2 + \langle\psi_-,[\hat{\nabla}_x,D]\psi_-\rangle)$$
$$- \int_{Y_0}\langle\psi_-,D\psi_-\rangle + \int_{Y_0}\langle\psi_0,D\psi_0\rangle \ . \tag{5.36}$$

Noting that the commutator $[\hat{\nabla}_x, D]$ is a vector bundle endomorphism which is bounded by $c_0R_*$, then the identity in (5.37) with (5.23) leads to the Hörmander inequality if the two boundary integrals sum to a non-negative number. And, this is precisely the case because the the sum of the two boundary integrals is equal to $(1 - (1\text{-}s)^2)$ times the integrals over $Y_0$ of $\langle\Pi_0\psi_0,D\Pi_0\psi_0\rangle$ -

$$- (1 - (1\text{-}s)^2)\int_{Y_0}\langle\psi_-,(1-\Pi_0)D\psi_-\rangle + (1 - (1\text{-}s)^2)\int_{Y_0}\langle\psi_0,\Pi_0D\psi_0\rangle \tag{5.37}$$

which is a sum of two non-negative integrals.

Of course, the analoguous argument works for the case when $\Psi = (0, \psi_0, \psi_+)$ with $\psi_0$ and $\psi_+$ having support very close to $Y_1$.

Granted that the Hörmander inequality holds for elements in $C$ that obey the boundary conditions in (5.28), then it must also hold for elements in $\mathcal{H}_s$ because both sides of the Hörmander inequality define bounded quadratic forms on $\mathcal{H}_s$. Then, given that the Hörmander inequality holds for elements in $\mathcal{H}_s$, it then follows by standard arguments (but for cosmetics, these are the same as those for any first order elliptic operator on a compact manifold--see [H]) that the operator $\mathcal{D}_C$ is closed on the domain $\mathcal{H}_s$, it is self-adjoint there, and its spectrum is purely point spectrum with finite multiplicities and no accumulations.

As for the final assertion of Lemma 5.3: This follows using bump functions as done previously to write any given element in $\mathcal{H}_s$ as a sum of an element with support very near to $Y_0$ and/or $Y_1$, and another element with support in the complement of $Y_0\cup Y_1$: Use (5.37) and (5.23) for the element with support near $Y_0$ (and the analog near $Y_1$); and for the element with compact support in the complement of $Y_0 \cup Y_1$, use the Bochner-Weitzenboch formula for $\mathcal{D}_\bullet$ with the knowledge that the right hand side of (5.16) is greater than $c_0^{-1}$ on the support.

***Proof of Lemma 5.4***: The lemma asserts a standard Weyl-type bound which can be proved by comparing the heat operator for $\mathcal{D}_\bullet^2$ (the operator $\exp(-\tau\mathcal{D}_\bullet^2)$ for $\tau \in (0, \infty)$) with the heat operator for the standard Laplacian (the operator $\exp(-\tau d^\dagger d)$) as is done in e.g. [CL]. To elaborate briefly: For any given $\tau \in (0,\infty)$, the operator $\exp(-\tau\mathcal{D}_\bullet^2)$ is a bounded, trace class operator on the $L^2$ completion of the space of smooth sections of $((\Lambda^+ \oplus \mathbb{R}) \oplus T^*X) \otimes \mathfrak{L}$. This operator, $K^\tau$, can be viewed as an endomorphism over $X \times X$ of the tensor product of the respective pull-backs of $((\Lambda^+ \oplus \mathbb{R}) \oplus T^*X) \otimes \mathfrak{L}$ via the two projections to X. The first important observation about $K^\tau$ is that its pointwise trace on the diagonal in $X \times X$ can be written written as $\sum_\lambda |\psi_\lambda|^2\, e^{-E_\lambda \tau}$ with this sum taken over a complete, $L^2$-orthonormal set of eigenvectors for $\mathcal{D}_\bullet^2$, and with $E_\lambda$ for any given $\lambda$ denoting the eigenvalue of the eigenvector. Meanwhile, the Bochner-Weitzenboch formula for $\mathcal{D}_\bullet^2$ can be used to prove (by comparising with that of $\exp(-\tau d^\dagger d)$) that the pointwise trace of $K^\tau$ on the diagonal obeys the bound below:

$$\mathrm{trace}(K^\tau(x, x)) \le c_0 \tfrac{1}{\tau^2}\, e^{\tau R_*} \,. \tag{5.38}$$

Taking the integral over X of both sides of (5.38) with $\mathrm{trace}(K^\tau(x,x)$ written as $\sum_\lambda |\psi_\lambda|^2\, e^{-E_\lambda \tau}$ gives an inequality that has the form

$$\sum_E n_E\, e^{-E\tau} \le \tfrac{1}{\tau^2}\, e^{\tau R_*} \tag{5.39}$$

with the sum on the left hand side being indexed by the eigenvalues of $\mathcal{D}_\bullet^2$, and with $n_E$ denoting the multiplicity of its label eigenvalue. Taking $\tau = 1/R_*$ in (5.39) now leads to a $c_0 R_*^2$ bound on the number of linearly independent eigenvectors of $\mathcal{D}_\bullet^2$ with eigenvalue less than $R_*$. (Just throw away terms in the left hand sum that correspond to eigenvectors with eigenvalue $R_*$ or greater.) Since each eigenvector of $\mathcal{D}_\bullet$ corresponds to at most two eigenvectors of $\mathcal{D}_*$, (these being $\psi \pm \sqrt{E}\mathcal{D}_\bullet\psi$ when $E \neq 0$), the bound on the number of $\mathcal{D}_\bullet^2$ eigenvectors with eigenvalue at most $R_*$ leads to a $c_0 R_*^2$ bound on the number of $\mathcal{D}_\bullet$ eigenvectors with eigenvalue having square at most $R_*$.

### f) Estimating the number of eigenvalues of $\mathcal{D}_\bullet$ near zero.

Returning to the context of Proposition 5.2: The applications to come require a useful bound on the number of linearly independent eigenvectors of the D and D´ versions of $\mathcal{D}_\bullet$ with eigenvalue near zero. In this regard, Proposition 5.2 says in effect that the number of linearly independent eigenvectors of $\mathcal{D}_\bullet$ with eigenvalue in some interval of width at most $e^{-\sqrt{R_*}/c_0}$ around zero is the same for the respective D and D´ versions of $\mathcal{D}_*$ when $R_*$ is large and when t is either $2\sqrt{2}c_0R_*$ or $2\sqrt{2}c_1R_*$; but Proposition 5.2 does not give an $R_*$-independent bound for this number. This subsection supplies a proposition that can be used to bound this number of small eigenvalues when $R_*$ is large. Looking ahead, the subsection has three parts: The first part states a lemma that elaborates on what is said about localization by Proposition 5.1. The second part of the subsection introduce a model operator that will be used to bound the eigenvalues of $\mathcal{D}_\bullet$ near zero. The third part of the subsection has the proposition that exploits the model operator to obtain desired bounds.

*Part 1*: The lemma that follows momentarily next says more about the localization of the eigenvectors that are described in Proposition 5.1. To set the notation for this lemma, reintroduce the endomorphism $\hat{\rho} \equiv \frac{1}{2} \frac{\varsigma^+_k}{|\varsigma^+|} \rho_k[\sigma, \cdot\,]$ from (5.19). (It is an endomorphism of $((\Lambda^+ \oplus \underline{\mathbb{R}}) \oplus T^*X) \otimes \mathfrak{L}$ that is defined over X–Z.) As noted subsequent to (5.19), this is symmetric, its square is the identity, and it commutes with $\Gamma$. These features imply in turn that $\hat{\rho}\,\Gamma$ is symmetric and it too has square 1. This last fact implies that $\hat{\rho}\,\Gamma$ has eigenvalues 1 and -1 at each point in X–Z. It is also the case that the corresponding eigenspaces have equal dimension (this is proved by noting that $\hat{\rho}\,\Gamma$ anti-commutes with the endomorphisms that have the form $\nu_k\rho_k$ with $\nu \in \Lambda^+$ being orthogonal in $\Lambda^+$ to $\varsigma^+$.) Granted these observations about $\hat{\rho}\,\Gamma$, the bundle $((\Lambda^+ \oplus \underline{\mathbb{R}}) \oplus T^*X) \otimes \mathfrak{L}$ can be written over X–Z as the orthogonal, direct sum of the +1 and -1 eigenbundles for $\hat{\rho}\,\Gamma$. These bundles are denoted in the lemma by $\mathcal{V}_+$ and $\mathcal{V}_-$.

**Lemma 5.6**: *There exists* $\kappa > 1$ *such that if* $R_*$ *is sufficiently large, then the following is true: Fix* $t \in [0, \infty)$ *to define the operator in (5.19). If* $\psi$ *is an eigenvector of this operator with eigenvalue* E *obeying* $|E| \le \frac{1}{\kappa^2} R_*$*, then the norm of the orthogonal projection of* $\psi$ *onto* $\mathcal{V}_-$ *is at most* $\kappa e^{-\sqrt{R_*}/\kappa} \int_X |\psi|^2$ *at all points where the distance to Z is greater than* $\frac{1}{25} r_1$.

***Proof of Lemma 5.6***: Because $\hat{\rho}^{\,2} = 1$, the right most term on the right hand side of (5.15) can be written where the distance to Z is greater than $\frac{1}{50} r_1$ as

$$2R_*^2(|\varsigma^+| - \tfrac{1}{2} \tfrac{t}{R_*} \hat{\rho}\,\Gamma)^2 \tag{5.40}$$

which is greater than $2R_*^2(|\varsigma^+| + \frac{1}{2}\frac{t}{R_*})^2$ on sections of $\mathcal{V}_-$. This being the case, the analysis leading to Proposition 5.1 can be repeated with $\chi_*$ replaced by $\chi(2 - \frac{50}{r_1}\text{dist}(\cdot,Z))$ to obtain the bounds in the lemma. (In this regard, the replacement cut-off function is equal to 1 where the distance to Z is greater than $\frac{1}{25}r_1$ and equal to zero where the distance to Z is less than $\frac{1}{50}r_1$.)

*Part 2*: The analysis of the $\mathcal{V}_+$ part of an eigenvector for $\mathcal{D}_\bullet$ requires the digression that follows to introduce a model version of $\mathcal{D}_\bullet$ for comparison purposes. To this end, let Y denote a compact, oriented 3-dimensional Riemannian manifold (the forthcoming applications have Y being $S^1\times S^2$). The version of the manifold X is $Y\times\mathbb{R}$ with its product metric. The Euclidean coordinate on the $\mathbb{R}$ factor is denoted by x.

As for the bundle $\mathcal{I}$, this is the pullback from Y via the projection of a real line bundle over Y. Likewise, the principal bundle P is pulled back from Y as is the isometric endomorphism $\sigma$: $\mathcal{I}\to \text{ad}(P)$. Finally, the connection on this pull-back of P is also pulled back from Y. This is the connection A; and as before, $\nabla_A\sigma = 0$. The bundle ad(P) then splits as the direct sum $\mathfrak{L}\oplus\mathcal{I}$ with $\mathfrak{L}$ being a real 2-plane bundle.

The definition of the model operator requires a choice of a smooth, nowhere zero function on Y to be denoted by $\mu$ in what follows. The discussion below assumes that $\mu$ is positive; the discussion when $\mu$ is negative is obtained from what is written below via some ± sign changes. Also needed is positive number to be denoted by R. A given $\mu$ and R model for (5.5)'s operator is denoted by $\mathcal{D}_\Diamond$, and it is this

$$\mathcal{D}_\Diamond = \gamma_k(\nabla_k + \tfrac{1}{2}(\nabla_k\hat{\rho})\hat{\rho}) + \text{Q} + \gamma_x\partial_x + \sqrt{2}R\mu\, x\hat{\rho} \tag{5.41}$$

with the notation as follows: First, the index k (which can have the values 1, 2, or 3) labels a local, orthonormal frame for TY; and it is summed implicitly over its possible values. What is denoted by $\nabla_k$ in (5.41) denotes the covariant directional derivative defined by the metric on Y and the given connection on $\mathfrak{L}$ along the k'th frame vector. As for $\hat{\rho}$, this denotes a symmetric endomorphism of $((\Lambda^+\oplus\underline{\mathbb{R}})\oplus T^*X)\otimes_{\mathbb{R}}\mathfrak{L}$ that can be written as $\frac{1}{2}\nu_k\rho_k[\sigma,\cdot]$ with $\nu$ denoting a unit normed section over Y of $\Lambda^+\otimes\mathcal{I}$. Note in particular that $\hat{\rho}\hat{\rho} = 1$ and that $\hat{\rho}$ anti-commutes with the $\gamma_k$ endorphisms. What is denoted by Q is a symmetric endomorphism of $((\Lambda^+\oplus\underline{\mathbb{R}})\oplus T^*X)\otimes_{\mathbb{R}}\mathfrak{L}$ which has the following properties: It anti-commutes with $\gamma_x$, $\hat{\rho}$ and $\Gamma$; and it is independent of the coordinate x. A dense domain for $\mathcal{D}_\Diamond$ is the space of compactly supported, sections over $Y\times\mathbb{R}$ of the $\hat{\rho}\Gamma = 1$ eigenbundle in $((\Lambda^+\oplus\underline{\mathbb{R}})\oplus T^*X)\otimes_{\mathbb{R}}\mathfrak{L}$. This is a $\mathbb{R}^8$-bundle over $\mathbb{R}\times Y$; it is denoted by $\mathcal{V}_+$ because it corresponds to the $\mathcal{V}_+$ from Part 1.

To understand $\mathcal{D}_\Diamond$, it is useful to first analyze the operator $\mathfrak{D}$ given below which acts on the space of maps from $\mathbb{R}$ to $\mathbb{R}^8$:

$$\mathfrak{D} = \gamma_x \partial_x + \sqrt{2}\, R_\mu\, x\, \hat{\rho} \ . \tag{5.42}$$

Here, $R_\mu$ is a positive constant. Meanwhile, $\gamma_x$ and $\hat{\rho}$ are (as before) anti-commuting endomorphisms of $\mathbb{R}^8$, respectively antisymmetric and symmetric, and with their respective squares being -1 and 1.

This operator is relevant because $\mathcal{D}_\Diamond$ can be written over any given point in Y as

$$\mathcal{D}_\Diamond = \gamma_k(\nabla_k + \tfrac{1}{2}(\nabla_k \hat{\rho})\hat{\rho}) + Q + \mathfrak{D}_\mu \tag{5.43}$$

where $\mathfrak{D}_\mu$ denotes the version of (5.42)'s operator $\mathfrak{D}$ that has $R_\mu$ being the value of $R\cdot\mu$ at the given point in Y. This identification of $\mathfrak{D}_\mu$ with the $R_\mu = R\cdot\mu$ version of $\mathfrak{D}$ also requires the choice of an isometric isomorphism between the fiber of $\mathcal{V}_+$ at the given point and $\mathbb{R}^8$. (Isomorphisms such as the latter will usually not be noted explicitly.)

Returning to $\mathfrak{D}$: The dense domain for $\mathfrak{D}$ is the space of smooth, compactly supported maps from $\mathbb{R}$ to $\mathbb{R}^8$. The lemma that follows summarizes what is needed regarding this operator. The lemma uses $\mathbb{L}$ to denote the Hilbert space completion of the space of compactly supported maps from $\mathbb{R}$ to $\mathbb{R}^8$ using the inner product that assigns the integral over $\mathbb{R}$ of $\psi^\dagger\eta$ to any given pair of maps. This inner product is denoted by $\langle \cdot , \cdot \rangle_{\mathbb{L}}$; the associated norm is denoted by $\| \cdot \|_{\mathbb{L}}$.

**Lemma 5.7**: *The operator $\mathfrak{D}$ extends to the Hilbert space $\mathbb{L}$ as an unbounded, self-adjoint operator with discrete spectrum having no accumulation points and finite multiplicities. In this regard,*

- *The kernel of $\mathfrak{D}$ is 4-dimensional over $\mathbb{R}$; and any kernel element can be written as*

$$\left(\frac{\sqrt{2}R_\mu}{\pi}\right)^{1/4} e^{-\frac{1}{\sqrt{2}}R_\mu x^2}\, \eta$$

  *where $\eta$ is an x-independent vector in the -1 eigenspace for the action of $\gamma_x\hat{\rho}$ on $\mathbb{R}^8$.*
- *An orthonormal basis of eigenvectors for the non-zero eigenvalues of $\mathfrak{D}_\mu$ is labeled by the data sets that have the form $(\varepsilon_0, n, \hat{\imath})$ with $\varepsilon_0$ from $\{1, -1\}$, with $\hat{\imath}$ from $\{1, 2, 3, 4\}$, and with n being a non-negative integer. The corresponding eigenvalue is the positive ($\varepsilon_0 = 1$) or negative ($\varepsilon_0 = -1$) square root of $2\sqrt{2}R_\mu n$.*
- *With regards to $\hat{\imath}$: This index labels the element of an orthonormal basis of the 4-dimensional space of eigenvectors with the same data set $(\varepsilon_0, n)$.*

This lemma is proved at the end of the subsection.

Returning now to $\mathcal{D}_\Diamond$, the upcoming Lemma 5.8 summarizes the story regarding this operator. What follows directly sets the stage for this lemma.

To start the stage setting: Let $P_+$ denote the principle O(8) bundle of orthonormal frames for the vector bundle $\mathcal{V}_+$. Then let $\mathbb{L}_+$ denote the associated Hilbert space vector bundle $P_+ \times_{O(8)} \mathbb{L}$. (By way of notation, the fiber metric is denoted by $\langle\, ,\rangle_{\mathbb{L}}$.) With $\mathbb{L}_+$ understood: Associate to each $y \in Y$ the $\mathbb{R}^4$-vector space given by the kernel of the corresponding version of $\mathfrak{D}_\mu$ (from (5.43)). This association determines a $\mathbb{R}^4$ vector subbundle over Y (denoted by $\mathcal{K}$) in the Hilbert space bundle $\mathbb{L}_+$. Let $\Pi$ now denote the linear map from $\mathbb{L}_+$ to itself with image $\mathcal{K}$ given by the fiberwise orthogonal projection in $\mathbb{L}_+$ to the kernel of $\mathfrak{D}_\mu$. Any given section of $\mathbb{L}_+$ that is fiberwise orthogonal in $\mathbb{L}_+$ to the kernel of $\mathfrak{D}_\mu$ will often be written as $\psi^\perp$ with the superscript $^\perp$ indicating this fiberwise orthogonality.

With regards to the bundle $\mathcal{K}$: As explained directly, the precise subspace in $\mathbb{L}_+$ at any given point in Y is determined solely by the value of $R\mu$ at that point. Even so, there is a canonical isomorphism between any two. This is so because: First, the endomorphism $\hat{\rho}\,\Gamma$ whose +1 eigenspace is $\mathcal{V}_+$ commutes with each operator from the set $\{(\nabla_k + \frac{1}{2}(\nabla_k \hat{\rho})\hat{\rho}\}$; and likewise, so does the endomorphism $\gamma_x \hat{\rho}$. Second, if $R_\mu$ and $R_\mu{}'$ are two positive numbers, and if $\psi$ is an element in the kernel of the $R_\mu$ version of the operator $\mathfrak{D}$ from Lemma 5.6, then

$$\frac{(R_\mu{}')^{1/4}}{(R_\mu)^{1/4}} e^{\frac{1}{\sqrt{2}}(R_\mu - R_\mu{}')x^2} \psi \tag{5.44}$$

is in the kernel of the $R_\mu{}'$ version of $\mathfrak{D}$. (See the upcoming (5.50).)

To continue with observations about $\mathcal{K}$: The space of sections of $\mathcal{K}$ has a covariant derivative (to be denoted by $\nabla^{\mathcal{K}}$) which is defined by the rule

$$\nabla^{\mathcal{K}}{}_k\psi = \Pi(\nabla_k + \tfrac{1}{2}(\nabla_k \hat{\rho})\hat{\rho})\psi \quad . \tag{5.45}$$

Of particular import is that this covariant derivative commutes with the isomorphism depicted in (5.44) that identifies a given $R\mu$ version of $\mathcal{K}$ with the corresponding $(R\mu)'$ version. Another important point is that $\mathcal{K}$ has the structure of a Clifford module for T*Y that comes via the action of the endomorphisms $\{\gamma_k\}_{k=1,2,3}$. Indeed, this action is well defined because these endomorphisms anti-commute with the fiberwise operator $\mathfrak{D}_\mu$. Thus, if $\mathfrak{D}_\mu\phi$ is zero at a point in Y, then any k = 1, 2, 3 version of $\mathfrak{D}_\mu(\gamma_k\phi)$ is also zero at that point.

Yet a third point regarding $\mathcal{K}$: The space of sections of $\mathcal{K}$ has an inner product that assigns to an ordered pair of sections $\phi$ and $\eta$ the integral over Y of $\phi^T\eta$. (This pointwise inner product is $\langle\phi,\eta\rangle_{\mathbb{L}}$ when $\phi$ and $\eta$ are viewed as sections of $\mathbb{L}_+$ that are pointwise annihilated by $\mathfrak{D}_\mu$).

With the preceding observations understood, let $D^{\mathcal{K}}$ denote the operator

$$D^{\mathcal{K}} = \gamma_k \nabla^{\mathcal{K}}{}_k + \Pi Q \Pi \; ; \tag{5.46}$$

it is a first order, symmetric, elliptic operator on the space of sections of the bundle $\mathcal{K}$. A crucial point now is that $D^{\mathcal{K}}$ can be viewed (via (5.44)) as being independent of the choices for R and the function $\mu$. (It does depend on $\hat{\rho}$ and, of course, Q.) In particular, the isomorphism described above between any $R\mu$ and $R'\mu'$ versions of $\mathcal{K}$ can be used to see that the eigenvalues of $D^{\mathcal{K}}$ and their eigenvectors do not depend on either R or $\mu$.

One final definition: The lemma uses $\underline{\mathbb{L}}$ to denote the completion of the space of sections of the Hilbert space bundle $\mathbb{L}_+$ using the norm whose square takes any given section $\psi$ to the integral over Y of the function $\|\psi\|_{\mathbb{L}}{}^2$. (This norm restricts to sections of $\mathcal{K}$ to give the norm for the Hermitian inner product that was described previously.)

**Lemma 5.8**: *The version of the operator $\mathcal{D}_\Diamond$ for any no-where zero function $\mu$ and endomorphism* Q *defines an unbounded, self-adjoint operator on $\underline{\mathbb{L}}$ with purely discrete spectrum of finite multiplicity and lacking accumulation points. Moreover, given that function $\mu$ and the endomorphism* Q *(and the connection on $\mathfrak{L}$), and given a positive number $\lambda$, there exists $\kappa_\lambda > 1$ such that if $R > \kappa_\lambda$, then what follows hold*:

- *Let $\psi$ denote an eigenvector of $\mathcal{D}_\Diamond$ with eigenvalue between $-\lambda$ and $\lambda$ and with $\underline{\mathbb{L}}$ norm equal to 1. Then, $\psi$ can be written as $\phi + \psi^\ddagger$ with $\phi$ being a section of $\mathcal{K}$ and an eigenvector of $D^{\mathcal{K}}$ whose eigenvalue differs from that of $\psi$ by at most $\kappa_\lambda \frac{1}{\sqrt{R}}$, and with $\psi^\ddagger$ obeying $\int_Y \|\psi^\ddagger\|_{\mathbb{L}}^2 \le \kappa_\lambda \frac{1}{R}$.*
- *Conversely, let $\phi$ denote a section of $\mathcal{K}$ and an eigenvector of $D^{\mathcal{K}}$ with eigenvalue between $-\lambda$ and $\lambda$ and with $\underline{\mathbb{L}}$ norm equal to 1. Then $\phi$ can be written as $\psi + \psi^\ddagger$ with $\psi$ denoting an eigenvector of $\mathcal{D}_\Diamond$ whose eigenvalue differs from that of $\phi$ by at most $\kappa_\lambda \frac{1}{\sqrt{R}}$, and with $\psi^\ddagger$ obeying $\int_Y \|\psi^\ddagger\|_{\mathbb{L}}^2 \le \kappa_\lambda \frac{1}{R}$.*

The rest of this part of the subsection is dedicated to first proving Lemma 5.7 and then proving Lemma 5.8.

***Proof of Lemma 5.7***: The proof is a dimension 1 version of the proof of Lemma 4.3. Indeed, by way of contrast with (4.6), the formula for $\mathfrak{D}^2$ is just

$$\mathfrak{D}^2 = -\frac{\partial^2}{\partial x^2} + \sqrt{2}\,R_\mu \gamma_x \hat{\rho} + 2R_\mu{}^2 x^2 \,. \tag{5.47}$$

This depiction implies via Lemma 4.5 the assertion in Lemma 5.7 to the effect that $\mathfrak{D}$ extends to the Hilbert space $\mathbb{L}$ as an unbounded, self-adjoint operator with discrete spectrum having no accumulation points and finite multiplicities.

To say more about the spectrum and eigenvectors of $\mathfrak{D}$ requires a short digression to describe those of $\mathfrak{D}^2$. The digression starts with the observation that $\mathfrak{D}^2$ commutes with $\gamma_x \hat{\rho}$ and so it preserves the $\gamma_x \hat{\rho}$ eigenspaces in $\mathbb{R}^8$. Letting $\iota$ denote a given eigenvalue of $\gamma_x \hat{\rho}$ (either +1 or -1), the right hand side of (5.47) on the corresponding eigenspace can be written as

$$\mathbf{a}^\dagger \mathbf{a} + (1+\iota)\sqrt{2}R \tag{5.48}$$

where $\mathbf{a}$ and its adjoint $\mathbf{a}^\dagger$ are the operators

$$\mathbf{a} = \frac{\partial}{\partial x} + \sqrt{2}\,R_\mu x \quad \textit{and} \quad \mathbf{a}^\dagger = -\frac{\partial}{\partial x} + \sqrt{2}\,R_\mu x \,. \tag{5.49}$$

Noting the commutation rule $[\mathbf{a}, \mathbf{a}^\dagger] = 2\sqrt{2}\,R_\mu$, it follows that the eigenvalues of (5.48) consist of the set $\{1+\iota+2n\}_{n=0,1,\ldots}$ with each eigenspace being 2-dimensional and with any integer n eigenvector having the form

$$\psi = (\mathbf{a}^\dagger)^n e^{-\frac{1}{\sqrt{2}}R_\mu x^2}\, \nu \tag{5.50}$$

with $\nu$ being a vector in $\mathbb{R}^8$ obeying $\gamma_x \hat{\rho}\, \nu = \iota\, \nu$.

The spectrum of $\mathfrak{D}^2$ leads to a description of the spectrum of $\mathfrak{D}$ as follows: The first point to note is that the kernels of $\mathfrak{D}^2$ and $\mathfrak{D}$ are the same, these being the n = 0 version of (5.50) with $\nu$ in the $\iota = -1$ eigenspace $\gamma_x \hat{\rho}$ ). Meanwhile, the non-zero eigenvalues are the ± square roots of the eigenvalues of $\mathfrak{D}^2$. To see about a corresponding eigenvector for $\mathfrak{D}$: Write a positive eigenvalue in $\mathfrak{D}^2$ as $\hat{E}^2$ with $\hat{E} > 0$ and let $\psi$ denote an eigenvector of $\mathfrak{D}^2$ in the $\gamma_x \hat{\rho} = -1$ eigenspace. Then $\psi \pm \frac{1}{\hat{E}} \mathfrak{D}\psi$ is an eigenvector for $\mathfrak{D}$ with eigenvalue eigenvalue $\pm\hat{E}$. Moreover, every eigenvector for $\mathfrak{D}$ with eigenvalue $\pm\hat{E}$ is obtained in this way from a unique eigenvalue $\hat{E}^2$ eigenvector of $\mathfrak{D}^2$ in the $\gamma_x \hat{\rho} = -1$ eigenspace.

***Proof of Lemma 5.8***: The proof is a replay of the proof of Lemma 4.4 with the operator $D^{\mathcal{K}}$ in (5.46) playing the role of $\gamma_s \frac{\partial}{\partial s}$. The details follow directly.

As was the case in that proof, this proof starts with the Bochner-Weitzenboch formula for $\mathcal{D}_\Diamond{}^2$ which is:

$$\mathcal{D}_\Diamond{}^2 = (D^{\mathcal{K}})^2 + \mathfrak{D}_\mu{}^2 + \sqrt{2}\, i\, R\, x\, \nabla_k \mu\, \gamma_k \hat{\rho} \ . \tag{5.51}$$

Since $D^{\mathcal{K}}$ is elliptic, this formula with (5.47) implies that $\mathcal{D}_*$ satisfies Lemma 4.5's requirements as an operator on the Hilbert space $\mathbb{L}$ and thus Lemma 4.5's conclusions.

To continue along the lines of the proof of Lemma 4.4, fix for the moment an open set in Y where the -1 eigenspace of $\gamma_x \hat{\rho}$ has an orthonormal frame, and then fix such a frame (the frame elements are denoted by $\{\nu_p\}_{p=1,2,3,4}$). There is a corresponding basis for $\mathcal{K}$ over this open set given by the sections $\{\phi_p\}_{p=1,2,3,4}$ depicted below:

$$\phi_p = \left(\frac{\sqrt{2}R\mu}{\pi}\right)^{1/4} e^{-\frac{1}{\sqrt{2}}R\mu x^2}\, \nu_p \ . \tag{5.52}$$

To say that this is a basis for $\mathcal{K}$ is to say that $\mathfrak{D}_\mu \phi_p = 0$ for $p \in \{1, 2, 3, 4\}$.

Differentiating the formulas in (5.52) leads to the following identity for each $\phi_p$:

$$(\nabla_k + \tfrac{1}{2}(\nabla_k \hat{\rho})\hat{\rho})\phi_p = (\Theta_k \cdot \phi)_p + \left(\tfrac{1}{4\mu} - \tfrac{1}{\sqrt{2}} R x^2\right)\nabla_k \mu\ \phi_p \tag{5.53}$$

where $\{\Theta_k\}_{k=1,2,3}$ are the components of an R-independent, SO(4)-Lie algebra valued 1-form on the open set (this 1-form is denoted subsequently by $\Theta$). Of particular note is that the indentity in (5.53) can be written as

$$(\nabla_k + \tfrac{1}{2}(\nabla_k \hat{\rho})\hat{\rho})\phi_p = (\Theta_k \cdot \phi)_p - \tfrac{1}{8\sqrt{2}R\mu^2}(\nabla_k \mu)\, \mathbf{a}^\dagger \mathbf{a}^\dagger \phi_p \tag{5.54}$$

which illustrates a key point: The right-most term on the right hand side of (5.53) and (5.54) is annihilated by $\Pi$, so it is in the image of $\mathfrak{D}_\mu$. A second key point is that the $\mathbb{L}$-norm of the right most term in (5.53) and (5.54) is bounded by $c_0 \frac{1}{|\mu|} |\nabla \mu|$. (The identity in (5.54) is the analog of (4.20) with $\Theta \cdot \phi$ playing the role that is played in (4.20) by $i\,\alpha_*$.)

Key points to keep in mind with regards to $\Theta$: As already noted, $\Theta$ is independent of R. And, it is the realization over the given open set of the connection on $\mathcal{K}$ that is implicit in the definition of the operator $\nabla^{\mathcal{K}}{}_k$. This is to say that the action of $\nabla^{\mathcal{K}}{}_k$ on a section $z = \sum_p z_p \phi_p$ of $\mathcal{K}$ sends that section to the section $\sum_p (\nabla_k z_p - (\Theta_k \cdot z)_p)\, \phi_p$ of $\mathcal{K}$.

Now fix $c > 1$ with a $c_0$ lower bound to be determined shortly and suppose that $\psi$ is an eigenvector for $\mathcal{D}_\Diamond$ with eigenvalue $\mathrm{E}$ obeying $|\mathrm{E}| \le \frac{1}{c}\sqrt{R}$. Write $\psi$ over this open set as $\psi = z + \psi^\perp$ with $z$ denoting $\Pi\psi$ and $\psi^\perp$ obeying $\Pi\psi^\perp = 0$ at each point. The respective $(1 - \Pi)$ and $\Pi$ projections of the equation $\mathcal{D}_\Diamond\psi = \mathrm{E}\psi$ when written in terms of $z$ and $\psi^\perp$ on the open set from (5.52) appear as follows:

- $(1 - \Pi)(\mathcal{D}_\Diamond - \mathrm{E})\,\psi^\perp = \frac{1}{8\sqrt{2}R\mu^2}\,(\nabla_k\mu)\gamma_k \sum_p z_p\, \mathbf{a}^\dagger\mathbf{a}^\dagger\phi_p$.
- $D^{\mathcal{K}}z - \mathrm{E}z = \frac{1}{8\sqrt{2}R\mu^2}\,\nabla_k\mu \sum_p \langle\gamma_k\, \mathbf{a}^\dagger\mathbf{a}^\dagger\phi_p,\, \psi^\perp\rangle_{\mathbb{L}}\phi_p$

(5.55)

As was the case in Lemma 4.4, if $|\mathrm{E}| \le c_0^{-1}\sqrt{R}$, then $(1 - \Pi)(\mathcal{D}_\Diamond - \mathrm{E})(1 - \Pi)$ is robustly invertible on the space of sections of $\mathbb{L}_+$ in the kernal of $\Pi$ at each point in Y. This is because there is the analog of (4.32) and (4.33) for such a section (denoted by $u$):

$$\int_Y \|\mathcal{D}_\Diamond u\|_{\mathbb{L}}^2 \ \ge c_0^{-1}\int_Y (\|\nabla u\|_{\mathbb{L}}^2 + \|\tfrac{\partial}{\partial x}u\|_{\mathbb{L}}^2 + R(1+R|\cdot|^2)\|u\|_{\mathbb{L}}^2)\ .$$

(5.56)

The invertibility of $(1 - \Pi)(\mathcal{D}_\Diamond - \mathrm{E})(1 - \Pi)$ can be stated formally as follows: Given a section of $\mathbb{L}_+$ (to be denoted by $q$) which obeys $\Pi q = 0$ at each point of Y, there is a unique section $u$ of $\mathbb{L}_+$ with $\Pi u = 0$ at each point and such that $(1 - \Pi)(\mathcal{D}_\Diamond - \mathrm{E})(1 - \Pi)u = q$. Moreover, the right hand side of (5.55) for this version of $u$ is bounded by $c_0 \int_Y \|q\|_{\mathbb{L}}^2$ .

What with (5.55), this invertibility of $(1 - \Pi)(\mathcal{D}_\Diamond - \mathrm{E})(1 - \Pi)$ has the following implication: Supposing that $z$ is any given section of $\mathcal{K}$, then the top bullet of (5.55) has a unique solution, $\psi^\perp$, obeying $\Pi\psi^\perp = 0$ at each point which obeys

$$\int_Y (\|\nabla\psi^\perp\|_{\mathbb{L}}^2 + \|\tfrac{\partial}{\partial x}\psi^\perp\|_{\mathbb{L}}^2 + R(1+R|\cdot|^2)\|\psi^\perp\|_{\mathbb{L}}^2) \ \le c_0\int_Y |z|^2\ .$$

(5.57)

With $\psi^\perp$ determined in this way by the section $z$ of $\mathcal{K}$, then the equation in the second bullet of (5.55) can be viewed as a linear equation for $z$ (and $\mathrm{E}$)

$$D^{\mathcal{K}}z - \mathrm{E}z = \tfrac{1}{8\sqrt{2}R\mu^2}\,\nabla_k\mu\sum_p\langle\gamma_k\,\mathbf{a}^\dagger\mathbf{a}^\dagger\phi_p, \big((1-\Pi)(\mathcal{D}_\Diamond - \mathrm{E})(1 - \Pi)\big)^{-1}(\tfrac{1}{8\sqrt{2}R\mu^2}\,\nabla_j\mu\sum_q z_q\gamma_j\,\mathbf{a}^\dagger\mathbf{a}^\dagger\phi_q)\rangle_{\mathbb{L}}\phi_p.$$

(5.58)

The right hand side of this equation when viewed as a function of $z$ defines a self-adjoint, bounded pseudo-differential operator on the Hilbert space completion of the space of sections of $\mathcal{K}$ using the norm whose square sends a section $z$ to the Y-integral of $|z|^2$.

This operator is denoted by L. With regards to L being a bounded operator: This is guaranteed when $|E| \le c_0^{-1}\sqrt{R}$. Assuming this bound for E, then the inequality in (5.57) leads directly to a $c_0 \frac{1}{R}$ bound on the Y-integral of $|L(z)|^2$.

With the preceding understood about L, fix for the moment a number δ between 0 and 1 and let $\Pi_{R,\delta}$ denote the orthogonal projection on $\mathcal{K}$ using the $\mathbb{L}$ inner product to the span of the eigenvectors of $D^{\mathcal{K}}$ with eigenvalue within δ of E. Act on both sides of (5.58) with $(1 - \Pi_{E,\delta})$ and then take the $\mathbb{L}$ norm of each side to see that

$$\int_Y \|(1-\Pi_{E,\delta})z\|_{\mathbb{L}}^2 \le c_0 \frac{1}{\delta R} . \tag{5.59}$$

This last bound leads directly to the assertion of Lemma 5.8's first bullet because for any given positive number λ, there are only finitely many eigenvalues of $D^{\mathcal{K}}$ between -λ and λ, and each has finite multiplicity. (Take δ to be less than $c_0^{-1}$ time the minimal gap between these eigenvalues.)

To see how the second bullet comes about, let $z$ now denote an eigenvector of $D^{\mathcal{K}}$ with eigenvalue $\hat{E}$ between -λ and λ, and with $\mathbb{L}$-norm equal to 1. Now it follows from (5.53) and (5.54) that if $\lambda < c_0^{-1}\sqrt{R}$, then $\mathcal{D}_\Diamond z - \hat{E}z$ is a section of $\mathbb{L}_+$ that is annihilated by Π and with a $c_0$ bound for its $\mathbb{L}$ norm. Then, as a consequence, there exists $\psi^\perp$, also annihilated by Π, that obeys the equation in the top bullet of (5.55) and (5.57). Then $\mathcal{D}_\Diamond(z+\psi^\perp) - \hat{E}(z+\psi^\perp)$, which is L(z), has $\mathbb{L}$ norm at most $c_0 \frac{1}{\sqrt{R}}$. Now suppose that $\delta > 0$ and that there is no $\mathcal{D}_\Diamond$ eigenvalue within δ of $\hat{E}$. Then, this $c_0 \frac{1}{\sqrt{R}}$ bound for the $\mathbb{L}$ norm of $\mathcal{D}_\Diamond(z+\psi^\perp) - \hat{E}(z+\psi^\perp)$ would lead directly to a $c_0 \frac{1}{\delta\sqrt{R}}$ bound on the norm of $z+\psi^\perp$ because the eigenvectors of $\mathcal{D}_\Diamond$ span the Hilbert space $\mathbb{L}$. Such a bound is nonsensical if $\delta < c_0 \frac{1}{\sqrt{R}}$ because $z$ has $\mathbb{L}$ norm equal to 1 and $\psi^\perp$ is $\mathbb{L}$-orthogonal to $z$.

*Part 3*: Return now to the context of (5.19) and the operator depicted there. Let $c$ denote a positive, regular value for the function $|\varsigma^+|$. Having fixed $R_*$ to be large (successively larger lower bounds are provided below), this last part of the subsection considers the spectrum of the version of (5.19) that is defined using the given value of $R_*$ with $t = 2\sqrt{2}R_* c$. The notation (again) has $\mathcal{D}_\bullet$ denoting this version of (5.19)'s operator.

By way of a look ahead: Proposition 5.1 says in effect that the eigenvectors of $\mathcal{D}_\bullet$ for those eigenvalues with absolute value at most $c_0^{-1}\sqrt{R_*}$ are concentrated very near the $|\varsigma^+| = c$ level set when $R_*$ is large. Now if $c$ is a regular value of $|\varsigma^+|$, then the $|\varsigma^+| = c$ level set is a smoothly embedded 3-manifold in X; and as explained below, the eigenvalues of $\mathcal{D}_\bullet$ with absolute value at most $c_0^{-1}\sqrt{R_*}$ and their eigenvectors are well approximated by those of an operator that is defined just on the normal bundle to this submanifold. This operator on the normal bundle is the version of the operator $\mathcal{D}_{\Diamond}$ in Lemma 5.8 with R being $R_*$, with $\mu$ being the norm of the 1-form $d|\varsigma^+|$ along the $|\varsigma^+| = c$ level set, and with Q manifesting the second fundamental form of the level set. Meanwhile, Lemma 5.8 says in effect that the eigenvalues of $\mathcal{D}_{\Diamond}$ with a given norm bound are very nearly those of (5.46)'s operator $D^{\mathcal{K}}$ when $R_*$ is large; and $D^{\mathcal{K}}$ is an $R_*$-independent operator that acts on sections over Y of the $R_*$-independent $\mathbb{R}^4$ bundle $\mathcal{K}$. Proposition 5.9 summarizes this with a statement which says in effect that the eigenvalues of $\mathcal{D}_{\bullet}$ with an a priori absolute value bound are (to a first approximation) independent of $R_*$ if $R_*$ is large.

To set the stage for the proposition: Let Y denote the $|\varsigma^+| = c$ level set in X. Use the 1-form $d|\varsigma^+|$ along Y to orient Y's normal bundle. Once oriented, Y's normal bundle has a canonical, isometric identification with $Y \times \mathbb{R}$ which will be implicitly used in what follows. The Euclidean coordinate for the $\mathbb{R}$ factor is again denoted by x. The Riemannian metric's exponential map along Y defines a smooth map from $Y \times \mathbb{R}$ to X which sends some positive $r_{\ddagger}$ version of $Y \times (-r_{\ddagger}, r_{\ddagger})$ diffeomorphically onto a tubular neighborhood of Y in X. (This tubular neighborhood is the set of points in X with distance less than $r_{\ddagger}$ from Y.) This diffeomorphism is used implicitly to identify the tubular neighborhood with the subset in $Y \times \mathbb{R}$ where $|x| < r_{\ddagger}$. With regards to $r_{\ddagger}$: An $R_*$-independent upper bound for this number is specified momentarily.

As for the bundle $\mathcal{V}_+$: Pull-back to Y via the exponential map defines $\mathcal{V}_+$ along $Y\times \{0\}$ in $Y\times(-r_{\ddagger},r_{\ddagger})$. Then, parallel transport from $Y \times\{0\}$ along the $(-r_{\ddagger},r_{\ddagger})$ fibers using the covariant derivative $\nabla_x + \frac{1}{2}(\nabla_x \hat{\rho})\hat{\rho}$ identifies the pull-back of $\mathcal{V}_+$ from X over $Y\times(-r_{\ddagger},r_{\ddagger})$ with its pull-back from $Y\times\{0\}$ via the projection map. (The covariant derivative $\nabla_x$ that is used here is defined using the pull-back from X of the connection on $\mathcal{V}_+$ which comes from the connection on $\mathfrak{L}$ and the Levi-Civita connection on TX.) These bundle identifications are used implicitly in what follows.

Proposition 5.9 refers to a version of the operator $D^{\mathcal{K}}$ from (5.46). In this regard: The restriction to Y of the pull-back connection from X for the X version of $\mathcal{V}_+$ is used for the covariant derivative along Y that is needed to define the operator $\gamma_k\nabla^{\mathcal{K}}{}_k$ on the space of sections of the bundle $\mathcal{K}$. (This connection comes from a given connection on $\mathfrak{L}$ and the Levi-Civita connection on TX.) As noted previously, the endomorphism Q for for the relevant version of (5.46) comes from the second fundamental form of the $|\varsigma^+| = c$ level set; it accounts for the difference between the respective Levi-Civita covariant derivatives along tangent vectors to Y as defined by the metric on Y and the metric on X.

With localization to a neighborhood of Y in mind, define a function $\chi_{\ddagger}$ on X by the rule whereby $\chi_{\ddagger}(\cdot) = \chi(\frac{100\mathrm{dist}(\cdot,Y)}{r_{\ddagger}} - 1)$, this being a function that is equal to 1 where the distance to Y is less than $\frac{1}{100}\mathrm{r}_{\ddagger}$ and equal to zero where the distance to Y is greater than twice that. Supposing that $\psi$ is a section over X of $\mathcal{V}_+$, then $\chi_{\ddagger}\psi$ with the identifications from the preceding paragraphs defines a section over $\mathrm{Y}\times\mathbb{R}$ of the $\mathrm{Y}\times\mathbb{R}$ version of $\mathcal{V}_+$. Moreover, since this section has compact support, it can (and will) be viewed equivalently as a section of the Hilbert space bundle $\mathbb{L}_+$ over Y. Conversely: Supposing that $\psi'$ denotes a section of the bundle $\mathbb{L}_+$ over Y, then $\chi_{\ddagger}\psi'$ with the preceding identifications defines a section of the X version of $\mathcal{V}_+$.

The preceding as background, what follows is the promised Proposition 5.9.

**Proposition 5.9**: *Let $c$ denote a non-zero regular value of the function $|\varsigma^+|$ and let* Y *denote the corresponding level set. Given $c$ to define the level set (and the connection on $\mathfrak{L}$), there exists $\kappa > 1$, and given also a positive number to be denoted by $\lambda$, there exists $\kappa_\lambda$; and these are such that what follows holds: Fix $\mathrm{R}_* > \kappa_\lambda$ and then use $\mathrm{t}=2\sqrt{2}\,\mathrm{R}_* c$ to define the operator $\mathcal{D}_\bullet$. Meanwhile, let* Y *denote the $|\varsigma^+| = c$ level set so as to consider the $\mathbb{R}^4$-bundle $\mathcal{K} \to \mathrm{Y}$ and the corresponding operator $\gamma_k\nabla^{\mathcal{K}}{}_k + \mathrm{Q}$ on its space of sections. Fix $\mathrm{r}_{\ddagger} = \frac{1}{\kappa}$ and use it define the function $\chi_{\ddagger}$.*

- *Let $\psi_+$ denote the $\mathcal{V}_+$ component of an eigenvector of $\mathcal{D}_\bullet$ with eigenvalue between $-\lambda$ and $\lambda$, and with the* X*-integral of the square of its norm being 1. Then $\chi_{\ddagger}\psi_+$ can be written as $\phi + \psi^{\ddagger}$ with $\phi$ being a section of $\mathcal{K}$ and an eigenvector of $\mathrm{D}^{\mathcal{K}}$ whose eigenvalue differs from that of $\psi$ by at most $\kappa_\lambda \frac{1}{\sqrt{\mathrm{R}}}$, and with $\psi^{\ddagger}$ being a section of $\mathbb{L}_+$ over* Y *that obeys $\int_{\mathrm{Y}} \|\psi^{\ddagger}\|_{\mathbb{L}}^2 \le \kappa_\lambda \frac{1}{\mathrm{R}}$.*
- *Conversely, let $\phi$ denote an eigenvector of $\mathrm{D}^{\mathcal{K}}$ with eigenvalue between $-\lambda$ and $\lambda$ and with the* Y*-integral of $|\phi|^2$ being 1. Then $\chi_{\ddagger}\phi$ can be written as $\psi + \psi^{\ddagger}$ with $\psi$ denoting an eigenvector of $\mathcal{D}_\bullet$ whose eigenvalue differs from that of $\phi$ by at most $\kappa_\lambda \frac{1}{\sqrt{\mathrm{R}}}$, and with $\psi^{\ddagger}$ obeying $\int_{\mathrm{X}} |\psi^{\ddagger}|^2 \le \kappa_\lambda \frac{1}{\mathrm{R}_*}$.*

The rest of this Part 3 proves this proposition.

***Proof of Proposition 5.9***: By way of a head's up for the proof: Except for the what is said in the subsequent two paragraphs, the proof is essentially a repetition of the proof of Lemma 5.8. Also by way of a head's up: What is denoted by $c_0$ in each instance depends on the value of $c$ and on the given connection on the X version of bundle $\mathfrak{L}$. However, the value of any incarnation of $c_0$ will not depend on $\mathrm{R}_*$ nor will it depend on any given eigenvector of $\mathcal{D}_\bullet$ or $\mathrm{D}^{\mathcal{K}}$.

The proof starts with a direct consequence of Proposition 5.1 and Lemma 5.6 which is this: Let $\psi$ denote an eigenvector for $\mathcal{D}_\bullet$ with the X integral of $|\psi|^2$ being 1, and with eigenvalue E obeying $|E| \le c_0^{-1}\sqrt{R_*}$. Let $\psi_+$ denote its $\mathcal{V}_+$ component. If $r_\ddagger \le c_0^{-1}$ and if $R_* \ge c_0$, then Proposition 5.1 and Lemma 5.6 say in effect that

$$|\psi - \chi_\ddagger\psi_+| \le e^{-\frac{1}{c_0}\sqrt{R_*}r_\ddagger} \quad \textit{and} \quad |\mathcal{D}_\bullet(\chi_\ddagger\psi_+) - E\chi_\ddagger\psi_+| \le e^{-\frac{1}{c_0}\sqrt{R_*}r_\ddagger} . \tag{5.60}$$

These bounds are exploited by viewing $\mathcal{D}_\bullet$ and $\chi_\ddagger\psi_+$ as a differential operator and section of $\mathcal{V}_+$ on the $|x| < r_\ddagger$ part of $Y\times\mathbb{R}$.

To view $\mathcal{D}_\bullet$ on the $|x| < r_\ddagger$ part of $Y\times\mathbb{R}$: The pull-back to $Y\times\mathbb{R}$ of the Riemannian metric on X can be compared with the product metric; and likewise the respective pull-backs to $Y\times\mathbb{R}$ from X of the connections on $\mathfrak{L}$ and $\mathcal{V}_+$ can be compared with their pull-backs from Y via the projection map. Do this to depict the operator $\mathcal{D}_\bullet$ on $Y \times (-r_\ddagger, r_\ddagger)$ acting on sections of $\mathcal{V}_+$. Here is the result (in schematic form):

$$\mathcal{D}_\bullet = \mathcal{D}_\Diamond + x\,(\mathfrak{P}_j\nabla_j + \mathfrak{P}_0) + \tfrac{1}{2}\,q_0\gamma_x + q_1R_*x^2\hat{\rho} \tag{5.61}$$

where the notation is as follows: First, regarding $\mathcal{D}_\Diamond$, this is the version from (5.41) that has $\mu$ being the norm of $d|\varsigma^+|$ along Y; and it has Q accounting for the difference between the respective Levi-Civita covariant derivatives along tangent vectors to Y as defined by the induced Riemannian metric on TY and by the given Riemannian metric on $TX|_Y$. Meanwhile, $\{\mathfrak{P}_j\}_{j=1,2,3}$ and $\mathfrak{P}_0$ denote smooth, $R_*$-independent endomorphisms of $\mathcal{V}_+$ defined over $Y\times[-r_\ddagger, r_\ddagger]$; and what are denoted by $q_0$ and $q_1$ denote smooth, $R_*$-independent functions on $Y\times[-r_\ddagger, r_\ddagger]$. Each of the endomorphisms and also $q_0$ and $q_1$ will, in general, depend on the x coordinate, and they will also vary over Y. (The function $q_0$ along Y is the trace of Y's second fundamental form; and $q_1$ along Y is half of the second x-derivative of $|\varsigma^+|$.) Some important points for the analysis follow directly:

- *Both* $\Pi\gamma_x\Pi$ *and* $\Pi\hat{\rho}\Pi$ *are zero (because both* $\gamma_x$ *and* $\hat{\rho}$ *anti-commute with* $\gamma_x\hat{\rho}$ *).*
- *The operator* $\Pi(x\,(\mathfrak{P}_j\nabla_j + \mathfrak{P}_0))\Pi$ *is a first order, differential operator on section of* $\mathcal{K}$ *whose coefficients are bounded by* $c_0\frac{1}{R_*}$. (*This is because the function* x on $Y\times\mathbb{R}$ *can be written using the operators in* (5.39) *as* $\frac{1}{2\sqrt{2}R\mu}(\mathbf{a} + \mathbf{a}^\dagger)$.)
- *In general, if z denotes any given section of* $\mathcal{K}$ *but viewed as a section of* $\mathbb{L}_+$ *(by writing it locally as* $z_a\phi_a$ *as in (5.55)), then*

$$\|x\,(\mathfrak{P}_j\nabla_j + \mathfrak{P}_0)z\|_{\mathbb{L}}^2 \le c_0\,\frac{1}{R_*}\Big(\int_Y (|\nabla^K z|^2 + |z|^2)\Big) .$$

- *If* $\psi$ *is any section of* $\mathbb{L}_+$ *along* Y *with* $\Pi\psi = 0$ *and with support where* $|x| \le r_\ddagger$*, then*

$$\|x(\mathfrak{P}_j\nabla_j + \mathfrak{P}_0)\psi\|_{\mathbb{L}}^2 \le c_0 r_\ddagger \Big(\int_Y (|\nabla\psi|^2 + |\psi|^2)\Big) .$$

(5.62)

Returning now to (5.60) and with (5.62) available, the proof of the proposition now repeats the arguments for the proof of Lemma 5.8. To elaborate very slightly: Start by writing $\chi_\ddagger\psi_+$ as $z + \psi^\perp$ with $z$ being the projection of $\chi_\ddagger\psi_+$ to $\mathcal{K}$. The argument for the proof of Lemma 5.8 leading to (5.57) are repeated almost verbatim using what is said in (5.62). Assuming $R_* > c_0$ and $r_\ddagger < c_0^{-1}$, then these arguments supply an analog of (5.57) that says this:

$$\int_Y (\|\nabla\psi^\perp\|_{\mathbb{L}}^2 + \|\tfrac{\partial}{\partial x}\psi^\perp\|_{\mathbb{L}}^2 + R(1 + R|\cdot|^2)\|\psi^\perp\|_{\mathbb{L}}^2) \le c_0 \int_Y (\tfrac{1}{R_*}|\nabla^K z|^2 + |z|^2) + e^{-\frac{1}{c_0}\sqrt{R_*}r_\ddagger} .$$

(5.63)

Continuing with the cosmetic modifications of the arguments for Lemma 5.8 then leads from (5.63) to an analog of (5.58) which can be written schematically as

$$D^{\mathcal{K}} z - \mathrm{E}\, z = \mathrm{L}_\ddagger(z) + \mathrm{B}$$

(5.64)

with $\mathrm{L}_\ddagger$ denoting a pseudo-differential operator and $\mathrm{B}$ denoting a $z$-independent section of $\mathcal{K}$, and with these obeying

$$\int_Y |\mathrm{L}_\ddagger(z)|^2 \le c_0 \tfrac{1}{R_*} \int_Y (\tfrac{1}{R_*}|\nabla^K z|^2 + |z|^2) \quad \textit{and} \quad \int_Y |\mathrm{B}|^2 \le e^{-\frac{1}{c_0}\sqrt{R_*}r_\ddagger} .$$

(5.65)

(These bounds follow directly using (5.63) and (5.62) and (5.60).) Granted (5.63) and (5.65), then the proof of Proposition 5.9 can be completed by directly copying what is said by the last two paragraphs of the proof of Lemma 5.8.

## g) STAGE 5 spectral flow for $\mathfrak{L}$-valued sections via comparison and excision

This subsection returns now to the original task of estimating the spectral flow along the STAGE 5 path that is depicted in (5.5). By way of a reminder, all operator on the path are versions of the operator $\mathcal{D}_\bullet$ from (5.19). In particular, all operators on the path are defined by the versions of the data set D that are identical except for the parameter t from D which varies from t = 0 to some very large T which will be much greater than $R_*$ (which in turn can be as large as desired, and in any event much greater than the number $r$ that characterizes the norm of the curvature of the original connection from the very beginning of the STAGE 1 path). The four parts of this subsection provide an estimate and bounds for the spectral flow along the STAGE 5 path by using Propositions 5.2 and 5.9 to compare the spectral flow along various parts of the STAGE 5 path with the spectral flow for paths of other versions of $\mathcal{D}_{\bullet}$.

*Part 1*: This part sets up the first comparison path. To this end, recall now that the line bundle $\mathcal{I}$ is isomorphic to the product line bundle on an $R_*$-independent tubular neighborhood of each component of the zero locus of $\varsigma^+$ (which is the set Z). Also, the Riemannian metric on the same tubular neighborhood of each component of Z is the pull-back via the tubular neighborhood projection map of the fiber metric on the component's normal bundle (which is a flat metric). With regards to the radius of these neighborhoods: The tubular neighborhood of any given component of Z contains the region where the distance to that component is less than a fixed positive number, this being the number that is denoted by $r_1$ in Section 5c.

Meanwhile, the section $\varsigma^+$ on each component's tubular neighborhood has one of two possible canonical forms. And an important point for what follows: The norm of either version of $\varsigma^+$ on the tubular neighborhood of any given component of Z is $\sqrt{2}$ times the distance to the given component. An observation now for what is to come: The non-zero values $|\varsigma^+|$ near each component are regular values of $|\varsigma^+|$ as they correspond to level sets of the distance function to the component.

As for the line bundle $\mathcal{I}$ and the principal bundle P: As explained in Part 1 of the proof of Proposition 5.8 both are isomorphic to product bundles over the radius $r_1$ tubular neighborhood of Z. Moreover, as noted in Section 5c, there exists an isomorphism of this sort that identifies the connection A with the product connection and identifies the section $\sigma$ with a constant element in $\mathfrak{su}(2)$. Since $\mathfrak{L}$ is the orthogonal complement in $\sigma$ in ad(P), this same isomorphism identifies $\mathfrak{L}$ with the product $\mathbb{R}^2$ bundle over the radius $r_1$ tubular neighorhood of Z.

What follows next describes a data set $D' = (X, P', \mathcal{I}', A', \sigma', \varsigma^+)$ that agrees with the original $D = (X, P, \mathcal{I}, A, \sigma, \varsigma^+)$ on the radius $r_1$ tubular neighborhood of Z. The description starts with the bundle ad(P′) which is extended as $\mathcal{I}\oplus(\mathcal{I}\oplus\mathbb{R})$ with the $(\mathcal{I}\oplus\mathbb{R})$ summand being the bundle $\mathfrak{L}'$. The bundle P′ is then the principal bundle of orthonormal frames in $\mathcal{I}\oplus(\mathcal{I}\oplus\mathbb{R})$. Meanwhile, the connection A′ is the flat connection that preserves all three summands, and $\sigma'$ is the isometry that sends $\mathcal{I}$ to the left most $\mathcal{I}$ summand in $\mathcal{I}\oplus(\mathcal{I}\oplus\mathbb{R})$. (The adjoint action of $\sigma'$ annihilates the first $\mathcal{I}$ summand in $\mathcal{I}\oplus(\mathcal{I}\oplus\mathbb{R})$ and it sends any given $(s_1, s_2)$ in the $\mathcal{I}\oplus\mathbb{R}$ summand to $(-s_2, s_1)$ in $\mathbb{R}\otimes\mathcal{I}$.)

The following lemma describes the spectral flow as the parameter t increases for the D′ version of the operator in (5.19).

**Lemma 5.10**: *There exists* $\kappa > 1$ *with the following significance: If* $R_* > \kappa$, *then the absolute value of the spectral flow for the* D´ *version of the family of operators depicted in (5.19) as* t *is increased from zero is independent of* t *when* $t > \kappa R_*$, *and then its absolute value is bounded by* $\kappa$. *Moreover, given a number* $c$ *which is a regular value on* X–Z *of the function* $|\varsigma^+|$, *there exists* $\kappa_c$ *such that if* $R_* > \kappa_c$, *then the absolute value of the spectral flow as* t *increases from* 0 *to* $2\sqrt{2}R_*c$ *is bounded by* $\kappa_c$.

This lemma is proved at the end of this subsection.

The next lemma is a direct corollary of Lemma 5.10 and Proposition 5.2: To set the stage, recall that the function $|\varsigma^+|$ on the radius $r_1$ tubular neighborhood of Z is equal to $\sqrt{2}$ times the distance to Z. Let $c$ denote a give nunber between $\frac{1}{8}r_1$ and $\frac{1}{4}r_1$. Any such choice for $c$ is a regular value of $|\varsigma^+|$; and the data sets D and D´ agree for any such choice where $|\varsigma^+| < 2c$.

**Lemma 5.11**: *There exists* $\kappa > 1$ *with the following significance: If* $R_* > \kappa$, *then the absolute value of the spectral flow for the* D *version of the family of operators depicted in* (5.19) *as* t *is inreased from* 0 *to* $2\sqrt{2}R_*c$ *is bounded by* $\kappa$.

Lemma 5.11 says in effect that the spectral flow for the STAGE 5 family of operators is significant only as t increases from $2\sqrt{2}R_*c$.

*Part 2*: This part of the proof describes the second and third comparison families of data sets. These data sets have manifold component $X_2$ which is the double of X across the boundary of a tubular neighborhood of Z. To say more about $X_2$, let r denote the $\frac{1}{1000}r_1$. Let Ñ denote the radius r tubular neighborhood of Z. The plan now is to take two copies of X–Ñ and glue them along their common boundary by the diffeomorphism that reverses the orientation of each component of Z. This new manifold $X_2$ is called 'the double' of X–Ñ.

To say more about $X_2$: Label the components of Z as $\{Z_1, \ldots, Z_p\}$ and for a given integer $k \in \{1, 2, \ldots, p\}$, let $N_{Zk}$ denote the radius $r_1$ tubular neighborhood of $Z_k$. Letting $\ell_k$ denote the length of $Z_k$, the radius $r_1$ tubular neighborhood $N_{Zk}$ has its Euclidean coordinate s for $\mathbb{R}/(\ell_k\mathbb{Z})$ and its $\mathbb{R}^3$ coordinate x with norm at most $r_1$. The metric appears with these coordinates as $ds \otimes ds + dx_k \otimes dx_k$ (see (5.10)). Use the coordinates (s, x) for the first copy of $N_k$ and let (s´, x´) denote the analogous coordinates for the second copy. The coordinate identification which is defined momentarily identifies the set where $|x| \in (\frac{1}{4}r, 4r)$ in the first copy of X with the set where $|x'| \in (\frac{1}{4}r, 4r)$ in the second copy. The diffeomorphism that implements the identification is this:

$$s' = -s \quad \textit{and} \quad x' = r^2 \frac{x}{|x|^2} \ .$$
(5.66)

The - sign for the identification makes this orientation preserving.

The Euclidean metric on the second version of $N_k$ pulls back to the metric below

$$ds \otimes ds + \frac{r^4}{|x|^4} dx_k \otimes dx_k \ ;$$
(5.67)

Although the gluing map in (5.66) is not isometric, the metric in (5.67) and the Euclidean metric on the $|x| \le r_1$ part of $N_k$ can be melded (as is done below in (5.68)) to obtain a metric that extends over the double of X. Here is the melded metric:

$$ds \otimes ds \ + \big(\chi(\tfrac{r}{|x|} -1) + \tfrac{r^4}{|x|^4} \chi(\tfrac{|x|}{r} -1)\big) \, dx_k \otimes dx_k \, .$$
(5.68)

which does extend over the doubled version of $N_k - (N_k \cap \tilde{N})$.

As for the principal bundle P and the line bundle $\mathcal{I}$: Both bundles have specific isomorphisms with respective product bundles over both the first and second version of $N_k$, and these isomorphisms are used to construct product bundle versions of P and $\mathcal{I}$ over the double of $N_k - (N_k \cap \tilde{N})$. Since the respective connections for P on the two versions of $N_k$ are the product connection, the product connection for the product bundle version of P over the double of $N_k - (N_k \cap \tilde{N})$ smoothly melds the two connections on the two versions of $N_k$. As for $\sigma$: Because this is the same constant map to $\mathfrak{su}(2)$ on each of the two versions of $N_k$, it is entirely consistent to define $\sigma$ on the double of $N_k - (N_k \cap \tilde{N})$ to be that same constant element in $\mathfrak{su}(2)$.

The last data element to extend over the double of $N_k - (N_k \cap \tilde{N})$ is $\varsigma^+$. To this end, remember that there are two canonical forms for $\varsigma^+$ on $N_k$. The first is

$$\varsigma^+ = x_k (ds \wedge dx_k + \tfrac{1}{2} \varepsilon_{kij} dx_i dx_j) \ .$$
(5.69)

The other canonical form is

$$\varsigma^+ = x_j \mathbb{M}_{jk}(s) \, (ds \wedge dx_k + \tfrac{1}{2} \varepsilon_{kij} dx_i dx_j)$$
(5.70)

with $\mathbb{M}_{jk}$ given in (5.12). Suppose first that $\varsigma^+$ has the form depicted in (5.69). In this case, the pull-back of the $(s', x')$ version of this 2-form by the diffeomorphism in (5.66) is

$$\frac{r^4}{|x|^4} x_k (ds \wedge dx_k + \tfrac{1}{2} \frac{r^2}{|x|^2} \varepsilon_{kij} dx_i dx_j)$$
(5.71)

whose wedge product with (5.69) is $2\frac{r^6}{|x|^6}\,ds\,dx_1dx_2dx_3$. This implies that the self dual projection on the double of $N_k-(N_k\cap\tilde{N})$ (as defined by the metric in (5.68)) of the 2-form

$$\chi(\tfrac{r}{|x|}-1)\,x_k(ds\wedge dx_k+\tfrac{1}{2}\,\varepsilon_{kij}dx_idx_j)\;+\;\frac{r^4}{|x|^4}\chi(\tfrac{|x|}{r}-1)\,x_k(ds\wedge dx_k+\tfrac{1}{2}\,\frac{r^2}{|x|^2}\,\varepsilon_{kij}dx_idx_j) \tag{5.72}$$

is nowhere zero. The self-dual projection of this 2-form (as defined by the metric in (5.68)) will play the role of $\varsigma^+$ on the double of $N_k-(N_k\cap\tilde{N})$ when $\varsigma^+$ on $N_k$ has the form depicted in (5.69).

Consider now the case when $\varsigma^+$ on $N_k$ has the form depicted in (5.70). In this case, the pull-back of the $(s', x')$ version of this 2-form by the diffeomorphism in (5.66) can be written as below:

$$\frac{r^4}{|x|^4}\,x_j\mathbb{M}_{jk}U_{kn}(ds\wedge dx_n+\tfrac{1}{2}\,\frac{r^2}{|x|^2}\,\varepsilon_{nij}dx_idx_j) \tag{5.73}$$

with $\{U_{jk}\}_{j,k=1,2,3}$ being functions on $\mathbb{R}^3-0$ defined by the rule below:

$$U_{jk} = -\,\delta_{jk} + 2\,\frac{x_jx_k}{|x|^2}\;. \tag{5.74}$$

These functions are the components of an x-dependent, $3\times 3$ special orthogonal matrix whose square is the identity matrix. This matrix (call it U) gives a map from the 2-sphere in $\mathbb{R}^3$ to SO(3). And, since $\pi_2(SO(3))$ is trivial, this matrix valued map is homotopic to the identity matrix. Let $\underline{U}$: $[0,1]\times S^2\to SO(3)$ denote a smooth homotopy that equals U near the 0 endpoint of the [0,1] factor and equals the identity matrix near the 1 endpoint of that factor. Use $\underline{U}$ to define a map to be denoted by $\hat{U}$ from $\mathbb{R}^3-0$ to SO(3) by the rule whereby $\hat{U}(x)$ is $\underline{U}(\frac{4|x|}{r}-1, x)$. Thus $\hat{U}$ equals U near where $|x|$ is $\frac{1}{4}r$, and $\hat{U}$ is the identity matrix near where $|x| = \frac{1}{2}r$. The desired $\varsigma^+$ is defined on the double of $N_k-(N_k\cap\tilde{N})$ to be the self-dual projection (as defined by the metric in (5.68)) of the following 2-form:

- *The 2-form is the* $(s',x')$ *version of (5.70) on the* $|x'|\geq 4r$ *part of the double of* $N_k-\tilde{N}$.
- *The 2-form is* $\frac{r^4}{|x|^4}\,x_j\mathbb{M}_{jk}\hat{U}_{kn}(ds\wedge dx_n+\tfrac{1}{2}\,\frac{r^2}{|x|^2}\,\varepsilon_{nij}dx_idx_j)$ *on the part of the double where* $|x|$ *is between* $\frac{1}{4}r$ *and* $\frac{1}{2}r$.
- *The 2-form is the one below on the part of the double where* $|x|$ *is greater than* $\frac{1}{2}r$.

$$\chi(\tfrac{r}{|x|}-1)x_j\mathbb{M}_{jk}(s)(ds\wedge dx_k+\tfrac{1}{2}\,\varepsilon_{kij}dx_idx_j)+\chi(\tfrac{|x|}{r}-1)\frac{r^4}{|x|^4}\,x_j\mathbb{M}_{jk}(ds\wedge dx_n+\tfrac{1}{2}\,\frac{r^2}{|x|^2}\,\varepsilon_{nij}dx_idx_j) \tag{5.75}$$

The square of this 2-form is everywhere positive and thus its self-dual projection is nowhere zero.

To summarize where things stand: A data set for the double of $X-\tilde{N}$ has been defined on the part of $X_2$ consisting of the doubles of the various $k \in \{1, \ldots, p\}$ versions of $N_k-(N_k \cap \tilde{N})$). The remaining part of $X_2$ is the disjoint union of two copies of $X-(\cup_{k=1,\ldots p} N_k)$. Two extensions of data set to the two copies of $X-(\cup_{k=1,\ldots p} N_k)$ in $X_2$ will be considered in what follows. To describe these, keep in mind that Part 1 considered two extensions of the data set on $\cup_k N_k$ to the whole of X. The interesting extension was denoted by D, and the rather boring extension was denoted by D´. The first of the two extensions to the whole of $X_2$ is denoted by $D_2$ and it is defined by using D on one of the copies of $X-(\cup_{k=1\ldots,p} N_k)$ and D´ on the other. The second extension to the whole of $X_2$ is denoted by $D_2´$ and it is defined by using D´ on both copies of $X-(\cup_{k=1\ldots,p} N_k)$.

Each of the data sets $D_2$ and $D_2´$ with the specficiation of a non-negative numbers $R_*$ and t have a corresponding version of $\mathcal{D}_\bullet$ depicted below:

$$\mathcal{D}_\bullet = \gamma^\alpha(\nabla_\alpha + \tfrac{1}{2}(\nabla_\alpha \hat{\rho})\hat{\rho}) + \sqrt{2}R_*|\varsigma^+|\hat{\rho} - \tfrac{1}{2}t\Gamma . \tag{5.76}$$

(This operator is defined over $X_2$; its self-adjoint there with a dense domain being the space of $\mathfrak{L}$ valued sections of $(\Lambda^+ \oplus \mathbb{R}) \otimes T^*X_2$.)

The lemma that follows bounds the absolute value of the spectral flow as a function of $t \in [0, \infty)$ for fixed large $R_*$ for the $D_2´$ version of the t-dependent family of operators that is depicted in (5.73).

**Lemma 5.12**: *There exists* $\kappa > 1$ *with the following significance: If* $R_* > \kappa$, *then the absolute value of the spectral flow for the* $D_2´$ *version of the 1-parameter family of operators depicted in (5.73) as* t *is increased from zero is independent of* t *when* $t > \kappa R_*$, *and then its absolute value is bounded by* $\kappa$. *Moreover, given a number* $c$ *which is a regular value on* $X_2$ *of the function* $|\varsigma^+|$, *there exists* $\kappa_c$ *such that if* $R_* > \kappa_c$, *then the absolute value of the spectral flow as* t *increases from* 0 *to* $2\sqrt{2}R_*c$ *is bounded by* $\kappa_c$.

This lemma is also is proved at the end of this subsection.

The next lemma is a direct corollary of Lemma 5.12 and Proposition 5.2: To set the stage, recall that the function $|\varsigma^+|$ on the radius $r_1$ tubular neighborhood of Z is equal to $\sqrt{2}$ times the distance to Z. Let $c$ denote a give nunber between $\frac{1}{8}r_1$ and $\frac{1}{4}r_1$. Any such choice for $c$ is a regular value of the $X_2$ version of $|\varsigma^+|$ and the data sets $D_2$ and $D_2´$ agree for any such choice where $|\varsigma^+| < 2c$.

**Lemma 5.13**: *There exists* $\kappa > 1$ *with the following significance: If* $R_* > \kappa$, *then the absolute value of the spectral flow for the* $D_2$ *version of the 1-parameter family of operators depicted in* (5.76) *as* t *is increased from* 0 *to* $2\sqrt{2}R_* c$ *is bounded by* $\kappa$.

The preceding lemma will be used with the lemma below (the latter is proved at the end of the subsection).

**Lemma 5.14**: *There exists* $\kappa > 1$ *with the following significance: If* $R_* > \kappa$, *then the absolute value of the spectral flow for the* $D_2$ *version of the 1-parameter family of operators depicted in (5.76) as* t *is increased to any value larger than* $\kappa R_*$ *differs by at most* $\kappa$ *from the absolute of the* D *version of that 1-parameter family as* t *is increased to that same value.*

The preceding lemma is also proved at the end of this subsection

*Part 3*: As is explained in the next paragraphs, the spectral flow for the $D_2$ version of the family in (5.76) as t increases from 0 and to any sufficiently large value can be calculated (when $R_*$ is large) using the strategy that was used in Part 4 of Section 3c.

The key point to keep in mind with regards to the upcoming calculation is that $|\varsigma^+|$ is nowhere zero on $X_2$. As a consequence of that, if $R_* > c_0$, and sufficiently large, then the t = 0 version of $\mathcal{D}_*$ in (5.76) has trivial kernel. In addition, the t = 0 version of (5.76) can be modified along a 1-parameter path of operators without incurring zero eigenvalues along the path so that the end member of the path has $|\varsigma^+|$ being constant. (Such a family is obtained from the t = 0 version of (5.76) as defined by the data set $D_2$ by changing only the functional form of $|\varsigma^+|$ according to the rule whereby the parameter $s \in [0, 1]$ member of the family has $|\varsigma^+|$ being the sum of s with (1 - s) times the original version of the function $|\varsigma^+|$.)

With it understood that $|\varsigma^+|$ is now constant and equal to 1, then for any $t \geq 0$, the square of the corresponding version of $\mathcal{D}_\bullet$ obeys the analog of (3.30) which is this:

$$\mathcal{D}_\bullet^2 = (\gamma_\alpha(\nabla_\alpha + \tfrac{1}{2}(\nabla_\alpha \hat{\rho})\hat{\rho})^2 + 2R_*^2(1 - \tfrac{t}{2\sqrt{2}R_*}\hat{\rho}\Gamma)^2 . \tag{5.77}$$

This implies in particular that $\mathcal{D}_\bullet$ has a non-zero kernel only in the event that $t = 2\sqrt{2}R_*$. Moreover, the kernel of $\mathcal{D}_\bullet$ for this value of t is the kernel of the operator

$$D = \gamma_\alpha(\nabla_\alpha + \tfrac{1}{2}(\nabla_\alpha \hat{\rho})\hat{\rho} \tag{5.78}$$

acting on the space of sections of the +1 eigenbundle of the endomorphism $\hat{\rho}\,\Gamma$. (The operator D commutes with $\hat{\rho}\,\Gamma$.) To determine the difference between the number of eigenvalues of D that cross 0 from below and the number that cross zero from above, keep in mind that both $\hat{\rho}$ and $\Gamma$ have square 1 and that they commute. Their joint eigenspaces decompose the bundle $((\Lambda^+\oplus\underline{\mathbb{R}})\otimes T^*X_2)\otimes\mathfrak{L}$ as a direct sum of 4 vector bundles each of real dimension 8. These are labeled below as $\mathcal{V}_{1,1}$, $\mathcal{V}_{1,-1}$, $\mathcal{V}_{-1,1}$ and $\mathcal{V}_{-1,-1}$ with the eigenvalue of $\hat{\rho}$ indicated by the left hand entry of the subscript and that of $\Gamma$ indicated by the right hand entry (the $\Gamma = 1$ eigenspace is the $(\Lambda^2\oplus\mathbb{R})\otimes\mathfrak{L}$ summand). The $\hat{\rho}\,\Gamma = 1$ eigenbundle is the 16 dimensional bundle given by the direct sum $\mathcal{V}_{1,1}\oplus\mathcal{V}_{-1,-1}$. Since D anti-commutes with both $\hat{\rho}$ and $\Gamma$, it maps the space of sections of $\mathcal{V}_{1,1}$ to the space of sections of $\mathcal{V}_{-1,-1}$ and vice-versa. Thus, it has the off diagonal form

$$D = \begin{pmatrix} 0 & \eth^* \\ \eth & 0 \end{pmatrix}, \tag{5.79}$$

with $\eth$ denoting the restriction of D to $C^\infty(X_2, \mathcal{V}_{1,1})$ and with $\eth^*$ denoting the formal adjoint of $\eth$. This decomposition of D is relevant to the question of the direction of the eigenvalues crossing zero because if $\psi$ denotes an eigenvector for $\mathcal{D}_\bullet$ with eigenvalue $\lambda$ for some given value of t, then

$$\tfrac{d}{dt}\langle\psi, \mathcal{D}_*\psi\rangle_2 = -\tfrac{1}{2}\langle\psi, \Gamma\psi\rangle_2 \tag{5.80}$$

where $\langle\ ,\ \rangle_2$ denotes the $L^2$ inner product on the space of sections of $((\Lambda^+\oplus\underline{\mathbb{R}})\otimes T^*X_2)\otimes\mathfrak{L}$. This formula implies directly that the elements in the kernel of $\mathcal{D}_\bullet$ from the kernel of D on the space of sections of $\mathcal{V}_{1,1}$ cross zero from above as t increases whereas those from the kernel of D on the space of sections of $\mathcal{V}_{-1,-1}$ cross zero from below. This is to say that the spectral flow for the family in (5.76) as t is increased from zero to any number greater than $2\sqrt{2}\,R_*$ is equal to -1 times the Fredholm index of the operator $\eth$ from (5.79).

As noted in Part 6 of Section 3c, if the bundle $\mathcal{I}$ is isomorphic to the product $\mathbb{R}$ bundle, then the index of $\eth$ has the form depicted in (3.33) with $t$ and [K] denoting cohomology classes on $X_2$, that of $t$ being Euler class of the line bundle $\mathfrak{L}$ and that of [K] being the Euler class of the orthogonal complement in $\Lambda^+$ of the span of $\varsigma^+$. Also, $b_1$ and $b_2^+$ are Betti numbers for $X_2$.

As noted in Part 7 of Section 3c, a version of the formula in (3.33) for the index of $\eth$ holds when $\mathcal{I}$ is not isomorphic to the product $\mathbb{R}$ bundle: This alternate version views both $t$ and [K] as 2-dimensional cohomology classes on the 2-fold cover of $X_2$ whose points are the set of norm 1 points in $\mathcal{I}$. Since both of these classes on this 2-fold cover change sign under the involution that interchanges the points in the fibers of the projection to $X_2$, the cup products of $t$ with [K] and $t$ with itself are well defined on $X_2$.

*Part 4*: This last part of the subsection has the proofs for Lemmas 5.10, 5.12, and 5.14.

***Proof of Lemma 5.10***: The assertion that the spectral flow is independent of t when t is greater than $c_0R_*$ follows from Proposition 5.1. Let $c_1$ denote this value of $c_1$

The fact that there is an $R_*$-independent upper bound for the spectral flow of the $\mathcal{D}_\bullet$ family (defined using D´) as t is increased from 0 to any $t > c_1R_*$ can be seen using the path depicted below (explanation follows):

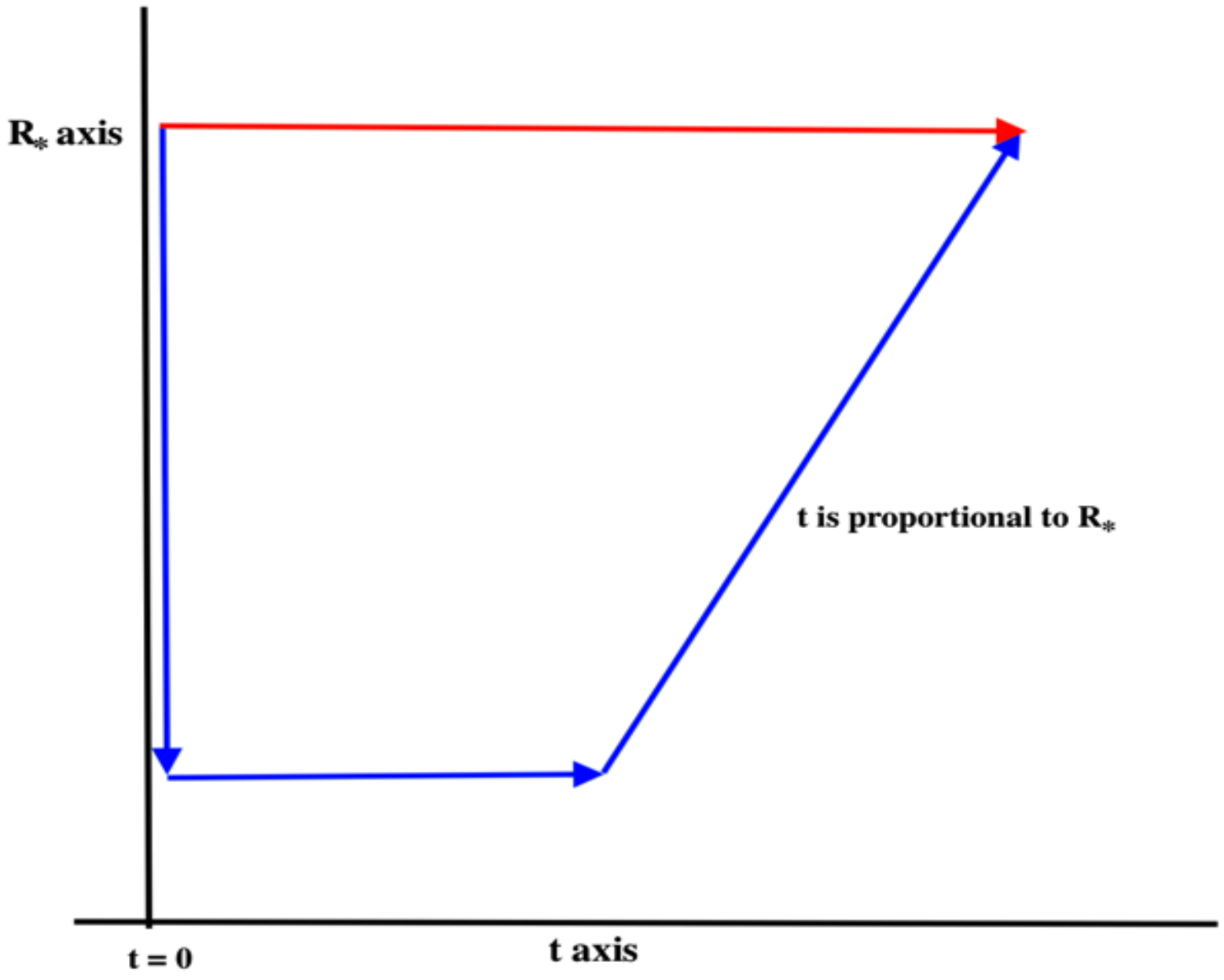

FIGURE 1

To explain: The spectral flow for a given value of $R_*$ along the red arrow path in the diagram is the spectral flow of interest. Use $R_{*1}$ to denote this given version of $R_*$. Supposing this $R_{*1}$ is large, there is (none-the-less) a $c_0$ version of $R_*$ (to be denoted by $R_{*0}$) which is such that Propositions 4.8 and 5.1 can be invoked for the D´ version of any $\mathcal{D}_\bullet$ path with any $R_*$ between $R_{*0}$ and $R_{*1}$. (In this regard: Remember that the connection A, the metrix and $\varsigma^+$ are all canonical for the D´ path which is why this $R_{*0}$ is bounded by $c_0$.) Now decrease $R_*$ from $R_{*1}$ to $R_{*0}$ with t = 0 and use Proposition 4.8 (see Section 4b) to conclude that there is an $R_{*1}$-independent bound for the absolute value of the spectral flow of the resulting path of operators when doing this. (This is the descending blue arrow at t = 0 in the diagram.) Next, note that when $t > c_1 R_{*0}$, there is no spectral flow by the first assertion of Lemma 5.10. And the absolute value of the spectral flow with $R_*$ is set at $R_{*0}$ and t increased from 0 to $t > c_1 R_{*0}$ is bounded by $c_0$ since all the parameters that are used for defining this path are independent of $R_{*1}$. This is the horizontal blue path in the diagram. Finally, with t set initially at $t_0 = c_1 R_{*0}$, increase $R_*$ from $R_{*0}$ to $R_{*1}$ and simultaneously increase t from $t_0$ according to the rule whereby $t(R_*) = c_{*1} R_*$. Then there is no spectral flow along this path either–invoke Proposition 5.1 at each point on the path. (This is the slanted blue arrow in the diagram.)

As for the third claim about the spectral flow from 0 to $2\sqrt{2}c R_*$: This follows using much the same strategy and diagram using Proposition 5.9 instead of Proposition 5.1 for the slanted path in the diagram. To elaborate: Make $R_{*0}$ larger if necessary (a lower bound will depend on $c$ now) so that Proposition 5.9 can be invoked for and $R_* \geq R_{*0}$. With this version of $R_{*0}$ understood, again use the descending blue path as before and then use the horizontal blue path as before, but stop the latter when $t = 2\sqrt{2}c R_{*0}$. The absolute value of the spectral flow along the descending blue path is at most $s_1$ and that along the horizontal blue path from 0 to $t = 2\sqrt{2}\, c\, R_{*0}$ is bounded by at most $s_2$ with both of these numbers being independent of $R_{*1}$ (but the latter depends on $c$). Next, increase $R_*$ from $R_{*0}$ to $R_{*1}$ while simultaneously increasing t by the rule $t(R_*) = 2\sqrt{2}c R_*$. By virtue of Proposition 5.9, the absolute value of the spectral flow along the latter path has a $c_0$ upper bound independent of $R_*$ that will be called $s_3$ It then follows directly that the absolute value of the spectral flow for the fixed $R_{*1}$ path as t is increased from 0 to $2\sqrt{2}\, c R_{*1}$ is no greater than $s_1 + s_2 + s_3$.

***Proof of Lemma 5.12***: The argument here is essentially identical to that for Lemma 5.10.

***Proof of Lemma 5.14***: Let $c = \frac{1}{6} r_1$. The absolute value of the spectral flow as t is increased from 0 to $2\sqrt{2}\, c R_*$ is bounded by $c_0$ if $R_* > c_0$, this being an instance of what is asserted by Lemma 5.13. With this understood, consider invoking Proposition 5.2 with $c_0 = c$ and with $c_1$ chosen to be larger than 1 plus the maximum of $|\varsigma^+|$. The comparison data sets are $\mathrm{D}_2$ for the first data set and a second data set to be denoted by $\mathrm{D}_2''$ which is as follows: The manifold component is the disjoint union of two copies of X with $\mathcal{D}_\bullet$ defined on one copy using the data set $\mathrm{D}$ and with $\mathcal{D}_\bullet$ defined on the other using the data set $\mathrm{D}'$. According to Propositions 5.1 and 5.2, if $R_* > c_0$ then the respective spectral flows for the $\mathrm{D}_2$ and $\mathrm{D}_2''$ versions of the family of operators as t is increased from $2\sqrt{2}\, c R_*$ to any $t > 2\sqrt{2}R_*c_1$ differ by at most $c_0$ from the number of eigenvalues of the $D_2$ version of the operator in (5.73) (defined for at $t = 2\sqrt{2}c$) between numbers $E_0$ and $-E_0$ with $E_0$ denoting a number that is at most $e^{-\sqrt{R_*}/c_0}$. And, meanwhile, an instance of Proposition 5.9 asserts that the latter number is also bounded by $c_0$.

### h) Summary for the STAGE 5 spectral flow and that of STAGES 6 and 7

This subsection has two purposes. The first purpose is to summarize the calculations in the previous subsection regarding the spectral flow for the STAGE 5 path depicted in (5.5). The second purpose is to comment on the STAGE 6 and 7 path spectral flows and then summarizes the analysis of all of the STAGES with a final proposition (Proposition 5.15). The two purposes are dealt with in Parts 1 and 2 of this subsection.

*Part 1*: Part 3 of the previous subsection found that the spectral flow for the STAGE 5 path depicted in (5.5) differs by at most $c_0$ from a number that can be written as $-1$ times the Fredholm index of the operator on the manifold $X_2$, this being the operator $\mathfrak{d}$ that appears in (5.76). As noted in that same part of the subsection, the index of $\mathfrak{d}$ can be written in terms of characteristic classes using the Atiyah-Singer index theorem. When $\mathcal{I}$ is isomorphic to the product line bundle, this formula is

$$\text{index}(\mathfrak{d}) = (1 + b_2^+ - b_1) + (t \cdot t - t \cdot [K]) . \tag{5.81}$$

with $t$ denoting the Euler class of the line bundle $\mathfrak{L}$ on the manifold $X_2$, and with [K] denoting the orthogonal complement in the bundle of self-dual 2-forms on $X_2$ to the nowhere zero section $\varsigma^+$. (The Betti numbers in (5.81) are those of $X_2$.) When the line bundle $\mathcal{I}$ is not isomorphic to the product $\mathbb{R}$-bundle, then the same formula holds with it understood that both $t$ and [K] should be interpreted as classes on the 2-fold cover of $X_2$ that change sign when pulled back by the action of the generator of the covering space $\mathbb{Z}/2$ action. The pairings in (5.81) should then be viewed as pairings of twisted cohomology classes on $X_2$.

The formula in (5.81) can be interpreted solely in terms of data on the original manifold X. This reinterpretation follows directly. To do this, start with the original connection A, the original version of the isometry $\sigma$, and the original version of $\varsigma^+$ on X. Let F denote for now the 2-form $2\langle\sigma F_A\rangle$, this being a closed, $\mathcal{I}$-valued 2-form on X. With regards to the constructions of Section 2, the 2-form F is $2mr\,\varsigma$ with $\varsigma$ denoting the harmonic, $\mathcal{I}$-valued 2-form whose self dual part is $\varsigma^+$.

If $\mathcal{I}$ is isomorphic to the product $\mathbb{R}$ bundle, then F is -i times the curvature 2-form for an orthogonal connection on the line bundle $\mathfrak{L}$. This product line bundle case will be considered first in detail. The story when $\mathcal{I}$ is not isomorphic to the product $\mathbb{R}$ bundle will be summarized at the end.

When $\mathcal{I}$ is isomorphic to the product $\mathbb{R}$ bundle, the bundle $\mathfrak{L}$ can be viewed as a complex line bundle over X. As such, it is isomorphic to the product $\mathbb{C}$ near Z and, with an isomorphism chosen, the connection on $\mathfrak{L}$ coming from A can be written near any given component of Z with respect to the corresponding product connection $\theta_0$ on $\mathfrak{L}$ as

$$\hat{A} = \theta_0 + i(2b_0 + \hat{p})\,ds + i\,\mathrm{A} \tag{5.82}$$

with $b_0$ being the constant that appears in (4.61), with $\hat{p}$ being an even integer multiple of $2\pi/\ell$ where $\ell$ here denotes the length of the relevant component of Z, and with the norm of $\mathrm{A}$ obeying the bound $|\mathrm{A}| \le c_0 m\, r\,\mathrm{dist}(\cdot, Z)$. The product structure near Z is further specified by the requiring that $\hat{p} = 0$. With this choice understood, then $\hat{A}$ is modified near Z to make it flat on a neighborhood of Z. This modification replaces (5.82) with

$$\hat{A}_2 = \theta_0 + i\chi_1(2b_0\,ds + \mathrm{A}) \tag{5.83}$$

with $\chi_1$ as described in the beginning of Section 4c. In this regard, $\chi_1$ is defined using a small, positive number $r_1$: It is equal to zero where the distance to Z is less than $r_1$ and it is equal to 1 where the distance to Z is greater than $2r_1$. This modification from $\hat{A}$ to $\hat{A}_2$ changes $\mathrm{F}$ to $\mathrm{F}_2$ with $\mathrm{F}_2$ having compact support well inside X-$\tilde{N}$ with $\tilde{N}$ denoting the radius $\frac{1}{1000} r_1$ tubular neighborhood of Z. Here is $\mathrm{F}_2$:

$$\mathrm{F}_2 = 2b_0\, d\chi_1 \wedge ds + d\chi_1 \wedge \mathrm{A} + \chi_1 \mathrm{F}\,. \tag{5.84}$$

Now remember that $X_2$ is the double of X–$\tilde{N}$. The operator $\eth$ on $X_2$ is defined using a connection on the bundle $\mathfrak{L}$ on $X_2$ whose curvature 2-form is $\mathrm{F}_2$ on the first copy of X–$\tilde{N}$ and zero on the second copy of X–$\tilde{N}$. As a consequence, the first Chern class of the bundle $\mathfrak{L}$ is represented in $H^2(X; \mathbb{Z})$ by the 2-form $\frac{1}{2\pi i}\mathrm{F}_2$.

As for the class [K] on $X_2$: The Poincare´ dual of this class can be represented by a piece-wise embedded submanifold in $X_2$ that is defined as follows: To start, note that this class [K] is the Euler class of the orthogonal complement to $\varsigma^+$ in the bundle of self-dual 2-forms on $X_2$. To represent this class, go first to X–Z where the orthogonal complement to the original version of $\varsigma^+$ on X–Z is an oriented 2-plane bundle there. Any such bundle has a section (denoted by $s$) which can be chosen so that its zero locus is a smooth, oriented properly embedded 2-dimensional submanifold X–Z whose closure in X can be viewed as a cycle (linear combination of 2-simplices) with boundary equal to 2Z. Let Σ denote the zero locus of such a submanifold. The X–$\tilde{N}$ part of this version of $\Sigma$ can then be glued to a second version across $\tilde{N}$ to obtain a 2-dimensional submanifold in $X_2$ that represents the Poincare´ dual of the class [K].

With the preceding understood, it then follows that what is denoted by $t\cdot[K]$ in (5.81) is equal to $\frac{1}{2\pi}$ times the integral of $F_2$ over $\Sigma$. That is an integer that differs by at most $c_0|b_0|$ from the integral of $F$ over $\Sigma$. In the context of the versions of (A, $\omega$) from Section 2, that integral is, in turn, proportional (by a universal factor) to the integral of $mr\varsigma$ over $\Sigma$. Thus, the $t\cdot[K]$ term can contribute a factor proportional to $mr$ to the right hand side of (5.81).

Meanwhile, what is denoted by $t\cdot t$ in (5.81) is proportional to the integral over $X-\tilde{N}$ of the 4-form $F_2 \wedge F_2$. This in turn is equal to the integral over X of $F_2 \wedge F_2$ which is equal to the integral of $F \wedge F$ over X. The choice of the number $r_1$ makes no difference in this regard, nor does the term with $b_0$ in (5.84). As for the X-integral of $F \wedge F$, this is proportional (by a universal factor) to the first Pontrjagin class of the bundle ad(P). Thus this term contributes at most $c_0$ to the right hand side of (5.81).

Now suppose that $\mathcal{I}$ is not isomorphic to the product $\mathbb{R}$ bundle. Even so, $\mathcal{I}$ is isomorphic to the product $\mathbb{R}$ bundle on $\tilde{N}$ and an isomorphism there should be chosen. Having done that, then $\mathcal{L}$ on $\tilde{N}$ can again be viewed as a complex line bundle and then the construction of $F_2$ can be made as before with it understood that it extends from $\tilde{N}$ into X as an $\mathcal{I}$-valued 2-form. In this regard, $F_2 \wedge F_2$ is as before, a real valued 4-form whose integral over X is universally proportional to the first Pontrjagin class of the bundle ad(P). This integral computes (up to a universal factor), what is denoted by $t\cdot t$ in (5.81).

Meanwhile, there is a version of $\Sigma$ in this case which is now a non-orientable surface, but orientable in a tubular neighborhood of Z. It is the zero locus of a section over $X-Z$ of the orthogonal compliment in $\mathrm{ad}(P) \otimes \mathcal{I}$ of $\varsigma^+$. (The bundle $\mathcal{L}$ over the complement of $\Sigma$ is isomorphic to $\mathbb{R} \oplus \mathcal{I}$ with the $\mathbb{R}$ factor being the span of the section.) An isomorphism between $\mathcal{I}$ and the product $\mathbb{R}$ bundle over $\tilde{N}$ orients $\Sigma$ on $\tilde{N}-Z$ and its closure there can be viewed as a 2-cycle with boundary 2Z. In any event, the same local calculations can be made to see that the integral of $F_2$ over $\Sigma$ makes sense. (This is because $\wedge^2 T\Sigma$ isomorphic to $\mathcal{I}$ and $F_2$ is $\mathcal{I}$ valued). The integral of $\frac{1}{2\pi} F_2$ over $\Sigma$ computes the $t\cdot[K]$ term in (5.78). This integral differs by at most $c_0$ from the integral of $\frac{1}{2\pi} F$ over $\Sigma$. In the context of the solutions from Section 2, the integral $F$ over $\Sigma$ is $mr$ times the integral of $\varsigma$ over $\Sigma$.

*Part 2*: The end member of the STAGE 5 spectral flow is the t = T version of the operator depicted in (5.5). If T is sufficiently large (and much greater than R), then it follows from the Bochner-Weitzenboch formula in (5.15) that there is no spectral flow in the STAGE 6 path (remember that this path takes $R_*$ to zero, changes A to a fiducial connection, and removes the T term in (5.5)).

Meanwhile, the partial Bochner-Weitzenboch formula below in (5.85) implies that there is no spectral flow in the STAGE 7 path.

$$(\gamma_\alpha \nabla_{A_0\alpha} + \tfrac{1}{2} t\Gamma)^2 = (\gamma_\alpha \nabla_{A_0\alpha})^2 + \tfrac{1}{4} t^2 . \tag{5.85}$$

The proposition below summarizes the spectral flow calculations for the 7 stages of the operator path that was described at the start of Section 5. This proposition refers to the surface $\Sigma$ in X–Z from Part 1 above. Remember also that in the context of a solution $(A,\omega)$ to (1.2) from Section 2, the $\mathcal{I}$-valued 2-form $\langle\sigma F_A\rangle$ is $mr\varsigma$.

**Proposition 5.15**: *Given* $\Xi > 1$ *and* $k \in \mathbb{Z}$, *there exists* $\kappa > 1$ *with the following significance: Let* $P \to X$ *denote a principle* SO(3) *bundle with* k *being the first Pontrjagin class of* ad(P). *Fix* $m < \frac{1}{\kappa}$ *and* $r > \kappa$; *and let* $(A,\omega)$ *denote a solution to (1.2) obeying (3.20) with* $\varsigma^+$ *obeying* $|\varsigma^+| \ge \frac{1}{\Xi}\mathrm{dis}(\cdot,Z)$. *Then the spectral flow for the operator*

$$\mathcal{D} = \gamma_\alpha \nabla_{A\alpha} + \tfrac{1}{\sqrt{2}} \varsigma^+{}_k \rho_k [\sigma, \cdot\,] - \tfrac{1}{2} m\Gamma - \tfrac{1}{2} m$$

*differs by at most* $\kappa$ *from the integral of* $\frac{1}{\pi}\langle\sigma F_A\rangle$ *over the surface* $\Sigma$.

The next two subsections remark on some instances of this proposition

**i) With regards to the constructions in Section 2**

The first remark concerns the $|\omega|$-diverging sequences of solutions to (1.2) that are described in Section 2: The sequences in that section for principle bundles with non-zero first Pontrjagin class have versions of the self-dual 2-form $\varsigma^+$ that are not independent of $r$. Rather, the corresponding versions of $\varsigma^+$ have $r$-dependence of the form $\varsigma^+ = \varsigma^+{}_0 + \mathcal{O}(r^{-\alpha})$ with $r$ increasing along the sequence and with $\alpha$ being a positive constant. The crucial point here is that Proposition 5.15 can be applied with the same version of $\kappa$ for all large $r$ solutions in an $|\omega|$ diverging sequence of (parametrized by an unbounded set of values for $r$ ) if the large $r$ versions of $\varsigma^+$ can be written as $\varsigma^+{}_0 + \mathcal{O}(r^{-\alpha})$ with $\varsigma^+{}_0$ obeying the restriction that $|\varsigma^+{}_0| \ge \frac{1}{\Xi}\mathrm{dis}(\cdot,Z)$.

By way of an example: When $t$ is defined as in (2.1), the sequence was initially parametrized by the value of the integer q with the integer p being $\frac{1}{4}(-q^2 + \varepsilon k)$. In this case, the sequence can be parametrized by $r$ which is proportional to p and thus to $q^2$.

With this parametrization, the self-dual 2-form $\varsigma^+_0$ will be a non-zero multiple of the self-dual part of the harmonic 2-form representing the class $x_2 + x_3$ (which is a class with zero self-cup product). Meanwhile, the $\mathcal{O}(r^{-\alpha})$ part will be the sum of an $\mathcal{O}(r^{-1/2})$ term proportional to the self dual part of the harmonic 2-form representing the class $x_1$ and an $\mathcal{O}(r^{-1})$ term proportional to the self dual part of the harmonic 2-form representing the class $x_2$ - $x_3$.

**j) With regards to the case of symplectic 4-manifolds**

Suppose now that X is a symplectic 4-manifold which is to say that it has a nowhere vanishing, closed 2-form which I will denote by w with $w \wedge w$ defining the orientation for X. It then has a Riemannian metric which makes w self-dual and normalized so that $w \wedge w$ is the volume form on X. This is to say that $|w| = 1$. I will use this metric in what follows. I will also assume that X has indefinite cup-product pairing on its second cohomology and that the self dual and anti-self dual Betti numbers $b_2^+$ and $b_2^-$ are at least 3. This is an interesting case because one might hope to mimic the constructions for Kähler manifolds in Section 3d to construct examples where the spectral flow diverges as $r$ gets ever larger and other examples where the spectral flow is bounded as $r$ gets ever larger.

The discussion below regarding this symplectic case has two parts. The first part concerns describes the surface $\Sigma$ in general terms and some the issues that arise that make the story interesting. The second part describes some specific instances of Proposition 5.15 in the symplectic case.

*Part 1*: Introduce by way of notation $K_X$ to denote the orthogonal complement to w in the bundle $\Lambda^+$. This is an oriented $\mathbb{R}^2$ bundle which, when viewed as a complex line bundle, is isomorphic to the 'canonical' bundle for any almost complex structure on X that is compatible with the symplectic form w.

Now suppose that $\varsigma$ is a harmonic 2-form that comes from a solution to (1.2) as described in (3.20) for an instance when the line bundle $\mathcal{I}$ is isomorphic to the product $\mathbb{R}$ bundle. The self-dual part of $\varsigma$, which is $\varsigma^+$, can be written as $\varsigma^+ = \mu w + c$ where $\mu$ is a non-negative number and $c$ is a closed, self-dual 2-form with the integral of $w \wedge c$ over X being zero. Note in particular that $c$ need not be point-wise orthogonal to w. In particular, $c$ can be written as $f_c w + c'$ with $f_c$ denoting a function on X and $c'$ denoting a section of $K_X$. With this understood, let S denote the zero locus of $c'$. If $c'$ vanishes transverally, then S is a smooth, oriented, 2-dimensional submanifold in X whose fundamental class when pushed into X is Poincare´ dual to the Euler class of $K_X$.

With regards to S (supposing it isn't empty): The zero locus of $\varsigma^+$ is the locus in S where the function $f_c+\mu$ is zero. That locus in S is the $\varsigma^+$ version of the set Z. The $\varsigma^+$ version of Σ can then be taken to be the disjoint union of the $f_c+\mu < 0$ and $f_c+\mu > 0$ parts of S with it understood that the former part, where $f_c+\mu < 0$, is oriented *opposite* to its orientation as an open set in S.

The orientation switch with regards to the $f_c+\mu < 0$ part of S has the following consequence: Whereas the integral of the closed form $rm\,\varsigma$ over S is proportional to the cup-product pairing between the respective Euler classes of $K_X$ and $\mathfrak{L}$, that pairing is not what is relevant for Proposition 5.15. Of concern in Proposition 5.15 is the integral of $rm\varsigma$ over Σ. In particular, the integral of $rm\varsigma$ over the $f_c+\mu < 0$ part of Σ is opposite in sign to how that integral contributes to the cup-product pairing between the two Euler classes. (The integral of ς over the $f_c+\mu > 0$ part of Σ is the same sign as its contribution to the cup-product pairing between those Euler classes.)

Interesting challenge is to understand how to compute these integrals and whether there is a way to choose ς so that the Σ integral is zero.

*Part 2*: This second part of the subsection describes examples where the complications with regards to Σ are not relevant; these are specifically cases where the Euler class of $K_X$ has non-positive self cup product. (The respective self-dual and anti-self dual second Betti numbers should be at least 3 also).

To start, assume that the self-cup product of the Euler class of $K_X$ is negative. This class can then be depicted as the class of a harmonic 2-form that can be written as $\alpha w + c + \beta y$ where the notation is as follows: What is denoted by α is a positive real number; what is denoted by $c$ is a closed, self-dual 2-form with the integral of $w \wedge c$ over X being zero; what is denoted by β is a positive real number; and what is denoted by y is a closed, anti-self dual 2-form with the integral over X of $y \wedge y$ being $-V$ where V denotes the volume of X. Let $c^2 V$ denote the integral over X of $c \wedge c$. The assertion that the self-cup product of the Euler class of $K_X$ is negative says in effect that

$$\alpha^2 + c^2 - \beta^2 < 0. \tag{5.86}$$

What follows describes a construction of ς for use in (3.20) and Proposition 5.15 with the following properties:

- *The self-dual part $\varsigma^+$ is nowhere zero and its orthogonal complement in $\Lambda^+$ is isomorphic to* $K_X$.
- *Its cohomology class has zero self-cup product.*

- *Its cohomology class has zero cup product with the Euler class of* $K_X$.

(5.87)

As a first attempt at $\varsigma$, consider the closed form $\beta w + \alpha y$. The cohomology class of this form has zero cup-product with the Euler class of $K_X$, but its self-cup product is positive (since $\beta^2 > \alpha^2$). If the anti-self dual Betti number of X is at least two, then there is a closed, anti-self dual 2-form $y'$ with the X-integral of $y' \wedge y$ being zero and with the integral of $y' \wedge y'$ being $-\nu$. Having chosen such a form, then second attempt at $\varsigma$ is $\beta w + \alpha y + (\beta^2 - \alpha^2)^{1/2} y'$ which is a harmonic 2-form obeying the conditions specified by the three bullets in (5.87).

A possible problem with taking $\varsigma$ as just described is that it need not define a rational cohomology class. To remedy this, one can add a very small harmonic 2-form to $\varsigma$ to make the result a rational class. This form can and should be chosen so that its cohomology class has zero cup product pairing with the Euler class of $K_X$. Let $e$ denote this harmonic 2-form. It's norm can be less than any given positive bound (but not in general zero). A third attempt at $\varsigma$ is then

$$\beta w + \alpha y + (\beta^2 - \alpha^2)^{1/2} y' + e.$$

(5.88)

This third attempt is now rational, it has zero cup product with the Euler class of $K_X$ but it's self-cup product need not be zero. As explained momentarily, if the self-cup product of the class in (5.88) is positive, and if $b^{2-} \geq 3$, then there is a very small (an upper bound on its norm gets ever smaller as the norm of what is denoted by $e$ tends to zero), harmonic 2-form whose cohomology class is rational, orthogonal to the class of the 2-form in (5.88) and to the Euler class of $K_X$, and has negative self-cup product which can be added to (5.87) to give the desired version of $\varsigma$. If the self cup product of the class in (5.88) is negative, then there is a very small harmonic 2-form whose cohomology class is rational, is orthogonal to the class in (5.88) and to the Euler class of $K_X$ and has positive self-cup product which can be added to (5.88) to give the desired version of $\varsigma$. In either case, if the norm of $e$ and this last added harmonic 2-form are sufficiently small, then the self-dual part of the resulting closed form will not vanish anywhere; and so the three bullets in (5.87) are obeyed.

To find the desired 2-form $e$: Note first that the manifold X cannot be a spin manifold because if $[K_X]\cdot[K_X]$ is negative, then there is a symplectic sphere with self-intersection -1. This being the case, the cup product pairing matrix on $H^2(X;\mathbb{Z})$ can be diagonalized. Let $\{P_j: j = 1, \ldots, b^{2+}\}$ denote a orthogonal set of classes with self-cup product equal to 1, and let $\{Q_\alpha: \alpha = 1, \ldots, b^{2-}\}$ denote an orthogonal set of classes with zero cup product against the P's and with self cup product equal to -1.

Write $[K_X]$ using this basis as $\sum_j n_j P_j + \sum_\alpha m_\alpha Q_\alpha$ with the n's and m's being integers that obey $\sum_j n_j^2 - \sum_\alpha m_\alpha^2 < 0$. Meanwhile, write the class depicted in (5.88) as $\sum_j a_j P_j + \sum_\alpha b_\alpha Q_\alpha$ with the a's and b's being rational numbers obeying $\sum_j a_j n_j - \sum_\alpha b_\alpha m_\alpha = 0$ and $\sum_j a_j^2 - \sum_\alpha b_\alpha^2 = \varepsilon$ with ε being a rational number whose norm is very small. With the preceding understood, fix two distinct P's to be denoted by $\{P_{j(1)}, P_{j(2)}\}$ and then two distinct Q's to be denoted by $\{Q_{\alpha(1)}, Q_{\alpha(2)}\}$. Fix a vector $(x_1, x_2)$ in $\mathbb{R}^2$ and numbers $\delta_1$ and $\delta_2$ from the set {-1, 1}. Use this data to define the class $s = \sum_{k=1,2} x_k (P_{j(k)} + \delta_k Q_{\alpha(k)})$. This class $s$ has zero self-cup product; and it is orthogonal to $[K_X]$ if $\sum_k x_k (n_{j(k)} - \delta_k m_{\alpha(k)})$ is zero. Meanwhile, its cup product with the class in (5.88) is $\sum_k x_k (a_{j(k)} - \delta_k b_{\alpha(k)})$. The latter won't be zero unless $a_{j(k)} = \mu\, n_{j(k)}$ and $b_{\alpha(k)} = \mu\, m_{\alpha(k)}$ for each k with $\mu$ being a real number. This can't happen for all choices of assignments $\{k \to (j(k), \alpha(k), \delta_k)\}_{k=1,2}$ because if that were the case, then the class in (5.88) would be a multiple of $[K_X]$ whereas it is, by construction, orthogonal to $K_X$. In particular, if ε is sufficiently small, there is an assignment $\{k \to (j(k), \alpha(k), \delta_k)\}_{k=1,2}$ such that system of equations

- $\sum_k x_k (n_{j(k)} - \delta_k m_{\alpha(k)}) = 0$ ,
- $\sum_k x_k (a_{j(k)} - \delta_k b_{\alpha(k)}) = -\frac{1}{2}\varepsilon$

(5.89)

can be solved with the solution having rational entries and obeying $x_1^2 + x_2^2 < c_0 \varepsilon^2$. (The latter bound comes by comparing the solution to (5.89) with the analogous equations for the case of (5.88) with $e = 0$. If $e = 0$, there is an assignment $\{k \to (j(k), \alpha(k), \delta_k)\}_{k=1,2}$ and then a solution $(x_1, x_2)$ to (5.89) with norm $\mathcal{O}(\varepsilon^2)$ but perhaps not with rational entries. The $e \neq 0$ version of (5.89) with $e$ having small norm using the same assignment will have a nearby solution (in $\mathbb{R}^2$) with rational entries.)

Adding the version of the class $s$ as just described to the class in (5.88) gives a non-zero, rational class with zero self cup product and zero cup product with $[K_X]$. As noted above, if ε is sufficiently small, then the corresponding harmonic 2-form will have nowhere zero self-dual part.

To find a version of ς that meets the requirements of the first two bullets in (5.87) but not the third, repeat the preceding construction with β in (5.88) replaced by α.

Now suppose that the Euler class of $K_X$ has zero self-cup product. The left hand side of (5.86) is zero in this case. If $c \neq 0$, then the construction just described for the negative self-cup product case can be repeated almost verbatim. If $c = 0$, then ς can be taken equal to the harmonic 2-form that represents the Euler class of $K_X$ unless that Euler class is a torsion class In the latter case, there are no $K_X$ constraints to speak of if the condition in the first bullet of (5.87) is met.